\documentclass[11pt,a4paper]{article}

\usepackage[utf8]{inputenc}
\usepackage[T1]{fontenc}
\usepackage{amsmath,amssymb,amsfonts,amsthm}
\usepackage{mathtools}
\usepackage[margin=1in]{geometry}
\usepackage{enumitem}
\usepackage{booktabs}
\usepackage{shuffle}     
\usepackage[numbers,sort&compress]{natbib}
\usepackage{xcolor}
\usepackage[colorlinks=true,linkcolor=blue,citecolor=blue,urlcolor=blue]{hyperref}
\usepackage{graphicx}

\allowdisplaybreaks[1]
\newcommand{\figorbox}[1]{%
  \IfFileExists{#1}{\includegraphics[width=\linewidth]{#1}}{%
    \fbox{\parbox[c][0.28\linewidth][c]{0.98\linewidth}{\centering\ttfamily
      [figure \detokenize{#1} not bundled with source; see companion code]}}}}

\theoremstyle{plain}
\newtheorem{theorem}{Theorem}[section]
\newtheorem{proposition}[theorem]{Proposition}
\newtheorem{lemma}[theorem]{Lemma}
\newtheorem{corollary}[theorem]{Corollary}
\theoremstyle{definition}
\newtheorem{definition}[theorem]{Definition}
\newtheorem{construction}[theorem]{Construction}
\theoremstyle{remark}
\newtheorem{remark}[theorem]{Remark}

\newcommand{\RR}{\mathbb{R}}
\newcommand{\ZZ}{\mathbb{Z}}
\newcommand{\NN}{\mathbb{N}}
\newcommand{\E}{\mathbb{E}}
\newcommand{\Prob}{\mathbb{P}}
\newcommand{\Sig}{\mathrm{Sig}}
\newcommand{\TM}{T^{(M)}(\RR^d)}
\newcommand{\Sblk}{\mathbb{S}}
\newcommand{\hSblk}{\widehat{\mathbb{S}}}
\newcommand{\dd}{\mathrm{d}}
\newcommand{\pvar}{p\text{-}\mathrm{var}}
\newcommand{\onevar}{1\text{-}\mathrm{var}}
\DeclareMathOperator{\Var}{Var}
\DeclareMathOperator{\Cov}{Cov}
\newcommand{\norm}[1]{\left\lVert #1 \right\rVert}
\newcommand{\bignormL}[2]{\bigl\lVert #1 \bigr\rVert_{#2}}
\newcommand{\hnorm}[1]{\left\vert\!\left\vert\!\left\vert #1 \right\vert\!\right\vert\!\right\vert}
\newcommand{\hs}{\mathrm{HS}}
\newcommand{\cP}{\mathcal{P}}
\newcommand{\cF}{\mathcal{F}}
\newcommand{\cN}{\mathcal{N}}
\newcommand{\cB}{\mathcal{B}}
\newcommand{\cW}{\mathcal{W}}
\newcommand{\cH}{\mathcal{H}}
\newcommand{\frH}{\mathfrak{H}}         
\newcommand{\Dsep}{\mathsf{D}}          
\newcommand{\Bop}{\mathsf{B}}           
\newcommand{\Cop}{\mathsf{C}}           
\newcommand{\cD}{\mathcal{D}}          
\newcommand{\one}{\mathbf 1}
\newcommand{\TA}{T\!\bigl((\RR^d)\bigr)}          
\newcommand{\alphabet}{\mathcal{A}}               
\newcommand{\emptyword}{\varnothing}               

\title{\bf Weak and strong discretization errors for signatures\\ of Gaussian processes}
\author{Florian Huber}
\date{October 2026}

\begin{document}
\maketitle

\begin{abstract}
Let $X$ be a continuous Gaussian process in $\RR^d$, $d\ge2$, with independent components of roughness
$H\in(\tfrac14,1)$, possibly with a linear drift, and let $X^n$ be its piecewise-linear interpolant on
a mesh of width $\delta$. We study how well the truncated signature of $X^n$ approximates that of the
natural lift of $X$, in mean (the weak error, the difference of the two expected signatures) and in
$L^2$ (the strong error). Over a class of drivers defined by quantitative regularity conditions, the
weak error at every fixed truncation level is $O(\delta^{\gamma(H)})$ with $\gamma(H)=\min(2H,4H-1)$,
uniformly over the class, and for truncation levels $M\ge4$ fractional Brownian motion attains this
rate. For centred drivers the weak error vanishes at levels $1$ to $3$, Gaussian or not. It first
appears at level~$4$, where we compute its leading constant for independent fractional Brownian
motion in $d=2$: on $(\tfrac12,1)$ it is $1/(2\sqrt3)$, independent of $H$, and on $(\tfrac14,\tfrac12)$
it is given by a $\zeta$-series. For fractional Brownian motion with drift, a level-$3$ error of the
same order appears for $H\ge\tfrac12$. Below $\tfrac12$ the drift reaches the rate only at levels $\ge5$ (proved for $d\ge3$). The strong error is of exact order $\delta^{2H-1/2}$ for
fractional Brownian motion with $H\in(\tfrac12,1)$, so in mean square the weak error lies a full power
of $\delta$ below it. For $H\le\tfrac12$ the strong error is $O(\delta^{\vartheta})$ for every
$\vartheta<2H-\tfrac12$, uniformly over the class, and the exponent $2H-\tfrac12$ cannot be improved.
Finally, we construct pairs of laws in the class whose mesh data have the same law but whose
expected signatures differ by order $\delta^{\gamma(H)}$. Hence no quantity computed from the mesh data
has mean within $o(\delta^{\gamma(H)})$ of the expected signature uniformly over the class.
\end{abstract}



\tableofcontents

\section{Introduction}\label{sec:intro}

Let $X=(X_t)_{t\in[0,1]}$ be a continuous process in $\RR^d$. Its signature $\Sig(X)$ is the
collection of all iterated integrals of $X$, and its truncation $\Sig^{\le M}(X)$ keeps those of order
at most $M$. When it has infinite radius of convergence, the full expected signature $\E[\Sig(X)]$ characterises
the law of $X$ \citep{ChevyrevLyons2016,ChevyrevOberhauser2022}, and its truncations are the features behind
many model-free methods for time series \citep{Lyons2014,Levin2013}. In practice $X$ is seen only on a
mesh $t_j=j\delta$, $\delta=1/n$, and one computes the signature of the piecewise-linear interpolant
$X^n$. We ask how much is lost by this step. There are two natural measures of the loss: the
weak error
\[
   \cW_{\le M}\ :=\ \E\bigl[\Sig^{\le M}(X^n)\bigr]-\E\bigl[\Sig^{\le M}(X)\bigr],
\]
the difference of the two expected signatures, and the strong error
\[
   \eta_\delta\ :=\ \bigl(\E\norm{\Sig^{\le M}(X^n)-\Sig^{\le M}(X)}_\hs^2\bigr)^{1/2}.
\]
By Jensen's inequality $\norm{\cW_{\le M}}_\hs\le\eta_\delta$. We show that for Gaussian drivers   $X$ (the
process whose signature is approximated)   with independent components, the weak error is much smaller than the strong one. We compute it to leading
order at level~$4$, and show that no procedure based on the mesh data removes it uniformly over the
class.

\paragraph{The setting.} We work over the class $\cP(H,L)$ of Definition~\ref{def:class}, with $d\ge2$
throughout. For $d=1$, the signature of a path is a function of its increment, so that
$\cW_{\le M}\equiv0$ and $\eta_\delta\equiv0$. The driver is $X_t=\theta te_1+G_t$, where $G$ is a
centred continuous Gaussian process with independent components of roughness $H\in(\tfrac14,1)$,
normalised by $\Var(G^i_1)=1$, and $|\theta|\le L$. Its signature is that of the natural lift,
the limit of the lifts of piecewise-linear approximations. The class is cut out by five quantitative
conditions (C1)--(C5), on the covariance and on the signatures of single cells, and by a moment bound
on the lift. These are, up to one branch of (C1), the facts the upper-bound proofs use. We verify them for three laws: fractional
Brownian motion (fBm) with independent components, its drifted version, and the tied-down
perturbation of fBm built in \S\ref{sec:mesh} (Appendix~\ref{app:verif}). Other drivers, such as fractional Ornstein--Uhlenbeck
processes, are natural candidates, but we do not verify the conditions for them. The centred
sub-class $\theta=0$ is $\cP_0(H,L)$. Throughout,
\[
   \gamma(H)\ =\ \min(2H,\,4H-1).
\]
The value $H=\tfrac12$ is the crossover between the Young range $(\tfrac12,1)$, where the rate is
$\delta^{2H}$, and the rough range $(\tfrac14,\tfrac12)$, where it is $\delta^{4H-1}$.

\paragraph{The mechanism.} Two effects of the interpolation produce the weak error. First, inside a
mesh cell $[t_a,t_{a+1})$ the interpolant has less variance than the process itself. For fBm, $X^n_t=(1-\lambda)X_{t_a}+\lambda
X_{t_{a+1}}$ with $\lambda=(t-t_a)/\delta$, and $\E(X^i_t)^2-\E(X^{n,i}_t)^2=\lambda(1-\lambda)\delta^{2H}$ up to
terms whose integral over $[0,1]$ is $O(\delta^{2H+1}+\delta^2)$. This deficit averages to
$\delta^{2H}/6$ over a cell. It moves the level-$4$ coordinates by order $\delta^{2H}$. Half of this cell average, $\tfrac1{12}$ in units of $\delta^{2H}$, enters with opposite signs
at the two ends of a pair of cells. These two contributions are the constants
$\Phi^{\rm s}=-\Phi^{\rm t}=-\tfrac1{12}$ of Lemma~\ref{lem:reflect}. The leading constants of
the level-$4$ coordinates are built from them. In the Young range this gives
$\cB=1/(2\sqrt3)$. Second, the interpolant has no L\'evy area inside a cell, while the
process has mean-zero cell areas of size $\delta^{2H}$. Their sum over the $n$ cells has a variance of
order $n\,\delta^{4H}=\delta^{4H-1}$, the square of the strong rate $\delta^{2H-1/2}$, and the
level-$4$ coordinates depend on the second moment of the areas at the same order $\delta^{4H-1}$. Hence
$\gamma(H)=\min(2H,4H-1)$: the variance deficit dominates above the crossover, while the cell areas dominate below.

\paragraph{Main result 1: the rate of the weak error.} For every fixed $M$ and every
$H\in(\tfrac14,1)$,
\[
   \sup_{P\in\cP(H,L)}\norm{\cW_{\le M}(P)}_\hs\ =\ O\bigl(\delta^{\gamma(H)}\bigr)
\]
(Propositions~\ref{prop:alllevels}, \ref{prop:alllevelsdrift}, \ref{prop:roughalllevels}
and~\ref{prop:drift5}), and for $M\ge4$ the rate is attained by fBm with independent components (Theorems~\ref{thm:2H} and~\ref{thm:roughrate}). On the centred class the weak error first appears at
level~$4$, where we compute its leading constant. For $d=2$ and fBm with independent components,
\[
   \norm{\cW_4}\ =\ \cB(H)\,\delta^{\gamma(H)}\bigl(1+o(1)\bigr),\qquad
   \cB(H)=\begin{cases}\dfrac1{2\sqrt3}, & H\in(\tfrac12,1),\\[6pt]
   \dfrac1{4\sqrt3}, & H=\tfrac12,\\[6pt]
   2\,|\kappa_{\rm alt}(H)|, & H\in(\tfrac14,\tfrac12),\end{cases}
\]
with the value at $H=\tfrac12$ exact for every mesh (Proposition~\ref{prop:C1exact},
Theorem~\ref{thm:2H}(i), Corollary~\ref{cor:kadj}). The even levels $6,\dots,M$ are bounded by
$O(\delta^{\gamma(H)})$. They may contribute at that order (numerically level~$6$ does below the
crossover, Remark~\ref{rem:M6}), and we do not evaluate their constants, so $\cB(H)$ is the constant
of $\norm{\cW_4}$. On the Young range, it is built from two first-moment constants
$\Phi^{\rm s}=-\Phi^{\rm t}=-\tfrac1{12}$, which do not depend on the word or on $H$
(Lemma~\ref{lem:reflect}). On the rough range $\kappa_{\rm alt}(H)<0$ is a $\zeta$-series
(Propositions~\ref{prop:kappaclosed} and~\ref{prop:kappasign}). The $\delta^{4H-1}$ term is carried by
the alternating word $1212$ and the nested word $1221$, with opposite coefficients
$\kappa_{\rm nest}=-\kappa_{\rm alt}$, while the adjacent word $1122$ has none
(Remark~\ref{find:kappa}). The constants for general $d$ follow
by counting pairs of letters (Remark~\ref{rem:dscope}, Corollary~\ref{cor:kadj}).

\paragraph{Main result 2: the level structure.} For a centred driver with independent components, the
weak error vanishes at levels $1$, $2$ and $3$, Gaussian or not, and at level~$4$ it depends on the law
only through the covariances of the components (Propositions~\ref{prop:level4}
and~\ref{prop:nongauss}). The reason is that the expected L\'evy area vanishes. A drift changes this.
Let $d=2$, $H\in(\tfrac12,1)$, and let the driver have stationary increments whose increment variance $\phi(u)=\E[(G^i_{t+u}-G^i_t)^2]$ satisfies Hypothesis~(S) of Lemma~\ref{lem:momentgen} (or let $G$ be
Brownian motion, $H=\tfrac12$). Then the drift creates the level-$3$ error
\[
   \cW_3\ =\ \theta\,\phi(\delta)\bigl(-\tfrac1{12},\ \tfrac16,\ -\tfrac1{12}\bigr)
   +O\bigl(\delta^2\bigr)
\]
on the words $122,212,221$, with the same two constants as at level~$4$
(Proposition~\ref{prop:level3drift}). When $\phi(\delta)\asymp\delta^{2H}$ this is of the order of the
level-$4$ error, so ``the weak error starts at level~$4$'' is a statement about $\cP_0$. Below the
crossover, the level-$3$ drift error is of the smaller order $\delta^{2H}$, and the level-$\le4$ constant
is unchanged (Proposition~\ref{prop:drift4}). The drift reaches the order $\delta^{4H-1}$ only at
levels $\ge5$ (proved for $d\ge3$; Propositions~\ref{prop:drift5} and~\ref{prop:attain5}).

\paragraph{Main result 3: the strong error.} For fBm with independent components and $M\ge2$,
\[
   \eta_\delta\ \asymp\ \delta^{2H-1/2},\qquad H\in(\tfrac12,1)
\]
(Proposition~\ref{prop:strong}). The lower bound comes from the L\'evy areas of the cells, and the
upper bound from a Wick expansion of the difference $\Sig^{\le M}(X^n)-\Sig^{\le M}(X)$ against a deterministic kernel
\citep{HuberStrong}. Hence on $(\tfrac12,1)$ and for $M\ge4$ the squared weak error
$\norm{\cW_{\le M}}_\hs^2\asymp\delta^{4H}$ lies a full power of $\delta$ below the mean-square strong
error $\eta_\delta^2\asymp\delta^{4H-1}$. For $H\in(\tfrac14,\tfrac12]$ the strong error is
$O(\delta^\vartheta)$ for every $\vartheta<2H-\tfrac12$, uniformly over the class
(Proposition~\ref{prop:strongunif}, a direct consequence of a theorem of \citet{FrizRiedel2014}), and
for fBm it is at least of order $\delta^{2H-1/2}$ (Lemma~\ref{lem:stronglower}).

\paragraph{Main result 4: the weak error cannot be removed using the mesh data.} To a fBm, we add a small
Gaussian ripple that vanishes at every mesh node (Construction~\ref{con:bridge}). The two laws $P_0$
(fBm) and $P_1$ (fBm plus ripple) lie in $\cP_0(H,L)$ for $L\ge L_0(H,\epsilon)$, $\epsilon$ denoting the size of the
ripple, and their mesh data have the same law
(Proposition~\ref{prop:adm}). Yet for $M\ge4$ their expected signatures differ by order
$\delta^{\gamma(H)}$ on the whole range $(\tfrac14,1)$ (Propositions~\ref{prop:sep},
\ref{prop:sepbrownian} and~\ref{prop:seprough}). On the rough range the sign of the separation comes
from two spectral identities in the eigenbasis of the covariance of tied-down fBm
(Lemma~\ref{lem:opsigns}). Consequently, for every quantity $T$ computed from the mesh data, with
$\mu(P):=\E_P[\Sig^{\le M}(X)]$ the expected signature under the law $P$,
\[
   \max_{j\in\{0,1\}}\ \bigl\lVert\E_{P_j}[T]-\mu(P_j)\bigr\rVert_\hs\ \ge\
   \tfrac12\norm{\mu(P_1)-\mu(P_0)}_\hs\ \ge\ c\,\delta^{\gamma(H)}
\]
(Corollary~\ref{cor:undetectable}).   So no quantity computed from the mesh data has a mean within
$o(\delta^{\gamma(H)})$ of the expected signature uniformly over the class, and the rate of Main
result~1 is the best possible for such procedures. This is a statement about the class: the second law is built for the mesh at
hand. For a fixed, known law the weak error may be computable. For fBm with known $H$, for instance,
the leading constants of \S\S\ref{sec:bias}--\ref{sec:rough} allow the leading term to be subtracted.

\paragraph{Motivation and companion paper.} When the expected signature is estimated from one long
trajectory cut into $K$ windows, each observed on a mesh, the natural estimator averages the
signatures of the interpolated windows. Its bias is exactly the weak error, and
its mean squared error splits into this bias and a fluctuation term. The statistical consequences of
our results, namely minimax rates under weak dependence and an ergodic central limit theorem, are
derived in the companion paper \citet{HuberII}. The present paper is self-contained and deals with a
single path on $[0,1]$.

\paragraph{Related work.} \citet{Schenck2026} proves finite-sample bounds for the block-averaging
estimator of the expected signature, with a discretization term controlled by the strong $L^2$
approximation error (his Assumption~2.6,
$\E\norm{S^{(M)}(\mathbf X)-S^{(M)}(\mathbf X^{(n)})}_\hs^2\le C_{\rm approx}\delta^{\gamma_{\rm S}}$). His
Proposition~2.7 states it for independent-coordinate fractional Ornstein--Uhlenbeck drivers,
with $\gamma_{\rm S}=2H$ for $H\in(\tfrac12,1)$ (his Remark~2.8) and $\gamma_{\rm S}=4H-1-\varepsilon$, any
$\varepsilon>0$, for $H\in(\tfrac14,\tfrac12)$. His part~(i), for Gaussian semimartingales, gives
$\gamma_{\rm S}=1-\varepsilon$ and includes $H=\tfrac12$. His Discussion asks whether ``direct $L^2$ techniques
can close the empirical gap for the Gaussian case''. For fBm drivers on $(\tfrac12,1)$,
Proposition~\ref{prop:strong} answers this: the mean-square strong error is of exact order
$\delta^{4H-1}$. However, his class and ours differ. We compare the two bounds at the end of
\S\ref{subsec:strongyoung}.

For the L\'evy area of fBm, \citet[Thm.~2]{NeuenkirchTindelUnterberger2010} prove the mean-square
bound $\delta^{4H-1}$ for the trapezoidal scheme, which is the L\'evy area of the piecewise-linear
interpolant, for all $H>\tfrac14$. The published version conjectures that this rate is exact, and the
arXiv version (arXiv:0902.0497, Theorem~1.2) states the exact asymptotics
$\alpha_4(H)\,n^{-4H+1}(1+o(1))$ on $[0,1]$ for $H>\tfrac12$. So exactness at level~$2$ on
$(\tfrac12,1)$ is originally their result, and our lower bound is a short independent proof of it. For $H\le\tfrac12$
\citet{FrizRiedel2014} prove the strong rate $2H-\tfrac12$ for Gaussian rough paths, up to an
arbitrarily small loss in the exponent. They note that the rate is sharp, as ``may be seen from a
precise result concerning L\'evy's stochastic area for fBM'', and that it had been conjectured in
\citet{DeyaNeuenkirchTindel2012}.
Our Proposition~\ref{prop:strongunif} is a direct consequence of their theorem, and
Lemma~\ref{lem:stronglower} is a short self-contained proof of the sharpness. What is new on the
strong side is the upper bound $\delta^{2H-1/2}$ at every signature level on $(\tfrac12,1)$
\citep{HuberStrong}. We are not aware of it in the literature. Discretization schemes
for differential equations driven by fBm are a related but different problem
\citep{HuLiuNualart2016,LiuTindel2019}, as is multilevel Monte Carlo for rough differential equations
driven by Gaussian processes \citep{BayerFrizRiedelSchoenmakers2016}.

For fBm with $H>\tfrac12$, \citet{Passeggeri2020} shows that every even-level coordinate
of the expected signature of the piecewise-linear approximation with $n$ segments converges at rate
$n^{-2H}=\delta^{2H}$. This is the fBm case of the upper bound of Main result~1 on $(\tfrac12,1)$. The rate was believed to be sharp. Theorem~\ref{thm:2H} confirms this at level~$4$. What is new here on $(\tfrac12,1)$ is the
uniformity over the class, the matching lower bound and the exact leading constants. The rough range
$(\tfrac14,\tfrac12)$ is new as well.

The convergence of the lifts of piecewise-linear approximations of fBm for $H>\tfrac14$ goes back to
\citet{CoutinQian2002}. For general Gaussian processes whose covariance has finite $\rho$-variation,
$\rho<2$, the limit is the natural lift of \citet[Thm.~15.33]{FrizVictoir2010}, and
\citet{FrizRiedel2014} give convergence rates in rough-path metrics. \citet{BaudoinCoutin2007}
computed the expected signature of fBm for $H>\tfrac12$, and \citet{CassFerrucci2024} gave Wiener-chaos
expansions of Gaussian signatures for $H\in(\tfrac14,1)$. Integrability of
Gaussian rough paths is studied by \citet{CassLittererLyons2013}. \citet{HagerPelizzari2026} compute
expected signatures of heterogeneous paths exactly: an invertible change of coordinates removes mixed
integration against the irregular component and turns the transformed expected signature into a
deterministic integral of increment correlators over a simplex, evaluated by quadrature. Their
Appendix~A shows that, for time-augmented one-dimensional fBm, the weak discretization error of the
piecewise-linear approximation is $O(\Delta t^{2H})$ at fixed truncation level. Remark~\ref{rem:hp}
compares this with our rough-regime rate: on $(\tfrac14,\tfrac12)$ the exponent $4H-1$ is smaller, so
the error is larger, and it is carried by a mechanism absent in one dimension.

\paragraph{Organisation.} Section~\ref{sec:prelim} recalls signatures and the shuffle algebra, and
\S\ref{sec:setup} introduces the class and proves the moment and Lipschitz lemmas. Section~\ref{sec:R3}
sets up the Gaussian Wick sub-model. Section~\ref{sec:bias} treats the weak error on the Young range
and the level structure, \S\ref{sec:rough} the rough range, and \S\ref{sec:strong} the strong error.
Section~\ref{sec:mesh} constructs the pairs of laws with the same mesh data, and \S\ref{sec:discussion}
discusses open questions. Appendix~\ref{app:verif} verifies the class conditions, and
Appendix~\ref{app:fp} introduces the
finite-part functional that identifies the two closed forms of the rough-regime constant.
Every numerical value quoted in the paper is reproduced by the code in the companion repository
\texttt{https://github.com/FlorianHuber-Math/expected-signature-minimax}, which accompanies both
papers. No proof uses a numerical value as an input.

\section{Preliminaries: the signature and the shuffle algebra}\label{sec:prelim}

Conventions follow \citet{LyonsCaruanaLevy2007,FrizVictoir2010} and \citet{Reutenauer1993}.
Readers familiar with rough paths may skip to \S\ref{subsec:notation}.

\subsection{Tensor algebra and the signature}\label{subsec:sig}
Let $\alphabet=\{1,\dots,d\}$ be the alphabet and $\alphabet^\ast$ the free monoid of words
$w=(i_1\cdots i_m)$, $i_j\in\alphabet$, of length $|w|=m$, with $\emptyword$ the empty word. The
tensor algebra $\TA=\bigoplus_{m\ge0}(\RR^d)^{\otimes m}$ has the words as a basis ($w$
identified with $e_{i_1}\otimes\cdots\otimes e_{i_m}$) and truncation
$\TM=\bigoplus_{m=0}^M(\RR^d)^{\otimes m}$; both carry the Hilbert--Schmidt inner product
$\langle\cdot,\cdot\rangle$ making the words orthonormal, and $\langle w,\xi\rangle$ is the
coefficient of $w$ in $\xi\in\TA$.

For $x\in C^{\onevar}([0,1];\RR^d)$ the signature is
$\Sig(x)=\bigl(\langle w,\Sig(x)\rangle\bigr)_{w\in\alphabet^\ast}\in\TA$ with
\[
   \langle w,\Sig(x)\rangle=\int_{0<t_1<\cdots<t_m<1}\dd x^{i_1}_{t_1}\cdots\dd x^{i_m}_{t_m},
   \qquad w=(i_1\cdots i_m),
\]
and $\langle\emptyword,\Sig(x)\rangle=1$. A geometric $p$-rough path $\mathbf x$ already prescribes
the levels $1,\dots,\lfloor p\rfloor$. Its signature $\Sig(\mathbf x)$ is the unique extension of these
levels to all levels with finite $p$-variation \citep{Lyons1998}, \citep[Thm.~9.5]{FrizVictoir2010},
and it is the limit of the signatures of smooth approximations of $\mathbf x$. For $p\ge2$ a path does
not determine its signature. For the driver $X$, ``the signature of $X$'' always means
$\Sig(\mathbf X)$ for the natural lift $\mathbf X$ of Definition~\ref{def:class}.

We call the $m$ positions of a word $w=(i_1\cdots i_m)$ the slots of the coordinate
$\langle w,\Sig(x)\rangle$. Slot $k$ carries the differential $\dd x^{i_k}_{t_k}$, and the coordinate
is linear in each slot separately. Two uses of this recur. If $x=h+g$ is a sum of two paths, then
expanding every slot gives $\langle w,\Sig(x)\rangle$ as a sum of $2^m$ iterated integrals, in each
of which every slot carries either $\dd h$ or $\dd g$. If $g$ is centred Gaussian, the expectation of
an iterated integral against $\dd g$ alone is a sum over pairings of its slots (Wick's formula),
and with independent components only pairs of slots with the same letter contribute.

The truncated signature $\Sig^{\le M}(x)$ is the image of $\Sig(x)$ in $\TM$. Writing
$\Sig(x)_{s,t}$ for the signature of the restriction of $x$ to $[s,t]$, signatures satisfy
Chen's identity $\Sig(x)_{s,u}=\Sig(x)_{s,t}\otimes\Sig(x)_{t,u}$, $s\le t\le u$, and are group-like:
they lie among the characters of the shuffle algebra below.

\subsection{The shuffle algebra}\label{subsec:shuffle}
The shuffle product $\shuffle:\alphabet^\ast\times\alphabet^\ast\to\ZZ\langle\alphabet\rangle$ is
defined recursively by $\emptyword\shuffle w=w\shuffle\emptyword=w$ and
\[
   (i\,u)\shuffle(j\,v)=i\,\bigl(u\shuffle(j\,v)\bigr)+j\,\bigl((i\,u)\shuffle v\bigr),
   \qquad i,j\in\alphabet,\ u,v\in\alphabet^\ast,
\]
extended bilinearly. $(\TA,\shuffle)$ is a commutative, associative, unital algebra, the
shuffle algebra. For a geometric signature the pointwise product of two coordinates is the
shuffle's coordinate, the shuffle identity:
\begin{equation}\label{eq:shuffle-id}
   \langle u,\Sig(x)\rangle\,\langle v,\Sig(x)\rangle=\langle u\shuffle v,\Sig(x)\rangle
   \qquad(u,v\in\alphabet^\ast).
\end{equation}
Two consequences recur.

\emph{(i) Products close one level up.} Since $u\shuffle v$ has length $|u|+|v|$, a product of two
level-$M$ coordinates is a linear combination of level-$2M$ coordinates. So covariances of level-$M$
coordinates lie within $T^{(2M)}(\RR^d)$, not within $T^{(M)}(\RR^d)$.

\emph{(ii) The area is the only new level-2 information.} With $u=(i)$, $v=(j)$, \eqref{eq:shuffle-id} gives
$\langle ij,\Sig\rangle+\langle ji,\Sig\rangle=\langle i,\Sig\rangle\langle j,\Sig\rangle$. The
symmetric part of level~$2$ is thus determined by the endpoint increment, and only the antisymmetric
part, the L\'evy area $A_{ij}=\tfrac12(\langle ij,\Sig\rangle-\langle ji,\Sig\rangle)$, is new. This is
the start of the level-4 weak-error mechanism (\S\ref{sec:bias}): for the centred independent-component
drivers of Definition~\ref{def:class}, $\E[A_{ij}]=0$, so the weak error cannot appear at level~$2$.

\subsection{The mesh data}\label{subsec:windows}
The driver $X=(X_t)_{t\in[0,1]}$ has roughness $H$: informally, its increments over an
interval of length $u$ are of size $u^H$, so its paths are $\alpha$-H\"older for $\alpha<H$ and of
finite $p$-variation for $p>1/H$, as for $\mathrm{fBm}(H)$. Formally, $H$ enters only through the
conditions of Definition~\ref{def:class}. It is not seen continuously:
it is sampled on the uniform mesh $t_j=j\delta$, $j=0,\dots,n$, $\delta=1/n$. The mesh data are
the nodal vector $\cD=(X_{t_j})_{j\le n}$, and anything computed from the samples is a function of
$\cD$. The object of interest is the expected signature
$\mu(P):=\E_P[\Sig^{\le M}(X)]\in\TM$ under the law $P$ of $X$.

\subsection{Notation}\label{subsec:notation}
\begin{center}
\renewcommand{\arraystretch}{1.15}
\begin{tabular}{@{}l p{0.5\linewidth}@{}}
\toprule
Symbol & Meaning\\
\midrule
$\alphabet=\{1,\dots,d\}$, $\alphabet^\ast$ & alphabet; words $w=(i_1\cdots i_m)$, length $|w|$\\
$\TA$, $\TM$ & tensor algebra; truncation to level $M$\\
$\langle w,\xi\rangle$, $\shuffle$ & word coefficient; shuffle product \eqref{eq:shuffle-id}\\
$\Sig(x),\ \Sig^{\le M}(x)$ & signature; truncated signature\\
$\Sblk=\Sig^{\le M}(X)$, $\mu=\E[\Sblk]$ & truncated signature of the driver; expected signature\\
$X^n$; $\hSblk_\delta=\Sig^{\le M}(X^n)$ & piecewise-linear interpolant; its signature\\
$\cW_{\le M}$, $\cW_m$, $\eta_\delta$ & weak error; its level-$m$ component; strong error (\S\ref{sec:setup})\\
$\delta=1/n$, $I_a$, $\cD$ & mesh width; mesh cells $I_a=[t_a,t_{a+1})$, $a=0,\dots,n-1$ ; mesh data\\
$H$ & roughness (Hurst) exponent\\
$\hnorm\cdot_{\pvar}$, $\rho_{\pvar}$ & homogeneous $p$-variation norm; inhomogeneous metric \eqref{eq:rhopvar}\\
$\cP(H,L)\supset\cP_0$ & model class; centred sub-class ($\theta=0$)\\
$\varrho(s,t)=c_H|t-s|^{2H-2}$ & fBm increment-covariance density, $c_H=H(2H-1)$\\
$r(k)$, $\widetilde r(u)$, $\varrho(u)$ & fGn autocovariance; its interpolant; $c_Hu^{2H-2}$\\
$\Delta E_\bullet$, $\cB(H)$ & word-differences at level~$4$; level-4 leading constant\\
$\Phi^{\rm s}=-\Phi^{\rm t}=-\tfrac1{12}$, $\Phi_1=\tfrac16$ & the two first-moment constants in cell units; their difference\\
$\kappa_{\rm alt}(H)$, $\kappa_{\rm sep}(H)$ & rough-regime constant; separation moment of $(P_0,P_1)$\\
$\rho_A(m)$, $\sigma^2_\infty$ & cross-cell L\'evy-area correlation; its long-run variance\\
$\varrho_{\max}(I,J)$ & canonical correlation between first chaoses (Proposition~\ref{prop:cancorr})\\
$D_a$ & within-cell area defect\\
$\Dsep_\bullet(\epsilon)$ & Brownian separation coefficients\\
$\Bop$, $\Cop$ & covariance operators of tied-down fBm and of fBm\\
$\frH=\cH\otimes\RR^d$ & chaos Hilbert space (Lemma~\ref{lem:cellproj})\\
$\Theta(H)$, $W(H)$ & the two functionals of tied-down fBm in Lemma~\ref{lem:opsigns}\\
$\cH_1(I_a)$ & first Wiener chaos of the increments on a cell\\
$\norm\cdot_\hs$, $A_{ij}$ & Hilbert--Schmidt norm; L\'evy area\\
\bottomrule
\end{tabular}
\end{center}

\section{Setup, the class, and moment lemmas}\label{sec:setup}

Throughout, $d\ge2$ (for $d=1$ see the Introduction), the truncation level $M\in\NN$ is fixed, and
$X=(X_t)_{t\in[0,1]}$ is $\RR^d$-valued with $X_0=0$. With the truncated tensor algebra $\TM$ of
\S\ref{sec:prelim}, the truncated signature of the driver, its expectation, and the truncated signature
of the interpolant are
\[
   \Sblk:=\Sig^{\le M}(X),\qquad \mu:=\E[\Sblk],\qquad \hSblk_\delta:=\Sig^{\le M}(X^n),
\]
where $X^n$ is the piecewise-linear interpolation of $X$ on the uniform mesh of width $\delta=1/n$.
For rough drivers $\Sig$ is the signature of the natural lift (Definition~\ref{def:class}). All
signatures are geometric (Stratonovich), so \eqref{eq:shuffle-id} holds pathwise.

\paragraph{The two errors.} The weak error and the strong error of the interpolant are
\[
   \cW_{\le M}:=\E\bigl[\hSblk_\delta\bigr]-\mu=\E[\Sig^{\le M}(X^n)]-\E[\Sig^{\le M}(X)],\qquad
   \eta_\delta(P):=\bigl(\E\norm{\hSblk_\delta-\Sblk}_\hs^2\bigr)^{1/2}.
\]
Both are functionals of the law $P$ of $X$, and $\norm{\cW_{\le M}}_\hs\le\eta_\delta(P)$ by Jensen's
inequality. We write $\cW_m$ for the level-$m$ component of $\cW_{\le M}$, so that
$\norm{\cW_{\le M}}_\hs^2=\sum_{m\le M}\norm{\cW_m}^2$, and we also call the weak error the
bias.

\paragraph{The $p$-variation exponent and the norms.} We fix once and for all the exponent at which
$p$-variation norms of the driver are taken in Definition~\ref{def:class}:
\begin{equation}\label{eq:pdelta}
   p_\delta(H)\ :=\ \begin{cases}\tfrac12\bigl(\tfrac1H+2\bigr), & H>\tfrac12,\\[2pt]
                                  \tfrac12\bigl(\tfrac1H+4\bigr), & H\le\tfrac12,\end{cases}
\end{equation}
so that $p_\delta(H)\in(1/H,2)$ above the crossover and $p_\delta(H)\in(1/H,4)$ at or below it. $\norm\cdot_{\pvar}$ is non-increasing in $p$ at fixed $\lfloor p\rfloor$ and
for fBm$(H)$ it blows up as $p\downarrow1/H$, since fBm$(H)$ has infinite $(1/H)$-variation.
For a path $\mathbf x$ with values in the step-$\lfloor p\rfloor$ nilpotent group,
$\hnorm{\mathbf x}_{\pvar;[s,t]}:=\sup_{(t_i)\subset[s,t]}(\sum_i\norm{\mathbf x_{t_i,t_{i+1}}}^p)^{1/p}$ is the
homogeneous $p$-variation norm, $\norm\cdot$ being the Carnot--Carath\'eodory norm
\citep[(8.1)]{FrizVictoir2010}. For $m\ge0$ we write $\pi_m:\TA\to(\RR^d)^{\otimes m}$ for the projection onto level~$m$,
$\pi_m(\xi)=\sum_{|w|=m}\langle w,\xi\rangle\,w$, and $\pi_{\le M}:=\sum_{m\le M}\pi_m:\TA\to\TM$
for the truncation. The same maps are used on $\TM$ and on the step-$N$ nilpotent group, which
sits inside $T^{(N)}(\RR^d)$. With this notation, by \citep[Prop.~7.49]{FrizVictoir2010},
$\norm g$ is comparable with $\max_{i\le\lfloor p\rfloor}|\pi_i(g)|^{1/i}$ for $g$ in the
step-$\lfloor p\rfloor$ group. Throughout, $\rho_{\pvar}$ is
the inhomogeneous rough-path metric
\begin{equation}\label{eq:rhopvar}
   \rho_{\pvar}(\mathbf x,\mathbf y)\ :=\ \max_{1\le i\le\lfloor p\rfloor}\,
   \bigl\lVert \mathbf x^{(i)}-\mathbf y^{(i)}\bigr\rVert_{(p/i)\text{-}\mathrm{var}} ,
\end{equation}
with no root taken at level $i$, and $\norm{\mathbf x}_{\pvar}:=\rho_{\pvar}(\mathbf x,\one)$ is the
inhomogeneous norm, $\one$ denoting the constant path at the unit element $(1,0,\dots,0)$, all of
whose increments have vanishing components at levels $i\ge1$. Here $\mathbf x^{(i)}_{s,t}:=\pi_i(\mathbf x_{s,t})$ is the level-$i$ component of the increment
$\mathbf x_{s,t}$, and $\rho_{\pvar}$ is the inhomogeneous $p$-variation metric of
\citet[Def.~8.6]{FrizVictoir2010} with $N=\lfloor p\rfloor$. For $\lambda>0$, the dilation
$\mathrm{dil}_\lambda$ multiplies the level-$i$ component of a rough path by $\lambda^i$.

\subsection{The class}

The class-uniform upper bounds of \S\S\ref{sec:bias}--\ref{sec:strong} hold over the following
class, and the lower bounds and the pairs of laws of \S\ref{sec:mesh} are exhibited inside it. Its
conditions are the quantitative inputs of those proofs (with one exception, Remark~\ref{rem:C1scope}).

 \begin{definition}[the class]\label{def:class}
Fix $H\in(\tfrac14,1)$ and $L<\infty$. A law $P$ of $X$ belongs to $\cP(H,L)$ if
\[
   X_t\ =\ \theta\,t\,e_1\ +\ G_t ,\qquad t\in[0,1],\qquad |\theta|\le L,
\]
where $e_1$ is the first standard basis vector of $\RR^d$ and $G$ is a centred continuous Gaussian
process on $[0,1]$ with $G_0=0$ and independent components sharing one scalar increment covariance
$R_G([s,t]\times[u,v]):=\E[G^i_{s,t}G^i_{u,v}]$, normalised by $\Var(G^i_1)=1$, such that
\emph{(C1)}--\emph{(C5)} below hold and the moment condition $\E\norm{\mathbf X}_{\pvar}^{r}\le L$
holds with $r:=2M$, at $p=p_\delta(H)$ of \eqref{eq:pdelta}, in the inhomogeneous norm. Here
$\mathbf X$ is the natural lift of $X$: for $H>\tfrac12$ the canonical (Young) lift; for
$H\le\tfrac12$ the translate by $t\mapsto\theta te_1$ of the natural lift $\mathbf G$ of $G$, the
limit of the lifts of the piecewise-linear interpolants of $G$ (it exists by (C5),
Remark~\ref{rem:classread}).

The covariance $R_G$ is represented by a kernel $\varrho$. For $H>\tfrac12$,
$R_G([s,t]\times[u,v])=\int_s^t\!\int_u^v\varrho$ with $\varrho\in L^1([0,1]^2)$. For $H=\tfrac12$,
$\varrho$ is a finite signed measure on $[0,1]^2$ with $R_G(A\times B)=\varrho(A\times B)$, and we
require that the first marginal of $|\varrho|$ be absolutely continuous, so that $\varrho$ has a
disintegration $\varrho(\dd s,\dd v)=\dd s\,\varrho(s,\dd v)$. In both
cases $\varrho(\dd s,\dd v)=\dd s\,\varrho(s,\dd v)$ is its disintegration and $|\varrho|(s,\cdot)$ the
total variation of $\varrho(s,\cdot)$. For $H<\tfrac12$, $\varrho$ is a locally integrable function
off the diagonal with $R_G(I\times J)=\int_I\!\int_J\varrho$ for intervals $I,J$ with disjoint
interiors. Write $I_a=[t_a,t_{a+1})$, $a=0,\dots,n-1$, for the mesh cells, $\varrho^n$ for the increment-covariance
density of the piecewise-linear interpolant $G^n$, $\cH_1(I_a)$ for the first Wiener chaos of $G$'s
increments on $I_a$, and
\[
   \varrho_{\max}(I_a,I_b):=\sup\Bigl\{\tfrac{|\Cov(\xi,\eta)|}{\norm\xi_{L^2}\norm\eta_{L^2}}:\
   \xi\in\cH_1(I_a)\!\setminus\!\{0\},\ \eta\in\cH_1(I_b)\!\setminus\!\{0\}\Bigr\}
\]
for the canonical correlation between two cells. We require, for every mesh width $\delta=1/n$:
\begin{itemize}[leftmargin=1.4em,topsep=2pt,itemsep=1pt]
\item[\emph{(C1)}] \emph{(interpolation error)} with $e:=\varrho^n-\varrho$,
$\sum_{|a-b|\ge1}\int_{I_a\times I_b}|e|\le L\,\delta^{2H-1}$. For $H\ge\tfrac12$ the same bound with
the sum over all cells, i.e.\ $\norm e_{L^1([0,1]^2)}\le L\delta^{2H-1}$ (total variation at
$H=\tfrac12$);
\item[\emph{(C2)}] \emph{(per-cell moments)} for every cell $I_a$ and every $q\le r$,
\[
   \E\bigl\|\Sig^{\le M}\bigl(\delta^{-H}(G_{t_a+\delta\,\cdot}-G_{t_a})\bigr)\bigr\|_\hs^{q}\ \le\ L ;
\]
\item[\emph{(C3)}] \emph{(kernel sections, Young regime)} for $H\ge\tfrac12$,
$\int_0^1|\varrho|(s,\dd v)\le L$ for almost every $s\in[0,1]$;
\item[\emph{(C4)}] \emph{(aggregated cross-cell bounds, rough regime)} for $H<\tfrac12$,
\[
   \sup_a\ \sum_{b:\,b\ne a}\int_{I_a\times I_b}|\varrho|\ \le\ L\,\delta^{2H}
   \qquad\text{and}\qquad
   \sup_a\ \sum_{b:\,b\ne a}\varrho_{\max}(I_a,I_b)\ \le\ L ;
\]
\item[\emph{(C5)}] \emph{(two-dimensional $\rho$-variation of the covariance)} with
\begin{equation}\label{eq:rhoH}
   \rho\ :=\ 1\vee\tfrac1{2H}\ \in\ [1,2)\qquad(H>\tfrac14),
\end{equation}
one has $V_\rho\bigl(R_G;[0,1]^2\bigr)\le L$, where $V_\rho$ is the two-dimensional $\rho$-variation
over rectangular grids.
\end{itemize}
Write $\cP_0(H,L)\subset\cP(H,L)$ for the centred sub-class $\theta=0$.
\end{definition}

\begin{remark}\label{rem:classread}
\emph{(a) The drift direction.} Taking the drift along $e_1$ is no restriction. The law of $G$ is
invariant under orthogonal maps $O$ of $\RR^d$ (its components are i.i.d.\ centred Gaussian
processes), and such a map acts isometrically on every level of the
signature (as $O^{\otimes m}$ on level $m$) and commutes with interpolation. So a drift $tv$, $v\in\RR^d$, gives the same
strong error as $\theta=|v|$, and a weak error that is the image under $(O^{-1})^{\otimes m}$, level by
level, of the one for $\theta=|v|$, where $Ov=|v|e_1$; in particular its norm is the same.

\emph{(b) The natural lift.} It exists by (C5): grid $\rho$-variation gives a two-dimensional control
at every $\rho'\in(\rho,2)$ \citep[Thm.~1]{FrizVictoir2011}, and then
\citet[Thm.~15.33(iv)]{FrizVictoir2010} applies. For fBm see also \citet{CoutinQian2002}. For
$p\ge2$ the same path has other geometric lifts \citep[Rem.~15.35]{FrizVictoir2010}, so the choice
matters. In this paper every signature is that of the natural lift.

\emph{(c) The kernel.} For Brownian motion $\varrho(\dd s,\dd t)=\dd s\,\delta_{\{s=t\}}$, where
$\delta_{\{s=t\}}$ is the unit point mass at $t=s$, so that $\varrho$ is Lebesgue measure carried by the
diagonal. The
interpolant's density $\varrho^n$ is constant on each cell product $I_a\times I_b$, where it equals
the cell average of $\varrho$. (C1) measures how far $\varrho$ is from these averages.
\end{remark}

\begin{remark}[features of the definition]\label{rem:classwording}
\emph{(a) (C2) is a moment bound, not exact self-similarity.} Requiring one law for the rescaled cell
path at every mesh and every cell would be much stronger: $n=1$ identifies that law with the law of
$G$, and general $n$ then implies $n^HG_{\cdot/n}\overset d=G$, hence $\Var(G_u)=u^{2H}$ and exact
self-similarity. For stationary increments this pins $G$ to $\mathrm{fBm}(H)$. It would exclude the
tied-down perturbation, and also natural candidates such as fractional Ornstein--Uhlenbeck
drivers.

\emph{(b) (C4) is aggregated over $b$.} Lag-wise bounds
$\int_{I_a\times I_b}|\varrho|\le\delta^{2H}\bar r(|a-b|)$ or
$\varrho_{\max}(I_a,I_b)\le\bar r(|a-b|)$, with one summable sequence $\bar r$ for all meshes, fail for
the tied-down perturbation. At a cell next to a block node the ripple kernel contributes about
$\epsilon^2\delta^{2H}\,\delta/\delta_0$ to $\int_{I_a\times I_b}|\varrho|$ at every lag
$|a-b|\le\delta_0/\delta$, so such an $\bar r$ would have to dominate $m\mapsto1/m$. The sums over $b$
stay bounded (Remark~\ref{rem:C1C2}), and they are all that the proofs use
(Proposition~\ref{prop:rough}, Lemmas~\ref{lem:ibp} and~\ref{lem:Vbound}).
\end{remark}

\begin{remark}[(C1) below the crossover]\label{rem:C1scope}
For $H<\tfrac12$, (C1) is restricted to off-diagonal cells because the unrestricted condition fails
whenever the cell-increment variances are of order $\delta^{2H}$, as for fBm. Indeed, if
$\varrho\in L^1([0,1]^2)$, then splitting $\Var(G_1)=1$ over cell pairs gives
$\sum_a\Var\bigl(G_{t_{a+1}}-G_{t_a}\bigr)=1-\sum_{a\ne b}\int_{I_a\times I_b}\varrho\to1-\int_{[0,1]^2}\varrho$
by dominated convergence, a finite limit, while the left side is of order $n\delta^{2H}=\delta^{2H-1}\to\infty$.
So $\norm\varrho_{L^1}=\infty$, and since $\varrho^n$ is a bounded step function, $\varrho^n-\varrho$
has infinite $L^1$ norm at every mesh. At $H=\tfrac12$, where $n\delta^{2H}=1$, the measure-valued
form of (C1) is used.

For the fBm kernel $\varrho=c_H|t-s|^{2H-2}$ the off-diagonal condition holds, with the same exponent, on all of $(\tfrac14,1)$. On an off-diagonal cell product
$\varrho^n$ is the cell average of $\varrho$, so $\int_{I_a\times I_b}|e|$ is $\delta^2$ times the mean
absolute deviation of $\varrho$ from its average there, which is at most half the oscillation of
$\varrho$ over the cell product (a function with values in an interval of length $\ell$ has mean
absolute deviation at most its standard deviation, hence at most $\ell/2$). For $m=|a-b|\ge2$ the cell product lies at distance
$\ge(m-1)\delta$ from the diagonal, where $|\partial_s\varrho|+|\partial_t\varrho|
=2|c_H|(2-2H)|t-s|^{2H-3}$. Since the cell product is convex with side $\delta$, the mean value theorem
bounds the oscillation by $\delta\cdot2|c_H|(2-2H)((m-1)\delta)^{2H-3}$, and therefore
\begin{equation}\label{eq:farcell}
   \int_{I_a\times I_b}|e|\ \le\ |c_H|(2-2H)\,\delta^{2H}\,(m-1)^{2H-3}\qquad(m=|a-b|\ge2).
\end{equation}
There are $2(n-m)\le2n$ ordered pairs $(a,b)$ at lag $m$, so with $n=1/\delta$ these cells contribute at
most $2|c_H|(2-2H)\,\delta^{2H-1}\sum_{m\ge2}(m-1)^{2H-3}=O(\delta^{2H-1})$, the lag sum $\zeta(3-2H)$ being finite for
every $H<1$. For $m=1$, the substitution $s=t_a+\delta u$, $t=t_a+\delta v$ gives $\int_{I_a\times I_b}|\varrho|=\delta^{2H}|c_H|\int_{[0,1]\times[1,2]}|u-v|^{2H-2}<\infty$, the only
singularity being a corner with $2H-2>-2$, and $\int_{I_a\times I_b}|e|\le2\int_{I_a\times
I_b}|\varrho|$ because the cell average satisfies $\int_{I_a\times I_b}|\varrho^n|\le\int_{I_a\times I_b}|\varrho|$.
The $2(n-1)$ such cells again give $O(\delta^{2H-1})$. Only
the diagonal band $m=0$ diverges. The proof of Proposition~\ref{prop:rough} therefore never treats it
through the kernel: it controls the diagonal cells through the  moments (C2) of the assembled cell signatures. Below the crossover the off-diagonal branch of (C1) is not used by any proof. We keep it so
that (C1) has one form on all of $(\tfrac14,1)$, at the price of a slightly smaller $\cP$.
\end{remark}

\emph{Which laws are in the class.} Since (C1)--(C5) are hypotheses, membership must be checked
law by law, at every mesh $\delta$. Remark~\ref{rem:C1C2} in Appendix~\ref{app:verif} does this, together with the moment condition, for three
laws: the independent-fBm law $P_0$, its tied-down perturbation $P_1^{(\delta_0)}$ at every block width
$\delta_0$ and ripple size $\epsilon$ (Construction~\ref{con:bridge} in \S\ref{sec:mesh}), and the drifted law $P_\theta$, the law of $\theta te_1$
plus a path of $P_0$, for $L$ large enough depending on $\theta$. Fractional Ornstein--Uhlenbeck drivers are
natural candidates, but we do not check the conditions for them.

\begin{remark}[the tied-down perturbation depends on the mesh]\label{rem:meshtied}
$P_1^{(\delta_0)}$ is built on a grid of blocks of width $\delta_0$ (Construction~\ref{con:bridge}), so it is a different law for each $\delta_0$, and the
construction at mesh $\delta$ uses $P_1^{(\delta)}$. This is unavoidable. A fixed $P_1\ne P_0$ whose
mesh data have the same law as those of $P_0$ at every mesh would agree with $P_0$ on a dense set of
times, hence everywhere by path continuity. The class
itself does not depend on the mesh: it contains $P_1^{(\delta_0)}$ for every $\delta_0$, with
(C1)--(C5) verified at every mesh with one $L$ (Remark~\ref{rem:C1C2}).
\end{remark}

\begin{remark}[the role of Gaussianity]\label{rem:whygauss}
Centring and mutual independence of the components already make the bias vanish at levels $2$ and $3$,
and they make the level-$4$ bias a functional of only the component covariances
(Proposition~\ref{prop:nongauss}). For example, independent Rosenblatt processes of Hurst index
$H\in(\tfrac12,1)$ have exactly the level-$\le4$ bias of independent fBm$(H)$. A third cumulant changes
the expected signature, $\mu_{iii}=\E(X^i_1)^3/6$, but it does not create a bias at level $3$.
Gaussianity is used elsewhere. Non-Gaussian features first enter
the bias at level $5$, through third-order moment functions of single components (words of
multiplicities $(3,2)$), and at level $6$ through fourth cumulants and products of third moments.
Their size for asymmetric laws is not controlled here, while laws with symmetric components have no
odd-level bias at all. The moment bounds, and through them the strong-error estimates of
\S\ref{sec:strong}, use Gaussian tools (Fernique, hypercontractivity, Lemma~\ref{lem:rhovarcontract}).

Independence cannot be weakened to uncorrelatedness: $Y=(B,B^2-t)$, $B$ Brownian, has centred
uncorrelated components, and
$\E[\Delta B_a\Delta B_b\Delta Y^2_c]=2\delta^2\one_{\{\max(a,b)=c\}}$ gives
$\cW_3=\delta(-\tfrac16,-\tfrac16,\tfrac13)$ on $112,121,211$.
\end{remark}

\subsection{Moment and Lipschitz lemmas}\label{sec:R2}

\paragraph{Signature bound.} For a geometric $p$-rough path $\mathbf x$ over $[0,1]$ and fixed
$M$, the extension $S_M(\mathbf x)$ of $\mathbf x$ to level $M$ (the Lyons extension, whose
increment over $[0,1]$ is $\Sig^{\le M}(\mathbf x)$) satisfies $\hnorm{S_M(\mathbf x)}_{\pvar}\le C_{M,p}\hnorm{\mathbf x}_{\pvar}$
\citep[Thm.~9.5]{FrizVictoir2010}. A group element of homogeneous norm at most $K$ has level-$m$
component of norm at most $C_mK^m$, and $\hnorm{\mathbf x}_{\pvar}\le C(1\vee\norm{\mathbf x}_{\pvar})$.
Hence
\begin{equation}\label{eq:sigbound}
   \norm{\Sig^{\le M}(\mathbf x)}_\hs\ \le\ c_{M,p}\bigl(1\vee\norm{\mathbf x}_{\pvar}\bigr)^{M}.
\end{equation}
The next lemma turns \eqref{eq:sigbound} into moment bounds for the signature of the driver. It
is used wherever an expectation of a signature coordinate must be finite, in particular in
Lemma~\ref{lem:sumzero} and in the proof of Proposition~\ref{prop:rough}.

  \begin{lemma}[$L^q$ moments]\label{lem:moments}
Let $p\ge1$, $q>0$, and let $\mathbf X$ be a geometric $p$-rough path. If
$\E\norm{\mathbf X}_{\pvar}^{Mq}<\infty$, then
$\E\norm{\Sblk}_\hs^q\le c_{M,p}^q(1+\E\norm{\mathbf X}_{\pvar}^{Mq})<\infty$. In particular
$\E\norm{\Sblk}_\hs^q<\infty$ for every $q$ (i) for a centred Gaussian driver with independent
components and $V_\rho(R)<\infty$, $\rho\in[1,2)$, (and its translates by $\theta te_1$) at every $p\in(2\rho,4)$, and (ii) for
$P\in\cP(H,L)$, $H>\tfrac12$, at $p=p_\delta(H)$, by Lemma~\ref{lem:youngmom}.
\end{lemma}
\begin{proof}
By \eqref{eq:sigbound}, $\norm{\Sblk}_\hs\le c_{M,p}\,(1\vee\norm{\mathbf X}_{\pvar})^{M}$
pathwise. Since $(1\vee r)^{Mq}\le1+r^{Mq}$,
\[
  \E\norm{\Sblk}_\hs^q\ \le\ c_{M,p}^q\,\E\bigl[(1\vee\norm{\mathbf X}_{\pvar})^{Mq}\bigr]
  \ \le\ c_{M,p}^q\bigl(1+\E\norm{\mathbf X}_{\pvar}^{Mq}\bigr).
\]
Now let the driver be centred Gaussian with independent components and covariance of finite
$\rho$-variation, $\rho\in[1,2)$ (for fBm, $\rho=1\vee\frac{1}{2H}$ as in \eqref{eq:rhoH}, $H>\tfrac14$;
see \citet[Prop.~15.5]{FrizVictoir2010}), and let $p\in(2\rho,4)$. By Lemma~\ref{lem:gaussmom} below the
homogeneous norm $\hnorm{\mathbf X}_{\pvar}$ of the natural lift has moments of every order. The
inhomogeneous norm is at most $C(1\vee\hnorm{\mathbf X}_{\pvar})^{\lfloor p\rfloor}$ by the ball--box
estimate, so it has moments of every order too. For a translate by $h=\theta te_1$, the
translation estimate $\hnorm{T_h(\mathbf x)}_{\pvar}\le C_1(\hnorm{\mathbf x}_{\pvar}+|\theta|)$
\citep[Thm.~9.33(i)]{FrizVictoir2010}, with $|h|_{\onevar}=|\theta|$, preserves all these moments. A $p<2$ arises only
at the exponent $p_\delta(H)$ of the Young range (case (ii)). There $\lfloor p\rfloor=1$, the covariance hypothesis
alone does not make $\norm{\mathbf X}_{\pvar}$ finite, and the moments over the class are those of
Lemma~\ref{lem:youngmom}. Hence $\E\norm{\Sblk}_\hs^q<\infty$ for every $q$. Class-uniform bounds are given in Lemma~\ref{lem:rhovarcontract}.
\end{proof}

The strong error is a difference of two signatures, so it needs a local Lipschitz bound for
$\Sig^{\le M}$ in the metric $\rho_{\pvar}$, with a constant polynomial in the norms of the two
paths. For $p<2$ this is the local Lipschitz continuity of the truncated signature on the Young
range, and for $p\ge2$ that of the Lyons extension \citep{Lyons1998,LyonsCaruanaLevy2007}. See
\citet[Thm.~9.10]{FrizVictoir2010} for the continuity of the Lyons lift $S_N$ in $\rho_{p,\omega}$.
Lemma~\ref{lem:liprough} below gives the form used in this paper.

The moments of the interpolant's lift must be bounded uniformly in the mesh. The next lemma
reduces this to the driver itself: interpolation does not increase $p$-variation. It is used in the
proofs of Lemmas~\ref{lem:rhovarcontract} and~\ref{lem:L1conv}.

\begin{lemma}[$p$-variation of a piecewise-linear interpolant is attained on its nodes]\label{lem:pvarcontract}
Let $x\in C([0,1];\RR^d)$, let $x^n$ be its piecewise-linear interpolation on any finite node set
$0=t_0<\cdots<t_n=1$, and let $p\ge1$. Then the supremum defining $\norm{x^n}_{\pvar}$ is attained
on
a partition contained in $\{t_0,\dots,t_n\}$, and consequently
\[
   \norm{x^n}_{\pvar}\ \le\ \norm{x}_{\pvar} .
\]
The same holds for $x-y$ interpolated on the same nodes, $(x-y)^n=x^n-y^n$.
\end{lemma}
\begin{proof}
The argument is by convexity: we show that moving a partition point that is not a node to a node,
or deleting it, never decreases the $p$-variation sum, so the supremum may be taken over partitions of
the node set, on which $x^n=x$. Let $\Pi$ be a partition and $s\in\Pi$ a non-node with$\Pi$-neighbours $u<s<v$, lying in the cell
$[t_j,t_{j+1}]$. On $J:=[\max(u,t_j),\min(v,t_{j+1})]$ the map $\sigma\mapsto x^n_\sigma$ is affine,
so $\sigma\mapsto|x^n_\sigma-x^n_u|^p+|x^n_v-x^n_\sigma|^p$ is a sum of increasing convex
functions of norms of affine functions, hence convex on $J$ and maximised at an endpoint of $J$.
Each endpoint is a node or one of $u,v$. In the latter case the corresponding term vanishes and $s$
may be deleted, giving a coarser partition without a smaller sum. Repeating this on the remaining non-nodes
takes at most $|\Pi|$ steps and never decreases the sum. Hence
$\norm{x^n}^p_{\pvar}=\sup\{\sum_i|x^n_{s_{i+1}}-x^n_{s_i}|^p:\ (s_i)\subseteq\{t_j\}\}$, and since
$x^n$ agrees with $x$ at the nodes this is the $p$-variation of $x$ restricted to the node set, at
most $\norm x_{\pvar}$.
\end{proof}

The next lemma is the two-dimensional analogue of Lemma~\ref{lem:pvarcontract} for covariances.
It gives lift moments bounded uniformly in the mesh, used in Lemmas~\ref{lem:L1conv}
and~\ref{lem:ibp} and Proposition~\ref{prop:attain5}.

  \begin{lemma}[the two-dimensional contraction, and uniform-in-mesh moments]\label{lem:rhovarcontract}
Let $Z$ be a continuous square-integrable process on $[0,1]$ with increment covariance
$R\bigl([s,t]\times[u,v]\bigr)=\E[Z_{s,t}Z_{u,v}]$, let $Z^n$ be its piecewise-linear interpolation
on a node set $T$, with increment covariance $R^n$, and write $V_\rho(R;[0,1]^2)$ for the two-dimensional $\rho$-variation of $R$ over
rectangular grids. Then for every $\rho\ge1$ the supremum defining $V_\rho(R^n)$ is attained on
grids contained in $T\times T$, whence
\[
   V_\rho\bigl(R^{n};[0,1]^2\bigr)\ \le\ V_\rho\bigl(R;[0,1]^2\bigr) ,
\]
and the same holds for the joint covariance of any finite family all of whose members are
interpolated on the same nodes, in particular for $Z-Z'$ and for the pair $(Z^n,Z'^n)$.
Consequently, for $P\in\cP(H,L)$, with $\rho=\rho(H)$ of \eqref{eq:rhoH}, every $p\in(2\rho,4)$ and every $r<\infty$,
\[
   \sup_{n}\ \E\bigl\|\mathbf X^n\bigr\|_{\pvar}^{\,r}\ \le\ C(r,p,\rho,d,L)\ <\ \infty ,
\]
$\mathbf X^n$ being the canonical lift of the interpolant $X^n$. The constant depends on $P$ only
through the class parameters, so the bound is uniform over $\cP(H,L)$.
\end{lemma}

\begin{proof}[Proof of Lemma~\ref{lem:rhovarcontract}]
The first claim follows from the convexity argument of Lemma~\ref{lem:pvarcontract}, applied in each
variable separately. The second combines it with the Gaussian moment bound of Lemma~\ref{lem:gaussmom}
and a translation estimate for the drift.

\emph{Step 1 (grids in $T\times T$ suffice).} For grids $(s_i)$, $(u_j)$ the
$\rho$-variation sum is
\[
   \textstyle\sum_{i,j}\bigl|R^n\bigl([s_i,s_{i+1}]\times[u_j,u_{j+1}]\bigr)\bigr|^\rho .
\]
  Fix the second-variable
grid and let $s$ be a first-variable non-node with neighbours $u<s<v$, in the cell $[t_a,t_{a+1}]$.
On $J:=[\max(u,t_a),\min(v,t_{a+1})]$ the maps $\sigma\mapsto Z^n_{u,\sigma}$ and
$\sigma\mapsto Z^n_{\sigma,v}$ are affine in $\sigma$, with random coefficients, so
$\sigma\mapsto R^n([u,\sigma]\times B)$ and $\sigma\mapsto R^n([\sigma,v]\times B)$ are affine for
every second-variable block $B$. The $\sigma$-dependent part of the double sum is then a sum of
functions $|\text{affine}|^\rho$, convex on $J$ (for a matrix-valued family, a norm of an affine map
composed with $x\mapsto x^\rho$, convex for the same reason), hence maximised at an endpoint of $J$:
a node, or one of $u,v$, where the corresponding terms vanish and $s$ may be deleted for a coarser
grid with no smaller sum. Iterating in each variable terminates. Since $Z^n=Z$ on $T$, $R^n=R$ on
rectangles with corners in $T\times T$, so the supremum over such grids is at most $V_\rho(R;[0,1]^2)$.

\emph{Step 2 (the Gaussian part).} For the second claim, write $X=\theta te_1+G$ as in Definition~\ref{def:class}. The linear drift is piecewise
linear, hence interpolated exactly, and $G^n$ is centred Gaussian with independent components.
Lemma~\ref{lem:gaussmom} below, applied to $G$ with $K=L$ (by (C5)), gives
$\bignormL{\hnorm{\mathbf G^n}_{\pvar}}{L^r}\le C(\rho,p,L,d)\sqrt r$ for $r\ge r_0$, uniformly in $n$ (for
$r<r_0$ use $\norm\cdot_{L^r}\le\norm\cdot_{L^{r_0}}$),
$\mathbf G^n$ the canonical lift of $G^n$.

\emph{Step 3 (adding the drift).} The linear drift is then added by translation. Both the drift and $G^n$ have bounded variation, so the
translation of $\mathbf G^n$ by $h=\theta te_1$ coincides with the canonical lift of the sum,
$T_h(S_{\lfloor p\rfloor}(G^n))=S_{\lfloor p\rfloor}(\theta te_1+G^n)$
\citep[Rem.~9.34]{FrizVictoir2010}. The translation estimate
\citep[\S9.4.6, Thm.~9.33(i)]{FrizVictoir2010} gives
$\hnorm{T_h(\mathbf x)}_{\pvar}\le C_1(\hnorm{\mathbf x}_{\pvar}+|h|_{q\text{-}\mathrm{var}})$, with
$q=1$ here, so the complementary Young condition $\tfrac1p+\tfrac1q>1$ holds for every $p$.
Converting to the inhomogeneous convention \eqref{eq:rhopvar} by the ball--box comparison
\citep[Prop.~8.9, from the ball--box estimate Prop.~7.49; (8.6) for the $p$-variation
form]{FrizVictoir2010},
\[
   \norm{\mathbf X^n}_{\pvar}\ \le\ C(p,L)\bigl(1+\hnorm{\mathbf G^n}_{\pvar}\bigr)^{\lfloor p\rfloor},
\]
so every moment transfers. The power $\lfloor p\rfloor$ comes from the conversion between the two
norms.
\end{proof}

The next lemma packages the Gaussian moment bound used above, for a general Gaussian process,
uniformly over its interpolants. Besides Lemmas~\ref{lem:moments} and~\ref{lem:rhovarcontract}, it
supplies the lift moments in Proposition~\ref{prop:rough}, Lemma~\ref{lem:subint},
Proposition~\ref{prop:strongunif} and Appendix~\ref{app:verif}.

\begin{lemma}[uniform moments of Gaussian lifts]\label{lem:gaussmom}
Let $Z$ be a centred continuous Gaussian process on $[0,1]$ with values in $\RR^d$ and independent
components, whose covariances satisfy $V_\rho(R_{Z^i};[0,1]^2)\le K$, $i\le d$, for some $\rho\in[1,2)$.
For a finite node set $T\ni0,1$ let $Z^T$ be the piecewise-linear interpolation of $Z$ on $T$, and let
$\mathbf Z$ be the natural lift of $Z$. Fix $p\in(2\rho,4)$. There are $r_0=r_0(\rho,p)$ and
$C=C(\rho,p,K,d)$ such that for every $r\ge r_0$ and every node set $T$,
\[
   \bignormL{\hnorm{S_3(Z^T)}_{\pvar;[0,1]}}{L^r}\ \le\ C\sqrt r,\qquad
   \bignormL{\hnorm{\mathbf Z}_{\pvar;[0,1]}}{L^r}\ \le\ C\sqrt r .
\]
The same holds, with $d$ replaced by $kd$, for a finite family $Z_1,\dots,Z_k$ of such processes, all
interpolated on the same nodes, whenever the $\RR^{kd}$-valued process $(Z_1,\dots,Z_k)$ is jointly
Gaussian with mutually independent scalar components: the bounds hold for the lift of
$(Z_1^T,\dots,Z_k^T)$ and for the natural lift of $(Z_1,\dots,Z_k)$.
\end{lemma}

\begin{proof}
The bound is the Gaussian moment estimate for natural lifts \citep[Thm.~15.33(ii)]{FrizVictoir2010}
combined with a Kolmogorov--Garsia-type criterion. The only point is that the constants depend on the
covariance through $K$ alone, which the contraction of Lemma~\ref{lem:rhovarcontract} makes uniform in
$T$.

\emph{Step 1 (contraction).} By the first claim of Lemma~\ref{lem:rhovarcontract},
$V_\rho(R_{Z^{T,i}})\le V_\rho(R_{Z^i})\le K$ for every component and every node set $T$.

\emph{Step 2 (a two-dimensional control).} Put $\rho':=\tfrac12(\rho+\tfrac p2)\in(\rho,\tfrac p2)$, so
that $\rho'<2$ and $2\rho'<p$. Grid $\rho$-variation gives controlled $\rho'$-variation with a constant
$c=c(\rho,\rho')$, and the $\rho'$-th power of the controlled $\rho'$-variation over rectangles is a
two-dimensional control \citep[Thm.~1]{FrizVictoir2011}. Let $\omega_T$ be the sum over the $d$
components of these controls for $Z^T$. It dominates the $\rho'$-variation of the covariance of $Z^T$
on every rectangle, and $\omega_T([0,1]^2)\le d(cK)^{\rho'}=:K'$, uniformly in $T$.

\emph{Step 3 (the moment bound).} $Z^T$ is a centred continuous Gaussian process with independent
components and bounded-variation paths, so its natural lift is the canonical lift $S_3(Z^T)$
\citep[Thm.~15.33(iv), Rem.~15.35]{FrizVictoir2010}, and \citet[Thm.~15.33(ii)]{FrizVictoir2010}, at
$\rho'$, gives for all $s<t$ and $q\ge1$
\[
   \bignormL{d\bigl(S_3(Z^T)_s,S_3(Z^T)_t\bigr)}{L^q}\ \le\ C_0(\rho')\sqrt q\;\omega_T([s,t]^2)^{1/(2\rho')} .
\]
The map $(s,t)\mapsto\omega_T([s,t]^2)$ is a one-dimensional control with total mass at most $K'$. As in
the proof of \citet[Prop.~15.22(ii)]{FrizVictoir2010}, reparametrise time by
$\tau(t):=\omega_T([0,t]^2)/\omega_T([0,1]^2)$ (if $\omega_T([0,1]^2)=0$ there is nothing to prove).
The reparametrised process $Y$, $Y_{\tau(t)}:=S_3(Z^T)_t$, satisfies
$\norm{d(Y_\sigma,Y_{\sigma'})}_{L^q}\le C_0\sqrt q\,K'^{1/(2\rho')}|\sigma'-\sigma|^{1/(2\rho')}$ for every $q$.
\citet[Thm.~A.12]{FrizVictoir2010}, with the parameter $r$ of that theorem equal to $2\rho'$ and H\"older exponent $1/p<1/(2\rho')$, gives
$q_0=q_0(\rho',p)$ and $C_1=C_1(\rho',p)$, not depending on $q$, such that
$\bignormL{\hnorm Y_{1/p\text{-H\"ol};[0,1]}}{L^q}\le C_1C_0K'^{1/(2\rho')}\sqrt q$ for $q\ge q_0$. As
$p$-variation is invariant under reparametrisation and dominated by the $1/p$-H\"older norm on $[0,1]$,
the first bound follows with $r_0=q_0$ and $C$ independent of $T$.

\emph{Step 4 (the natural lift and families).} The covariance of $Z$ itself has the same bounds, so
Steps~2--3 apply verbatim to $\mathbf Z$ \citep[Thm.~15.33(ii)]{FrizVictoir2010}. A family as in the
statement is a single process in $\RR^{kd}$ with independent components, each with grid
$\rho$-variation at most $K$, and its interpolation on common nodes is the interpolation of that
process.
\end{proof}

On the Young range the moment condition of Definition~\ref{def:class} is imposed at the single
order $2M$ and at an exponent $p<2$, where Lemma~\ref{lem:gaussmom} does not apply. The next lemma
upgrades it to moments of every order. It is used in Lemma~\ref{lem:moments}(ii) and in the
dominated-convergence step of Lemma~\ref{lem:L1conv}.

\begin{lemma}[all moments on the strict Young range]\label{lem:youngmom}
Let $H\in(\tfrac12,1)$ and $P\in\cP(H,L)$, and let $p=p_\delta(H)\in(1/H,2)$ be the exponent
\eqref{eq:pdelta} of the moment condition in Definition~\ref{def:class}. There is an absolute constant $c$ such that for every $q\ge2$
\[
   \bignormL{\norm{\mathbf X}_{\pvar}}{L^q}\ \le\ C_L\sqrt q,\qquad
   C_L:=(1+c)\bigl((1\vee L)^{1/(2M)}+L\bigr)+L .
\]
\end{lemma}

\begin{proof}
The $p$-variation of the Gaussian part is the supremum of a countable family of centred Gaussian
variables, so the Borell--TIS inequality gives it Gaussian tails, with a scale controlled by its second
moment. The moment condition bounds that second moment.

\emph{Step 1 (second moment of the Gaussian part).} As $\lfloor p\rfloor=1$, $\norm{\mathbf X}_{\pvar}$ is the $p$-variation of the path, and
$\norm{\theta te_1}_{\pvar}=|\theta|\le L$. So $G=X-\theta te_1$ satisfies
\[
   \bignormL{\norm G_{\pvar}}{L^2}\ \le\ \bignormL{\norm{\mathbf X}_{\pvar}}{L^{2M}}+L
   \ \le\ L^{1/(2M)}+L\ \le\ (1\vee L)^{1/(2M)}+L,
\]
by the moment condition of Definition~\ref{def:class} and Jensen's inequality ($2M\ge2$).

\emph{Step 2 (Gaussian concentration).}
 $G$ is centred Gaussian with continuous paths, and $\Lambda:=\norm G_{\pvar}=\sup_{f\in F}f(G)$ for a
countable symmetric family $F$ of continuous linear functionals (rational partitions, and rational
dual sequences in the unit ball of $\ell^{p/(p-1)}(\RR^d)$), finite almost surely. By the
Borell--TIS inequality \citep[Thm.~2.1.1]{AdlerTaylor2007}, $\E\Lambda<\infty$ and
$\Prob(\Lambda-\E\Lambda>u)\le e^{-u^2/(2\sigma^2)}$ with $\sigma^2=\sup_{f\in F}\E f(G)^2\le\E\Lambda^2$ , $F$ being
symmetric (so $|f(G)|\le\Lambda$). Since $0\le\Lambda\le\E\Lambda+(\Lambda-\E\Lambda)_+$ and a variable with
such a Gaussian tail has $\norm{(\Lambda-\E\Lambda)_+}_{L^q}\le c\sqrt q\,\sigma$ , for $q\ge2$
\[
   \norm\Lambda_{L^q}\ \le\ \E\Lambda+c\sqrt q\,\sigma\ \le\ (1+c)\sqrt q\,\norm\Lambda_{L^2}.
\]

\emph{Step 3 (adding the drift back).} Finally $\bignormL{\norm{\mathbf X}_{\pvar}}{L^q}\le \norm\Lambda_{L^q} +|\theta|\le
(1+c)\sqrt q\bigl((1\vee L)^{1/(2M)}+L\bigr)+L\le C_L\sqrt q$, as $\sqrt q\ge1$.
\end{proof}

The next estimate makes the polynomial dependence on the norms explicit for rough paths. It is
used for the class-uniform strong rate (Proposition~\ref{prop:strongunif}). The homogeneous norm
$\hnorm\cdot_{\pvar}$ and the dilation $\mathrm{dil}_\lambda$ are those of \S\ref{sec:setup}.

\begin{lemma}[Lipschitz estimate with dilation]\label{lem:liprough}
Let $p\ge1$ and $M\ge1$. There is $C=C(p,M)$ such that for all geometric $p$-rough paths
$\mathbf x,\mathbf y$ over $[0,1]$,
\[
   \norm{\Sig^{\le M}(\mathbf x)-\Sig^{\le M}(\mathbf y)}_\hs\ \le\ C\,N^{M-1}\rho_{\pvar}(\mathbf x,\mathbf y),
   \qquad N:=1\vee\hnorm{\mathbf x}_{\pvar}\vee\hnorm{\mathbf y}_{\pvar} .
\]
\end{lemma}

\begin{proof}
We first prove the bound with $N=1$, for paths in the unit ball of the homogeneous norm, from the
local Lipschitz continuity of the Lyons extension. Dilation by $\lambda=N^{-1}$ then moves any pair into
that ball, and undoing it costs the factor $N^{m-1}$ at level $m$.

\emph{Step 1 (the unit ball).} The Lyons extension is Lipschitz in $\rho_{\pvar}$ on the unit ball of the homogeneous
norm, with a constant $C=C(p,M)$. Indeed, for $\hnorm{\mathbf x}_{\pvar},\hnorm{\mathbf y}_{\pvar}\le1$ and $\varepsilon:=\rho_{\pvar}(\mathbf x,\mathbf y)>0$,
the control $\omega(s,t):=\hnorm{\mathbf x}^p_{\pvar;[s,t]}+\hnorm{\mathbf y}^p_{\pvar;[s,t]}+\sum_{i\le\lfloor p\rfloor}
\bigl(\norm{\mathbf x^{(i)}-\mathbf y^{(i)}}_{(p/i)\text{-var};[s,t]}/\varepsilon\bigr)^{p/i}$ has $\omega(0,1)\le2+\lfloor p\rfloor$ and
$\rho_{p,\omega}(\mathbf x,\mathbf y)\le\varepsilon$, so \citet[Thm.~9.10]{FrizVictoir2010} bounds the level-$m$ difference by
$C\varepsilon\,\omega(0,1)^{m/p}\le C\varepsilon$:
\[
   \bigl\lVert\pi_m\bigl(\Sig(\mathbf x)-\Sig(\mathbf y)\bigr)\bigr\rVert_\hs\ \le\ C\,\rho_{\pvar}(\mathbf x,\mathbf y)
   \qquad(m\le M,\ \hnorm{\mathbf x}_{\pvar},\hnorm{\mathbf y}_{\pvar}\le1).
\]
On $\{\varepsilon=0\}$ there is nothing to prove.

\emph{Step 2 (dilation).} For $\lambda:=N^{-1}\le1$ the
dilated paths lie in that ball, $\rho_{\pvar}(\mathrm{dil}_\lambda\mathbf x,\mathrm{dil}_\lambda\mathbf y)\le\lambda\,\rho_{\pvar}(\mathbf x,\mathbf y)$ (level $i$ scales by
$\lambda^i\le\lambda$), and level $m$ of the signature scales by $\lambda^m$. So the level-$m$ difference is
$\le C\lambda^{1-m}\rho_{\pvar}(\mathbf x,\mathbf y)=C\,N^{m-1}\rho_{\pvar}(\mathbf x,\mathbf y)$, and,
summing the squares over $m\le M$ and using $N\ge1$, $\norm{\Sig^{\le M}(\mathbf x)-\Sig^{\le M}(\mathbf y)}_\hs\le C\,N^{M-1}\rho_{\pvar}(\mathbf x,\mathbf y)$, with
$C$ enlarged by the factor $\sqrt M$.
\end{proof}

The proofs of \S\S\ref{sec:bias}--\ref{sec:strong} compute the expected signature of the natural
lift as the limit of the expected signatures of finer interpolants. The next lemma justifies this.

 \begin{lemma}[convergence of the interpolants' signatures]\label{lem:L1conv}
Let $P\in\cP(H,L)$, and let $n'\to\infty$ along any sequence of uniform meshes. Then
$\Sig^{\le M}(X^{n'})\to\Sig^{\le M}(X)$ in $L^q$ for every $q<\infty$. In particular, for every mesh
$\delta=1/n$,
\[
   \cW_{\le M}\ =\ \lim_{k\to\infty}\Bigl(\E\bigl[\Sig^{\le M}(X^n)\bigr]-\E\bigl[\Sig^{\le M}(X^{kn})\bigr]\Bigr).
\]
\end{lemma}

\begin{proof}
We prove almost sure convergence (Young range) or convergence in probability (rough range), and
upgrade it to $L^q$ by uniform integrability, which comes from moment bounds uniform in the mesh.

\emph{Case $H>\tfrac12$.} Let $p=p_\delta(H)\in(1/H,2)$ and $p'\in(p,2)$. By Lemma~\ref{lem:pvarcontract},
$\norm{X^{n'}}_{\pvar}\le\norm X_{\pvar}$, and $X^{n'}\to X$ uniformly, hence in $p'$-variation by
interpolation, $\norm{y}_{p'\text{-var}}\le\norm y_{\pvar}^{p/p'}|y|_0^{1-p/p'}$ with
$|y|_0:=\sup_{u,v}|y_u-y_v|$
\citep[Prop.~5.5(i)]{FrizVictoir2010}, applied to $y=X^{n'}-X$. The signature is continuous in
$p'$-variation on the Young range \citep[Cor.~9.11(i)]{FrizVictoir2010}, so the convergence holds almost
surely, and it is dominated by $c_{M,p}(1\vee\norm X_{\pvar})^M$ (by \eqref{eq:sigbound} and
$\norm{X^{n'}}_{\pvar}\le\norm X_{\pvar})$, which lies in every $L^q$ by
Lemma~\ref{lem:youngmom}. Dominated convergence gives the $L^q$ convergence.

\emph{Case $H\le\tfrac12$.} Let $p=p_\delta(H)\in(2\rho,4)$ (with $\rho=\frac{1}{2H}$ of
\eqref{eq:rhoH}, $2\rho=\frac{1}{H}<p_\delta(H)$). The lifts $S_3(G^{n'})$
converge in probability, in homogeneous $p$-variation, to the natural lift $\mathbf G$
\citep[Thm.~15.33(iv)]{FrizVictoir2010}, after passing to a two-dimensional control at some
$\rho'\in(\rho,p/2)$ \citep[Thm.~1]{FrizVictoir2011}. The translation by the drift
\citep[Cor.~9.35]{FrizVictoir2010} and the Lyons extension to level $M$
\citep[Cor.~9.11(i)]{FrizVictoir2010} are continuous, so
$\Sig^{\le M}(X^{n'})\to\Sig^{\le M}(X)$ in probability. By the signature bound \eqref{eq:sigbound} and
Lemma~\ref{lem:rhovarcontract} the family $\norm{\Sig^{\le M}(X^{n'})}_\hs$ is bounded in every $L^q$,
hence its $q$-th powers are uniformly integrable for every $q$, and the convergence holds in $L^q$ (Vitali).

\emph{The claim.} Convergence in $L^1$ gives $\E[\Sig^{\le M}(X^{kn})]\to\E[\Sig^{\le M}(X)]$
as $k\to\infty$. Subtracting from $\E[\Sig^{\le M}(X^n)]$ gives the claim.
\end{proof}

\section{The Gaussian Wick sub-model}\label{sec:R3}

Let $X$ be a continuous centred $\RR^d$-valued Gaussian process with its natural lift. By Wick's
(Isserlis') theorem \citep{Janson1997}, $\mu_w=\E\langle w,\Sblk\rangle=0$ for $|w|$ odd,
while for $|w|=2\ell$ it sums over the $(2\ell-1)!!$ pairings of the $2\ell$ integration variables
the corresponding ordered-simplex integral of the increment-covariance measure. For a piecewise-linear
Gaussian path, whose increment covariance is a bounded density, this is Isserlis' theorem for the
finite Gaussian vector of cell increments. For the natural lift of a member of the class we use it
through Lemma~\ref{lem:L1conv}, as the limit of the same expressions for the interpolants. For
$H>\tfrac12$ the limit can be taken inside: by (C1), $\varrho^{n'}\to\varrho$ in $L^1$, and the simplex
integrals of products of kernels on disjoint pairs of variables are continuous in $L^1$, so the formula
holds with $\varrho$ itself.
For fBm, \citet{BaudoinCoutin2007} computed the expected signature for $H>\tfrac12$, and
\citet{CassFerrucci2024} gave Wiener-chaos expansions of the signature for $H\in(\tfrac14,1)$.

The next lemma shows that the level-$2$ diagonal coordinates are fixed by the increments. The
proof of Proposition~\ref{prop:level4} uses it at level~$2$.

\begin{lemma}[level-2 diagonal]\label{lem:diag2}
Let $M\ge2$, let $X$ with $X_0=0$ have a geometric $p$-rough lift $\mathbf X$, with $\Sblk=\Sig^{\le M}(\mathbf X)$,
and let $\E[(X^i_1)^2]<\infty$ for some $i\in\alphabet$. Then
$\mu_{ii}=\tfrac12\E[(X^i_1)^2]=\tfrac12\Var(X^i_1)+\tfrac12(\E
X^i_1)^2$.
\end{lemma}
\begin{proof}
The signature of a geometric $p$-rough path is a limit of signatures of bounded
variation paths, so it satisfies the shuffle relation
$\langle u,\Sblk\rangle\langle v,\Sblk\rangle=\langle u\shuffle v,\Sblk\rangle$ pathwise, by
integration by parts \citep{Reutenauer1993}, \citep[Ch.~7 and Ch.~9]{FrizVictoir2010}. With $u=v=(i)$ and
$(i)\shuffle(i)=2\,(ii)$, this gives $\langle i,\Sblk\rangle^2=2\langle ii,\Sblk\rangle$. Since
the level-$1$ signature is the increment, $\langle i,\Sblk\rangle=X^i_1$, so
$(X^i_1)^2=2\langle ii,\Sblk\rangle$. Taking expectations gives the claim. No Gaussianity is used.
\end{proof}

\noindent So the expected signature is not group-like: $2\mu_{ii}-\mu_i\mu_i=\Var(X^i_1)$ is the defect of
the shuffle relation, and in the centred class $\mu_{ii}=\tfrac12\Var(X^i_1)$.

\begin{remark}[centring and independent components annihilate the expected area]\label{rem:revarea}
Let $A^{ij}=\tfrac12(\langle ij,\Sblk\rangle-\langle ji,\Sblk\rangle)$ be the L\'evy area. For
$i=j$ it vanishes pathwise. For $i\ne j$ and a centred Gaussian driver with independent
components, $\E\langle ij,\Sblk\rangle=\E\int_{0<s<t<1}\dd X^i_s\,\dd X^j_t$ has a single Wick
pairing, joining a slot of component $i$ to a slot of component $j$. The two components are
independent and centred, so the pairing vanishes, and likewise for $\langle ji,\Sblk\rangle$. Hence
\begin{equation}\label{eq:EA0}
   \E[A^{ij}]=0 .
\end{equation}
The same argument works on every sub-interval and for any process whose components are independent
and centred, in particular the piecewise-linear interpolant $X^n$. For the natural lift of a rough $X$,
$\E[A^{ij}]$ is the limit of the same expectation for the interpolants (Lemma~\ref{lem:L1conv}).
A linear drift changes nothing: for $i\ne j$ at most one of the two slots can be a drift slot, since
the drift has only the component $1$, and the other slot then carries $\dd G$, whose mean vanishes. So
\eqref{eq:EA0} holds throughout the class of Definition~\ref{def:class}. It drives the level-$4$
bias analysis of \S\ref{sec:bias}.

Independence of the components, not Gaussianity, is what matters. With
correlated components a centred stationary Gaussian driver can have $\E[A^{12}]\ne0$: for the
stationary $2$-D Ornstein--Uhlenbeck process $\dd X_t=-BX_t\,\dd t+\dd W_t$ with
$B=\begin{psmallmatrix}\lambda&-\omega\\ \omega&\lambda\end{psmallmatrix}$, $\lambda>0$, $\omega\ne0$, the
expected area over a short interval of length $\delta$ is $-\omega\delta^2/2+O(\delta^3)$.
\end{remark}

\section{The weak error on the Young range, and the level structure}\label{sec:bias}

This section studies the weak error $\cW_{\le M}$ on the Young range. Unless a statement says
otherwise, the driver is the independent-fBm$(H)$. Proposition~\ref{prop:level4} and the bounds of
\S\ref{subsec:alllevels} hold for every (centred, resp.\ drifted) member of the class;
Proposition~\ref{prop:level3drift} needs stationary increments and Hypothesis~(S); and
Proposition~\ref{prop:nongauss} concerns non-Gaussian drivers outside the class.

\subsection{The bias is a level-4 effect}

The next proposition locates the first level that carries a weak error in the centred class.
Propositions~\ref{prop:alllevels} and~\ref{prop:adm} use it for the levels $m\le3$, and it gives
the Gaussian part of Main result~2.

  \begin{proposition}[the bias begins at level 4 in the centred class]\label{prop:level4}
For a centred Gaussian driver with independent components ($P\in\cP_0$) the weak error
$\cW_m:=\E\langle\cdot,\hSblk_\delta\rangle-\E\langle\cdot,\Sblk\rangle$ vanishes for $m\le3$.
For the independent-fBm with $d\ge2$ and $H\in(\tfrac14,1)$, $\cW_4\ne0$ on all sufficiently fine
meshes (on every mesh at $H=\tfrac12$), so there the weak error first appears at level~$4$.
\end{proposition}

Centring cannot be dropped: for $\theta\ne0$ a level-$3$ word with an odd number of drift slots has
an even number of Gaussian slots, and the level-$3$ bias is nonzero and of order $\delta^{2H}$
(Proposition~\ref{prop:level3drift}; Proposition~\ref{prop:drift4} below the crossover, where it is
of lower order than level~4). Without Gaussianity, centring together with independence of the
components still kills levels $1$--$3$ (Proposition~\ref{prop:nongauss}(i)), but centring alone does
not (Remark~\ref{rem:whygauss}).

\begin{proof}
We compare the expected signatures of $X$ and of its interpolant level by level. Levels $1$
and~$2$ are fixed by the endpoint and by the vanishing expected area, level~$3$ by Wick's theorem for
odd moments. For level~$4$ we cite the exact constants proved later.
The interpolant $X^n$ matches $X$ at every node, in particular $X^n_1=X_1$. We use one fact twice.

\emph{(F1)} $X^n$ is again centred Gaussian with independent components, being a componentwise
linear image of the Gaussian nodal vector.

\emph{Level $1$.} $\langle i,\Sblk\rangle=X^{i}_1$ depends only on the endpoint, which $X^n$
reproduces, and $\E X^i_1=0$, so $\cW_1=0$.

\emph{Level $2$.} The symmetric part is determined by the endpoint through the shuffle relation
$\langle ij,\Sblk\rangle+\langle ji,\Sblk\rangle=X^i_1X^j_1$ (Lemma~\ref{lem:diag2} for $i=j$, the
same shuffle for $i\ne j$), hence is the same for $X$ and $X^n$. The antisymmetric part is the area
$A^{ij}$. By \eqref{eq:EA0} of Remark~\ref{rem:revarea}, $\E[A^{ij}(X)]=0$, and
$\E[A^{ij}(X^n)]=0$ by the same computation, using (F1). So $\cW_2=0$.

\emph{Level $3$.} Every level-$3$ word has odd length, so both means vanish, by Wick's theorem for
the centred Gaussian $X$ (through its interpolants, Lemma~\ref{lem:L1conv}) and by (F1) for $X^n$;
$\cW_3=0$.

\emph{Level $4$.} Here the two centred Gaussian laws of $X$ and $X^n$ can first have different
even-order Wick sums, since the increment covariance enters quadratically and, off the symmetric
part, is not determined by the endpoint. That $\cW_4\ne0$ in the stated cases follows from the exact constants of Proposition~\ref{prop:C1exact} on $(\tfrac12,1)$,
from Theorem~\ref{thm:2H}(i) at $H=\tfrac12$, and from $\kappa_{\rm alt}\ne0$ on $(\tfrac14,\tfrac12)$
(Propositions~\ref{prop:Lexists} and~\ref{prop:kappasign}).
\end{proof}

Gaussianity enters the proof of Proposition~\ref{prop:level4} only through Wick's theorem. The
next proposition shows that independence, centring and moments suffice for levels $1$--$4$. It gives
the non-Gaussian half of the level structure stated in the Introduction.

\begin{proposition}[levels $3$ and $4$ need no Gaussianity]\label{prop:nongauss}
Let $Y=(Y^1,\dots,Y^d)$ have mutually independent components with continuous paths,
$Y^i_0=0$, $\E Y^i_t=0$ and $\sup_t\E|Y^i_t|^M<\infty$ (so that, by H\"older's inequality, products
of at most $M$ increments of one component are integrable), and suppose
$\Sig^{\le M}(Y^{n'})\to\Sig^{\le M}(Y)$ in $L^1$ along the uniform meshes, $Y^{n'}$ the
interpolant. Then, for $M\ge4$:
\begin{enumerate}[label=(\roman*),leftmargin=2em,itemsep=1pt,topsep=2pt]
\item $\cW_1(Y)=\cW_2(Y)=\cW_3(Y)=0$;
\item $\cW_4(Y)$ depends on the law of $Y$ only through the covariance functions of its components,
and equals $\cW_4$ of the centred Gaussian process with independent components and the same
covariances whenever that process's natural lift is the $L^1$ limit of its interpolants;
\item if every component is symmetric, $Y^i\overset d=-Y^i$, then $\cW_m(Y)=0$ for every odd $m\le M$;
\item in general, for $m\le M$, $\cW_m(Y)$ depends on the law only through the moment functions of the single
components of orders $2,\dots,m-2$.
\end{enumerate}
For level~$3$ alone, keeping the continuity, moment and $L^1$ hypotheses, mutual independence
and centring can be replaced by the condition that all
mixed third moment functions $\E[Y^i_sY^j_tY^k_u]$ ($i,j,k$ not all equal) vanish: under it
$\cW_3(Y)=0$.
 \end{proposition}

The third-moment condition is sufficient but not necessary for $\cW_3(Y)=0$. The condition does not imply $\cW_2(Y)=0$. For instance, let
$Y_t=(\xi t,\xi t^2)$ with $\xi\sim N(0,1)$. This process is centred and Gaussian, so all its
third moments vanish and the condition holds. On a cell of width $\delta$, the curve
$t\mapsto(t,t^2)$ has L\'evy area $\delta^3/6$, so $\E[A_c]=\delta^3/6$, where $A_c$ is the
L\'evy area of $Y$ over the cell $I_c$. By \eqref{eq:areaChen}, the L\'evy-area component of
$\cW_2$ is then $-\sum_c\E[A_c]=-\delta^2/6\ne0$.

  The $L^1$ convergence hypothesis holds whenever $\E\norm Y^M_{\pvar}<\infty$ for some $p<2$. Indeed,
by \citet[Prop.~5.20]{FrizVictoir2010} the interpolants are bounded in $p$-variation by
$3^{1-1/p}\norm Y_{\pvar}$ and converge in $p'$-variation, $p'\in(p,2)$, so the Young signatures
converge a.s., and $|\Sig^{(m)}(Y^{n'})|\le C\norm Y^m_{\pvar}$ dominates. The lift condition in (ii)
holds, for instance, if the covariance has finite $\rho$-variation, $\rho<2$, dominated by a
two-dimensional control. This domination holds at every $\rho'\in(\rho,2)$ by
\citet[Thm.~1(ii)]{FrizVictoir2011}, applied to the interpolants' covariances after
\citet[Prop.~15.11]{FrizVictoir2010} with a constant uniform in the mesh. The convergence then
follows from \citet[Thm.~15.33(iv) and its proof, Cor.~15.29(ii), Thm.~9.5(ii),
Cor.~9.11(i)]{FrizVictoir2010}.

\begin{proof}
\emph{Strategy.} Fix a mesh $\delta=1/n$. By the $L^1$ hypothesis,
$\E\langle w,\Sig(Y)\rangle=\lim_{n'\to\infty}\E\langle w,\Sig(Y^{n'})\rangle$ for every word $w$ with
$|w|\le M$, so
\[
   \langle w,\cW_{\le M}(Y)\rangle=\E\langle w,\Sig(Y^n)\rangle-\lim_{n'\to\infty}\E\langle w,\Sig(Y^{n'})\rangle .
\]
Both terms involve only signatures of piecewise-linear paths. It therefore suffices to compute
$\E\langle w,\Sig(Z)\rangle$ for a piecewise-linear interpolant $Z$ of $Y$ on an arbitrary uniform mesh,
and to check that the claimed property holds for every mesh. It then passes to the limit $n'\to\infty$.

\emph{Step 1: the signature of a piecewise-linear path.} Let $Z$ be linear on each cell $I_a$, with
increment $\Delta_a=Z_{t_{a+1}}-Z_{t_a}\in\RR^d$. The signature of a straight segment with increment
$\Delta$ is $\exp(\Delta)=\sum_{k\ge0}\Delta^{\otimes k}/k!$, so by Chen's identity
$\Sig(Z)=\exp(\Delta_0)\otimes\cdots\otimes\exp(\Delta_{n-1})$. Taking the coefficient of
$w=(w_1\cdots w_m)$,
\begin{equation}\label{eq:chenPL}
   \langle w,\Sig(Z)\rangle=\sum_{0\le a_1\le\cdots\le a_m\le n-1}c_{\mathbf a}\prod_{l=1}^m\Delta^{w_l}_{a_l},
   \qquad c_{\mathbf a}=\prod_{b=0}^{n-1}\frac{1}{k_b(\mathbf a)!},\quad k_b(\mathbf a):=\#\{l:a_l=b\}.
\end{equation}
 Each slot $l$ of $w$ is assigned a cell $a_l$, the assignment is non-decreasing, and the factor
$1/k!$ comes from $k$ slots sharing one cell.

\emph{Step 2: factorisation over letters.} Let $m_i(w):=\#\{l:w_l=i\}$ be the multiplicity of the
letter $i$ in $w$. Grouping the factors of each product in \eqref{eq:chenPL} by letter and using the
independence of the components,
\begin{equation}\label{eq:factorletters}
   \E\prod_{l=1}^m\Delta^{w_l}_{a_l}=\prod_{i=1}^d\E\Bigl[\prod_{l:\,w_l=i}\Delta^i_{a_l}\Bigr].
\end{equation}
The $i$-th factor is a joint moment of order $m_i(w)$ of increments of the single component $Y^i$,
finite by the moment assumption. Expanding each increment into a difference of values shows that it
is a linear combination of the moment functions $\E[Y^i_{s_1}\cdots Y^i_{s_k}]$, $k=m_i(w)$.

\emph{Step 3: two kinds of words with zero weak error.}
(a) If some letter $i$ has $m_i(w)=1$, its factor in \eqref{eq:factorletters} is
$\E\Delta^i_a=0$, since $Y^i$ is centred. Then $\E\langle w,\Sig(Z)\rangle=0$ for every interpolant,
so both terms of $\langle w,\cW_{\le M}(Y)\rangle$ vanish.
(b) If $w=i^m$ is a single-letter word, then for every path $z$,
$\langle i^m,\Sig(z)\rangle=(z^i_1-z^i_0)^m/m!$, the signature of a scalar path. The interpolant has
the same endpoints as $Y$, so $\langle i^m,\Sig(Y^n)\rangle=\langle i^m,\Sig(Y)\rangle$ pathwise and
$\langle i^m,\cW_{\le M}(Y)\rangle=0$.

\emph{Step 4: proof of (i).} A word with a single distinct letter is of type (b), and a word
of length at most $3$ with at least two distinct letters has a letter of multiplicity one, so it is
of type (a). Hence
$\cW_1(Y)=\cW_2(Y)=\cW_3(Y)=0$.

\emph{Step 5: proof of (ii).} At length $4$ the multiplicity patterns are $(4)$, $(3,1)$, $(2,1,1)$,
$(1,1,1,1)$ and $(2,2)$. All but $(2,2)$ are of type (a) or (b). A $(2,2)$-word, such as $iijj$ or
$ijij$ with $i\ne j$, has by \eqref{eq:factorletters}
\[
   \E\prod_{l=1}^4\Delta^{w_l}_{a_l}=\E\bigl[\Delta^i_{a}\Delta^i_{a'}\bigr]\,
   \E\bigl[\Delta^j_{b}\Delta^j_{b'}\bigr],
\]
where $a,a'$ are the cells of the two $i$-slots and $b,b'$ those of the two $j$-slots. This is a
product of increment covariances of $Y^i$ and $Y^j$. So $\E\langle w,\Sig(Z)\rangle$, for $Z=Y^n$ and
$Z=Y^{n'}$, depends on the law of $Y$ only through the covariance functions of its components, and so
does the limit $\E\langle w,\Sig(Y)\rangle$. The centred Gaussian process with independent components
and the same covariances gives the same value at every mesh. If its natural lift is the $L^1$ limit of
its interpolants, its limit is the same as well, which proves the comparison in (ii).

\emph{Step 6: proof of (iii).} If $m$ is odd, then $\sum_im_i(w)=m$ implies that some $m_i(w)$ is odd.
If $Y^i\overset d=-Y^i$, every joint moment of odd order of increments of $Y^i$ vanishes, so the
$i$-th factor in \eqref{eq:factorletters} is zero. Hence $\E\langle w,\Sig(Z)\rangle=0$ for every
interpolant and every word of length $m$, and $\cW_m(Y)=0$.

\emph{Step 7: proof of (iv).} After Step 3, only words in which every letter appears at least twice
(not type (a)) and at least two distinct letters appear (not type (b)) remain. In such a word
$2\le m_i(w)\le m-2$ for every letter $i$ present. By Steps 1 and 2, $\E\langle w,\Sig(Z)\rangle$ is a
linear combination of products of moment functions of the single components of orders between $2$ and
$m-2$, and this persists in the limit.

\emph{The level-$3$ refinement.} Now drop independence and centring, and assume only that
$\E[Y^i_sY^j_tY^k_u]=0$ whenever $i,j,k$ are not all equal. Single-letter words $iii$ have zero weak
error by Step 3(b), which uses neither assumption. For any other word $w=(ijk)$, \eqref{eq:chenPL}
writes $\E\langle w,\Sig(Z)\rangle$ as a combination of mixed third moments
$\E[\Delta^i_a\Delta^j_b\Delta^k_c]$. Expanding the increments into values turns these into
combinations of $\E[Y^i_sY^j_tY^k_u]$, which vanish. So $\E\langle w,\Sig(Y^{n'})\rangle=0$ for every
$n'$, hence $\E\langle w,\Sig(Y)\rangle=0$, and $\cW_3(Y)=0$.

\emph{The condition is not necessary.} Let $Y^i_t=B^i_t+\zeta_it$, $i=1,2$, where $B$ is a
two-dimensional Brownian motion and $(\zeta_1,\zeta_2)$ is independent of $B$, centred, with moments
of order $M$ and $\E[\zeta_1^2\zeta_2]\ne0$. This $Y$ satisfies the standing hypotheses:
$\sup_t\E|Y^i_t|^M<\infty$, and the level-$m$ coordinates of the interpolants are sums of products of
monomials in $\zeta$ of degree at most $m\le M$ with iterated integrals of the Brownian interpolants
against $\dd B$ and $\dd t$. The latter converge in every $L^q$ to their Stratonovich limits, and $\zeta$
is independent of $B$, so $\Sig^{\le M}(Y^{n'})$ converges in $L^1$. Expanding the product and using that $B$ is centred,
Gaussian and independent of $\zeta$,
\[
   \E[Y^1_sY^1_tY^2_u]=(s\wedge t)\,u\,\E\zeta_2+stu\,\E[\zeta_1^2\zeta_2]=stu\,\E[\zeta_1^2\zeta_2]\ne0,
\]
so the condition fails. Still $\cW_3(Y)=0$. For a word $w$ of length $3$, the same expansion gives
$\E[\Delta^{w_1}_a\Delta^{w_2}_b\Delta^{w_3}_c]=\delta^3\,\E[\zeta_{w_1}\zeta_{w_2}\zeta_{w_3}]$ on cells of
width $\delta$: every term containing an increment of $B$ has mean zero (one or three increments
of $B$: an odd Gaussian moment; two increments: the factor $\E\zeta_i=0$). By \eqref{eq:chenPL},
\[
   \E\langle w,\Sig(Y^n)\rangle=\E[\zeta_{w_1}\zeta_{w_2}\zeta_{w_3}]\sum_{a_1\le a_2\le a_3}c_{\mathbf a}\,\delta^3
   =\tfrac16\,\E[\zeta_{w_1}\zeta_{w_2}\zeta_{w_3}],
\]
since the sum is the level-$3$ coordinate of the scalar path $t\mapsto t$, which equals $1/3!$ at every
mesh. The value does not depend on $n$, so it is also $\E\langle w,\Sig(Y)\rangle$, and the weak error
vanishes. The reason is that the third moment function $stu\,\E[\zeta_1^2\zeta_2]$ is linear in each
time variable, so linear interpolation reproduces it exactly.
\end{proof}

For centred independent components, a level-$4$ word in which some letter appears exactly once
has mean zero for $X$ and for $X^n$, and the words $iiii$ are determined by the endpoint. So the
level-$4$ weak error lives in the $6\binom d2$ words $iijj,ijij,ijji,jiij,jiji,jjii$, $\{i,j\}$ a
pair of distinct letters. By the single surviving Wick pairing and the symmetry
$i\!\leftrightarrow\!j$, each pair reduces to the same three values
\[
   E_{\rm adj}:=\E\langle1122,\Sig(X)\rangle,\quad
   E_{\rm alt}:=\E\langle1212,\Sig(X)\rangle,\quad
   E_{\rm nest}:=\E\langle1221,\Sig(X)\rangle,
\]
and $E^n_\bullet$ for the same coordinates of the interpolant $X^n$.
For the interpolant, they are sums over ordered cell assignments of products of two increment
covariances, with the coefficients $c_{\mathbf a}$ of \eqref{eq:chenPL}.
The expectations $E_\bullet$ are integrals over the ordered simplex of products of the increment
covariance, one factor for each Wick pair of slots. When the covariance measure charges the diagonal,
as for Brownian motion, two adjacent slots of a pair can coincide, and we must decide how much of this
diagonal mass the ordered simplex contains. We count one half: for each Wick pair of adjacent slots,
\begin{equation}\label{eq:coinc}
   \int_{\{s\le t\}}\mu\ :=\ \int_{\{s<t\}}\mu\ +\ \tfrac12\int_{\{s=t\}}\mu .
\end{equation}
For symmetric $\mu$ this equals $\tfrac12\int_{[0,1]^2}\mu$, in line with
$\langle ii,\Sig(X)\rangle=(X^i_1)^2/2$. It is also the limit of the factor $1/2!$ that the
piecewise-linear formula \eqref{eq:chenPL} gives to two slots in one cell. Two non-adjacent slots can
coincide only if the slots between them coincide as well, and in every case used here such sets
carry no mass.
For $H>\tfrac12$ the diagonal is $\varrho$-null and \eqref{eq:coinc} is the ordinary integral. At
$H=\tfrac12$, where $\varrho(\dd s,\dd t)=\dd s\,\delta_{\{s=t\}}$ is carried by the diagonal, the
strict convention, which gives the diagonal the weight $0$,  would return $0$ for every even-level word, and \eqref{eq:coinc} is what makes the
Brownian evaluations agree with $\E\Sig=\exp(\tfrac12\sum_ie_i^{\otimes2})$. We use it only for
Brownian computations, and the general bounds of \S\ref{subsec:alllevels} avoid it by working with
interpolants, whose kernels are bounded densities. Below the crossover
$\varrho$ is not a measure and \eqref{eq:coinc} is not used: there the level-$4$ coordinates are computed from the covariance (Lemma~\ref{lem:cov4}). With $\cW_4$ the level-4 bias
and $\Delta E_\bullet=E^n_\bullet-E_\bullet$,
\begin{equation}\label{eq:W4norm}
   \norm{\cW_4}=\Bigl(\tbinom d2\cdot2\sum_\bullet(\Delta E_\bullet)^2\Bigr)^{1/2},
\end{equation}
and $\sum_\bullet\Delta E_\bullet=0$ by the endpoint identity (Lemma~\ref{lem:sumzero}).
Throughout \S\ref{sec:bias} we display the case $d=2$, where $\binom d2=1$. General-$d$
constants are in Remark~\ref{rem:dscope}.

\subsection{The exponent \texorpdfstring{$2H$}{2H} on the Young range}\label{subsec:gamma2H}

By Proposition~\ref{prop:level4}, level~$4$ is the first level with a weak error. The next
theorem gives its order on the Young range, which is the rate $\delta^{2H}$ of Main result~1 there,
and its proof supplies the double-error bound \eqref{eq:T3bound} used in
Proposition~\ref{prop:C1exact}.

\begin{theorem}[level-4 weak-error exponent]\label{thm:2H}
Let $d\ge2$ and let $X$ have independent $\mathrm{fBm}(H)$ components. For $H\in[\tfrac12,1)$,
$\norm{\cW_4}\asymp\delta^{2H}$. Precisely:
\begin{enumerate}[leftmargin=1.6em,topsep=2pt]
\item[(i)] At $H=\tfrac12$ (Brownian), $\Delta E_{\rm adj}=-\delta/12$, $\Delta E_{\rm
alt}=\Delta E_{\rm nest}=+\delta/24$, so $\norm{\cW_4}=\delta/(4\sqrt3)$ exactly at $d=2$; for general $d$ multiply by $\sqrt{\binom d2}$.
\item[(ii)] For $H\in(\tfrac12,1)$, $\norm{\cW_4}=O(\delta^{2H})$, and $\delta^{-2H}\norm{\cW_4}$
converges to a positive limit as $\delta\to0$, so $\norm{\cW_4}\asymp\delta^{2H}$.
\end{enumerate}
\end{theorem}

The upper bound in (ii) is proved below. That the limit exists and is positive follows from its
exact value, computed in Proposition~\ref{prop:C1exact}.

\begin{proof}[Proof of (i)]
We compute the three interpolant values $E^n_\bullet$ exactly from \eqref{eq:chenPL} and
subtract the Brownian values. Write $\Gamma_{jk}:=\E[\Delta^i_j\Delta^i_k]$ for the covariance of the increments of one component
over cells $j$ and $k$. For Brownian motion, increments over disjoint cells are independent, so
$\Gamma_{jk}=\delta\,\one_{j=k}$. In \eqref{eq:chenPL}, with the two letters Wick-paired, a term is
therefore nonzero only if paired slots lie in the same cell.
For $1122$ this means that slots $1,2$ lie in a cell $a$ and slots $3,4$ in a cell $b\ge a$. If $a<b$
the coefficient is $\tfrac1{2!}\cdot\tfrac1{2!}$; if $a=b$ it is $\tfrac1{4!}$, and
$\E[(\Delta^1_a)^2(\Delta^2_a)^2]=\delta^2$. There are $\binom n2$ pairs $a<b$, each
contributing $\tfrac14\delta\cdot\delta$, and $n$ cells $a=b$, each contributing
$\tfrac1{24}\delta^2$. Hence, with $n\delta=1$,
\[
   E^n_{\rm adj}=\tfrac{\delta^2}4\tbinom n2+\tfrac{\delta^2n}{24}
   =\tfrac{1-\delta}8+\tfrac{\delta}{24} =\tfrac18-\tfrac{\delta}{12}.
\]
For $1212$ and $1221$ the paired slots are not adjacent, and placing each pair in one cell puts all four slots into the same cell: for $1212$, $a_1=a_3$ and $a_1\le a_2\le a_3$ give
$a_2=a_1$, and then $a_4=a_2$; for $1221$, $a_1=a_4$ gives $a_1=a_2=a_3=a_4$ directly, so $E^n_{\rm alt}=E^n_{\rm nest}=n\cdot\tfrac{\delta^2}{24}=\tfrac{\delta}{24}$.
For Brownian motion itself, $\E\Sig(B)=\exp(\tfrac12\sum_ie_i^{\otimes2})$ has level-$4$ part
$\tfrac18\sum_{i,j}e_i\otimes e_i\otimes e_j\otimes e_j$, so $(E_{\rm adj},E_{\rm alt},E_{\rm nest})
=(\tfrac18,0,0)$, the limits of the interpolant values as $\delta\to0$. Subtracting,
$\Delta E_{\rm adj}=-\tfrac{\delta}{12}$ and $\Delta E_{\rm alt}=\Delta E_{\rm nest}=\tfrac{\delta}{24}$,
and \eqref{eq:W4norm} gives $\norm{\cW_4}=\bigl(2(\tfrac{\delta^2}{144}+\tfrac{2\delta^2}{576})\bigr)^{1/2}
=\delta/(4\sqrt3)$.
\end{proof}

\begin{proof}[Proof of (ii), upper bound]
Writing $\varrho^n=\varrho+e$ splits each $\Delta E_\bullet$ into two single-error terms and
one double-error term (Step~2). A single-error term integrates $e$, which has zero mass on every
cell product, against a Lipschitz weight, so it is at most a constant times
$\delta\norm e_{L^1}=O(\delta^{2H})$ (Steps~1 and~3). The double-error term is $O(\delta^{4H-1})$
(Step~3$'$). Here $\varrho(s,t)=c_H|s-t|^{2H-2}$, $c_H=H(2H-1)>0$, is the increment-covariance
density of fBm$(H)$, and $\varrho^n$ is its cell average, so $e:=\varrho^n-\varrho$ satisfies
$\int_{I_j\times I_k}e=0$ on every cell product.

\emph{Step 1 (the $L^1$ size of the error kernel).} We show $\norm e_{L^1}=O(\delta^{2H-1})$.
Split cells by lag $m=|j-k|$.
\begin{itemize}[leftmargin=1.4em,topsep=1pt,itemsep=1pt]
\item \emph{Diagonal and adjacent cells ($m\le1$, the singular band).} Here
$\int_{I_j\times I_k}|e|\le2\int_{I_j\times I_k}|\varrho|$, since $\varrho^n$ is constant on the
cell and $\int_{I_j\times I_k}|\varrho^n|=\bigl|\int_{I_j\times I_k}\varrho\bigr|\le\int_{I_j\times I_k}|\varrho|$, and by the self-similarity
$\varrho(\delta u,\delta v)=\delta^{2H-2}\varrho(u,v)$,
\[
  \int_{I_j\times I_k}|\varrho|\ \le\ c_H\,\delta^{2H}\!\!\int_{[0,2]^2}\!\!|u-v|^{2H-2}\,\dd u\,\dd v
  =\ (2\delta)^{2H},
\]
finite because $2H-2>-1$; here $c_H\int_{[0,2]^2}|u-v|^{2H-2}$ is the variance $2^{2H}$ of an fBm increment over an
interval of length~$2$. There are $n+2(n-1)<3n=3/\delta$ such cells, contributing at most
$6\cdot2^{2H}\delta^{2H-1}$.
\item \emph{Far cells ($m\ge2$).} Each such cell satisfies the far-cell bound
\eqref{eq:farcell} of Remark~\ref{rem:C1scope} (mean value theorem on the cell), with
$|c_H|=c_H$. There are
$2(n-m)\le2n$ cells at lag $m$, so, with $n=1/\delta$,
\[
  \begin{aligned}
  \sum_{m\ge2}\ \sum_{|j-k|=m}\int_{I_j\times I_k}|e|
  &\ \le\ 2n\,c_H(2-2H)\,\delta^{2H}\sum_{m\ge2}(m-1)^{2H-3}\\
  &\ =\ 2c_H(2-2H)\,\zeta(3-2H)\,\delta^{2H-1},
  \end{aligned}
\]
the sum $\sum_{m\ge1}m^{2H-3}=\zeta(3-2H)$ converging since $3-2H>1$.
\end{itemize}
Both bands are $O(\delta^{2H-1})$, proving Step~1:
\begin{equation}\label{eq:eL1fbm}
   \norm e_{L^1}\ \le\ \bigl(6\cdot2^{2H}+2c_H(2-2H)\,\zeta(3-2H)\bigr)\,\delta^{2H-1}.
\end{equation}

\emph{Step 2 (decomposition of $\Delta E_{\rm adj}$).} After the single surviving Wick pairing, the
level-$4$ $(2,2)$-word integral is a $4$-fold ordered-simplex integral of a product of two kernel
factors, and discretizing replaces $\varrho$ by $\varrho^n$ in each. Writing $\varrho^n=\varrho+e$
in both factors splits the difference into two single-error terms and one double-error term:
\begin{equation}\label{eq:eexpand}
  \varrho^n(t_1,t_2)\varrho^n(t_3,t_4)-\varrho(t_1,t_2)\varrho(t_3,t_4)
  =e(t_1,t_2)\varrho(t_3,t_4)+\varrho(t_1,t_2)e(t_3,t_4)+e(t_1,t_2)e(t_3,t_4).
\end{equation}
Integrating over $0<t_1<t_2<t_3<t_4<1$, and in the two single-error terms integrating out the exact
factor over its ordered pair, gives, exactly,
\[
  \Delta E_{\rm adj}=T_1+T_2+T_3,\qquad
  T_1=\!\!\int_{t_1<t_2}\!\!e(t_1,t_2)\,w_\sharp(t_2),\qquad
  T_2=\!\!\int_{t_3<t_4}\!\!e(t_3,t_4)\,w_\flat(t_3),
\]
with $T_3=\int_{t_1<t_2<t_3<t_4}e(t_1,t_2)\,e(t_3,t_4)$. The weights are the partial integrals of
the exact continuum factor, $w_\sharp(t_2)=\int_{t_2<t_3<t_4<1}\varrho(t_3,t_4)$ and
$w_\flat(t_3)=\int_{0<t_1<t_2<t_3}\varrho(t_1,t_2)$, that is, half the variance of an fBm increment:
\[
  w_\sharp(a)=\tfrac12(1-a)^{2H},\qquad w_\flat(b)=\tfrac12\,b^{2H},
\]
both $C^1$ on $[0,1]$ with $w_\sharp'(a)=-H(1-a)^{2H-1}$, $w_\flat'(b)=H\,b^{2H-1}$, hence (using $2H-1>0$)
\[
  \norm{w_\sharp'}_\infty=\norm{w_\flat'}_\infty=H .
\]

\emph{Step 3 (bounding $T_1,T_2$).} Since $e$ is symmetric and the diagonal is null for
$H>\tfrac12$, $e$ also has zero mass on each half $\{t_1<t_2\}\cap(I_j\times I_j)$ of a diagonal
cell, and the domain $\{t_1<t_2\}$ of $T_1$ is a union of such pieces and of off-diagonal cell
products. On a piece with $t_2\in I_k$ fix any $t_2^\ast\in I_k$. Subtracting the constant
$w_\sharp(t_2^\ast)$ changes nothing (cell cancellation), and the remaining weight is at most
(mesh)$\times$(derivative), which gains the factor $\delta$:
\[
   T_1=\sum_{\text{pieces}}\int e(t_1,t_2)\bigl(w_\sharp(t_2)-w_\sharp(t_2^\ast)\bigr).
\]
As $t_2$ and $t_2^\ast$ lie in one cell, $|w_\sharp(t_2)-w_\sharp(t_2^\ast)|\le\norm{w_\sharp'}_\infty\delta=H\delta$,
and by Step~1
\[
  |T_1|\ \le\ H\delta\int_{t_1<t_2}|e|\ \le\ H\delta\,\norm e_{L^1}\ =\ H\,\delta\cdot O(\delta^{2H-1})\ =\ O(\delta^{2H}).
\]
The same argument, with $t_3$ in place of $t_2$ and $w_\flat$ in place of $w_\sharp$, gives
$|T_2|=O(\delta^{2H})$.

\emph{Step 3$'$ (the other two words, and the double error).} For $E_{\rm alt}$ and
$E_{\rm nest}$ Step~2 gives the same three terms, the exact factor now integrating to
the corresponding weight $w$ of Lemma~\ref{lem:weights} instead of $w_\sharp,w_\flat$. Each of the six
weights is a combination of $x\mapsto x^{2H}$ at $s,t,1-s,1-t,t-s$ whose two partial derivatives
are bounded by $H$ in modulus. For example, for the alternating word with the error on slots
$\{1,3\}$, $\partial_sw=H[(t-s)^{2H-1}-(1-s)^{2H-1}]$ is a difference of two numbers in $[0,1]$. Each
piece $(I_j\times I_k)\cap\{s<t\}$ is convex, so the mean value theorem along a segment applies, and
for any point $(s^\ast,t^\ast)$ of the closed piece
$|w(s,t)-w(s^\ast,t^\ast)|\le H|s-s^\ast|+H|t-t^\ast|\le2H\delta$.
Hence the bound for
$T_1,T_2$ holds with $2H$ in place of $H$. The double error needs a separate argument, which we
give for all three words at once. Fix the Wick pairing. In $T_3$ both covariance factors are replaced by the error $e$. Let $t_A$ be
the two time variables of the factor that contains $t_1$, namely $(t_1,t_2)$, $(t_1,t_3)$,
$(t_1,t_4)$ for the adjacent, alternating and nested words, and $t_B$ the two time variables of the
other factor, and integrate over $t_B$ first:
\[
   T_3^\bullet=\int_{\{t_A\ \text{ordered}\}}\!\! e(t_A)\,h_\bullet(t_A)\,\dd t_A ,
   \qquad
   h_\bullet(t_A)=\int_{R_\bullet(t_A)\cap\{s<t\}}\!\! e(s,t)\,\dd s\,\dd t ,
\]
where the simplex ordering cuts the $t_B$-variables out of the rectangle
\[
   R_{\rm adj}(t_1,t_2)=(t_2,1)^2,\qquad
   R_{\rm alt}(t_1,t_3)=(t_1,t_3)\times(t_3,1),\qquad
   R_{\rm nest}(t_1,t_4)=(t_1,t_4)^2 .
\]
The edges of $R_\bullet(t_A)$ need not lie on cell boundaries (for the alternating word they
are $t_2=t_1$, $t_2=t_3$, $t_4=t_3$), and on cells cut by an edge the zero cell averages of $e$ are
lost. The following two estimates control these boundary cells, for all three words. The cells
$I_a=[t_a,t_{a+1})$ are half-open, so each endpoint of a subinterval of $[0,1]$ lies in exactly one
of them.

\emph{(a) Row mass.} Step~1, applied to a single row $a$, gives
\begin{equation}\label{eq:rowmass}
   \max_{0\le a<n}\ \int_{I_a\times[0,1]}|e|\ \le\ C^{\rm row}_H\,\delta^{2H},
   \qquad C^{\rm row}_H:=2\cdot3^{2H}+2c_H(2-2H)\,\zeta(3-2H):
\end{equation}
the three cells $|b-a|\le1$ lie in a square of side $3\delta$ meeting the diagonal, on which
$\int|\varrho|=(3\delta)^{2H}$ and $\int|e|\le2\int|\varrho|$. The row contains two far cells
of each lag $m\ge2$, each bounded by $c_H(2-2H)\delta^{2H}(m-1)^{2H-3}$ by \eqref{eq:farcell}, and
$\sum_{m\ge2}(m-1)^{2H-3}=\zeta(3-2H)$. Since $e(s,t)=e(t,s)$, the same bound holds for the column masses
$\int_{[0,1]\times I_b}|e|$.

\emph{(b) Boundary strips.} For every rectangle $R=J\times J'\subseteq[0,1]^2$,
\begin{equation}\label{eq:striplemma}
   \Bigl|\int_{R\cap\{s<t\}}e\Bigr|\ \le\ 4\,C^{\rm row}_H\,\delta^{2H}.
\end{equation}
Indeed, split $[0,1]^2$ into cells. A cell $I_a\times I_b\subseteq R$ contributes $0$: for
$a<b$ the constraint $s<t$ holds on the whole cell and $\int_{I_a\times I_b}e=0$. For $a>b$
the cell does not meet $\{s<t\}$. For $a=b$ the symmetry of $e$
gives $\int_{(I_a\times I_a)\cap\{s<t\}}e=\tfrac12\int_{I_a\times I_a}e=0$. A cell disjoint from $R$
contributes $0$. Any other cell meeting $R$ has $I_a\not\subseteq J$ or $I_b\not\subseteq J'$, hence
contains an endpoint of $J$ in its first factor or an endpoint of $J'$ in its second. As the cells
are half-open, each of those four endpoints lies in exactly one cell, so these cells lie in at most
two columns and two rows, whose $|e|$-mass is at most $4C^{\rm row}_H\delta^{2H}$ by \eqref{eq:rowmass} and
its column form.

Applying \eqref{eq:striplemma} to $R_\bullet(t_A)$ gives $\norm{h_\bullet}_\infty\le4C^{\rm row}_H\delta^{2H}$
for each of the three words, and $\int_{s<t}|e|=\tfrac12\norm e_{L^1}$ by symmetry, so
\begin{equation}\label{eq:T3bound}
   |T_3^\bullet|\ \le\ \norm{h_\bullet}_\infty\int_{s<t}|e|
   \ \le\ 2C^{\rm row}_H\,\delta^{2H}\norm e_{L^1}\ =\ O(\delta^{4H-1})\ =\ o(\delta^{2H})
   \qquad(\bullet\in\{\rm adj,alt,nest\}),
\end{equation}
as $4H-1>2H\iff H>\tfrac12$.

\emph{Conclusion.} For each word, $\Delta E_\bullet$ is the sum of two single-error terms, each
$O(\delta^{2H})$ by Steps~3 and~3$'$, and the double-error term $T_3^\bullet$, which is
$O(\delta^{4H-1})$ by \eqref{eq:T3bound} and is absorbed into $O(\delta^{2H})$ because $4H-1>2H$.
Hence $\Delta E_\bullet=O(\delta^{2H})$ for all three words, and $\norm{\cW_4}=O(\delta^{2H})$ by
\eqref{eq:W4norm}.
\end{proof}

The adjacent word $1122$ reduces to one-dimensional integrals of the variance function, so its
interpolation error has a closed form for every $H\in(0,1)$, used in
Proposition~\ref{prop:C1exact}, Remark~\ref{rem:merger}, Remark~\ref{find:kappa}(ii) and
Corollary~\ref{cor:kadj}.

\begin{proposition}[the adjacent coordinate in closed form]\label{prop:adjexact}
Let $X$ have independent fBm$(H)$ components, $H\in(0,1)$, and for a function $g$ on $[0,1]$ write
$I_ng$ for its piecewise-linear interpolant at the nodes $t_a=a\delta$ and $\epsilon_g:=I_ng-g$; for
$t\in I_a=[t_a,t_{a+1})$ put $\lambda(t)=(t-t_a)/\delta$. With $f(t)=t^{2H}$ and
$\psi(t)=\frac{1}{2}(1-(1-t)^{2H})$,
\[
   E_{\rm adj}(X^n)=\tfrac12\sum_{a=0}^{n-1}\bigl(\psi(t_{a+1})-\psi(t_a)\bigr)\frac{f(t_a)+f(t_{a+1})}{2}
   -\frac{\delta^{2H}}{12}.
\]
For $H>\tfrac14$, $E_{\rm adj}(X)=\frac{1}{2}\int_0^1f\psi'$ and
$\Delta E_{\rm adj}=-\frac{\delta^{2H}}{12}+\frac{1}{2}\mathfrak D_n$ with
$\mathfrak D_n:=\sum_a\int_{I_a}\bigl[\frac{f(t_a)+f(t_{a+1})}{2}-f\bigr]\psi'$, and
$|\mathfrak D_n|\le\frac{1}{2}\omega_f(\delta)$ ($\omega_f$ the modulus of continuity of $f$),
$|\mathfrak D_n|\le C_H(\delta^{2H+1}+\delta^2)$. In particular
$\Delta E_{\rm adj}=-\frac{1}{12}\delta^{2H}(1+o(1))$ for every $H\in(\tfrac14,1)$, and
$\Delta E_{\rm adj}=o(\delta^{4H-1})$ for $H<\tfrac12$.
\end{proposition}

\begin{proof}
We write $E_{\rm adj}$ as $\frac{1}{2}\int v\,\dd\psi$, with $v$ the variance of the first component
and $\psi(t)=\E[X^2_1X^2_t]-\tfrac12\E(X^2_t)^2$ for the second, both for the interpolant and for the driver, and compare the two
integrals. The interpolant's $v_n$ and $\psi_n$ are linear interpolants plus an explicit quadratic
correction on each cell. The corrections produce the term $-\delta^{2H}/12$, and the linear
interpolants produce $\frac{1}{2}\mathfrak D_n$.

\emph{Interpolant.} For continuous piecewise-linear $Y$ with $Y_0=0$, integrating out $t_1<t_2<t$
(which gives $\frac{1}{2}(Y^1_t)^2$) and $t_4>t$ (which gives $Y^2_1-Y^2_t$) yields $\langle1122,\Sig(Y)\rangle=\int_0^1\frac{1}{2}(Y^1_t)^2(Y^2_1-Y^2_t)\,\dd Y^2_t$ (Riemann--Stieltjes), so by
$X^1\perp X^2$, $E_{\rm adj}(X^n)=\frac{1}{2}\int_0^1v_n\psi_n'$ with $v_n(t)=\E(X^n_t)^2$ and
$\psi_n(t)=\E[X^n_1X^n_t]-\frac{1}{2}v_n(t)$, since $\E[(Y_1-Y_t)\dot Y_t]=\frac{\dd}{\dd t}(\E[Y_1Y_t]-\frac{1}{2}\E
Y_t^2)$ inside each cell. As $X^n_t=(1-\lambda)X_{t_a}+\lambda X_{t_{a+1}}$ and $X^n_1=X_1$,
$v_n=I_nf-\lambda(1-\lambda)\delta^{2H}$ and $\E[X^n_1X^n_t]=I_n[R(1,\cdot)]$; since
$R(1,\cdot)-\frac{1}{2}f=\psi$, $\psi_n=I_n\psi+\frac{1}{2}\lambda(1-\lambda)\delta^{2H}$ and
$\psi_n'=(I_n\psi)'+\frac{1}{2}(1-2\lambda)\delta^{2H-1}$. Here
$R(s,t)=\frac{1}{2}(f(s)+f(t)-f(|t-s|))$ is the fBm covariance, and $v_n$ follows from
$\Var((1-\lambda)X_{t_a}+\lambda X_{t_{a+1}})=(1-\lambda)f(t_a)+\lambda f(t_{a+1})-\lambda(1-\lambda)f(\delta)$.
Expanding the product,
\[
   \begin{aligned}
   \tfrac12\int v_n\psi_n'
   &=\tfrac12\!\int\! I_nf\,(I_n\psi)'+\tfrac14\delta^{2H-1}\!\!\int\! I_nf\,(1-2\lambda)\\
   &\quad-\tfrac12\delta^{2H}\!\!\int\!\lambda(1-\lambda)(I_n\psi)'
   -\tfrac14\delta^{4H-1}\!\!\int\!\lambda(1-\lambda)(1-2\lambda).
   \end{aligned}
\]
The first term is the displayed sum, since $(I_n\psi)'$ is constant on $I_a$ and
$\int_{I_a}I_nf=\delta\frac{f(t_a)+f(t_{a+1})}{2}$. On one cell, $\int_{I_a}(1-2\lambda)=0$,
$\int_{I_a}\lambda(1-2\lambda)=-\frac{\delta}{6}$ and $\int_{I_a}\lambda(1-\lambda)=\frac{\delta}{6}$, so,
summing over $a$, $\frac{1}{4}\delta^{2H-1}\int I_nf\,(1-2\lambda)=-\frac{\delta^{2H}}{24}(f(1)-f(0))$ and
$-\frac{1}{2}\delta^{2H}\int\lambda(1-\lambda)(I_n\psi)'=-\frac{\delta^{2H}}{12}(\psi(1)-\psi(0))$ are both
$-\frac{\delta^{2H}}{24}$, because $f(1)-f(0)=1$ and $\psi(1)-\psi(0)=\frac{1}{2}$, and
$\int\lambda(1-\lambda)(1-2\lambda)=0$.

\emph{Continuum.} $E_{\rm adj}(X)=\frac{1}{2}\int v\,\dd\psi$ with $v=f$ is Lemma~\ref{lem:cov4}. It is
also the limit, as $\mathfrak D_N\to0$ and $E_{\rm adj}(X^N)\to E_{\rm adj}(X)$
(natural lift, Lemma~\ref{lem:L1conv}, $H>\tfrac14$). Then $\frac{1}{2}\int I_nf(I_n\psi)'-\frac{1}{2}\int f\psi'
=\frac{1}{2}\mathfrak D_n$, because $\int_{I_a}(I_n\psi)'=\psi(t_{a+1})-\psi(t_a)=\int_{I_a}\psi'$.
Subtracting, $\Delta E_{\rm adj}=E_{\rm adj}(X^n)-E_{\rm adj}(X)=-\frac{\delta^{2H}}{12}+\frac{1}{2}\mathfrak D_n$.

\emph{Bounds.} The first is immediate: on $I_a$,
$|\frac{f(t_a)+f(t_{a+1})}{2}-f|\le\omega_f(\delta)$, while $\psi'>0$ and $\int_0^1\psi'=\frac{1}{2}$. For the second write
$\frac{f(t_a)+f(t_{a+1})}{2}-f=(\frac{1}{2}-\lambda)(f(t_{a+1})-f(t_a))+\epsilon_f$ on $I_a$, so
$\mathfrak D_n=\int\epsilon_f\psi'-\sum_a(f(t_{a+1})-f(t_a))\tau_a$ with
$\tau_a=\frac{\psi(t_a)+\psi(t_{a+1})}{2}-\frac{1}{\delta}\int_{I_a}\psi=\frac{\delta^2}{12}\psi''(\xi_a)$,
some $\xi_a\in I_a$, for $a\le n-2$, where $\psi$ is $C^2$ on the closed cell. Indeed, integration by parts on $I_a$ gives $\int_{I_a}(\frac{1}{2}-\lambda)\psi'=-\tau_a$,
and $\delta\tau_a$ is the error of the trapezoidal rule for $\int_{I_a}\psi$. Use
$|\epsilon_f|\le\delta^{2H}$ on $I_0$ and $\le C\delta^2t_a^{2H-2}$ on $I_a$, $a\ge1$;
$|\tau_a|\le C\delta^2(1-t_{a+1})^{2H-2}$ for $a\le n-2$ and $\le C\delta^{2H}$ for $a=n-1$;
$|f(t_{a+1})-f(t_a)|\le C\delta t_a^{2H-1}$ for $a\ge1$; $0<\psi'\le C(1-t)^{2H-1}$; and
$\delta\sum_{a=1}^n(a\delta)^\alpha=O(1)$ for $\alpha>-1$, $O(\delta^{1+\alpha})$ for $\alpha<-1$ (the case
$\alpha=-1$ arises only at $H=\frac{1}{2}$, where $\mathfrak D_n\equiv0$). Splitting at $t=\frac{1}{2}$, both
parts are $O(\delta^{2H+1}+\delta^2)$.
\end{proof}

The same variance functions give the adjacent double-error term in closed form. The next corollary
shows that it is $o(\delta^{4H-1})$, so that the alternating and nested double errors cancel at
order $\delta^{4H-1}$. It is used in Remark~\ref{rem:T3hyp}.

  \begin{corollary}[the adjacent double error]\label{cor:T3adj}
Let $H\in(\tfrac12,1)$, let $e=\varrho^n-\varrho$ be the interpolation error of the fBm kernel
$\varrho(s,t)=c_H|t-s|^{2H-2}$, and let $\epsilon_f$ and $\lambda$ be as in
Proposition~\ref{prop:adjexact}. Define the double-error terms
\[
   T_3^{\rm adj}:=\int_{\mathcal T_4}e(t_1,t_2)e(t_3,t_4),\quad
   T_3^{\rm alt}:=\int_{\mathcal T_4}e(t_1,t_3)e(t_2,t_4),\quad
   T_3^{\rm nest}:=\int_{\mathcal T_4}e(t_1,t_4)e(t_2,t_3),
\]
where $\mathcal T_4=\{0<t_1<t_2<t_3<t_4<1\}$. Then
\[
   T_3^{\rm adj}=-\tfrac14\Bigl[\int_0^1\epsilon_f\,\epsilon_{\tilde f}'\,\dd t
   -2\delta^{2H-1}\int_0^1(1-2\lambda)\,\epsilon_f\,\dd t\Bigr]=O(\delta^{4H}+\delta^3),
   \qquad\tilde f(t)=(1-t)^{2H}.
\]
Moreover $\sum_\bullet T_3^\bullet=0$, so $T_3^{\rm alt}+T_3^{\rm nest}=-T_3^{\rm adj}=o(\delta^{4H-1})$.
\end{corollary}

\begin{proof}
For the adjacent pairing, the four-fold integral splits at the gap between $t_2$ and $t_3$, so
it is an integral of a product of two one-variable functions, the half variances of
Proposition~\ref{prop:adjexact}. We expand that product and bound the two surviving terms.

\emph{The formula.} For kernels $k_1,k_2$ the bilinear form
$\mathcal E(k_1,k_2)=\int_{\mathcal T_4}k_1(t_1,t_2)k_2(t_3,t_4)$
equals $-\int_0^1F_1G_2'$ with $F_1(u)=\int_{t_1<t_2<u}k_1$, $G_2(u)=\int_{u<t_3<t_4}k_2$
(integrate $F_1'(t_2)G_2(t_2)$ by parts; $F_1(0)=G_2(1)=0$). For
$\varrho$ and $\varrho^n$ these are half variances, $F=\frac{1}{2}f$, $F_n=\frac{1}{2}(I_nf-\lambda(1-\lambda)\delta^{2H})$,
$G=\frac{1}{2}\tilde f$, $G_n=\frac{1}{2}(I_n\tilde f-\lambda(1-\lambda)\delta^{2H})$, so
$T_3^{\rm adj}=\mathcal E(e,e)=-\int(F_n-F)(G_n-G)'$. Expanding, the $\delta^{4H-1}$ term integrates
$\lambda(1-\lambda)(1-2\lambda)$ to $0$. Integrating by parts on each cell,
$\int_{I_a}\lambda(1-\lambda)\epsilon_{\tilde f}'=-\delta^{-1}\int_{I_a}(1-2\lambda)\epsilon_{\tilde f}$, and the
reflection $t\mapsto1-t$ of the uniform mesh gives $\int(1-2\lambda)\epsilon_{\tilde f}=-\int(1-2\lambda)\epsilon_f$.
This is the formula.

\emph{The size.} On $I_a$, $a\ge1$, $\epsilon_f=+\frac{1}{2}f''(\xi_t)(t-t_a)(t_{a+1}-t)$
(linear interpolation overestimates the convex $f$), and replacing
$f''(\xi_t)$ by its midpoint value costs $\delta\sup_{I_a}|f'''|$ while
$\int_0^1\lambda(1-\lambda)(1-2\lambda)=0$, so $|\int(1-2\lambda)\epsilon_f|\le\delta^{2H+1}+C\delta^4\sum_a
(t_a^{2H-3}+1)=O(\delta^{2H+1})$. For the first integral, $|\epsilon_f|\le C\delta^2t_a^{2H-2}$ and
$|\epsilon_{\tilde f}'|\le\delta\sup_{I_a}|\tilde f''|$ on interior cells, and the end cells contribute
$O(\delta^{2H+2})$; as $2H-2>-1$ the total is $O(\delta^3)$. Finally $3>4H-1$, since $H<1$, and
$4H>4H-1$, so $T_3^{\rm adj}=o(\delta^{4H-1})$. The three pairings exhaust the interleavings of two
ordered pairs, so $\sum_\bullet T_3^\bullet=\tfrac18(\int_{[0,1]^2}e)^2=0$ by cell cancellation. This gives the last claim.
\end{proof}

\begin{remark}[the inputs of the double-error bound \eqref{eq:T3bound}, and sharpness]\label{rem:T3hyp}
The double-error bound uses cell cancellation and the symmetry $e(s,t)=e(t,s)$, the global bound
$\norm e_{L^1}\le L\delta^{2H-1}$ (the all-cells branch of \emph{(C1)} of
Definition~\ref{def:class} for $H\ge\tfrac12$), and the row-wise bound \eqref{eq:rowmass}. The
last is not a class hypothesis: (C1) controls only the sum over rows, and \emph{(C3)} gives only
$\int_{I_a\times[0,1]}|e|\le2L\delta$. From (C1), (C3) and cell cancellation one still gets
$|T_3^\bullet|\le4L^2\delta^{2H}$, which suffices for Theorem~\ref{thm:2H}(ii) but not for the
$o(\delta^{2H})$ used in Proposition~\ref{prop:C1exact}.

By Corollary~\ref{cor:T3adj}, $\sum_\bullet T_3^\bullet=0$ and
$T_3^{\rm alt}+T_3^{\rm nest}=o(\delta^{4H-1})$. Numerically, the alternating and nested double
errors are each of exact order $\delta^{4H-1}$ (with opposite leading coefficients, by the
corollary), so \eqref{eq:T3bound} is, numerically, sharp in order (see also Remark~\ref{rem:nonuniform}).
\end{remark}

\noindent\emph{The leading constant.} At $H=\tfrac12$, Theorem~\ref{thm:2H}(i) gives
$\cB(\tfrac12):=\lim_{\delta\to0}\norm{\cW_4}/\delta^{2H}=1/(4\sqrt3)$. For $H\in(\tfrac12,1)$ the
leading constant comes from an exact first-moment identity (Lemma~\ref{lem:momentid}), fed into a
single-error asymptotic (Lemma~\ref{lem:singleerror}). Assembling the words gives the three
constants exactly (Proposition~\ref{prop:C1exact}). Every word constant is a rational multiple of
the constant $\Phi_1=\tfrac16$ of Lemma~\ref{lem:momentid}, so no numerical input is needed.

\begin{lemma}[first-moment identity]\label{lem:momentid}
Let $r(k)=\tfrac12[(k+1)^{2H}-2k^{2H}+|k-1|^{2H}]$, $k\ge0$, be the
fGn autocovariance, so $r(0)=1$. Equivalently, the cell-averaged interpolant density is
$\varrho^n=\delta^{2H-2}r(|j-k|)$ on $I_j\times I_k$. Let
 $\widetilde r$ be its continuous piecewise-linear interpolant ($\widetilde r(k)=r(k)$), and $\varrho(u)=c_Hu^{2H-2}$,
$c_H=H(2H-1)$. The first diagonal moment of the averaged error is
\[
   \Phi_1:=\int_0^\infty\! u\,\bigl(\widetilde r(u)-\varrho(u)\bigr)\,\dd u=\tfrac16
   \qquad\text{for every }H\in(0,1),
\]
the integral converging absolutely.
\end{lemma}

For $H\ne\tfrac12$ neither of the two separate first moments, of $\widetilde r$ and of $\varrho$, converges. Only their
difference does. The restriction $H<1$ is needed: at $H=1$ one has $r\equiv1\equiv\varrho$, so
$\Phi_1=0$, and the constant is not continuous up to the endpoint.

\begin{proof}
The lemma is the case $\phi(u)=u^{2H}$, $\delta=1$ of Lemma~\ref{lem:momentgen} below, whose
proof does not use the present lemma. Hypothesis~(S) holds with $\beta=2H$, since
$|\phi^{(j)}(u)|$ is a constant multiple of $u^{2H-j}$. In this case $\hat r=r$,
$\hat\varrho=\tfrac12\phi''=c_Hu^{2H-2}=\varrho$ and $\hat r(0)=1$, so $\Phi_1=\hat\Phi_1=\tfrac16$,
and the integral converges absolutely.
\end{proof}

The level-$3$ drift bias of Proposition~\ref{prop:level3drift} involves the same first moments
for a stationary-increment driver that need not be self-similar. The next lemma extends
Lemma~\ref{lem:momentid} to such drivers, at a fixed mesh and in cell units, and supplies the lag
bound \eqref{eq:lagbound} used in Proposition~\ref{prop:level3drift} and
Remark~\ref{rem:singleerrorscope}.

\begin{lemma}[the first moment for a general stationary-increment driver]\label{lem:momentgen}
Let $G$ have stationary increments with $\phi(u)=\E[(G^i_{t+u}-G^i_t)^2]$ and assume

\smallskip
\noindent\emph{Hypothesis (S).} $\phi(0)=0$, $\phi\in C^4\bigl((0,\infty)\bigr)$, and there are $c_\phi$
and $\beta<2$ with, for $j=1,2,3,4$,
\[
   |\phi^{(j)}(u)|\le c_\phi\,u^{2H-j}\quad(0<u\le1),\qquad
   |\phi^{(j)}(u)|\le c_\phi\,u^{\beta-j}\quad(u\ge1).
\]

\smallskip
\noindent Fix a mesh width $\delta$, write $\hat\phi(u):=\phi(\delta u)$, $\hat r(k):=\tfrac12\bigl[\hat\phi(k+1)-2\hat\phi(k)+\hat\phi(|k-1|)\bigr]=\E[\Delta_a\Delta_{a+k}]$, $k\ge0$,
 for the cell-increment autocovariance ($\Delta_a$ the increment of one component of $G$ over the
cell $I_a$), so that
$\hat r(0)=\hat\phi(1)=\Var(\Delta_a)=\phi(\delta)$, let $\widetilde{\hat r}$ be the continuous
piecewise-linear interpolant of $\hat r$ and $\hat\varrho(u):=\tfrac12\hat\phi''(u)$ the
increment-covariance density in the same units. Then
\[
    \boxed{\ \hat\Phi_1\ :=\ \int_0^\infty\! u\,\bigl(\widetilde{\hat r}(u)-\hat\varrho(u)\bigr)\,\dd u
   \ =\ \tfrac16\,\hat r(0)\ =\ \tfrac16\,\phi(\delta)\ }
\]
for every $H\in(0,1)$, the integral converging absolutely. Let now $H\in(\tfrac12,1)$, let
$\hat e(u,v):=\hat r(m)-\hat\varrho(v-u)$ on the cell $C_m:=[0,1]\times[m,m+1]$ (intersected with
$\{u<v\}$ for $m=0$) be the cell-averaging error in cell units, and put
\[
   \hat\phi^{\rm s}(m):=\int_{C_m}\hat e(u,v)\,(u-\tfrac12)\,\dd u\,\dd v,\qquad
   \hat\phi^{\rm t}(m):=\int_{C_m}\hat e(u,v)\,(v-m-\tfrac12)\,\dd u\,\dd v .
\]
Then $\sum_m|\hat\phi^\bullet(m)|<\infty$, and the constants
$\hat\Phi^\bullet:=\sum_{m\ge0}\hat\phi^\bullet(m)$ satisfy
\[
   \hat\Phi^{\rm t}-\hat\Phi^{\rm s}=\hat\Phi_1=\tfrac16\phi(\delta),
   \qquad
   \hat\Phi^{\rm s}=-\hat\Phi^{\rm t}=-\tfrac1{12}\,\phi(\delta) .
\]
\end{lemma}

Since $G$ lives on $[0,1]$, its law determines $\phi(u)$ only for lags $u\in[0,1]$. The values of
$\phi$ on $(1,\infty)$ in Hypothesis~(S) are an analytic convenience: they enter only through cell
lags $m\ge n$, that is, time lags $m\delta\ge1$, which do not occur on $[0,1]$, and whose
contribution is bounded separately wherever the lemma is used. Any $C^4$ extension with the stated
decay may therefore be used, for instance one that is affine on $[2,\infty)$. For fBm, $\phi(u)=u^{2H}$ and
$\beta=2H$. For a sum $G^i=\sum_k\sqrt{a_k}\,B^{H_k}$ of independent fBms with $a_k>0$ and $2H\le2H_k<2$, used
in the numerical check of Proposition~\ref{prop:level3drift}, the incremental variance is
$\phi(u)=\sum_ka_ku^{2H_k}$, and (S) holds with $\beta=\max_k2H_k$: for $j\le4$ each term satisfies
$u^{2H_k-j}\le u^{2H-j}$ on $(0,1]$ and $u^{2H_k-j}\le u^{\beta-j}$ on $[1,\infty)$. In original units $\varrho^{\rm s}_1(a,b)=\delta\,\hat\phi^{\rm s}(b-a)$ and
$\varrho^{\rm t}_1(a,b)=\delta\,\hat\phi^{\rm t}(b-a)$ in the notation of Lemma~\ref{lem:reflect}.

\begin{proof}[Proof of Lemma~\ref{lem:momentgen}]
We compute $\int_0^Nu(\widetilde{\hat r}-\hat\varrho)$ by summation by parts and Taylor
expansions of $\hat\phi$ at $N$, and let $N\to\infty$. The divergent terms cancel. For the lag moments we use that the densities of $v-u$ on the cells $C_m$ are the hat
functions of linear interpolation, so the differences $\hat\phi^{\rm t}(m)-\hat\phi^{\rm s}(m)$ sum to
$\hat\Phi_1$; Lemma~\ref{lem:reflect} then separates the two constants.

\emph{The identity for $\hat\Phi_1$.} On $[k,k+1]$ the interpolant is $\widetilde{\hat r}(u)=\hat r(k)(k+1-u)+\hat r(k+1)(u-k)$, and the
substitution $u=k+s$ gives
\[
   \int_k^{k+1}u\,\widetilde{\hat r}(u)\,\dd u
   =\hat r(k)\bigl(\tfrac k2+\tfrac16\bigr)+\hat r(k+1)\bigl(\tfrac k2+\tfrac13\bigr).
\]
Sum over $k=0,\dots,N-1$. An interior node $1\le j\le N-1$ is the left end of $[j,j+1]$ and the
right end of $[j-1,j]$, so $\hat r(j)$ receives $\tfrac j2+\tfrac16$ from the first interval and
$\tfrac{j-1}2+\tfrac13=\tfrac j2-\tfrac16$ from the second, in total $j$. The end nodes receive
only one contribution, $\tfrac16$ at $j=0$ and $\tfrac N2-\tfrac16$ at $j=N$. Hence
\begin{equation}\label{eq:momgen1}
   \int_0^N\! u\,\widetilde{\hat r}(u)\,\dd u
   \ =\ \tfrac16\hat r(0)+\sum_{k=1}^{N-1}k\,\hat r(k)+\hat r(N)\bigl(\tfrac N2-\tfrac16\bigr).
\end{equation}
With $V(k):=\hat\phi(k+1)-\hat\phi(k)$ one has $2\hat r(k)=V(k)-V(k-1)$, and summing by parts
$\sum_{k=1}^{M}k\bigl(V(k)-V(k-1)\bigr)=M\,V(M)-\sum_{k=0}^{M-1}V(k)=M\,V(M)-\hat\phi(M)$, so
\begin{equation}\label{eq:momgen2}
   \sum_{k=1}^{N-1}k\,\hat r(k)\ =\ \tfrac12\bigl[(N-1)V(N-1)-\hat\phi(N-1)\bigr] .
\end{equation}
On the continuous side, $u\hat\phi'(u)\to0$ as $u\downarrow0$ by (S) with $j=1$, so
$\int_0^Nu\hat\varrho=\tfrac12\bigl[N\hat\phi'(N)-\hat\phi(N)\bigr]$. Expanding
$V(N-1)=\hat\phi'-\tfrac12\hat\phi''+\tfrac16\hat\phi'''+O(\hat\phi^{(4)})$ and
$\hat\phi(N-1)=\hat\phi-\hat\phi'+\tfrac12\hat\phi''-\tfrac16\hat\phi'''+O(\hat\phi^{(4)})$ at $N$,
with $O(\hat\phi^{(4)})$ standing for $O(\sup_{[N-1,N+1]}|\hat\phi^{(4)}|)$, the $\hat\phi'$,
$\hat\phi''$ and $\hat\phi'''$ terms without a factor $N$ cancel in pairs and \eqref{eq:momgen2}
becomes $\tfrac12[N\hat\phi'-\hat\phi]-\tfrac N4\hat\phi''+\tfrac N{12}\hat\phi'''+O(N\hat\phi^{(4)})$,
while $\hat r(N)(\tfrac N2-\tfrac16)=\tfrac N4\hat\phi''-\tfrac1{12}\hat\phi''+O(N\hat\phi^{(4)})$.
Substituting \eqref{eq:momgen2} into \eqref{eq:momgen1}, the terms $\pm\frac{N}{4}\hat\phi''$ cancel and
\[
   \int_0^N\! u\,\widetilde{\hat r}\,\dd u
   =\tfrac16\hat r(0)+\tfrac12\bigl[N\hat\phi'-\hat\phi\bigr]+\tfrac N{12}\hat\phi'''-\tfrac1{12}\hat\phi''
   +O\bigl(N\hat\phi^{(4)}\bigr).
\]
The divergent part $\frac{1}{2}[N\hat\phi'-\hat\phi]$ is exactly $\int_0^Nu\hat\varrho$, leaving
\begin{equation}\label{eq:momgen3}
   \int_0^N\! u\bigl(\widetilde{\hat r}-\hat\varrho\bigr)\,\dd u
   \ =\ \tfrac16\hat r(0)\ +\ \tfrac N{12}\hat\phi'''(N)-\tfrac1{12}\hat\phi''(N)
   \ +\ O\bigl(N\hat\phi^{(4)}(N)\bigr).
\end{equation}
Since $\hat\phi^{(j)}(u)=\delta^j\phi^{(j)}(\delta u)$, (S) at infinity bounds the three terms
after $\frac{1}{6}\hat r(0)$ by
$O(\delta^2(\delta N)^{\beta-2})$ once $\delta N\ge2$, and they vanish as $N\to\infty$ because
$\beta<2$. This gives $\hat\Phi_1=\frac{1}{6}\hat r(0)$. Absolute convergence: at infinity $\widetilde{\hat r}-\hat\varrho=O(\hat\phi^{(4)})$, both
from the second central difference and from the linear interpolation error, so
$u(\widetilde{\hat r}-\hat\varrho)=O(\delta^4u(\delta u)^{\beta-4})$ is integrable. At the origin
$u\hat\varrho(u)=O(\delta^{2H}u^{2H-1})$ is integrable for $H>0$.

\emph{The lag moments.} Let $H\in(\tfrac12,1)$, so that $\hat\varrho$ is integrable at $0$ and every
$\hat\phi^\bullet(m)$ is finite. For $m\ge2$, $\int_{C_m}\hat e=0$ and
$\int_{C_m}\hat r(m)(u-\tfrac12)=0$ give
$\hat\phi^{\rm s}(m)=-\int_{C_m}[\hat\varrho(v-u)-\hat\varrho(m)](u-\tfrac12)$, so
$|\hat\phi^{\rm s}(m)|\le\sup_{|x-m|\le1}|\hat\varrho'(x)|$, and likewise for $\hat\phi^{\rm t}$. By (S),
since $\hat\varrho'(x)=\frac{1}{2}\delta^3\phi'''(\delta x)$ and $x\ge m-1\ge\frac{m}{2}$ on the cell, for $m\ge2$,
\begin{equation}\label{eq:lagbound}
   |\hat\phi^\bullet(m)|\ \le\ \sup_{|x-m|\le1}|\hat\varrho'(x)|\ =\
   \begin{cases}O\bigl(\delta^{2H}m^{2H-3}\bigr), & m\delta\le\tfrac12,\\[2pt]
   O\bigl(\delta^3(m\delta)^{\beta-3}\bigr), & m\delta\ge\tfrac12,\end{cases}
\end{equation}
with constants depending only on $c_\phi$, $H$ and $\beta$. The right-hand side is summable in $m$
(as $2H-3<-1$ and $\beta<2$), and $\hat\phi^\bullet(0),\hat\phi^\bullet(1)$ are finite. So
$\sum_m|\hat\phi^\bullet(m)|<\infty$. Next, $(v-m-\tfrac12)-(u-\tfrac12)=(v-u)-m$ and
$\int_{C_m}\hat e=0$ give the per-cell identity
\[
   \hat\phi^{\rm t}(m)-\hat\phi^{\rm s}(m)=\int_{C_m}\hat e(u,v)\,(v-u)\,\dd u\,\dd v
   =\int_0^\infty x\,\bigl[\hat r(m)-\hat\varrho(x)\bigr]\,T_m(x)\,\dd x,
\]
where $T_m(x)=(1-|x-m|)_+$ for $m\ge1$ and $T_0(x)=(1-x)_+$ is the density of $v-u$ on $C_m$. As
$\sum_{m=0}^NT_m=1$ on $[0,N]$, $T_N(x)=N+1-x$ on $[N,N+1]$, and
$\sum_{m=0}^N\hat r(m)T_m=\widetilde{\hat r}$ on $[0,N]$ (the $T_m$ are the hat functions of linear
interpolation),
\[
   \sum_{m=0}^{N}\bigl(\hat\phi^{\rm t}(m)-\hat\phi^{\rm s}(m)\bigr)
   =\int_0^N x\bigl(\widetilde{\hat r}-\hat\varrho\bigr)\dd x
   +\int_N^{N+1}x\bigl[\hat r(N)-\hat\varrho(x)\bigr](N+1-x)\,\dd x .
\]
The last integral is at most $(N+1)\sup_{[N-1,N+1]}|\hat\phi'''|=O(\delta^2(\delta N)^{\beta-2})\to0$, and
the first tends to $\hat\Phi_1$. Hence $\hat\Phi^{\rm t}-\hat\Phi^{\rm s}=\hat\Phi_1=\tfrac16\phi(\delta)$.
Finally $\hat\phi^{\rm s}(m)=-\hat\phi^{\rm t}(m)$ for every $m$ by Lemma~\ref{lem:reflect}(i), applied
to the stationary kernel $\hat\varrho$, so $\hat\Phi^{\rm s}=-\hat\Phi^{\rm t}=-\tfrac1{12}\phi(\delta)$.
\end{proof}

Lemmas~\ref{lem:momentid} and~\ref{lem:momentgen} determine only the difference of the two
first-moment constants. Lemma~\ref{lem:singleerror} (hypothesis (K2)), Proposition~\ref{prop:level3drift}
and Proposition~\ref{prop:sep} need them separately. The next lemma supplies the missing relation
$\Phi^{\rm s}=-\Phi^{\rm t}$, by a reflection of each cell.

\begin{lemma}[reflection symmetry, the two moment constants are $\mp\tfrac1{12}$]\label{lem:reflect}
Let $e$ be either \emph{(i)} the cell-averaging error $\varrho^n-\varrho$ of a driver with stationary
increments, or \emph{(ii)} a symmetric cell-mean-zero kernel invariant under the cell
reflection $(u,v)\mapsto(1-v,1-u)$, such as the ripple density $b_H$ of
Proposition~\ref{prop:sep}. For $a\le b$, $m=b-a$, set
\[
   \varrho^{\rm s}_1(a,b)=\!\!\int_{(I_a\times I_b)\cap\{s<t\}}\!\!\!\! e(s,t)\,(s-c_a)\,\dd s\,\dd t,
   \qquad
   \varrho^{\rm t}_1(a,b)=\!\!\int_{(I_a\times I_b)\cap\{s<t\}}\!\!\!\! e(s,t)\,(t-c_b)\,\dd s\,\dd t,
\]
$c_a,c_b$ the cell centres and $I_a=[t_a,t_{a+1})$. Then $\varrho^{\rm s}_1(a,b)=-\varrho^{\rm t}_1(a,b)$ for every $a\le b$. Hence,
whenever $\varrho^\bullet_1(a,b)=\delta^{2H+1}\phi^\bullet(b-a)$ for some sequences $\phi^\bullet$
(the scaling form in (K2) of Lemma~\ref{lem:singleerror}), $\phi^{\rm s}(m)=-\phi^{\rm t}(m)$ for every $m$, and the constants
$\Phi^\bullet=\sum_m\phi^\bullet(m)$ satisfy
\begin{equation}\label{eq:reflectanti}
   \Phi^{\rm s}=-\Phi^{\rm t}.
 \end{equation}
For the fGn cell-averaging error with $H\in(\tfrac12,1)$ and $\Var(X_1)=1$,
$\Phi^{\rm t}-\Phi^{\rm s}=\Phi_1=\tfrac16$ (Lemma~\ref{lem:momentgen} with $\phi(u)=u^{2H}$ at $\delta=1$,
and Lemma~\ref{lem:momentid}), so $\Phi^{\rm s}=-\Phi^{\rm t}=-\tfrac12\Phi_1$, that is,
\[
   \boxed{\ \Phi^{\rm s}=-\Phi^{\rm t}=-\tfrac1{12}
   \qquad\text{(fGn cell-averaging error, }H\in(\tfrac12,1),\ \Var(X_1)=1\text{)}.\ }
\]
\end{lemma}

\begin{proof}[Proof of Lemma~\ref{lem:reflect}]
Rescale by $s=c_a+\delta u$, $t=c_b+\delta v$ with $u,v\in(-\tfrac12,\tfrac12)$, and let
$\iota(u,v):=(-v,-u)$. This involution maps $(-\tfrac12,\tfrac12)^2$ onto itself with unit
Jacobian, preserves the constraint $u<v$ (active only for $m=0$), and exchanges the first-moment
weights, $(s-c_a)\mapsto-(t-c_b)$. So we only need $\iota$-invariance of $e$, i.e.\
$e(s,t)=e(c_a+c_b-t,\,c_a+c_b-s)$ in the original variables, then applying $\iota$ to the integral
defining $\varrho^{\rm s}_1$ turns it into $-\varrho^{\rm t}_1$.

In case (i) the map $(s,t)\mapsto(c_a+c_b-t,\,c_a+c_b-s)$ sends the cell $I_a\times I_b$ onto
itself and preserves $t-s$. The cell average $\varrho^n$ is constant on the cell, and $\varrho$ depends
on $(s,t)$ only through $t-s$ by stationarity of the increments, so both are invariant, and so is $e$.
(The error $e$ itself is not a function of $t-s$, since $\varrho^n$ is constant on cells.) The
argument is cell by cell, so it applies equally to the cells $C_m$ of Lemma~\ref{lem:momentgen}. In
case (ii), the ripple density $b_H$ is not a function of $v-u$, but it is block diagonal (so only
$m=0$ occurs), and on the diagonal cell
$b_H(u,v)=c_H|u-v|^{2H-2}-f(u)-f(v)+1$ with $f(u)=f(1-u)$ (Proposition~\ref{prop:sep}), so
$b_H(u,v)=b_H(1-v,1-u)$, i.e.\ $\iota$-invariance. Under the scaling form, $\varrho^{\rm s}_1(a,b)=-\varrho^{\rm t}_1(a,b)$
gives $\phi^{\rm s}(m)=-\phi^{\rm t}(m)$, and summing  over $m$ gives
$\Phi^{\rm s}=-\Phi^{\rm t}$. With $\Phi^{\rm t}-\Phi^{\rm s}=\Phi_1=\tfrac16$ for the fGn error
(Lemma~\ref{lem:momentgen} with $\phi(u)=u^{2H}$, in cell units) this gives $2\Phi^{\rm
t}=\frac{1}{6}$, the boxed values.
\end{proof}

The proofs of Theorem~\ref{thm:2H}(ii) and of Proposition~\ref{prop:C1exact} reduce each
level-$4$ word to integrals $\int_{s<t}e\,w$ of the error kernel against a weight $w$. The next lemma
gives the six weights in closed form and checks the regularity that Lemma~\ref{lem:singleerror}
requires of them.

\begin{lemma}[the word weights, in closed form]\label{lem:weights}
Let $H\in(\tfrac12,1)$ and $\varrho(u,v)=c_H|u-v|^{2H-2}$, $c_H=H(2H-1)$. Consider the three
$(2,2)$-words $1122$, $1212$, $1221$ (adjacent, alternating, nested) with their Wick pairs of slots
$\{1,2\},\{3,4\}$, resp.\ $\{1,3\},\{2,4\}$, resp.\ $\{1,4\},\{2,3\}$. For one of the two pairs, $P$,
write $P'$ for the other one. The weight with the error on $P$ (``err $P$'') is
\[
   w(s,t):=\int\varrho(t_{P'})\,\dd t_{P'},\qquad 0\le s\le t\le1,
\]
the integral running over the two times of $P'$ such that, with $s<t$ placed at the slots of $P$, the
four times $t_1<t_2<t_3<t_4$ are increasing in $(0,1)$. The six weights are, on
$T:=\{0\le s\le t\le1\}$,
 \[
\begin{array}{ll@{\qquad}ll}
\text{adj, err }\{1,2\}: & w=\tfrac12(1-t)^{2H}, &
\text{adj, err }\{3,4\}: & w=\tfrac12 s^{2H},\\[2pt]
\text{nest, err }\{1,4\}: & w=\tfrac12(t-s)^{2H}, &
\text{alt, err }\{2,4\}: & w=\tfrac12\bigl[t^{2H}-(t-s)^{2H}-s^{2H}\bigr],\\[4pt]
\multicolumn{4}{l}{\text{alt, err }\{1,3\}:\quad w=\tfrac12\bigl[(1-s)^{2H}-(1-t)^{2H}-(t-s)^{2H}\bigr],}\\[2pt]
\multicolumn{4}{l}{\text{nest, err }\{2,3\}:\quad w=\tfrac12\bigl[1-(1-s)^{2H}-t^{2H}+(t-s)^{2H}\bigr].}
\end{array}
\]

Each $w$ is continuous on $T$ and $C^2$ in its interior, $\nabla w$ is bounded and
$(2H-1)$-H\"older on $T$, and there is a constant $C_w=C_w(H)$ with
\begin{equation}\label{eq:W2corner}
   |D^2w(s,t)|\ \le\ C_w\bigl[|t-s|^{2H-2}+s^{2H-2}+(1-t)^{2H-2}\bigr],\qquad 0<s<t<1.
\end{equation}
So each $w$ satisfies Hypothesis~\emph{(W)} of Lemma~\ref{lem:singleerror} below.
\end{lemma}
\begin{proof}
\emph{The closed forms.} Direct integration of $\varrho$, using
$\int_b^{c}(x-a)^{2H-2}\dd x=\tfrac1{2H-1}[(c-a)^{2H-1}-(b-a)^{2H-1}]$ twice. For example, for the
alternating word with the error on slots $\{1,3\}$ (at times $s=t_1$, $t=t_3$), the exact factor
sits on slots $\{2,4\}$, and the ordering $t_1<t_2<t_3<t_4$ gives $t_2\in(s,t)$ and $t_4\in(t,1)$.
Hence
\[
\begin{aligned}
   w(s,t)&=\int_s^t\!\dd t_2\int_t^1\varrho(t_2,t_4)\,\dd t_4
   =H\int_s^t\bigl[(1-t_2)^{2H-1}-(t-t_2)^{2H-1}\bigr]\dd t_2\\
   &=\tfrac12\bigl[(1-s)^{2H}-(1-t)^{2H}-(t-s)^{2H}\bigr].
\end{aligned}
\]

The other five weights follow in the same way.

\emph{Regularity, Hypothesis~(W).} Each $w$ is a finite combination of
$x\mapsto x^{2H}$ at $s,t,1-s,1-t,t-s$, hence continuous on $T$ and $C^2$ in its interior. Since
$x\mapsto2Hx^{2H-1}$ is bounded and $(2H-1)$-H\"older on $[0,1]$, $\nabla w$ is bounded and $(2H-1)$-H\"older on
the closed triangle, which gives (W1) and (W3). The second derivatives blow up on the diagonal and also at the corners: differentiating
$s^{2H}$ and $(1-t)^{2H}$ twice gives $s^{2H-2}$ and $(1-t)^{2H-2}$. So a bound by
$|t-s|^{2H-2}$ alone fails for the four weights containing one of these terms. The terms $t^{2H}$ and $(1-s)^{2H}$ cause no problem, since
$t\ge t-s$ and $1-s\ge t-s$ and $2H-2<0$. Hence all six weights satisfy \eqref{eq:W2corner}, which is (W2).
\end{proof}

The next lemma evaluates such integrals $\int_{s<t}e\,w$ to leading order: they equal $\delta^{2H}$
times the two first-moment constants of $e$ paired with the diagonal integrals of $\nabla w$. It gives
the leading constants in Proposition~\ref{prop:C1exact} and, with the ripple density in place of $e$,
in Proposition~\ref{prop:sep}.

\begin{lemma}[single-error asymptotic]\label{lem:singleerror}
Let $H\in(\tfrac12,1)$, $\delta=1/n$, and let $e$ be a symmetric kernel on $[0,1]^2$ with the following
two properties, for a constant $C_e$.

\smallskip
\noindent\emph{(K1)} $e$ is integrable on every cell product $I_a\times I_b$, with
$\int_{I_a\times I_b}e=0$.

\noindent\emph{(K2)} There are numbers $\phi^{\rm s}(m)=-\phi^{\rm t}(m)$ and $\psi(m)$, $m\ge0$
(they may depend on $\delta$), such that for all cells $a\le b$, with $\varrho^{\rm s}_1,\varrho^{\rm t}_1$
as in Lemma~\ref{lem:reflect},
\begin{gather*}
   \varrho^{\bullet}_1(a,b)=\delta^{2H+1}\phi^\bullet(b-a),\qquad
   \norm e_{L^1((I_a\times I_b)\cap\{s<t\})}\le\delta^{2H}\psi(b-a),\\
   |\phi^\bullet(m)|+\psi(m)\le C_e(1+m)^{2H-3}.
\end{gather*}

\smallskip
\noindent Let $w$ satisfy

\smallskip
\noindent\emph{Hypothesis (W).} $w$ is continuous on $T=\{0\le s\le t\le1\}$ and $C^2$ in its
interior, with a constant $C_w$ such that \emph{(W1)} $|\nabla w|\le C_w$ on $T$. \emph{(W2)}
$|D^2w(s,t)|\le C_w[|t-s|^{2H-2}+s^{2H-2}+(1-t)^{2H-2}]$ on the interior. \emph{(W3)} $\nabla w$
extends continuously to $T$ and is $(2H-1)$-H\"older there.

\smallskip
\noindent Then, with $\Phi^\bullet:=\sum_{m\ge0}\phi^\bullet(m)$, so that $\Phi^{\rm s}=-\Phi^{\rm t}$,
\[
   \int_{0<s<t<1}\!e(s,t)\,w(s,t)\,\dd s\,\dd t
   =\delta^{2H}\Bigl[\Phi^{\rm s}\!\int_0^1\!\partial_s w(x,x)\,\dd x
   +\Phi^{\rm t}\!\int_0^1\!\partial_t w(x,x)\,\dd x\Bigr]
   +O\bigl(\delta^{4H-1}+\delta^{2}\log\tfrac1\delta\bigr),
\]
the implied constant depending only on $C_w$, $C_e$ and $H$. For the fGn cell-averaging error at
$\Var(X_1)=1$, $\Phi^{\rm s}=-\Phi^{\rm t}=-\tfrac1{12}$, so the main term is
$\tfrac1{12}\delta^{2H}\int_0^1(\partial_t w-\partial_s w)(x,x)\,\dd x$.
\end{lemma}

\begin{remark}[the kernels to which the lemma is applied]
\label{rem:singleerrorscope}
\emph{(i) The fGn error} $e=\varrho^n-\varrho$, $\varrho=c_H|t-s|^{2H-2}$, used in
Proposition~\ref{prop:C1exact}. (K1) holds because $\varrho^n$ is the cell average of $\varrho$.
Stationarity of the increments and the exact scaling $\varrho(\delta\cdot,\delta\cdot)=\delta^{2H-2}\varrho$
give (K2) with $\phi^\bullet(m)$, $\psi(m)$ the lag moments and $L^1$ norms of the unit-mesh error,
independent of $\delta$. $\phi^{\rm s}(m)=-\phi^{\rm t}(m)$ by Lemma~\ref{lem:reflect}(i),
and $|\phi^\bullet(m)|\le\sup_{|x-m|\le1}|\varrho'(x)|=O(m^{2H-3})$ for $m\ge2$ by
\eqref{eq:lagbound} with $\phi(u)=u^{2H}$, $\delta=1$. Further $\psi(m)\le2\int_{C_m}|\varrho|<\infty$ for $m\le1$, while
$\psi(m)\le2\sup_{|x-m|\le1}|\varrho'(x)|$ for $m\ge2$. Lemma~\ref{lem:momentgen} with $\phi(u)=u^{2H}$
gives $\Phi^{\rm t}-\Phi^{\rm s}=\Phi_1=\tfrac16$, hence $\Phi^{\rm s}=-\Phi^{\rm t}=-\tfrac1{12}$.
\emph{(ii) The ripple density} $b_H$ of Proposition~\ref{prop:sep}, used there. It is block
diagonal, so (K2) holds with $\phi^\bullet(m)=\psi(m)=0$ for $m\ge1$. On the diagonal cell the
ripple's exact scaling gives $\varrho^\bullet_1(a,a)=\delta^{2H+1}\phi^\bullet(0)$ with
$\phi^{\rm s}(0)=\int_{0<u<v<1}b_H(u,v)(u-\tfrac12)$, $\phi^{\rm t}(0)=-\phi^{\rm s}(0)$ by
Lemma~\ref{lem:reflect}(ii), and $\psi(0)=\int_{[0,1]^2}|b_H|<\infty$. (K1) is checked in the proof of
Proposition~\ref{prop:sep}. \emph{(iii) Kernels of a stationary-increment driver under (S)}, which
are not self-similar, enter only through the affine weights of Proposition~\ref{prop:level3drift}.
There the argument simplifies and is given in full, with the constants of Lemma~\ref{lem:momentgen}.
The decay $(1+m)^{2H-3}$ in (K2) is needed in Step~4 below: summability alone would leave its second
remainder at $\delta^{2H}$, the size of the main term.
\end{remark}

\begin{proof}[Proof of Lemma~\ref{lem:singleerror}]
The idea is cell cancellation. The kernel $e$ has zero mass on every piece of the cell
decomposition, so on each piece only the first-order Taylor term of $w$ survives, and it produces the
first moments $\varrho^\bullet_1$ multiplied by $\nabla w$ at the centre. Summed over the pieces, these
products are Riemann sums for the diagonal integrals $\int_0^1\partial_\bullet w(x,x)\,\dd x$.

\emph{The cell decomposition.} Let $I_a=[a\delta,(a+1)\delta)$, $a=0,\dots,n-1$, with centres
$c_a=(a+\tfrac12)\delta$. If $s<t$ with $s\in I_a$ and $t\in I_b$, then $a\le b$. Hence, up to a
null set, the triangle $\{0<s<t<1\}$ is the disjoint union of the pieces
\[
   Q_{ab}:=(I_a\times I_b)\cap\{s<t\},\qquad 0\le a\le b\le n-1 .
\]
For $a<b$, $Q_{ab}$ is the whole square $I_a\times I_b$, since every point of $I_a$ lies to the left
of every point of $I_b$. For $a=b$, $Q_{aa}$ is the half of the diagonal square above the diagonal.
Grouping the pieces by their lag $m=b-a$, of which there are $n-m$ for each $m$,
\begin{equation}\label{eq:celldecomp}
   \int_{0<s<t<1}e\,w\ =\ \sum_{m=0}^{n-1}\ \sum_{a=0}^{n-1-m}\ \int_{Q_{a,a+m}}e\,w .
\end{equation}
Each piece carries zero mass of $e$. For $a<b$ this is (K1). For $a=b$, $e$ is symmetric and the
diagonal is null, so the two halves of $I_a\times I_a$ carry equal integrals, each equal to
$\tfrac12\int_{I_a\times I_a}e=0$.

\emph{Step 1 (expansion on one piece).} On $Q_{ab}$ expand $w$ to first order about the centre
$(c_a,c_b)$:
\[
   w(s,t)=w(c_a,c_b)+\partial_sw(c_a,c_b)\,(s-c_a)+\partial_tw(c_a,c_b)\,(t-c_b)+r_{ab}(s,t).
\]
Multiply by $e$ and integrate over $Q_{ab}$. The constant term contributes
$w(c_a,c_b)\int_{Q_{ab}}e=0$. The two linear terms give, by definition, the first moments
$\varrho^{\rm s}_1(a,b)$ and $\varrho^{\rm t}_1(a,b)$ of Lemma~\ref{lem:reflect}. Hence
\[
   \int_{Q_{ab}}e\,w=\partial_sw(c_a,c_b)\,\varrho^{\rm s}_1(a,b)+\partial_tw(c_a,c_b)\,\varrho^{\rm t}_1(a,b)
   +R_{ab},\qquad R_{ab}:=\int_{Q_{ab}}e\,r_{ab}.
\]
$Q_{ab}$ is convex, contains $(c_a,c_b)$ in its closure, and $\nabla w$ is continuous up to the
boundary by (W3). The mean value theorem therefore gives
$|r_{ab}(s,t)|\le\sqrt2\,\delta\,\mathrm{osc}_{Q_{ab}}(\nabla w)$, and so
\[
   |R_{ab}|\ \le\ \sqrt2\,\delta\ \mathrm{osc}_{Q_{ab}}(\nabla w)\ \norm{e}_{L^1(Q_{ab})} .
\]
The remainder must be bounded piece by piece, because by (W2) the second derivatives of the word
weights blow up on the diagonal (for the alternating word,
$\partial_s\partial_tw(s,t)=-\varrho(s,t)\to\infty$ as $t\downarrow s$). Inserting the expansion
into \eqref{eq:celldecomp} and using (K2), $\varrho^\bullet_1(a,a+m)=\delta^{2H+1}\phi^\bullet(m)$, gives
\begin{equation}\label{eq:lagsum}
   \int_{0<s<t<1}e\,w
   =\delta^{2H+1}\sum_{m=0}^{n-1}\sum_{a=0}^{n-1-m}
   \Bigl[\phi^{\rm s}(m)\,\partial_sw(c_a,c_{a+m})+\phi^{\rm t}(m)\,\partial_tw(c_a,c_{a+m})\Bigr]
   +\sum_{0\le a\le b\le n-1}R_{ab}.
\end{equation}

\emph{Step 2 (summability of the lag sequences).} By (K2), $\sum_m|\phi^\bullet(m)|<\infty$ and
$\sum_m\psi(m)<\infty$, as $2H-3<-1$, while $\sum_m m\,|\phi^\bullet(m)|$ may diverge; below we
therefore use only sums truncated at $m<n$ and bound the tail separately.

\emph{Step 3 (the Taylor remainders ).} We bound $\mathrm{osc}_{Q_{ab}}(\nabla w)$ in one of
two ways: on every piece by the H\"older bound (W3). Or, on a cell $I_a\times I_b$ whose closure lies in
the interior of $T$, by the mean value theorem and (W2), $\mathrm{osc}_{\rm cell}(\nabla w)\le\sqrt2\,\delta
\sup_{\rm cell}|D^2w|$. The first is used near the diagonal and near the edges $s=0$, $t=1$ of $T$, where
$D^2w$ may blow up, the second elsewhere.

\emph{(a) Pieces near the diagonal ($m\le1$).} A piece $Q_{ab}$ has diameter at most
$\sqrt2\,\delta\le2\delta$, so (W3) gives $\mathrm{osc}_{Q_{ab}}(\nabla w)\le C_w(2\delta)^{2H-1}$.
By (K2), $\norm e_{L^1(Q_{ab})}\le\delta^{2H}\psi(m)\le C_e\,\delta^{2H}$. Hence each such
piece has $|R_{ab}|\le\sqrt2\,\delta\cdot C_w(2\delta)^{2H-1}\cdot C_e\,\delta^{2H}\lesssim\delta^{4H}$.
There are $n$ pieces with $m=0$ and $n-1$ with $m=1$, so fewer than $2n=2/\delta$ in total, and
\[
   \sum_{0\le b-a\le1}|R_{ab}|\ \lesssim\ 2n\,\delta^{4H}\ =\ 2\,\delta^{4H-1}.
\]
\emph{(b) Edge cells ($m\ge2$, and $a\le1$ or $b\ge n-2$).} These lie in the rows $a=0,1$ or in the
columns $b=n-2,n-1$ (row $a$ is $I_a\times[0,1]$, column $b$ is $[0,1]\times I_b$). By (W3), as in (a),
$|R_{ab}|\lesssim\delta\cdot\delta^{2H-1}\cdot\delta^{2H}\psi(m)=\delta^{4H}\psi(m)$. Each of these four
rows and columns contains at most one cell of each lag $m$, and $\sum_{m\ge0}\psi(m)<\infty$, so
\[
   \sum_{\text{edge cells}}|R_{ab}|\ \lesssim\ 4\,\delta^{4H}\sum_{m\ge0}\psi(m)\ =\ O(\delta^{4H}).
\]

\emph{(c) Interior cells ($m\ge2$, $a\ge2$, $b\le n-3$).} On such a cell $t-s\ge(m-1)\delta$,
$s\ge a\delta$ and $1-t\ge(n-1-b)\delta\ge2\delta$, so its closure lies in the interior of $T$ and (W2) gives
\[
   \mathrm{osc}_{\rm cell}(\nabla w)\ \le\ \sqrt2\,\delta\,C_w\bigl[((m-1)\delta)^{2H-2}+(a\delta)^{2H-2}
   +((n-1-b)\delta)^{2H-2}\bigr],
\]
while $|R_{ab}|\le\sqrt2\,\delta\,\mathrm{osc}_{\rm cell}(\nabla w)\,\delta^{2H}\psi(m)$ by Step~1 and (K2).
We sum the three terms separately. \emph{The diagonal term.} There are at most $n$ cells of each lag $m$,
and $\psi(m)\le C_e(1+m)^{2H-3}$, so these contributions add up to at most a constant times
\[
   n\,\delta^{2}\,\delta^{2H-2}\,\delta^{2H}\sum_{m\ge2}(m-1)^{2H-2}m^{2H-3}\ =\ O(\delta^{4H-1}),
\]
the sum converging because $4H-5<-1$. \emph{The corner term $s^{2H-2}$.} Each row $a$ contains at
most one cell of each lag, so these contributions add up to at most a constant times
\[
   \delta^{2}\,\delta^{2H-2}\,\delta^{2H}\Bigl(\sum_{a=2}^{n-1}a^{2H-2}\Bigr)\Bigl(\sum_{m\ge0}\psi(m)\Bigr)
   \ =\ \delta^{4H}\,O(n^{2H-1})\ =\ O(\delta^{2H+1}),
\]
using $\sum_{a\le n}a^{2H-2}=O(n^{2H-1})$ for $H>\tfrac12$ and $n=1/\delta$. \emph{The corner term
$(1-t)^{2H-2}$} is the mirror image of the previous one under $(s,t)\mapsto(1-t,1-s)$, which maps
row $a$ to column $n-1-a$, and gives the same bound. As $\delta^{2H+1}=o(\delta^{4H-1})$ for $H<1$, the
corner terms are negligible.

\emph{Collecting (a)--(c),}
\[
   \sum_{0\le a\le b\le n-1}|R_{ab}|\ =\ O(\delta^{4H-1}).
\]

\emph{Step 4 (main term, and the tail $m\ge n$).} The remaining sum is
$\delta^{2H}\sum_{m=0}^{n-1}\bigl[\phi^{\rm s}(m)\,S^{\rm s}_m+\phi^{\rm t}(m)\,S^{\rm t}_m\bigr]$ with
$S^\bullet_m:=\delta\sum_{a\le n-1-m}\partial_\bullet w(c_a,c_a+m\delta)$. By (W3), $\partial_\bullet w$ is
$(2H-1)$-H\"older, so the midpoint rule has error $O(\delta^{2H-1})$, and the diagonal offset and the
$m$ missing cells cost $O((m\delta)^{2H-1})$:
\[
   S^\bullet_m=\int_0^1\partial_\bullet w(x,x)\,\dd x+\varepsilon^\bullet(m,\delta),\qquad
   |\varepsilon^\bullet(m,\delta)|\le C\bigl(\delta^{2H-1}+(m\delta)^{2H-1}\bigr).
\]
The induced error is
\[
   \delta^{2H}\!\!\sum_{m<n}\!|\phi^\bullet(m)|\,|\varepsilon^\bullet|
   \ \lesssim\ \delta^{4H-1}\!\sum_{m\ge0}\!|\phi^\bullet(m)|
   \ +\ \delta^{4H-1}\!\!\sum_{1\le m<n}\!m^{2H-3}m^{2H-1}
   \ =\ O\bigl(\delta^{4H-1}+\delta^{2}\log\tfrac1\delta\bigr),
\]
since $\sum_{m<n}m^{4H-4}$ is $O(1)$ for $H<\tfrac34$, $O(\log n)$ at $H=\tfrac34$, and
$O(n^{4H-3})=O(\delta^{3-4H})$ for $H>\tfrac34$, where $\delta^{4H-1}\delta^{3-4H}=\delta^2$. Finally,
replacing $\sum_{m<n}\phi^\bullet(m)$ by $\Phi^\bullet$ costs
$\delta^{2H}C_w\sum_{m\ge n}|\phi^\bullet(m)|\le C_wC_e\,\delta^{2H}\sum_{m\ge n}m^{2H-3}=O(\delta^{2H}n^{2H-2})
=O(\delta^2)$.
\emph{Conclusion.} Inserting Steps~3 and~4 into \eqref{eq:lagsum} gives the main term of the
claim plus the errors $O(\delta^{4H-1})$ (Step~3), $O(\delta^{4H-1}+\delta^2\log\tfrac1\delta)$
(Riemann sums, Step~4) and $O(\delta^2)$ (tail $m\ge n$), the last contained in the second. All
constants depend only on $C_w$, $C_e$ and $H$. For the fGn error,
Remark~\ref{rem:singleerrorscope}(i) gives $\Phi^{\rm s}=-\Phi^{\rm t}=-\tfrac1{12}$.
\end{proof}

\begin{remark}[the leading constant is a non-uniform limit as $H\downarrow\tfrac12$]\label{rem:nonuniform}
Since $4H-1>2H$ exactly when $H>\tfrac12$, the error term of Lemma~\ref{lem:singleerror} is
$o(\delta^{2H})$ on $(\tfrac12,1)$, but not uniformly as $H\downarrow\tfrac12$. Relative to the
$\delta^{2H}$ main term it has size
\[
   O\bigl(\delta^{2H-1}+\delta^{2-2H}\log\tfrac1\delta\bigr),
\]
and the second term dominates exactly when $H>\tfrac34$, so the governing relative exponent is
$\min(2H-1,\,2-2H)$. Numerically, the first term appears only for words with a $\delta^{4H-1}$
coefficient, that is, the alternating and the nested word. The adjacent word and the level-$3$ drift
words of Proposition~\ref{prop:level3drift}, which are linear in $e$, show only the second (the
values are in the companion repository). Two consequences follow. \emph{(i)} The constants below hold as
$\delta\to0$ at fixed $H>\tfrac12$, not uniformly on $(\tfrac12,1)$. Near the crossover no
accessible mesh exhibits them for the alternating and nested words, so a Richardson extrapolation
of numerical values must use the ladder $\{1,\delta^{2H-1},\delta^{2-2H},\dots\}$.
\emph{(ii)} The value at $H=\tfrac12$ is not the limit of the values on $(\tfrac12,1)$ and is
computed separately (Theorem~\ref{thm:2H}(i) for the bias, Proposition~\ref{prop:sepbrownian} for
the separation of the pair $(P_0,P_1)$, see Remark~\ref{rem:merger}).
\end{remark}

The next proposition combines Lemmas~\ref{lem:weights} and~\ref{lem:singleerror} into the exact
level-$4$ leading constants. It shows that the limit in Theorem~\ref{thm:2H}(ii) exists and is positive,
and it gives the value $\cB(H)$ of Main result~1.

\begin{proposition}[exact leading constants, centred independent-fBm driver]\label{prop:C1exact}
Let $H\in(\tfrac12,1)$, $d=2$, and let the driver be the centred independent-fBm law with
$\Var(G^i_1)=1$. Then the three level-$4$ leading constants $C_1^\bullet(H):=\lim_{\delta\to0}\delta^{-2H}\Delta E_\bullet$ are
\[
   C_1^{\rm adj}=-\tfrac1{12},\qquad C_1^{\rm alt}=\tfrac16,\qquad C_1^{\rm nest}=-\tfrac1{12},
\]
and hence
\[
   \boxed{\ \cB(H)=\lim_{\delta\to0}\frac{\norm{\cW_4}}{\delta^{2H}}
   =\Bigl(2\!\sum_\bullet (C_1^\bullet)^2\Bigr)^{1/2}=\sqrt{\tfrac1{12}}=\frac{1}{2\sqrt3}\qquad(H\in(\tfrac12,1),\ d=2).\ }
\]
\end{proposition}

The constants do not depend on $H$, and each is a rational multiple of
$\Phi_1=\Phi^{\rm t}-\Phi^{\rm s}=\tfrac16$. At $H=\tfrac12$, $\cB$ jumps to $1/(4\sqrt3)$
(Remark~\ref{rem:merger}).

\begin{proof}
We split each word-difference into two single-error terms and one double-error term, show that
the double-error term is of lower order, and evaluate the single-error terms by
Lemma~\ref{lem:singleerror} with the six weights of Lemma~\ref{lem:weights}.

\emph{Reduction to single-error integrals.} By the expansion \eqref{eq:eexpand}, with the time
pairs $t_A,t_B$ of the two Wick pairs in place of $(t_1,t_2),(t_3,t_4)$, $\Delta E_\bullet$ is the
sum of two single-error terms, in which one covariance factor is the error $e$ and the other is
exact, and one double-error term $\int e(t_A)e(t_B)$. The
double-error term is $O(\delta^{4H-1})$ for each of the three words, by \eqref{eq:T3bound}, and
$4H-1>2H$ for $H>\tfrac12$, so it does not contribute at order $\delta^{2H}$. In a single-error term,
integrating out the two times of the exact factor leaves
$\int_{0<s<t<1}e(s,t)\,w(s,t)\,\dd s\,\dd t$, where $s<t$ are the two times of the pair carrying the
error and $w$ is the weight of Lemma~\ref{lem:weights}. Lemma~\ref{lem:singleerror} therefore applies
(Remark~\ref{rem:singleerrorscope}(i)). Throughout,
$\int_0^1 H(1-x)^{2H-1}\dd x=\int_0^1Hx^{2H-1}\dd x=\tfrac12$.

\emph{The six diagonal integrals.} By Lemma~\ref{lem:singleerror}, with
$\Phi^{\rm s}=-\Phi^{\rm t}=-\tfrac1{12}$, it suffices to compute the two diagonal integrals
$\int_0^1\partial_sw(x,x)\,\dd x$ and $\int_0^1\partial_tw(x,x)\,\dd x$ for each of the six
weights.

\emph{Adjacent word, Wick pairs $\{1,2\}$ and $\{3,4\}$.} Error on slots $\{1,2\}$:
$w=\tfrac12(1-t)^{2H}$ depends only on the larger time $t$, so $\partial_sw=0$ and
$\int_0^1\partial_tw(x,x)\,\dd x=-\int_0^1H(1-x)^{2H-1}\dd x=-\tfrac12$; contribution
$-\tfrac12\Phi^{\rm t}=-\tfrac1{24}$. Error on slots $\{3,4\}$: $w=\tfrac12s^{2H}$ depends only on the
smaller time $s$, so $\partial_tw=0$ and $\int_0^1\partial_sw(x,x)\,\dd x=\tfrac12$; contribution
$\tfrac12\Phi^{\rm s}=-\tfrac1{24}$. Total
$C_1^{\rm adj}=\tfrac12(\Phi^{\rm s}-\Phi^{\rm t})=-\tfrac12\Phi_1=-\tfrac1{12}$.

\emph{Alternating word, Wick pairs $\{1,3\}$ and $\{2,4\}$.} Error on slots $\{1,3\}$:
$w=\tfrac12\bigl[(1-s)^{2H}-(1-t)^{2H}-(t-s)^{2H}\bigr]$ vanishes on the diagonal, so differentiating
$x\mapsto w(x,x)=0$ gives $\partial_sw(x,x)=-\partial_tw(x,x)$. Moreover
$\partial_tw(x,x)=H(1-x)^{2H-1}$, since the term $(t-s)^{2H-1}$ vanishes on the diagonal for
$H>\tfrac12$. Hence $\int_0^1\partial_tw=\tfrac12$ and $\int_0^1\partial_sw=-\tfrac12$; contribution
$\tfrac12(\Phi^{\rm t}-\Phi^{\rm s})=\tfrac12\Phi_1=\tfrac1{12}$. Error on slots $\{2,4\}$:
$w=\tfrac12\bigl[t^{2H}-(t-s)^{2H}-s^{2H}\bigr]$ also vanishes on the diagonal, and
$\partial_tw(x,x)=Hx^{2H-1}$; again the contribution is $\tfrac12\Phi_1=\tfrac1{12}$. Total
$C_1^{\rm alt}=\Phi_1=\tfrac16$.

\emph{Nested word, Wick pairs $\{1,4\}$ and $\{2,3\}$.} Error on the outer pair $\{1,4\}$:
$w=\tfrac12(t-s)^{2H}$ and its gradient vanish on the diagonal, since $2H>1$; contribution $0$.
Error on the inner pair $\{2,3\}$: $w=\tfrac12\bigl[1-(1-s)^{2H}-t^{2H}+(t-s)^{2H}\bigr]$ has
$\partial_sw(x,x)=H(1-x)^{2H-1}$ and $\partial_tw(x,x)=-Hx^{2H-1}$, with integrals $\tfrac12$ and
$-\tfrac12$; contribution $\tfrac12(\Phi^{\rm s}-\Phi^{\rm t})=-\tfrac12\Phi_1=-\tfrac1{12}$. Total
$C_1^{\rm nest}=-\tfrac12\Phi_1=-\tfrac1{12}$.

\emph{Conclusion.} For each word, $\Delta E_\bullet$ is $\delta^{2H}C_1^\bullet$ plus three errors:
$O(\delta^{4H-1})$ from the double-error term and $O(\delta^{4H-1}+\delta^2\log\frac{1}{\delta})$ from each of
the two single-error terms. All are $o(\delta^{2H})$, because $4H-1>2H$ for $H>\tfrac12$ and $2>2H$ for
$H<1$, so the limits $C_1^\bullet$ exist and take the values above. The adjacent value agrees with the
closed form $\Delta E_{\rm adj}=-\tfrac1{12}\delta^{2H}(1+o(1))$ of Proposition~\ref{prop:adjexact}.

The three constants sum to $0$, as forced by the shuffle identity
$\langle 11,\Sblk\rangle\langle22,\Sblk\rangle=\langle(11)\shuffle(22),\Sblk\rangle$ of
Lemma~\ref{lem:sumzero}, and
$\cB(H)=\sqrt{2(\tfrac1{144}+\tfrac1{36}+\tfrac1{144})}=\sqrt{\tfrac1{12}}=1/(2\sqrt3)$.
\end{proof}

\begin{remark}[the scope of the constants: $d$, the normalisation, and the driver]\label{rem:dscope}
\emph{(a) Dimension.} By \eqref{eq:W4norm}, over an alphabet of size $d$ and for
$H\in(\tfrac12,1)$, $\cB(H,d)=\sqrt{\binom d2}\big/(2\sqrt3)$. At $H=\tfrac12$ the value $\delta/(4\sqrt3)$ of
Theorem~\ref{thm:2H}(i) is multiplied by the same factor $\sqrt{\binom d2}$.
For the drifted independent-$\mathrm{fBm}$ law $P_\theta$ and general $d$ the surviving level-$3$
words are $1kk,k1k,kk1$, $k\ge2$, so
$\norm{\cW_3(P_\theta)}^2=(d-1)\theta^2\delta^{4H}/24\,(1+o(1))$ by Proposition~\ref{prop:level3drift}.
On the same law the level-$4$ bias changes too. By Proposition~\ref{prop:driftconst} and the remark
after it, it gains the terms
$\theta^2\delta^{2H}(-\tfrac1{24},\tfrac1{12},-\tfrac1{12},0,\tfrac1{12},-\tfrac1{24})$ on the words
$11kk,1k1k,1kk1,k11k,k1k1,kk11$, which add to the centred values
$\delta^{2H}(-\tfrac1{12},\tfrac16,-\tfrac1{12},-\tfrac1{12},\tfrac16,-\tfrac1{12})$ of
Proposition~\ref{prop:C1exact}. On the words $1k1k$ and $k1k1$ both terms are positive, so
$(\tfrac16+\tfrac{\theta^2}{12})^2\ge4\cdot\tfrac16\cdot\tfrac{\theta^2}{12}$ gives
$\norm{\cW_4(P_\theta)}^2\ge(d-1)\theta^2\delta^{4H}/9\,(1+o(1))$, which exceeds
$\norm{\cW_3(P_\theta)}^2$. Compared on one law, the level-$4$ bias is therefore always the larger,
for every $\theta$ and $d$.

\emph{(b) Normalisation.} $\Delta E_\bullet$ is quadratic in $\varrho$, so a driver with
$\Var(G^i_1)=\sigma^2$ scales the three word constants, and $\cB$, by $\sigma^4$. This is why
Definition~\ref{def:class} imposes $\Var(G^i_1)=1$, which also fixes $\sum_\bullet E^\infty_\bullet
=\tfrac18$ (Lemma~\ref{lem:sumzero}).

\emph{(c) Driver.} The constants are computed for the independent-component fBm:
Lemma~\ref{lem:weights} comes from integrating $\varrho(u,v)=c_H|u-v|^{2H-2}$ and
Lemma~\ref{lem:momentid} from the fractional-Gaussian-noise autocovariance. The second dependence
can be removed for stationary-increment kernels: Lemma~\ref{lem:momentgen} gives $\hat\Phi_1=\tfrac16\Var(G^i_\delta)$ for every such kernel, which
is what extends Proposition~\ref{prop:level3drift}. The first cannot, and neither survives for a
Gaussian driver with correlated components, which has extra Wick pairings
and cross-component kernels. Correlated fractional drivers are open
(\S\ref{sec:discussion}).
\end{remark}

The next lemma shows that the three word-differences $\Delta E_\bullet$ sum to zero. It checks
Proposition~\ref{prop:C1exact} and reduces the level-$4$ computations of \S\ref{sec:rough} to two of
the three words.

  \begin{lemma}[the endpoint identity]\label{lem:sumzero}
Let the components $X^1,X^2$ be independent and identically distributed, with integrable level-$4$
coordinates of $\Sig(X)$ and $\Sig(X^n)$. Then
$\sum_\bullet\Delta E_\bullet=0$ at every mesh.
\end{lemma}

The hypotheses hold for every centred member of the class: its components are i.i.d.\ by
Definition~\ref{def:class}, and the level-$4$ coordinates are integrable by Lemma~\ref{lem:moments}.
\begin{proof}
By \eqref{eq:shuffle-id},
$\langle 11,\Sblk\rangle\langle22,\Sblk\rangle=\langle(11)\shuffle(22),\Sblk\rangle$. The six words on
the right, $1122,1212,1221,2112,2121,2211$, each appear once. The swap of the letters
$1\!\leftrightarrow\!2$ preserves the law of $(X^1,X^2)$ and maps $2211,2121,2112$ to $1122,1212,1221$,
so the expectation of the right-hand side is $2(E_{\rm adj}+E_{\rm alt}+E_{\rm nest})=2\sum_\bullet E_\bullet$,
while by the letter-shuffles $(i)\shuffle(i)=2(ii)$ the left-hand side
equals $\tfrac14 (X^1_1)^2(X^2_1)^2$ and is determined by the endpoint. Hence, by independence,
$2\sum_\bullet E_\bullet=\tfrac14\E[(X^1_1)^2]\,\E[(X^2_1)^2]$ is the same for $X$ and $X^n$, which have
equal nodal endpoints. In particular $\sum_\bullet E^\infty_\bullet=\tfrac18$ under the
normalisation $\Var(X^i_1)=1$ of Definition~\ref{def:class}, in agreement with
Proposition~\ref{prop:altcontinuum}.
\end{proof}

\begin{remark}[the $H=\tfrac12$ merger]\label{rem:merger}
The two one-sided limits add to the exact Brownian coefficients:
\[
   \underbrace{\bigl(-\tfrac1{12},\ \tfrac16,\ -\tfrac1{12}\bigr)}_{\lim_{H\downarrow1/2}C_1^\bullet(H)}
   \ +\ \underbrace{\bigl(0,\ -\tfrac18,\ \tfrac18\bigr)}_{\lim_{H\uparrow1/2}(\kappa_{\rm adj},\kappa_{\rm alt},\kappa_{\rm nest})}
   \ =\ \underbrace{\bigl(-\tfrac1{12},\ \tfrac1{24},\ \tfrac1{24}\bigr)}_{\text{exact Brownian, Theorem~\ref{thm:2H}(i)}}.
\]
Consequently $\cB$ jumps from $1/(2\sqrt3)$ (the limit from above) to $1/(4\sqrt3)$ at $H=\tfrac12$.

The second summand is a limit as $H\uparrow\tfrac12$, not a value at $H=\tfrac12$: with
$\kappa_\bullet=\lim_{\delta\to0}\delta^{-(4H-1)}\Delta E_\bullet$, Theorem~\ref{thm:2H}(i) gives the
value $(-\tfrac1{12},+\tfrac1{24},+\tfrac1{24})$ at $H=\tfrac12$. Indeed the closed form
$T(k)=-\tfrac12[k^{2H}-(k-1)^{2H}]$ behind $\lim_{H\uparrow1/2}\kappa_{\rm alt}=-\tfrac18$ is valid
only for $H<\tfrac12$: $\sum_{m\ge1}r(m)=-\tfrac12$ for every $H<\tfrac12$ but $r(m)\equiv0$ for
$m\ge1$ at $H=\tfrac12$ (Remark~\ref{rem:kappahalf}). The outer entries add nothing new:
$\kappa_{\rm adj}=0$ (Proposition~\ref{prop:adjexact}) and $\kappa_{\rm nest}=-\kappa_{\rm alt}$
(Corollary~\ref{cor:kadj}).
\end{remark}

In the centred class the level-$3$ bias vanishes (Proposition~\ref{prop:level4}). The next
proposition shows that a drift restores it at the order of the level-$4$ bias, so that on $\cP$
the bias need not begin at level~$4$. It is used in Remark~\ref{rem:dscope}(a).

  \begin{proposition}[the level-3 drift bias]\label{prop:level3drift}
Let $d=2$ and let $P\in\cP(H,L)$ have drift $\theta\ne0$ in coordinate $1$ and independent centred
Gaussian components, $X^1=\theta t+G^1$, $X^2=G^2$. Assume either that $H\in(\tfrac12,1)$ and $G$ has
stationary increments whose increment variance $\phi(u)=\E[(G^i_{t+u}-G^i_t)^2]$ satisfies Hypothesis~\emph{(S)} of Lemma~\ref{lem:momentgen}, or that $H=\tfrac12$ and $G$ is Brownian motion
($\phi(u)=u$). Then the
only level-$3$ words with nonvanishing bias are $122,212,221$, and
\[
   \boxed{\ \cW_3=\theta\,\phi(\delta)\bigl(-\tfrac1{12},\ \tfrac16,\ -\tfrac1{12}\bigr)
   +O\bigl(\delta^2\bigr),\qquad
   \norm{\cW_3}_\hs=\frac{|\theta|}{2\sqrt6}\,\phi(\delta)\,(1+o(1)),\ }
\]
with no remainder in the Brownian case, and otherwise with $|O(\delta^2)|\le|\theta|\,C\delta^2$,
$C$ depending only on $H$, $c_\phi$ and $\beta$. The bias is nonzero, and the second identity follows from the
first, as soon as $\phi(\delta)/\delta^2\to\infty$, which holds for every $H<1$ when
$\phi(\delta)\asymp\delta^{2H}$.
\end{proposition}

Here $\phi(\delta)=\Var(G^i_\delta)$ is the cell-increment variance.

\begin{proof}[Proof of Proposition~\ref{prop:level3drift}]
By Wick's theorem only the words with a single letter $1$ carry a bias, and the bias of each is a
single-error integral $\theta\int e\,w_\ast$ against an affine weight $w_\ast$. For $H>\tfrac12$ we
evaluate these integrals by cell cancellation, which is exact up to the lag sums. At $H=\tfrac12$ we
compute them directly.

\emph{Which words survive.} Let $w$ have $k$ letters equal to $1$ and $3-k$ equal to $2$. The
$G^2$-slots must be Wick-paired among themselves ($G^2\perp G^1$, both centred), so $3-k$ is even,
i.e.\ $k\in\{1,3\}$. For $k=3$ the shuffle gives $\langle111,\Sblk\rangle=\tfrac1{3!}(X^1_1)^3$,
endpoint-determined, hence zero bias. So $k=1$: the words $122,212,221$. In each, the single
$1$-slot contributes the deterministic $\theta\,\dd t$ or $\dd G^1$. The latter leaves one
$G^1$-slot and two $G^2$-slots, odd in $G^1$, hence zero mean. Only the drift term survives:
\[
   \E\langle w,\Sblk\rangle=\theta\!\!\int_{0<t_1<t_2<t_3<1}\!\!\varrho(t_i,t_j)\,\dd t ,
\]
$\{i,j\}$ the two $2$-slots.

\emph{The constants.} The drift is piecewise linear, hence interpolated exactly, so
discretization replaces $\varrho$ by $\varrho^n$:
$\cW_w=\theta\int_{0<s<t<1}e(s,t)\,w_\ast(s,t)$, where $e:=\varrho^n-\varrho$ (the sign matches
$\cW=\E\Sig(X^n)-\E\Sig(X)$), with $w_\ast$ obtained by integrating out the
single free variable:
\[
   w_\ast=s\ (\text{word }122),\qquad w_\ast=t-s\ (\text{word }212),\qquad w_\ast=1-t\ (\text{word }221),
\]
each affine. As $w_\ast$ is affine, its first-order expansion about the centre $(c_a,c_b)$ of
each piece $Q_{ab}=(I_a\times I_b)\cap\{s<t\}$, $0\le a\le b\le n-1$, is exact. The constant term
integrates to zero against $e$ (cell cancellation; for $a=b$ by the symmetry of $e$), and the linear
terms give the first moments of Lemma~\ref{lem:reflect}, which by stationarity depend only on the lag
$m=b-a$: in original units $\varrho^\bullet_1(a,b)=\delta\hat\phi^\bullet(b-a)$
(Lemma~\ref{lem:momentgen}). There are $n-m$
pieces of lag $m$, and $(n-m)\delta=1-m\delta$, so, exactly,
\[
   \int_{0<s<t<1}e\,w_\ast=\sum_{m<n}(1-m\delta)\bigl[\hat\phi^{\rm s}(m)\,\partial_sw_\ast
   +\hat\phi^{\rm t}(m)\,\partial_tw_\ast\bigr].
\]
Replacing $\sum_{m<n}(1-m\delta)\hat\phi^\bullet(m)$ by $\hat\Phi^\bullet=\sum_{m\ge0}\hat\phi^\bullet(m)$
costs at most
 \[
   \Bigl|\sum_{m<n}\hat\phi^\bullet(m)\,m\delta\Bigr|+\sum_{m\ge n}|\hat\phi^\bullet(m)|
   \ \lesssim\ \delta^{2H+1}\!\!\sum_{1\le m<n}\!m^{2H-2}+\delta^2\!\int_{1/2}^\infty\!u^{\beta-3}\,\dd u
   \ =\ O\bigl(\delta^{2}\bigr),
\]
by the bounds $|\hat\phi^\bullet(m)|=O(\delta^{2H}m^{2H-3})$ for $1\le m\le n/2$ and
$O(\delta^3(m\delta)^{\beta-3})$ for $m\ge n/2$, and by $\sum_{m<n}m^{2H-2}=O(n^{2H-1})$ for
$H>\tfrac12$. For $m\ge2$ these are the bounds \eqref{eq:lagbound}. For $m=1$, $|\hat\phi^\bullet(1)|\le\tfrac12\int_{C_1}|\hat e|\le\int_{C_1}|\hat\varrho|=O(\delta^{2H})$
by (S) with $j=2$. The implied constants depend only on $H$, $c_\phi$ and $\beta$. This is the $\delta^{2-2H}$ relative error recorded
for these words in Remark~\ref{rem:nonuniform}. What survives is
\[
   \int_{0<s<t<1}\!e\,w_\ast\ =\ \hat\Phi^{\rm s}\,\partial_sw_\ast+\hat\Phi^{\rm t}\,\partial_tw_\ast
   \ +\ O\bigl(\delta^2\bigr),
\]
and Lemma~\ref{lem:momentgen} gives
$\hat\Phi^{\rm s}=-\hat\Phi^{\rm t}=-\tfrac1{12}\phi(\delta)$, whence
\[
   \cW_{122}=\theta\,\hat\Phi^{\rm s}=-\tfrac{\theta}{12}\phi(\delta),\quad
   \cW_{212}=\theta(\hat\Phi^{\rm t}-\hat\Phi^{\rm s})=\tfrac{\theta}{6}\phi(\delta),\quad
   \cW_{221}=-\theta\,\hat\Phi^{\rm t}=-\tfrac{\theta}{12}\phi(\delta),
\]
up to the stated remainder. Unlike at level~$4$, $\hat\Phi^{\rm s}$ and $\hat\Phi^{\rm t}$ appear
separately here, not only through their difference; this is why Lemma~\ref{lem:reflect} is needed.
The three values sum to zero, as forced by the shuffle identity
$\langle1,\Sblk\rangle\langle22,\Sblk\rangle=\langle(1)\shuffle(22),\Sblk\rangle=
\langle122+212+221,\Sblk\rangle$, whose left-hand side is determined by the endpoint. Finally
$\norm{\cW_3}_\hs=|\theta|\phi(\delta)\bigl(\tfrac1{144}+\tfrac1{36}+\tfrac1{144}\bigr)^{1/2}
=|\theta|\phi(\delta)/(2\sqrt6)$. For the $\mathrm{fBm}$ law $\phi(\delta)=\delta^{2H}$.

\emph{The Brownian case $H=\tfrac12$.} A level-$3$ word with one drift slot contains a single
covariance factor, so $\cW_3$ is linear in $e$ and, unlike $\cW_4$ (Remark~\ref{rem:merger}), has no
double-error term. The formula above therefore extends to $H=\tfrac12$ unchanged, as we now check.

For Brownian motion the covariance measure $\varrho(\dd s,\dd t)=\delta_{\{s=t\}}\, \dd s$ lives on
the diagonal. Its cell average is $\varrho^n=\delta^{-1}$ on each diagonal cell $I_a\times I_a$ and
$0$ elsewhere. We integrate over the ordered simplex in the sense \eqref{eq:coinc}, so $\varrho$
contributes half of its diagonal mass, $\int_{\{s\le t\}}\varrho\,w=\tfrac12\int_0^1w(x,x)\,\dd x$,
while the bounded density $\varrho^n$ gives
$\delta^{-1}\sum_a\int_{(I_a\times I_a)\cap\{s<t\}}w$. For example, for $w=s$ each cell gives
$\tfrac{a\delta^2}2+\tfrac{\delta^2}6$, and summing over $a$ gives
$\tfrac{\delta^2}2\binom n2+\tfrac{n\delta^2}6=\tfrac14-\tfrac\delta{12}$. For the three affine
weights we obtain, exactly and for every $n$,
\[
\begin{array}{lccc}
   \text{word} & \int\varrho^n w_\ast & \int\varrho\,w_\ast & \int e\,w_\ast\\[2pt]
   122\ (w_\ast=s)   & \tfrac14-\tfrac{\delta}{12} & \tfrac14 & -\tfrac{\delta}{12}\\[2pt]
   212\ (w_\ast=t-s) & \tfrac{\delta}{6}           & 0        & \tfrac{\delta}{6}\\[2pt]
   221\ (w_\ast=1-t) & \tfrac14-\tfrac{\delta}{12} & \tfrac14 & -\tfrac{\delta}{12}
\end{array}
\]
Since $\phi(\delta)=\delta$, the last column is $(-\tfrac1{12},\tfrac16,-\tfrac1{12})\,\phi(\delta)$,
so the formula for $\cW_3$ holds at $H=\tfrac12$ with no remainder. The convention matters: with
the strict ordering $\{s<t\}$ the diagonal mass of $\varrho$ would be lost, the middle column would
vanish, and the first and last entries of the error column would be $\tfrac14-\tfrac\delta{12}$,
of order~$1$.
\end{proof}

\noindent\emph{A numerical check.} Exact evaluations of $\int_{0<s<t<1}e\,w_\ast$ on the mesh,
for $H=0.7$, $\phi(u)=u^{1.4}$ and two sums $\phi(u)=\sum_ka_ku^{2H_k}$ as after
Lemma~\ref{lem:momentgen}, approach
$\phi(\delta)\bigl(-\tfrac1{12},\tfrac16,-\tfrac1{12}\bigr)$ at the predicted relative rate
$\delta^{2-2H}$, and not the same multiples of $\delta^{2H}$. The values are reproduced by the code
in the companion repository.

\smallskip\noindent\textbf{The role of centring and independent components.} In the centred
class $\cP_0$ the bias first appears at level~4 (Proposition~\ref{prop:level4}), and
$\cB(H)=1/(2\sqrt3)$ (independent fBm, $d=2$) is a statement about $\cP_0$. Over the full
class $\cP$ we prove only the upper bound
$\bigl(\norm{\cW_3}^2+\norm{\cW_4}^2+\cdots\bigr)^{1/2}=O(\delta^{2H})$
(Proposition~\ref{prop:alllevelsdrift}). A $\theta$-dependent leading constant is only computed for
stationary-increment members at level~$3$ (Proposition~\ref{prop:level3drift}) and for drifted fBm at
level~$4$ (Proposition~\ref{prop:driftconst}). Even on $\cP_0$, $\cB(H)$ is the leading constant of
$\norm{\cW_4}$. For $M\ge6$ the constant of $\norm{\cW_{\le M}}$ also collects the even levels
$6,\dots,M$, each $O(\delta^{2H})$ (Proposition~\ref{prop:alllevels}), whose constants we do not
evaluate. The constants are fixed-$H$ limits (Remark~\ref{rem:nonuniform}).

\subsection{All truncation levels}\label{subsec:alllevels}

The weak error of the truncated signature is $\cW_{\le M}=\sum_{m\le M}\cW_m$, so the rate
$\delta^{2H}$ needs a bound at every level $m\le M$. The next proposition gives it on the centred
class, and with it the upper bound of Main result~1 there for $H\in[\tfrac12,1)$.

\begin{proposition}[full-truncation weak error]\label{prop:alllevels}
Let $H\in[\tfrac12,1)$, let $P\in\cP_0(H,L)$, let $X^n$ be its piecewise-linear
interpolant on the uniform mesh, and let $M$ be fixed. Then $\cW_m=0$ for
$m\le3$ (Proposition~\ref{prop:level4}), $\cW_m=0$ for all odd $m$, and $\norm{\cW_m}=O(\delta^{2H})$
for every even $m$ with $4\le m\le M$. Hence $\norm{\cW_{\le M}}_\hs=O(\delta^{2H})$, uniformly over
$\cP_0(H,L)$, with an implied constant depending only on $M$, $d$ and $L$.
\end{proposition}

\begin{proof}
At an even level we compare the interpolants on the
meshes $1/n$ and $1/(kn)$ and let $k\to\infty$. Isserlis' theorem and a telescoping sum write the
difference as a finite sum of single-error integrals $\int e\,w$. Cell cancellation against a weight $w$
that is Lipschitz on each piece of the cell decomposition then gains the factor $\delta$ that turns
$\norm e_{L^1}=O(\delta^{2H-1})$ into $O(\delta^{2H})$.

\emph{Odd levels.} $X^n$ is centred Gaussian, being a linear image of the Gaussian nodal vector, and the law of $G$ is
invariant under $G\mapsto-G$, which commutes with interpolation and with the natural lift. So
odd-length coordinates have zero mean for $X$ and for $X^n$, and $\cW_m=0$ for odd $m$. Fix an even
$m=2\ell$ with $4\le m\le M$ and a word $w$ of length $m$.

\emph{Refinement.} Fix $k\ge1$ and put $n':=kn$. The nodes are nested, so the interpolant of $X^{n'}$
on the mesh $1/n$ is $X^n$. By Lemma~\ref{lem:L1conv},
$\langle w,\cW_{\le M}\rangle=\lim_k\bigl(\E\langle w,\Sig(X^n)\rangle-\E\langle w,\Sig(X^{n'})\rangle\bigr)$,
so it suffices to bound the bracket by $C\delta^{2H}$ uniformly in $k$. Both paths are piecewise
linear and centred Gaussian, with bounded step densities as increment covariances: $\varrho^{n'}$,
constant on the cell products of the mesh $1/n'$, and $\varrho^n$, constant on those of the mesh
$1/n$, where it is the cell average of $\varrho$, hence also of $\varrho^{n'}$. Isserlis' theorem for the
finite Gaussian vector of cell increments gives
\[
   \E\langle w,\Sig(X^{n'})\rangle=\sum_{\mathfrak p}\ \int_{0<t_1<\cdots<t_{2\ell}<1}\
   \prod_{\{l,l'\}\in\mathfrak p}\varrho^{n'}(t_l,t_{l'})\ \dd t,
\]
summed over the finitely many pairings $\mathfrak p$  of the $2\ell$ slots compatible with $w$
(at most $(2\ell-1)!!$ of them), and the same
with $\varrho^n$ for $X^n$. No diagonal mass and no coincidence convention enter, also at
$H=\tfrac12$.

\emph{Single-error integrals.} Fix a pairing
$\mathfrak p=\{\{l_1,l_1'\},\dots,\{l_\ell,l_\ell'\}\}$, with $l_j<l_j'$, and write
$\varrho^n_j:=\varrho^n(t_{l_j},t_{l_j'})$ and $\varrho^{n'}_j:=\varrho^{n'}(t_{l_j},t_{l_j'})$. Telescoping the
product difference, which replaces one factor at a time,
\[
   \prod_{j=1}^{\ell}\varrho^n_j-\prod_{j=1}^{\ell}\varrho^{n'}_j
   =\sum_{r=1}^{\ell}\Bigl(\prod_{j<r}\varrho^n_j\Bigr)\,e_r\,\Bigl(\prod_{j>r}\varrho^{n'}_j\Bigr),
   \qquad e_r:=(\varrho^n-\varrho^{n'})(t_{l_r},t_{l_r'}),
\]
writes the bracket as a finite sum, over the pairings $\mathfrak p$ and $r=1,\dots,\ell$, of
single-error integrals
\[
   I_{\mathfrak p,r}=\int_{\rm simplex}e(t_{l_r},t_{l_r'})\prod_{j\ne r}(\text{kernel})_j\,\dd t,\qquad
   e:=\varrho^n-\varrho^{n'},
\]
each kernel a $\varrho^{n'}$ or a $\varrho^n$. Since they are paired with each other, the variables
$s:=t_{l_r}$, $t:=t_{l_r'}$ occur in no other kernel and enter only through the ordering
constraints. Integrating out the other $2\ell-2$ variables gives
\begin{gather*}
   I_{\mathfrak p,r}=\int_{0<s<t<1}e(s,t)\,w_{\mathfrak p,r}(s,t)\,\dd s\,\dd t,\\
   w_{\mathfrak p,r}(s,t)=\int_{\substack{\text{other }2\ell-2\text{ vars,}\\ \text{ordered around }s,t}}\
   \prod_{j\ne r}(\text{kernel})_j ,
\end{gather*}
No positivity
of the kernels is used: the bound below only needs the total-variation hypothesis (C3), which also
covers the sign-changing kernel of Construction~\ref{con:bridge}.

\emph{Uniform bound on the partial derivatives of the weight.} We claim that
$|\partial_sw_{\mathfrak p,r}|$ and $|\partial_tw_{\mathfrak p,r}|$ are at most $\Lambda:=2L^{\ell-1}$,
uniformly in $\delta$ and $k$. Since $(s,t)$ enter $w_{\mathfrak p,r}$ only as integration limits,
$\partial_s w_{\mathfrak p,r}$ is a sum of at most two boundary terms, each anchoring at $s$
a variable adjacent to $s$ in the ordering and integrating the rest (no term arises when the
neighbour of $s$ is $t$, $0$ or $1$). The anchored factor contributes
$\int_0^1|\varrho^{n'}(s,v)|\dd v$ or $\int_0^1|\varrho^n(s,v)|\dd v$. Each is at most $L$, for every
$s$: for $s$ in a
cell $I$ of the relevant mesh, it equals $|I|^{-1}\sum_J|\varrho(I\times J)|\le
|I|^{-1}\int_I|\varrho|(u,[0,1])\,\dd u\le L$ by the disintegration of $\varrho$ and (C3). The
remaining $\ell-2$ kernels sit on disjoint pairs of variables, so their absolute values integrate to
at most $\norm\varrho_{\rm TV}^{\ell-2}\le L^{\ell-2}$, again by (C3) (cell averaging does not
increase the total variation). So each boundary term is at most $L^{\ell-1}$, and
$|\partial_s w_{\mathfrak p,r}|\le2L^{\ell-1}=\Lambda$. The same holds for $\partial_t$.

\emph{From the derivative bound to a Lipschitz bound.} The step kernels make $w_{\mathfrak p,r}$
only piecewise $C^1$, with kinks along mesh lines, so the bound on $\nabla w_{\mathfrak p,r}$ holds only
almost everywhere, and an almost-everywhere bound alone does not control oscillations. Cell
cancellation needs a Lipschitz bound, and we obtain it from absolute continuity in each variable.
For fixed $t$ and $s<s'$, Fubini's theorem in the variables adjacent to $s$ gives
$w_{\mathfrak p,r}(s',t)-w_{\mathfrak p,r}(s,t)=\int_s^{s'}B_t(\sigma)\,\dd\sigma$, where $B_t(\sigma)$ is the
sum of the boundary terms at $\sigma$, and $|B_t(\sigma)|\le\Lambda$ for every $\sigma$ by the bounds
above. So $s\mapsto w_{\mathfrak p,r}(s,t)$ is absolutely continuous and $\Lambda$-Lipschitz, and
likewise $t\mapsto w_{\mathfrak p,r}(s,t)$. Now let $Q_{ab}=(I_a\times I_b)\cap\{s<t\}$ be the pieces of
\eqref{eq:celldecomp}. Any two points $(s,t),(s',t')$ of $Q_{ab}$ are joined inside
$Q_{ab}$ by one horizontal and one vertical segment: for $a<b$ because $Q_{ab}$ is a product; for $a=b$
go horizontally first if $s'<t$, and vertically first otherwise, in which case $s<t\le s'<t'$.
Along the two segments the one-variable bounds give
\begin{equation}\label{eq:pieceLip}
   |w_{\mathfrak p,r}(s,t)-w_{\mathfrak p,r}(s',t')|\le\Lambda\bigl(|s-s'|+|t-t'|\bigr),
   \qquad\text{hence}\qquad \mathrm{osc}_{Q_{ab}}(w_{\mathfrak p,r})\le2\delta\Lambda .
\end{equation}

\emph{Cell cancellation.} $\int_{Q_{ab}}e=0$ on every piece: for $a<b$ because
$\varrho^n$ and $\varrho^{n'}$ have the same integral $\varrho(I_a\times I_b)$ there, and for $a=b$
because $e$ is symmetric and bounded, so that the two halves of $I_a\times I_a$ carry equal integrals,
each $\tfrac12\int_{I_a\times I_a}e=0$ (cells are taken half-open). Subtracting from
$w_{\mathfrak p,r}$ its mean $\bar w_{ab}$ on $Q_{ab}$ therefore changes nothing, and what remains is at
most the oscillation \eqref{eq:pieceLip}, which gains the factor $\delta$:
\[
   |I_{\mathfrak p,r}|=\Bigl|\sum_{a\le b}\int_{Q_{ab}}e\,(w_{\mathfrak p,r}-\bar w_{ab})\Bigr|
   \le\sum_{a\le b}\norm e_{L^1(Q_{ab})}\,\mathrm{osc}_{Q_{ab}}(w_{\mathfrak p,r})
   \le2\delta\Lambda\,\norm e_{L^1}.
\]
By (C1) at the meshes $1/n$ and $1/n'$, with the total-variation norm at $H=\tfrac12$,
$\norm e_{L^1}\le\norm{\varrho^n-\varrho}+\norm{\varrho^{n'}-\varrho}\le L\delta^{2H-1}+L(\delta/k)^{2H-1}
\le2L\delta^{2H-1}$, as $2H-1\ge0$. So $|I_{\mathfrak p,r}|\le2\delta\cdot2L^{\ell-1}\cdot2L\delta^{2H-1}
=8L^{\ell}\,\delta^{2H}$ uniformly in $k$. For the
independent-fBm the $L^1$ bound is \eqref{eq:eL1fbm}.

\emph{Conclusion.} Summing over
the finitely many pairings $\mathfrak p$, indices $r$ and words $w$ at each even level and letting
$k\to\infty$ gives $\norm{\cW_m}=O(\delta^{2H})$, with a constant depending only on $M$, $d$ and
$L$, and
$\norm{\cW_{\le M}}_\hs^2=\sum_{m=4}^M\norm{\cW_m}^2=O(\delta^{4H})$, since the levels $m\le3$
and the odd levels contribute nothing.
\end{proof}

The next proposition extends the upper bound of
Proposition~\ref{prop:alllevels} to the drifted class, where level~$3$ need not vanish
(Proposition~\ref{prop:level3drift}). Together they give Main result~1 on all of
$\cP(H,L)$ for $H\in[\tfrac12,1)$.

\begin{proposition}[full-truncation weak error on the drifted class]\label{prop:alllevelsdrift}
Let $H\in[\tfrac12,1)$, let $P\in\cP(H,L)$ with drift $\theta$, $|\theta|\le L$, and let $M$ be
fixed. Then $\norm{\cW_{\le M}}_\hs=O(\delta^{2H})$, uniformly over $\cP(H,L)$. The implied constant
is $C(M,d)\max(1,L)^M$, with $C(M,d)$ depending only on $M$ and $d$.
\end{proposition}

\begin{proof}
The proof follows that of Proposition~\ref{prop:alllevels}. What changes is that each slot with
letter $1$ is expanded into its drift and its Gaussian part, so every term carries a set $D$ of drift
slots, and the drift times enter the weights with the factor $1$ in place of a kernel. Only the
derivative bound on the weights has to be redone.

\emph{Reduction to interpolants.} Write $X=\theta te_1+G$. The drift is linear, hence
interpolated exactly, so $X^{n'}=\theta te_1+G^{n'}$ for every $n'$. Fix a word $w=(w_1\cdots w_m)$
with $m\le M$, fix $k\ge1$ and put $n':=kn$. As in the proof of Proposition~\ref{prop:alllevels}, by
Lemma~\ref{lem:L1conv} it suffices to bound the bracket
$\E\langle w,\Sig(X^n)\rangle-\E\langle w,\Sig(X^{n'})\rangle$ by $C\delta^{2H}$, uniformly in $k$.

\emph{Expansion in the drift.} For a piecewise-linear path, $\langle w,\Sig\rangle$ is an iterated
Riemann--Stieltjes integral, linear in each slot, and
$\dd X^{n',i}_t=\theta\,\one_{\{i=1\}}\dd t+\dd G^{n',i}_t$. Expanding every slot, a slot $l$ with
$w_l=1$ carries either the drift $\theta\,\dd t_l$ or the Gaussian increment $\dd G^{n',1}_{t_l}$, and
every other slot carries $\dd G^{n',w_l}_{t_l}$. With $D$ the set of drift slots,
\[
   \langle w,\Sig(X^{n'})\rangle=\sum_{D\subseteq\{l:\,w_l=1\}}\theta^{|D|}
   \int_{0<t_1<\cdots<t_m<1}\ \prod_{l\in D}\dd t_l\ \prod_{l\notin D}\dd G^{n',w_l}_{t_l}.
\]
For fixed times, the integrand of each term is a polynomial of degree $m-|D|$ in the finite centred
Gaussian vector of cell increments of $G^{n'}$. If $m-|D|$ is odd, the term has mean zero. If
$m-|D|=2\ell$ is even, Fubini and Isserlis' theorem give a sum over the pairings $\mathfrak p$  of the
Gaussian slots. Only pairs of slots with the same letter survive, because the components are
independent, and every surviving pair carries the same kernel, because they share one scalar
increment covariance:
\[
   \E\langle w,\Sig(X^{n'})\rangle=\sum_{D:\ m-|D|\ \text{even}}\theta^{|D|}\sum_{\mathfrak p} \
   \int_{0<t_1<\cdots<t_m<1}\ \prod_{\{l,l'\}\in\mathfrak p} \varrho^{n'}(t_l,t_{l'})\,\dd t ,
\]
the inner sum over the pairings of the $2\ell=m-|D|$ Gaussian slots that pair equal letters, and the
same with $\varrho^n$ for $X^n$. A drift slot contributes the factor $1$ against Lebesgue measure. The
term $D=\emptyset$ is the centred expression of Proposition~\ref{prop:alllevels}. If every slot is a
drift slot ($\ell=0$, possible only for $w=1^m$), the term is $\theta^m/m!$ on both meshes and
cancels in the bracket.

\emph{Single-error integrals.} Fix $D$ and a pairing
$\mathfrak p=\{\{l_1,l_1'\},\dots,\{l_\ell,l_\ell'\}\}$ of the Gaussian slots, $l_j<l_j'$, with $\ell\ge1$. Telescoping the difference of the two
kernel products as in the proof of Proposition~\ref{prop:alllevels}, with $e:=\varrho^n-\varrho^{n'}$,
writes the contribution of $(D,\mathfrak p)$ to the bracket as
$\theta^{|D|}\sum_{r=1}^{\ell}I_{D,\mathfrak p,r}$, where
\[
   I_{D,\mathfrak p,r}=\int_{0<s<t<1}e(s,t)\,w_{D,\mathfrak p,r}(s,t)\,\dd s\,\dd t .
\]
Here $s=t_{l_r}$ and $t=t_{l_r'}$ are the two times of the pair carrying the error, and
$w_{D,\mathfrak p,r}(s,t)$ is the integral over the other $m-2$ times, ordered around $s$ and $t$, of the
remaining $\ell-1$ kernels, each a $\varrho^n$ or a $\varrho^{n'}$. The $|D|$ drift times enter with
weight $1$.

\emph{The derivative bound.} As in the proof of Proposition~\ref{prop:alllevels},
$w_{D,\mathfrak p,r}$ is absolutely continuous in each variable, and $\partial_sw_{D,\mathfrak p,r}$ is
a sum of at most two boundary terms, in each of which the time next to $s$ in the ordering is set
equal to $s$. A Gaussian neighbour contributes at most $L$ by (C3), as there, and a drift
neighbour, anchored against Lebesgue measure, contributes $1$. Every remaining kernel integrates in
absolute value to at most $\norm\varrho_{\rm TV}\le L$, and every remaining drift time to at most
$1$. Hence $|\partial_sw_{D,\mathfrak p,r}|\le2\max(1,L)^{\ell-1}$
for every $s$, likewise for $\partial_t$. The
derivation of \eqref{eq:pieceLip} uses only these one-variable bounds and absolute continuity in each
variable, so \eqref{eq:pieceLip} holds for $w_{D,\mathfrak p,r}$ with $\Lambda:=2\max(1,L)^{\ell-1}$,
uniformly in $\delta$ and $k$.

\emph{Conclusion.} The cell cancellation in the proof of Proposition~\ref{prop:alllevels} uses only
three facts, none of which involves the drift:$\int_{Q_{ab}}e=0$ for all $a\le b$
(for $a=b$ by the symmetry of $e$), $\norm e_{L^1}\le2L\delta^{2H-1}$ by (C1), with the total-variation
norm at $H=\tfrac12$, and the Lipschitz bound \eqref{eq:pieceLip} on the weight. Hence
\[
   |I_{D,\mathfrak p,r}|\ \le\ 2\delta\Lambda\,\norm e_{L^1}
   \ \le\ 2\delta\cdot2\max(1,L)^{\ell-1}\cdot2L\delta^{2H-1}\ \le\ 8\max(1,L)^{\ell}\,\delta^{2H}.
\]
With $|\theta|^{|D|}\le\max(1,L)^{|D|}$ and $|D|+\ell\le m\le M$, each contribution is at most
$8\max(1,L)^{M}\delta^{2H}$, uniformly in $k$. Summing over the finitely many words, drift sets $D$,
pairings $\mathfrak p$ and indices $r$, whose number depends only on $M$ and $d$, and letting $k\to\infty$, gives $\norm{\cW_{\le M}}_\hs=O(\delta^{2H})$
with constant $C(M,d)\max(1,L)^M$, uniformly over $\cP(H,L)$ and also at $H=\tfrac12$. Only first
derivatives of the weights enter, so Hypothesis~(W) is not needed.
\end{proof}

\section{The rough regime \texorpdfstring{$H\in(\tfrac14,\tfrac12)$}{H in (1/4,1/2)}: the sharp rate \texorpdfstring{$\delta^{4H-1}$}{delta to the power 4H-1}}
\label{sec:rough}

Below $H=\tfrac12$ the analytic proof of Theorem~\ref{thm:2H}(ii) breaks down. The
increment-covariance kernel $c_H|s-t|^{2H-2}$ has $c_H=H(2H-1)<0$ and exponent $2H-2<-1$, so it is
not integrable. Hence the error kernel $e=\varrho^n-\varrho$ has $\norm e_{L^1}=\infty$, the $L^1$
bound behind Theorem~\ref{thm:2H}(ii) is unavailable, and the level-$4$ expectations
are not Lebesgue integrals against $\varrho\otimes\varrho$. They are still well defined, through the natural
lift and without renormalisation. We take a route that avoids the kernel near the diagonal and gives
the sharp exponent and its constant.

\begin{remark}[scaling heuristic]
Sort the level-$4$ weak error by the cells that carry its four tensor slots. Two facts persist for
$H>\tfrac14$: self-similarity, $X_{\delta\,\cdot}\stackrel{d}{=}\delta^H X$, so a single-cell area
scales as $\delta^{2H}$; and the per-cell area defect $D_a=A_a(X)-A_a(X^n)=A_a(X)$ has
$\norm{D_a}_{L^2}\asymp\delta^{2H}$ and mean zero, by centring and component independence
(\eqref{eq:EA0} of Remark~\ref{rem:revarea}, applied on the cell $I_a$). Two contributions arise. \emph{(I) Diagonal}, all four
slots in one cell: a mean level-$4$ defect $\asymp\delta^{4H}$ per cell over $n=1/\delta$ cells
gives $\delta^{4H-1}$. \emph{(II) Cross-cell}, a level-$2$ defect in cell $a$ against one in cell
$b$: $\sum_{a\ne b}\Cov(D_a,D_b)\asymp\delta^{4H}\sum_m(n-m)\rho_A(m)$ with
$\rho_A(m)\asymp m^{4H-6}$ (Lemma~\ref{lem:rhoAasym}; antisymmetry of the L\'evy area gains
two powers over the naive $m^{4H-4}$). This is summable for all $H<\tfrac54$, so again
$n\,\delta^{4H}=\delta^{4H-1}$. The mixed contribution, a local defect against a smooth global
weight, stays at the smaller order $\delta^{2H}$ below the crossover. As $4H-1<2H$ for
$H<\tfrac12$, (I) and (II) dominate.
\end{remark}

Theorem~\ref{thm:roughrate} confirms this exponent: for $H\in(\tfrac14,1)$,
\begin{equation}\label{eq:roughrate}
  \norm{\cW_4}\ \asymp\ \delta^{\gamma(H)},\qquad
  \gamma(H)=\min(2H,\,4H-1)=\begin{cases}2H,& H\ge\tfrac12,\\[2pt] 4H-1,& \tfrac14<H\le\tfrac12,\end{cases}
\end{equation}
The exponent is continuous at the crossover
$H=\tfrac12$ (both branches give $1$, matching Theorem~\ref{thm:2H}(i)), and $\gamma(H)\downarrow0$ as
$H\downarrow\tfrac14$, so the bias stops vanishing at $H=\tfrac14$. This is the boundary of the range
on which \citet{CoutinQian2002} construct the geometric lift as an $L^2$ limit of piecewise-linear
approximations. The proof has two halves: a renormalisation-free
cell-defect decomposition for the upper bound $O(\delta^{4H-1})$ (Proposition~\ref{prop:rough}),
and a reduction of the lower bound to one coordinate, whose leading constant $\kappa_{\rm alt}(H)$
is a strictly negative Riemann-zeta series (Lemma~\ref{lem:altreduction},
Propositions~\ref{prop:altcontinuum}--\ref{prop:kappaclosed}).

\paragraph{Cell telescoping.} Let $X$ be a driver of the class of Definition~\ref{def:class}
(so $H>\tfrac14$ and $\Sig(X)$ is the signature of the natural lift), let $I_a=[t_a,t_{a+1})$,
$a=0,\dots,n-1$, be the mesh cells, $G_a:=\Sig(X)_{t_a,t_{a+1}}$ the signature of $X$ on $I_a$ and
$\Delta_a:=X_{t_{a+1}}-X_{t_a}$ the cell increment. Chen's identity gives
$\Sig(X)=G_0\otimes\cdots\otimes G_{n-1}$ and, since the signature of a linear segment is the tensor
exponential of its increment, $\Sig(X^n)=\exp\Delta_0\otimes\cdots\otimes\exp\Delta_{n-1}$. Put
$T_a:=G_0\otimes\cdots\otimes G_{a-1}\otimes\exp\Delta_a\otimes\cdots\otimes\exp\Delta_{n-1}$, so
that $T_0=\Sig(X^n)$ and $T_n=\Sig(X)$. Replacing the factors $\exp\Delta_b$ by $G_b$ one at a
time, from left to right, and telescoping, $\Sig(X)-\Sig(X^n)=\sum_{a=0}^{n-1}(T_{a+1}-T_a)$, that is,
\begin{equation}\label{eq:celltelescope}
  \Sig(X)-\Sig(X^n)=\sum_{a=0}^{n-1} G^{<a}\otimes\bigl(G_a-\exp\Delta_a\bigr)\otimes\exp\Delta^{>a},
\end{equation}
with $G^{<a}:=\bigotimes_{b<a}G_b$ and $\exp\Delta^{>a}:=\bigotimes_{b>a}\exp\Delta_b$ (empty
products $=\one$). The $a$-th term follows $X$ before cell $a$ ($G^{<a}$), the defect in cell $a$,
and the interpolant after it ($\exp\Delta^{>a}$). The identity is also used in
Propositions~\ref{prop:roughalllevels} and~\ref{prop:drift5}.

\begin{proposition}[cell-defect decomposition and the $O(\delta^{4H-1})$ upper bound]\label{prop:rough}
Let $H\in(\tfrac14,\tfrac12)$, $M\ge4$ and $P\in\cP_0(H,L)$. Then
$\norm{\cW_4}\le C\,\delta^{4H-1}$ for every mesh $\delta=1/n$, with $C$ depending only on $d$
and $L$.
\end{proposition}

\begin{proof}[Proof of Proposition~\ref{prop:rough}]
\emph{Idea and plan.} In the $a$-th term of \eqref{eq:celltelescope} only cell $a$ carries
the defect, which starts at level $2$, has level-$\ell$ part of size $\delta^{\ell H}$ and centred
level-$2$ part. Each of the $n$ cells therefore contributes $\delta^{4H}$ times a sum over the other
cells that the decorrelation clause of (C4) keeps bounded, in total $n\cdot\delta^{4H}=\delta^{4H-1}$.

\emph{Step 1 (telescoping over the cells).} For $H>\tfrac14$ the natural lift exists
(Definition~\ref{def:class}), so \\eqref{eq:celltelescope} holds and no
renormalisation is needed.
Since $\Sig(X)$, $\Sig(X^n)$ and the cell signatures have finite moments of all orders (Lemma~\ref{lem:moments}, through Lemma~\ref{lem:gaussmom}), each term of
\eqref{eq:celltelescope} is integrable by H\"older's inequality.

\emph{Step 2 (the level-$4$ part).} Write $E_a:=G_a-\exp\Delta_a$ for the defect of cell $a$: the
signature of $X$ on the cell minus that of its linear interpolant. The tensor algebra is graded, so the level-$4$ part of a
product $A\otimes B\otimes C$ is $\sum_{p+q+r=4}A^{(p)}\otimes B^{(q)}\otimes C^{(r)}$. The defect has
no levels $0$ and $1$: $E_a^{(0)}=1-1=0$, and $E_a^{(1)}=0$ because the interpolant passes through
the nodes, so both cell increments equal $\Delta_a$. Its level~$2$ is the cell L\'evy area,
$E_a^{(2)}=A_a$, since the symmetric parts $\tfrac12\Delta_a^{\otimes2}$ cancel. Hence, writing
$p,q,r$ for levels, taking expectations in \eqref{eq:celltelescope} and recalling
$\cW_4=\E\,\pi_4\Sig(X^n)-\E\,\pi_4\Sig(X)$, which is minus the expectation of the level-$4$ part of
\eqref{eq:celltelescope},
\[
   \cW_4=-\sum_{a=0}^{n-1}\ \sum_{\substack{p+q+r=4\\ q\ge2}}\
   \E\bigl[(G^{<a})^{(p)}\otimes E_a^{(q)}\otimes(\exp\Delta^{>a})^{(r)}\bigr].
\]
The factors $(G^{<a})^{(p)}$ and $(\exp\Delta^{>a})^{(r)}$ are built from the cells $b\ne a$. We
call them the macroscopic factors, and each increment $\Delta_b$, $b\ne a$, occurring in them a
macroscopic leg. Since $q\ge2$, there are $p+r\le2$ macroscopic levels. We
bound these groups using two inputs.

\emph{(I) Per-cell moments.} By (C2) of Definition~\ref{def:class},
$\norm{E_a^{(\ell)}}_{L^2}\lesssim\delta^{\ell H}$ and $\norm{G_a^{(\ell)}}_{L^2}\lesssim\delta^{\ell H}$
uniformly in $a,\delta$ for $\ell\le4$, with constants depending on $d$ and $L$ (note
$(\exp\Delta_a)^{(\ell)}=\Delta_a^{\otimes\ell}/\ell!$ and $\Delta_a=G_a^{(1)}$. The eighth moment of
$\Delta_a$ is covered because $8\le2M$, the largest moment order allowed in (C2)). Also $\E[E_a^{(2)}]=\E[A_a]=0$, by \eqref{eq:EA0}
of Remark~\ref{rem:revarea} applied on the cell $I_a$. And $\E[E_a^{(3)}]=0$ by parity: $G$ and $-G$ have the
same law, and the map $G\mapsto-G$ multiplies every level-$3$ signature coordinate by $-1$.

\emph{(II) Cross-cell decorrelation.} The estimates below pair a cell-$a$ factor with factors on
cells $b\ne a$, and use only the decorrelation condition of (C4),
$\sup_a\sum_{b\ne a}\varrho_{\max}(I_a,I_b)\le L$. For the independent-fBm law it follows from
Proposition~\ref{prop:cancorr}: $\varrho_{\max}(I_a,I_b)\le C(H)\bar r(|a-b|)$ with
$\bar r(m)\asymp m^{2H-2}$ defined there, hence
summable because $H<\tfrac12$. Within a cell only the moments (I) are used, never the kernel.

\emph{Step 3 (the three groups)}
\begin{itemize}[leftmargin=1.4em,topsep=2pt,itemsep=1pt]
\item $q=4$ (\emph{diagonal}, no macroscopic legs): everything lives in cell $a$, so (I) alone gives
$\bigl|\sum_a\E E_a^{(4)}\bigr|\le\sum_a\norm{E_a^{(4)}}_{L^1}\lesssim n\,\delta^{4H}=\delta^{4H-1}$.
\item $q=2$ (\emph{cross-cell and mixed}): the excess is the mean-zero area $E_a^{(2)}=A_a$, and the
two macroscopic levels sit on cells $\ne a$. Here $E_a^{(2)}$ is measurable for the increment
$\sigma$-field $\cF_{I_a}=\sigma(G_{s,t}:s,t\in I_a)$, the level-$2$ lift on $I_a$ being an a.s.\
limit of piecewise-linear approximations built from $G|_{I_a}$ \citep{CoutinQian2002}. By Chen,
\[
  (G^{<a})^{(2)}=\sum_{b<a}G_b^{(2)}+\sum_{b_1<b_2<a}\Delta_{b_1}\otimes\Delta_{b_2},\qquad
  (\exp\Delta^{>a})^{(2)}=\sum_{b>a}\tfrac12\Delta_b^{\otimes2}
  +\sum_{a<b_1<b_2}\Delta_{b_1}\otimes\Delta_{b_2},
\]
while the split $p=r=1$ contributes $\Delta_{b_1}\otimes\Delta_{b_2}$ with $b_1<a<b_2$, because
$(G^{<a})^{(1)}=\sum_{b<a}\Delta_b$ and $(\exp\Delta^{>a})^{(1)}=\sum_{b>a}\Delta_b$. So every
scalar word coefficient of this group is a sum of terms of one of two types:
\begin{itemize}[leftmargin=1.4em,topsep=1pt,itemsep=1pt]
\item[(i)] $\E[E_a^{(2),\bullet}\Delta^\bullet_{b_1}\Delta^\bullet_{b_2}]$ with $b_1,b_2\ne a$,
possibly $b_1=b_2$, which Lemma~\ref{lem:multileg} with $k=2$ and (C2) bound by
$\sqrt2\,\varrho_{\max}(I_a,I_{b_1})\varrho_{\max}(I_a,I_{b_2})\norm{E_a^{(2)}}_{L^2}
\norm{\Delta_{b_1}}_{L^2}\norm{\Delta_{b_2}}_{L^2}
\lesssim\delta^{4H}\varrho_{\max}(I_a,I_{b_1})\varrho_{\max}(I_a,I_{b_2})$;
\item[(ii)] $\E[E_a^{(2),\bullet}G_b^{(2),\bullet}]$ with $b<a$, arising only from the prefix. As
$\E E_a^{(2)}=0$ this equals $\Cov(E_a^{(2),\bullet},G_b^{(2),\bullet})$, and
Lemma~\ref{lem:twocell} with (C2) bounds it by
$\varrho_{\max}(I_a,I_b)\norm{E_a^{(2)}}_{L^2}\norm{G_b^{(2)}}_{L^2}
\lesssim\delta^{4H}\varrho_{\max}(I_a,I_b)$.
\end{itemize}
Summing (i) over $b_1,b_2\ne a$ and (ii) over $b<a$, and then over $a$, the decorrelation condition of
(C4) gives, up to a constant depending on $d$ and $L$,
\[
   \sum_{a}\Bigl[\sqrt2\Bigl(\sum_{b\ne a}\varrho_{\max}(I_a,I_b)\Bigr)^2+\sum_{b\ne a}\varrho_{\max}(I_a,I_b)\Bigr]\delta^{4H}
   \ \le\ n\,\delta^{4H}\bigl(\sqrt2L^2+L\bigr)=O(\delta^{4H-1}).
\]
  \item $q=3$ (\emph{level-$3\times$level-$1$}): since $\E E_a^{(3)}=0$,
$\E[E_a^{(3)}\otimes M^{(1)}]=\Cov(E_a^{(3)},M^{(1)})$, where $M^{(1)}=\sum_{b<a}\Delta_b$ (the
triple $(1,3,0)$, with the factors in the order $M^{(1)}\otimes E_a^{(3)}$) or
$M^{(1)}=\sum_{b>a}\Delta_b$ (the triple $(0,3,1)$) is a macroscopic level-$1$ factor. The case $k=1$ of Lemma~\ref{lem:multileg} (for the
independent-fBm law, equivalently Proposition~\ref{prop:cancorr}) gives
$|\Cov(E_a^{(3)},\Delta_b)|
\le\varrho_{\max}(I_a,I_b)\norm{E_a^{(3)}}_{L^2}\norm{\Delta_b}_{L^2}
\lesssim\delta^{4H}\,\varrho_{\max}(I_a,I_b)$ at lag $|a-b|\ge1$. Summing over $b\ne a$ with the
decorrelation clause of (C4), and then over $a$, gives $\lesssim n\,\delta^{4H}L=\delta^{4H-1}$.
\end{itemize}
\emph{Conclusion.} As $q\ge2$ there are at most two macroscopic legs, so no lag-sum beyond $L^2$ arises. Summing the
three groups over the $d^4$ words,
\[
   \norm{\cW_4}\ \lesssim\ \underbrace{n\,\delta^{4H}}_{q=4}+\underbrace{n\,\delta^{4H}(L^2+L)}_{q=2}
   +\underbrace{n\,\delta^{4H}L}_{q=3}\ \lesssim\ n\,\delta^{4H}(1+L+L^2)=O(\delta^{4H-1}),
\]
using $n=1/\delta$. Only the
per-cell moments (C2) and the decorrelation clause of (C4) enter, so the implied constant
depends only on $d$ and $L$.
\end{proof}

\paragraph{The chaos frame.}\label{par:chaosframe}
The decorrelation lemmas used above rest on a Wiener-chaos description of the cell
$\sigma$-fields, which we now set up. Let $G$ be as in Definition~\ref{def:class}: $d$ independent components with one scalar increment
covariance. Let $\cH$ be the Hilbert space obtained by completing
$\mathrm{span}\{\one_{[s,t]}:0\le s\le t\le1\}$ for
$\langle\one_{[s,t]},\one_{[u,v]}\rangle_\cH:=\E[G^1_{s,t}G^1_{u,v}]$ (a positive semidefinite form,
quotient by its null space, then complete), and set $\frH:=\cH\otimes\RR^d$. Then
$W(\one_{[s,t]}\otimes e_i):=G^i_{s,t}$ extends to an isonormal Gaussian process $W$ on $\frH$,
because the components are independent with the same $\varrho$. Write $I_q$ for the $q$-th
multiple Wiener--It\^o integral, $\frH^{\odot q}$ for the symmetric tensor power. For an interval
$I$ put
\[
  \cH(I):=\overline{\mathrm{span}}\{\one_{[s,t]}:s,t\in I\}\subset\cH,\quad \frH(I):=\cH(I)\otimes\RR^d,\quad
  P_I:=\text{orth.\ projection }\frH\to\frH(I),
\]
and let $\cF_I:=\sigma(G_{s,t}:s,t\in I)$ be the increment $\sigma$-field of the cell.
Then $\cH_1(I)=\{W(k):k\in\frH(I)\}$ is exactly Definition~\ref{def:class}'s first chaos of
increments on $I$, $\norm{W(k)}_{L^2(\Omega)}=\norm k_\frH$, and
\[
  \varrho_{\max}(I,J)=\sup\Bigl\{\tfrac{|\langle k,l\rangle_\frH|}{\norm k\norm l}:
  k\in\frH(I)\setminus\{0\},\ l\in\frH(J)\setminus\{0\}\Bigr\}.
\]
We never realise $\frH$ as a function space, so no pointwise kernel is needed. This matters for
$H<\tfrac12$, where for fBm the increment covariance has no kernel that is a measure and the L\'evy area is only an
$L^2$-limit.

Since $\cF_I=\sigma(W(k):k\in\frH(I))$ up to null sets and $\frH(I)$ is a closed subspace of
$\frH$, the Wiener chaos decomposition over this Gaussian Hilbert space \citep[Ch.~2]{Janson1997}
gives
\begin{equation}\label{eq:cellchaos}
\begin{gathered}
   L^2(\Omega,\cF_I,\Prob)=\bigoplus_{q\ge0}I_q\bigl(\frH(I)^{\odot q}\bigr),\quad\text{that is,}\\
   F=\sum_{q\ge0}I_q(f_q),\quad f_q\in\frH(I)^{\odot q},\quad f_0=\E F,\quad
   \norm F_{L^2}^2=\sum_{q\ge0}q!\,\norm{f_q}^2 .
\end{gathered}
\end{equation}

The scalar and vector definitions of $\varrho_{\max}$ agree. For $k=\sum_ik_i\otimes e_i\in\frH(I)$
and $l=\sum_il_i\otimes e_i\in\frH(J)$, Cauchy--Schwarz gives
$\langle k,l\rangle=\sum_i\langle k_i,l_i\rangle_\cH
\le\varrho^{\rm sc}_{\max}(I,J)\sum_i\norm{k_i}\norm{l_i}\le\varrho^{\rm sc}_{\max}(I,J)\norm k\norm l$,
where $\varrho^{\rm sc}_{\max}$ is the supremum over the scalar chaoses. The reverse inequality is
trivial. So (C4) may be checked with either definition.

The next lemma bounds the projection of a first-chaos element on one cell onto another cell by
the canonical correlation. It is the only place where $\varrho_{\max}$ enters the proofs of
Lemmas~\ref{lem:multileg} and~\ref{lem:twocell}.

  \begin{lemma}[cell projection]\label{lem:cellproj}
For all intervals $I,J$ and every $l\in\frH(J)$, $\norm{P_Il}_\frH\le\varrho_{\max}(I,J)\norm l_\frH$.
Consequently $\norm{P_I^{\otimes q}l^{(q)}}\le\varrho_{\max}(I,J)^q\norm{l^{(q)}}$ for every
$l^{(q)}\in\frH(J)^{\otimes q}$.
\end{lemma}
\begin{proof}
$P_Il\in\frH(I)$ and $\langle h,P_Il\rangle=\langle h,l\rangle$ for $h\in\frH(I)$, so
$\norm{P_Il}=\sup\{\langle h,l\rangle:h\in\frH(I),\norm h\le1\}\le\varrho_{\max}(I,J)\norm l$ by the
definition of $\varrho_{\max}$. The tensor statement is the multiplicativity of operator norms on
Hilbert tensor products applied to $P_I|_{\frH(J)}$.
\end{proof}

The next lemma bounds the terms (i) of the group $q=2$ ($k=2$) and the group $q=3$ ($k=1$)
in the proof of Proposition~\ref{prop:rough}.

\begin{lemma}[multi-leg separated-support contraction]\label{lem:multileg}
Let $I_a$ be a cell and $I_{b_1},\dots,I_{b_k}$ cells with $b_i\ne a$ for every $i$ (the $b_i$ need not be
distinct). Let $F\in L^2(\Omega,\cF_{I_a},\Prob)$ with $\E F=0$, and let
$\eta_i\in\cH_1(I_{b_i})$. Write $\rho_i:=\varrho_{\max}(I_a,I_{b_i})$ and let $\mathfrak M_k^\circ$
be the set of partial matchings $\sigma$ of $\{1,\dots,k\}$ that are not perfect, and $V(\sigma)$
the set of matched indices. Then
\begin{equation}\label{eq:multileg}
  \Bigl|\E\Bigl[F\prod_{i=1}^k\eta_i\Bigr]\Bigr|\ \le\ \norm F_{L^2}
  \sum_{\sigma\in\mathfrak M_k^\circ}\sqrt{(k-2|\sigma|)!}\
  \prod_{\{i,j\}\in\sigma}\bigl|\E[\eta_i\eta_j]\bigr|
  \prod_{i\notin V(\sigma)}\rho_i\norm{\eta_i}_{L^2}.
\end{equation}
In particular
\[
  k=1:\ |\E[F\eta_1]|\le\rho_1\norm F\norm{\eta_1};\qquad
  k=2:\ |\E[F\eta_1\eta_2]|\le\sqrt2\,\rho_1\rho_2\norm F\norm{\eta_1}\norm{\eta_2},
\]
and the constant $\sqrt2$ at $k=2$ is approached: take $F=W(h)^2-1$ and $\eta_1=\eta_2=W(g)$ with unit
vectors $h\in\frH(I_a)$, $g\in\frH(I_{b_1})$ whose correlation $\langle h,g\rangle$ approaches
$\varrho_{\max}(I_a,I_{b_1})$.
\end{lemma}
\begin{proof}
\emph{Idea and plan.} We expand $F$ in Wiener chaos on the cell $I_a$ and the product of the
legs by the Wick product formula. Chaos orthogonality pairs the order-$q$ part of $F$ with the
terms of the product that leave $q$ legs unmatched, and each unmatched leg can be replaced by its
projection onto the cell, which costs the factor $\rho_i$ (Lemma~\ref{lem:cellproj}).

\emph{Step 1 (chaos frame).} By \eqref{eq:cellchaos} with $I=I_a$, $F=\sum_{q\ge1}I_q(f_q)$ with
$f_q\in\frH(I_a)^{\odot q}$ and $f_0=\E F=0$, and $\sum_qq!\norm{f_q}^2=\norm F_{L^2}^2$. Write
$\eta_i=W(g_i)$ with $g_i\in\frH(I_{b_i})$, so $\norm{\eta_i}_{L^2}=\norm{g_i}$ and
$\E[\eta_i\eta_j]=\langle g_i,g_j\rangle$.

\emph{Step 2 (Wick expansion of the legs).} For centred jointly Gaussian $W(g_1),\dots,W(g_k)$,
\begin{equation}\label{eq:wickprod}
  \prod_{i=1}^kW(g_i)=\sum_{\sigma\in\mathfrak M_k}
  \Bigl(\prod_{\{i,j\}\in\sigma}\langle g_i,g_j\rangle\Bigr)\,
  I_{k-2|\sigma|}\Bigl(\widetilde{\bigotimes_{i\notin V(\sigma)}}\,g_i\Bigr),
\end{equation}
the sum over all partial matchings $\mathfrak M_k$. This is the classical inversion of the Wick
product for jointly Gaussian variables \citep[Ch.~3]{Janson1997}. For completeness, induct on $k$ using the
multiplication formula $W(g)I_q(h)=I_{q+1}(g\tilde\otimes h)+q\,I_{q-1}(h\otimes_1g)$. The case
$k=1$ is $W(g_1)=I_1(g_1)$. If \eqref{eq:wickprod} holds for $k$, multiply by $W(g_{k+1})$: for a
symmetric $h=\widetilde\bigotimes_{i\in S}g_i$ with $|S|=q$ one has
$q\,(h\otimes_1g_{k+1})=\sum_{l\in S}\langle g_l,g_{k+1}\rangle\widetilde\bigotimes_{i\in S\setminus\{l\}}g_i$,
so the first term of the multiplication formula produces the matchings of $\{1,\dots,k+1\}$ leaving
$k+1$ unmatched and the second produces those matching $k+1$ to some $l$, each exactly once.

\emph{Step 3 (pairing off the chaoses).} By \eqref{eq:wickprod} and
$\E[I_q(f_q)I_{q'}(h)]=\delta_{qq'}q!\langle f_q,\tilde h\rangle$,
\[
  \E\Bigl[F\prod_i\eta_i\Bigr]=\sum_{\sigma\in\mathfrak M_k}
  \Bigl(\prod_{\{i,j\}\in\sigma}\langle g_i,g_j\rangle\Bigr)\,q_\sigma!\,
  \Bigl\langle f_{q_\sigma},\ \bigotimes_{i\notin V(\sigma)}g_i\Bigr\rangle,
  \qquad q_\sigma:=k-2|\sigma|,
\]
where the symmetrisation may be dropped because $f_{q_\sigma}$ is symmetric. Perfect matchings
have $q_\sigma=0$, hence carry the factor $f_0=\E F=0$, and drop out.

\emph{Step 4 (the projection).} $f_{q}\in\frH(I_a)^{\otimes q}$, so
$\langle f_{q},\bigotimes_ig_i\rangle=\langle f_{q},\bigotimes_iP_{I_a}g_i\rangle$, and
Lemma~\ref{lem:cellproj} with Cauchy--Schwarz gives
$|\langle f_q,\bigotimes_ig_i\rangle|\le\norm{f_q}\prod_{i\notin V(\sigma)}\rho_i\norm{g_i}$.
Finally $q!\norm{f_q}\le\sqrt{q!}\,\norm F_{L^2}$, because $q!\norm{f_q}^2\le\norm F_{L^2}^2$ by
\eqref{eq:cellchaos}. Hence, for every non-perfect $\sigma$,
\[
   q_\sigma!\,\Bigl|\Bigl\langle f_{q_\sigma},\bigotimes_{i\notin V(\sigma)}g_i\Bigr\rangle\Bigr|
   \ \le\ \sqrt{q_\sigma!}\;\norm F_{L^2}\prod_{i\notin V(\sigma)}\rho_i\norm{\eta_i}_{L^2},
\]
and summing over $\sigma$ with the factors $|\langle g_i,g_j\rangle|=|\E[\eta_i\eta_j]|$ of the matched
pairs gives \eqref{eq:multileg}. For the sharpness claim, $F=W(h)^2-1$ has $\norm F=\sqrt2$, $f_2=h\otimes h$,
and $\E[F\,W(g)^2]=2\langle h,g\rangle^2$ for unit vectors $h,g$; as $\langle h,g\rangle\to\rho_1$ this
tends to $2\rho_1^2=\sqrt2\cdot\sqrt2\cdot\rho_1^2$, the bound.
\end{proof}

In the group $q=2$ of the proof of Proposition~\ref{prop:rough}, the area of cell $a$ is also
paired with the level-$2$ signature of an earlier cell $b$, which is not in the first chaos. The
next lemma bounds such covariances of arbitrary square-integrable cell functionals.

  \begin{lemma}[two-cell decorrelation]\label{lem:twocell}
Let $F_1\in L^2(\cF_{I_a})$ and $F_2\in L^2(\cF_{I_b})$. Then
$|\Cov(F_1,F_2)|\le\varrho_{\max}(I_a,I_b)\norm{F_1}_{L^2}\norm{F_2}_{L^2}$.
\end{lemma}
\begin{proof}
With $F_1=\sum_{q\ge0}I_q(f^{1}_q)$, $F_2=\sum_{q\ge0}I_q(f^{2}_q)$, $f^{1}_q\in\frH(I_a)^{\odot q}$,
$f^{2}_q\in\frH(I_b)^{\odot q}$ the chaos expansions \eqref{eq:cellchaos} on the two cells,
$\Cov(F_1,F_2)=\sum_{q\ge1}q!\langle f^{1}_q,f^{2}_q\rangle=\sum_{q\ge1}q!\langle f^{1}_q,P_{I_a}^{\otimes q}f^{2}_q\rangle$,
so by Lemma~\ref{lem:cellproj} and $\varrho_{\max}\le1$,
$|\Cov(F_1,F_2)|\le\varrho_{\max}\sum_{q\ge1}q!\norm{f^{1}_q}\norm{f^{2}_q}\le\varrho_{\max}\norm{F_1}\norm{F_2}$
by Cauchy--Schwarz.
\end{proof}

Lemma~\ref{lem:twocell} is the Gaussian maximal-correlation identity of
\citet{KolmogorovRozanov1960}, specialised to a pair of cells.

The next proposition verifies the decorrelation condition of (C4) for the independent-fBm law
(input (II) of the proof of Proposition~\ref{prop:rough}).

\begin{proposition}[separated-support decorrelation for fBm cells]\label{prop:cancorr}
Let $H\in(0,1)$, $c_H=H(2H-1)$, and let $\cH_1(I)$ be the first Wiener chaos generated by the increments of
$\mathrm{fBm}(H)$ on an interval $I$. For integers $m\ge1$ put
\[
   \bar r(m):=|c_H|\int_{[0,1]\times[m,m+1]}|u-v|^{2H-2}\,\dd u\,\dd v ,
\]
which is finite and, for $H\ne\tfrac12$, satisfies $\bar r(m)\asymp m^{2H-2}$.
  For unit cells $I_a=[a,a+1]$, $I_b=[b,b+1]$ at integer lag
$m=b-a\ge1$, the canonical correlation
\[
   \varrho_{\max}(I_a,I_b):=\sup\Bigl\{\tfrac{|\Cov(\xi,\eta)|}{\norm\xi_{L^2}\norm\eta_{L^2}}:\
   \xi\in\cH_1(I_a)\!\setminus\!\{0\},\ \eta\in\cH_1(I_b)\!\setminus\!\{0\}\Bigr\}
\]
satisfies $\varrho_{\max}(I_a,I_b)\le C(H)\,\bar r(m)$, with $C(H)$ depending only on $H$. The same bound, with the
same $\bar r$, holds for cells of width $\delta$ at integer lag $m$. For
$F\in L^2\bigl(\sigma(X_{s,t}:s,t\in I_a)\bigr)$ with $\E F=0$ and $\eta\in\cH_1(I_b)$,
$|\Cov(F,\eta)|\le C(H)\bar r(m)\norm F_{L^2}\norm\eta_{L^2}$.
\end{proposition}
\begin{proof}
\emph{Idea and plan.} Off the diagonal the covariance of $\xi=\int_{I_a}\varphi\,\dd X$ and
$\eta=\int_{I_b}\varphi'\,\dd X$ is $\int\!\int\varphi(s)\varphi'(t)\,c_H|t-s|^{2H-2}\,\dd s\,\dd t$
(the prime is a label, not a derivative). We
expand the kernel in powers of the in-cell offsets, so the covariance becomes a series in the
moments $\int\varphi(s)(s-c)^j\,\dd s$ of the integrands, which the spectral isometry bounds
by $\norm\xi_{L^2}$ (Step~1). The series then gives the decay $m^{2H-2}$ (Step~2).

The width-$\delta$ case reduces to the unit case: $\varrho_{\max}$ is a correlation, hence invariant
under the self-similar rescaling $X_{\delta\,\cdot}\stackrel d=\delta^HX$, which maps the width-$\delta$
configuration bijectively onto the unit one. Every $\xi\in\cH_1(I_a)$ is an $L^2$-limit of
$\int_{I_a}\varphi\,\dd X$ with $\varphi$ a step function, and for such $\varphi$, with
$\widehat\varphi(\lambda)=\int\varphi(s)e^{\mathrm i\lambda s}\dd s$, the fBm spectral isometry reads
\begin{equation}\label{eq:cancorrspec}
   \Bigl\lVert\int\varphi\,\dd X\Bigr\rVert^2_{L^2(\Omega)}
   =\kappa_H\int_\RR|\lambda|^{1-2H}|\widehat\varphi(\lambda)|^2\,\dd\lambda
   =\kappa_H\norm{\varphi}^2_{\dot H^{\sigma}},\qquad \sigma:=\tfrac12-H,
\end{equation}
with a constant $\kappa_H>0$ depending only on $H$.

\emph{Step 1 (moment bound).} Let $c$ be the centre of $I_a$ and $\psi_j(s):=(s-c)^j\one_{I_a}(s)$.
Then $|\psi_j|\le2^{-j}$ and $\psi_j$ has bounded variation with $\mathrm{TV}(\psi_j)\le2^{2-j}$, so
$|\widehat\psi_j|\le2^{-j}$ and $|\widehat\psi_j(\lambda)|\le\mathrm{TV}(\psi_j)/|\lambda|$, whence
\[
   \norm{\psi_j}^2_{\dot H^{-\sigma}}=\int_\RR|\lambda|^{2H-1}|\widehat\psi_j|^2\,\dd\lambda
   \ \le\ \frac{2^{-2j}}{H}+\frac{2^{4-2j}}{1-H}\ <\ \infty ,
\]
the two pieces converging because $2H-1>-1$ and $2H-3<-1$. Cauchy--Schwarz in the weighted
Fourier pairing then bounds the moments of $\varphi$:
\begin{equation}\label{eq:momentbound}
   |\mathfrak m_j(\varphi)|:=\Bigl|\int_{I_a}\varphi(s)(s-c)^j\,\dd s\Bigr|
   =\tfrac1{2\pi}\bigl|\langle\widehat\varphi,\widehat{\psi_j}\rangle\bigr|
   \ \le\ A_H\,(1+j)\,2^{-j}\,\norm\xi_{L^2(\Omega)} .
\end{equation}
The decay in $j$ is geometric. Only compact support of $\varphi$ and $|1-2H|<1$ are used.

\emph{Step 2 ($m\ge2$).} On separated cells $\varrho(s,t)=c_H|s-t|^{2H-2}$ is real-analytic. With $s=c_a+u$, $t=c_b+v$, $|u|,|v|\le\tfrac12$ and $\beta:=2-2H\in(0,2]$, so that $t-s=m+v-u\ge m-1\ge1$ and,
for step functions $\varphi$ on $I_a$ and $\varphi'$ on $I_b$,
\[
   \Cov(\xi,\eta)=c_H\int_{I_a}\!\int_{I_b}\varphi(s)\,\varphi'(t)\,(m+v-u)^{-\beta}\,\dd t\,\dd s ,
\]
 the binomial series
\[
   (m+v-u)^{-\beta}=m^{-\beta}\Bigl(1+\frac{v-u}{m}\Bigr)^{-\beta}
   =m^{-\beta}\sum_{k\ge0}\binom{-\beta}{k}m^{-k}(v-u)^k,
 \]
with the generalised binomial coefficient
$\binom{-\beta}{k}=(-1)^k\beta(\beta+1)\cdots(\beta+k-1)/k!$,
converges absolutely, since $|v-u|/m\le\tfrac12$ for $m\ge2$ and
$|\binom{-\beta}{k}|\le k+1$ for $\beta\le2$. Expanding $(v-u)^k$ and applying \eqref{eq:momentbound} to $\varphi$ on
$I_a$ and to $\varphi'$ on $I_b$,
\[
   |\Cov(\xi,\eta)|\le|c_H|\,m^{2H-2}A_H^2\,\norm\xi_{L^2}\norm\eta_{L^2}
   \sum_{k\ge0}\Bigl|\tbinom{-\beta}{k}\Bigr|(1+k)^2\,m^{-k}
   \ \le\ C(H)\,m^{2H-2}\,\norm\xi_{L^2}\norm\eta_{L^2},
\]
using $\sum_l\binom kl=2^k$ against the $2^{-k}$ from \eqref{eq:momentbound},
$|\binom{-\beta}{k}|\lesssim(1+k)^{\beta-1}$, and $m^{-k}\le2^{-k}$ for $m\ge2$, which makes the
$k$-sum converge. For $H\ne\tfrac12$ and $m\ge2$ the integrand of $\bar r(m)$ is at least
$(m+1)^{2H-2}\ge(\tfrac32)^{2H-2}m^{2H-2}$, so $m^{2H-2}\le C(H)\,\bar r(m)$ and the bound is of the
stated form. At $H=\tfrac12$, $c_H=0$ and the covariance vanishes.

\emph{Step 3 ($m=1$).} Here the cells share an endpoint and Step~2 does not apply. But
$\varrho_{\max}\le1$ always, while $\bar r(1)=|c_H|\int_{[0,1]\times[1,2]}|u-v|^{2H-2}$ is
finite (the only singularity is a corner, integrable since $2H-2>-2$) and, for $H\ne\tfrac12$, strictly
positive, so the claim holds at $m=1$ with $C=1/\bar r(1)$. At $H=\tfrac12$ disjoint cells are independent and
$\varrho_{\max}=0$.

\emph{Step 4 (the consequence).} If $F$ is centred and measurable for the increment $\sigma$-field
$\sigma(X_{s,t}:s,t\in I_a)$, its first-chaos projection $\Pi_1F$ lies in $\cH_1(I_a)$ with $L^2$
norm $\le\norm F_{L^2}$. By chaos orthogonality $\Cov(F,\eta)=\Cov(\Pi_1F,\eta)$ for
$\eta\in\cH_1(I_b)$, and the bound follows. \qedhere
\end{proof}

\begin{remark}[sharpness]\label{rem:cancorrsharp}
Numerically, Proposition~\ref{prop:cancorr} is sharp in order. Discretizing both cells into $N$
subintervals, the largest singular value of $G_{aa}^{-1/2}G_{ab}G_{bb}^{-1/2}$ has log--log slopes
in $m$ within $10^{-3}$ of $2H-2$ at $H=0.30,0.35,0.40,0.45,0.60,0.70$, with
$\varrho_{\max}/\bar r(m)\to C_\infty(H)\in[1.00,1.06]$ and no drift in $N$. We have not found this
bound in the literature.
\end{remark}

\subsection*{The lower bound: reduction to a single coordinate}

The next lemma reduces the lower bound matching Proposition~\ref{prop:rough} to showing that
one scalar is nonzero.

\begin{lemma}[coordinate reduction of the rough lower bound]\label{lem:altreduction}
Let $d\ge2$, let $X$ have independent $\mathrm{fBm}(H)$ components with $H\in(\tfrac14,\tfrac12)$,
and let $M\ge4$. For every mesh, $\norm{\cW_4}\ge\sqrt2\,|\Delta E_{\rm alt}|$. The limit
$\kappa_{\rm alt}(H):=\lim_{\delta\to0}\delta^{-(4H-1)}\Delta E_{\rm alt}$ exists
(Proposition~\ref{prop:Lexists}). If it is nonzero, the rough-regime lower bound
$\norm{\cW_4}\gtrsim\delta^{4H-1}$ holds.
\end{lemma}
\begin{proof}
The norm dominates each coordinate: by \eqref{eq:W4norm} and $\binom d2\ge1$,
\[
   \norm{\cW_4}=\Bigl(\tbinom d2\cdot2\sum_\bullet(\Delta E_\bullet)^2\Bigr)^{1/2}
   \ge\bigl(2(\Delta E_{\rm alt})^2\bigr)^{1/2}=\sqrt2\,|\Delta E_{\rm alt}| .
\]
If $\kappa_{\rm alt}\ne0$ then
$|\Delta E_{\rm alt}|\ge\tfrac12|\kappa_{\rm alt}|\,\delta^{4H-1}$ for all sufficiently small $\delta$, whence $\norm{\cW_4}\gtrsim\delta^{4H-1}$.
\end{proof}

\begin{remark}[the $\delta^{4H-1}$ coefficients of the three words]\label{find:kappa}
Write $\kappa_\bullet(H):=\lim_{\delta\to0}\delta^{-(4H-1)}\Delta E_\bullet$ for the three
$\delta^{4H-1}$ coefficients. \emph{(i)} By Lemma~\ref{lem:sumzero},
$\sum_\bullet\Delta E_\bullet=0$ at every mesh, hence also
$\kappa_{\rm adj}+\kappa_{\rm alt}+\kappa_{\rm nest}=0$ whenever the three limits exist. This
improves the coordinate reduction at every mesh: minimising $x^2+y^2$ subject to $x+y=-z$ gives
$x^2+y^2\ge\tfrac12z^2$. With $x=\Delta E_{\rm adj}$, $y=\Delta E_{\rm nest}$, $z=\Delta E_{\rm alt}$
and $d=2$ this gives $2\sum_\bullet(\Delta E_\bullet)^2\ge3(\Delta E_{\rm alt})^2$, i.e.\
$\norm{\cW_4}\ge\sqrt3\,|\Delta E_{\rm alt}|$ by \eqref{eq:W4norm}. With
$\Delta E_{\rm alt}=\kappa_{\rm alt}(H)\delta^{4H-1}+O(\delta^{2H})$
(Proposition~\ref{prop:Lexists}),
\begin{equation}\label{eq:sqrt3}
   \norm{\cW_4}\ \ge\ \sqrt3\,|\kappa_{\rm alt}(H)|\,\delta^{4H-1}+O(\delta^{2H}).
\end{equation}
\emph{(ii)} By Proposition~\ref{prop:adjexact}, $\Delta E_{\rm adj}=-\frac{1}{12}\delta^{2H}+O(\delta^{2H+1}+\delta^2)$
on all of $(\tfrac14,1)$: the Young-regime constant $C^{\rm adj}_1=-\tfrac1{12}$ of
Proposition~\ref{prop:C1exact} persists, and $\kappa_{\rm adj}=0$ as $2H>4H-1$ below the crossover. \emph{(iii)} Hence
$\kappa_{\rm nest}=-\kappa_{\rm alt}$, and \eqref{eq:sqrt3} sharpens to
$\norm{\cW_4}=2|\kappa_{\rm alt}|\delta^{4H-1}(1+o(1))$ at $d=2$ (Corollary~\ref{cor:kadj}). Propositions~\ref{prop:kappaclosed} and~\ref{prop:kappasign} show
that $\kappa_{\rm alt}<0$ on $(\tfrac14,\tfrac12)$.
\end{remark}

The proof of Proposition~\ref{prop:Lexists} subtracts the continuum value $E^\infty_{\rm alt}$
from the interpolant value $E^n_{\rm alt}$. The next proposition gives the three continuum
coordinates in closed form.

\begin{proposition}[closed form of the continuum level-4 coordinates]\label{prop:altcontinuum}
For $H\in(\tfrac14,1)$ the continuum expected-signature
$(2,2)$-coordinates of fBm$(H)$ are
\[
   E^\infty_{\rm adj}=\frac{H}{4(4H-1)}\,\frac{\Gamma(2H)^2}{\Gamma(4H-1)},\quad
   E^\infty_{\rm alt}=\frac{H}{4(4H-1)}\Bigl(1-\frac{\Gamma(2H)^2}{\Gamma(4H-1)}\Bigr),\quad
   E^\infty_{\rm nest}=\frac{2H-1}{8(4H-1)},
\]
with $E^\infty_{\rm adj}+E^\infty_{\rm alt}+E^\infty_{\rm nest}=\tfrac18$.
\end{proposition}
Here the signature is that of the natural lift. For $H>\tfrac12$ the expected signature of fBm was computed by \citet{BaudoinCoutin2007}, and
\citet{Passeggeri2020} shows that every even-level coordinate of the expected signature of
its piecewise-linear approximation with $n$ segments converges at rate $n^{-2H}=\delta^{2H}$. This
contains the upper bound of Theorem~\ref{thm:2H}(ii) and the fBm case of
Proposition~\ref{prop:alllevels} for $H\in(\tfrac12,1)$.
\begin{proof}
\emph{Idea and plan.} For $H>\tfrac12$ each coordinate is an absolutely convergent integral
over the ordered simplex, which we evaluate by integrating out the outer variables first.

\emph{The case $H>\tfrac12$.}
  Each coordinate is $\int_{0<t_1<t_2<t_3<t_4<1}$ of a product of two factors $\varrho(t_i,t_j)=
c_H|t_i-t_j|^{2H-2}$, $c_H=H(2H-1)$, one for each pair $\{i,j\}$ of the Wick pairing,
which pairs the two slots carrying letter $1$ with each other and the two slots carrying letter $2$
with each other (here $t_1,\dots,t_4$
are integration variables, not mesh nodes). Put $\vartheta=2H-1$, so $c_H^2=H^2\vartheta^2$.
\emph{Adjacent} $\{1,2\}\{3,4\}$: integrating $t_1\in(0,t_2)$, $t_4\in(t_3,1)$ gives
$c_H^2\vartheta^{-2}\!\int_{0<t_2<t_3<1}t_2^{\vartheta}(1-t_3)^{\vartheta}
=c_H^2\vartheta^{-2}(\vartheta+1)^{-1}\mathrm{B}(\vartheta+1,\vartheta+2)$, and
$\mathrm{B}(\vartheta+1,\vartheta+2)=\Gamma(\vartheta+1)^2/[2(2\vartheta+1)\Gamma(2\vartheta+1)]$ gives $E^\infty_{\rm adj}$.
\emph{Nested} $\{1,4\}\{2,3\}$: the inner pair $\int(t_3-t_2)^{2H-2}$, then the outer, gives
$c_H^2/[2\vartheta(\vartheta+1)^2(2\vartheta+1)]=(2H-1)/[8(4H-1)]$.
\emph{Alternating} $\{1,3\}\{2,4\}$: the outer variables $t_1,t_4$, then the remaining
$2\times2$ triangle, give $E^\infty_{\rm alt}=H^2/[2(\vartheta+1)(2\vartheta+1)]\,
[1-\Gamma(\vartheta+1)^2/\Gamma(2\vartheta+1)]$. Substituting $\vartheta+1=2H$, $2\vartheta+1=4H-1$ gives the three
formulas, of sum $\tfrac{H}{4(4H-1)}+\tfrac{2H-1}{8(4H-1)}=\tfrac18$, consistent
with Lemma~\ref{lem:sumzero}.

For $H\in(\tfrac14,\tfrac12)$ the simplex integrals diverge, and analytic continuation in $H$
would not identify the stochastic object $\E\langle w,\Sig(\mathbf X^H)\rangle$. Instead, Lemma~\ref{lem:cov4} writes all three as absolutely
convergent integrals of $C$ and its increment variance for every $H>\tfrac14$, and Remark~\ref{rem:noAC}
evaluates them, reproducing the three formulas. The same values follow from the Wiener-chaos
expansion of \citet{CassFerrucci2024}, proved for fBm$(H)$ on all of $(\tfrac14,1)$ with the signature
defined as above. Its level-$0$ term, the expected signature, with the pairing integrals written via
$\partial_{12}R(s,t)=H(2H-1)(t-s)^{2H-2}$ and $\tfrac12R'(t)-\partial_2R(s,t)=H(t-s)^{2H-1}$, gives the
three displayed expressions at the $(2,2)$-words.
\end{proof}

The sign of $\kappa_{\rm alt}$ is decided by the value of $\zeta$ on $(0,1)$, where $\zeta$ is
negative. The next lemma gives the two elementary bounds used for it in the proofs of
Propositions~\ref{prop:kappaclosed} and~\ref{prop:kappasign}.

\begin{lemma}[two bounds for $\zeta$ on $(0,1)$]\label{lem:zetabound}
For $s\in(0,1)$ the function $f(s):=\zeta(s)-\frac{1}{s-1}$ satisfies
\[
   \tfrac12\ \le\ f(s)\ \le\ \tfrac12+\tfrac{s(s+1)}{12}\,\zeta(s+2)\ \le\ \tfrac12+\tfrac{s(s+2)}{12}\ <\ \tfrac34 .
\]
Consequently $\zeta(s)<-\tfrac12$ and $|\zeta(s)|=-\zeta(s)\ge\frac{1}{1-s}-\frac{3}{4}$ for $s\in(0,1)$.
\end{lemma}

\begin{proof}
\emph{Idea.} We compare $\sum_nn^{-s}$ with $\int x^{-s}\,\dd x$ by the trapezoid rule on each
interval $[n,n+1]$. This writes $f(s)$ as $\tfrac12$ plus the sum of the trapezoid errors, which are
nonnegative by convexity and small by the error formula. Analytic continuation carries the identity
from $s>1$ to $(0,1)$.

  For $s>0$ let $g(x)=x^{-s}$, convex with $g''(x)=s(s+1)x^{-s-2}$ decreasing, and let
$\tau_n(s):=\tfrac12(g(n)+g(n+1))-\int_n^{n+1}g$ be the trapezoid errors. As $g$ is convex and the
trapezoid error equals $\tfrac1{12}g''(\xi_n)$ for some $\xi_n\in(n,n+1)$,
$0\le\tau_n(s)\le\tfrac1{12}g''(n)$. For $s>1$,
$\sum_{n<N}\tau_n(s)=\sum_{n\le N}n^{-s}-\tfrac12-\tfrac12N^{-s}-\int_1^Nx^{-s}\dd x\to\zeta(s)-\tfrac12-\frac{1}{s-1}$,
so 
\[
   f(s)=\tfrac12+\sum_{n\ge1}\tau_n(s)\qquad(s>1).
\]
Both sides are analytic on $\{\Re s>-1\}$: $f$ is entire,
and $\sum_n\tau_n(s)$ converges locally uniformly there, as
$|\tau_n(s)|\le\tfrac1{12}|s(s+1)|\,n^{-\Re s-2}$. So the identity holds on $(0,1)$, and
$0\le f(s)-\tfrac12\le\tfrac{s(s+1)}{12}\sum_nn^{-s-2}=\tfrac{s(s+1)}{12}\zeta(s+2)$. Comparing the sum
with an integral, $\zeta(s+2)\le1+\int_1^\infty x^{-s-2}\dd x=\frac{s+2}{s+1}$, which gives the third
inequality, and $s(s+2)<3$ on $(0,1)$ gives the last. Finally $\zeta(s)=f(s)-\frac{1}{1-s}$, so
$-\zeta(s)\ge\frac{1}{1-s}-\frac{3}{4}$, and $\zeta(s)<-\tfrac12$ because
$\frac{1}{1-s}-1=\frac{s}{1-s}>\frac{s(s+2)}{12}$, i.e.\ $(s+2)(1-s)<12$.
\end{proof}

Near $H=\tfrac14$ the sign of $\kappa_{\rm alt}$ follows from the pole of $\zeta$ at $1$, and
the constant blows up as $\gamma(H)=4H-1\downarrow0$.

  \begin{proposition}[rough lower bound near the wall $H=\tfrac14$]\label{prop:roughwall}
As $H\downarrow\tfrac14$,
$\kappa_{\rm alt}(H)=-\dfrac{1}{16(4H-1)}+O(1)\to-\infty$. Hence there is $\eta>0$ with
$\kappa_{\rm alt}(H)<0$, and therefore $\norm{\cW_4}\gtrsim\delta^{4H-1}$
(Lemma~\ref{lem:altreduction}), for all $H\in(\tfrac14,\tfrac14+\eta)$.
\end{proposition}
\begin{proof}
By Proposition~\ref{prop:kappaclosed} and Corollary~\ref{cor:geomkappa},
$\kappa_{\rm alt}=-\tfrac1{16}+\tfrac18\sum_{j\ge1}c_j\zeta(2j-4H)$ with $c_j=|\binom{2H}j|>0$,
$\sum_jc_j=1$ and $c_1=2H$. Only the $j=1$ term is singular at the wall: $\zeta$ has its pole at $1$,
$\zeta(s)=\tfrac1{s-1}+O(1)$, so $\tfrac18c_1\zeta(2-4H)=-\tfrac{2H}{8(4H-1)}+O(1)=
-\tfrac1{16(4H-1)}+O(1)$ (write $2H=\tfrac12+\tfrac12(4H-1)$), and the remainder is bounded, $0<\sum_{j\ge2}c_j\zeta(2j-4H)\le
(1-2H)\zeta(4-4H)\le\zeta(2)$ for $H\in(\tfrac14,\tfrac12)$, since $\sum_{j\ge2}c_j=1-2H$, $\zeta$
decreases on $(1,\infty)$ and $2j-4H\ge4-4H>2$ for $j\ge2$. This gives the asymptotic, and
Lemma~\ref{lem:altreduction} gives $\norm{\cW_4}\gtrsim\delta^{4H-1}$ near the wall.
\end{proof}

On the whole interval $\kappa_{\rm alt}$ is a Riemann-zeta series. The next proposition gives
it, with the limit $-\tfrac18$ as $H\uparrow\tfrac12$ used in Remark~\ref{rem:merger} and a proof
of the sign independent of Proposition~\ref{prop:kappasign}.

\begin{proposition}[closed form and sign of $\kappa_{\rm alt}$]\label{prop:kappaclosed}
For $H\in(\tfrac14,\tfrac12)$,
\[
   \kappa_{\rm alt}(H)\ =\ -\tfrac1{16}+\tfrac18\sum_{j\ge1}\bigl|\tbinom{2H}{j}\bigr|\,\zeta(2j-4H)
   \ =\ \tfrac1{16}+\tfrac18\sum_{k\ge2}\bigl[k^{4H}-(k^2-1)^{2H}\bigr],
\]
where $\zeta$ is the Riemann zeta function and the divergent second series is understood as its finite part
(Definition~\ref{def:fp}). In particular
$\lim_{H\uparrow1/2}\kappa_{\rm alt}(H)=-\tfrac18$,
$\kappa_{\rm alt}(H)\sim-\tfrac{1}{16(4H-1)}$ as $H\downarrow\tfrac14$, and
\[
   \kappa_{\rm alt}(H)<0\quad\text{for all }H\in(\tfrac14,\tfrac12).
\]
\end{proposition}
\begin{proof}
\emph{Closed form.} By Proposition~\ref{prop:Lexists} , $\kappa_{\rm alt}$ exists
and equals $H^2\zeta(2-4H)+\sum_{h\ge1}[U(h)^2-H^2h^{4H-2}]$, which is the finite part of $\sum_hU(h)^2$, $U(h)=-\tfrac14[(h+1)^{2H}-(h-1)^{2H}]$, and
Corollary~\ref{cor:regagree} identifies this finite part with the $\zeta$-series. The second form is the same finite part with the square expanded,
$16\,U(h)^2=(h+1)^{4H}+(h-1)^{4H}-2(h^2-1)^{2H}$, and each power sum continued
($\sum_{k\ge1}k^{4H}=\zeta(-4H)$). With $c_j:=|\binom{2H}{j}|$, which satisfy
$\sum_{j\ge1}c_j=1$ and $(-1)^j\binom{2H}{j}=-c_j$ for $j\ge1$ as $2H\in(0,1)$, the binomial
expansion of $(1-k^{-2})^{2H}$ gives (Proposition~\ref{prop:regagree})
\[
   k^{4H}-(k^2-1)^{2H}=-k^{4H}\sum_{j\ge1}\tbinom{2H}{j}(-1)^jk^{-2j}=\sum_{j\ge1}c_j\,k^{4H-2j},
\]
whose finite part over $k\ge2$ is $\sum_{j\ge1}c_j(\zeta(2j-4H)-1)=\sum_{j\ge1}c_j\zeta(2j-4H)-1$; so
$\tfrac1{16}+\tfrac18(\sum_jc_j\zeta(2j-4H)-1)$ is the first form.

\emph{Limits.} As $H\uparrow\tfrac12$ only $j=1$ survives ($\binom1j=0$ for $j\ge2$), giving
$-\tfrac1{16}+\tfrac18\zeta(0)=-\tfrac18$, a limit of values on the open interval, not the value
$\kappa_{\rm alt}=+\tfrac1{24}$ at $H=\tfrac12$ (Remark~\ref{rem:kappahalf}). The pole of $\zeta(2-4H)$ at $H=\tfrac14$ gives
$\kappa_{\rm alt}\sim-\tfrac1{16(4H-1)}$, as in Proposition~\ref{prop:roughwall}.

\emph{Sign.} In the first form, $c_j>0$ and $c_1=2H$. Since $\zeta$ decreases on $(1,\infty)$ and
$2j-4H\ge4-4H>2$ for $j\ge2$, $\zeta(2j-4H)\le\zeta(4-4H)$, so
$\sum_j c_j\zeta(2j-4H)\le g(H):=2H\,\zeta(2-4H)+(1-2H)\,\zeta(4-4H)$. Three facts about
$\zeta$ finish it: (i) $\zeta(4-4H)\le\zeta(2)=\pi^2/6$ (as $4-4H\ge2$); (ii) $\zeta(s)<-\tfrac12$ for
$s\in(0,1)$; (iii) $\zeta(s)\le\tfrac34-\tfrac1{1-s}$ for $s\in(0,1)$; (ii) and (iii) are
Lemma~\ref{lem:zetabound}. For $H\in[0.28,\tfrac12)$, (ii) gives $2H\zeta(2-4H)<-H$, so
$g(H)<-H+(1-2H)\tfrac{\pi^2}6\le-0.28+0.44\cdot1.645<\tfrac12$. For $H\in(\tfrac14,0.28)$, $u:=4H-1<0.12$
and (iii) gives $2H\zeta(2-4H)\le2H(\tfrac34-\tfrac1u)<\tfrac12(\tfrac34-8)$, while
$(1-2H)\zeta(4-4H)<\tfrac12\zeta(2)<0.83$, so $g(H)<-2.7$. Hence $g(H)<\tfrac12$ on $(\tfrac14,\tfrac12)$, so
$\sum_j c_j\zeta(2j-4H)<\tfrac12$ and $\kappa_{\rm alt}<-\tfrac1{16}+\tfrac18\cdot\tfrac12=0$.
\end{proof}

\subsection*{Existence of the limit, and the value of $\kappa_{\rm alt}$}

Proposition~\ref{prop:rough} only gives $\Delta E_{\rm alt}=O(\delta^{4H-1})$. We now prove that
$\delta^{-(4H-1)}\Delta E_{\rm alt}$ converges and give an absolutely convergent formula for the limit,
which also yields the values of Remark~\ref{rem:kappanumbers}. Two summation lemmas come first.
Throughout them
\[
   r(m)=\tfrac12\bigl[|m+1|^{2H}-2|m|^{2H}+|m-1|^{2H}\bigr],\qquad
   T(h)=-\tfrac12\bigl[h^{2H}-(h-1)^{2H}\bigr],
\]
so that $r(h)=T(h)-T(h+1)$ identically, and $T(h)=\sum_{m\ge h}r(m)$ when $H<\tfrac12$.

The proof of Proposition~\ref{prop:Lexists} reduces $E^n_{\rm alt}$ to sums of the triangular
sums $\Phi(N)$ below. The next lemma writes $\Phi(N)$ exactly as $\sum_{h\le N}T(h)^2$ minus an
explicit defect.

\begin{lemma}[the triangular truncation defect is explicit]\label{lem:trunc}
Let $H\in(0,1)$ and set
\[
   D(k):=(k+1)^{2H}-k^{2H}\quad(k\ge0),\qquad
   \mathcal D(N):=\sum_{a=1}^{N-1}D(a)\,D(N-a)\quad(N\ge1),
\]
with $\mathcal D(1):=0$. Then for every integer $N\ge1$,
\[
   \Phi(N):=\!\!\sum_{\substack{h\ge1,\;p,q\ge0\\ h+p+q\le N}}\!\!r(h+p)\,r(h+q)
   \;=\;\sum_{h\le N}T(h)^2\;-\;\Bigl[\tfrac12D(N)-\tfrac14D(N)^2+\tfrac14\mathcal D(N)\Bigr].
\]
\end{lemma}
\begin{proof}
For $H<\tfrac12$, $T(h)=\sum_{m\ge h}r(m)$, so $\sum_{p,q\ge0}r(h+p)r(h+q)=T(h)^2$: $\Phi(N)$ is
$\sum_hT(h)^2$ with the indices truncated to the triangle $h+p+q\le N$. The proof uses only finite
sums and the following identities, each a telescoping sum of the closed forms of $r$, $T$ and
$D$, so it holds for every $H\in(0,1)$:
\begin{gather}
   \sum_{q=i}^{j}r(q)=T(i)-T(j+1)\quad(1\le i\le j),\qquad
   \sum_{h=1}^{J}T(h)=-\tfrac12J^{2H}\quad(J\ge0),\label{eq:telA}\\
   T(N{+}1)=-\tfrac12D(N),\qquad 2r(a)=D(a)-D(a{-}1)\ (a\ge1),\qquad
   \Sigma(N):=\sum_{a=1}^Nr(a)=\tfrac12\bigl[D(N)-1\bigr],\label{eq:telB}
\end{gather}
the last because $D(0)=1$.

\emph{Step 1 (sum over $q$).} Since $p,q\ge0$, the constraint $h+p+q\le N$ implies $h\le N$, and for
fixed $h$ this yields $p+q\le N-h$. Hence, $\Phi(N)=\sum_{h=1}^N\Psi_h$ with
$\Psi_h:=\sum_{p=0}^{N-h}r(h+p)\sum_{q=0}^{N-h-p}r(h+q)$. By \eqref{eq:telA} with $i=h\le
j=N-p$, the inner sum is $T(h)-T(N+1-p)$. Hence,
\[
   \Psi_h=T(h)\sum_{p=0}^{N-h}r(h+p)-\sum_{p=0}^{N-h}r(h+p)\,T(N{+}1{-}p)
   =T(h)\bigl[T(h)-T(N{+}1)\bigr]-\sum_{a=h}^{N}r(a)\,T(N{+}1{-}a{+}h),
\]
using \eqref{eq:telA} again and substituting $a=h+p$. Equivalently,
\begin{equation}\label{eq:Ehexact}
   T(h)^2-\Psi_h=T(h)\,T(N{+}1)+\sum_{a=h}^{N}r(a)\,T(N{+}1{-}a{+}h).
\end{equation}

\emph{Step 2 (sum over $h$).} Sum \eqref{eq:Ehexact} over $1\le h\le N$. The first term gives
$T(N{+}1)\sum_{h\le N}T(h)=-\tfrac12N^{2H}T(N{+}1)$ by \eqref{eq:telA}. In the double sum, exchange
the order: for fixed $a$, $h$ runs over $1,\dots,a$, so $j:=N{+}1{-}a{+}h$ runs over
$N{+}2{-}a,\dots,N{+}1$, and by \eqref{eq:telA}
\[
   \sum_{j=N+2-a}^{N+1}T(j)=-\tfrac12\bigl[(N{+}1)^{2H}-(N{+}1{-}a)^{2H}\bigr].
\]
Hence, with $\mathcal S(N):=\sum_{a=1}^{N}r(a)\,(N{+}1{-}a)^{2H}$,
\begin{equation}\label{eq:defect1}
   \sum_{h\le N}T(h)^2-\Phi(N)
   =-\tfrac12N^{2H}T(N{+}1)-\tfrac12(N{+}1)^{2H}\,\Sigma(N)+\tfrac12\,\mathcal S(N).
\end{equation}

\emph{Step 3 (summation by parts in $\mathcal S$).} Put $v(a):=(N{+}1{-}a)^{2H}$, so that $v(N)=1$,
$v(1)=N^{2H}$ and $v(a)-v(a{+}1)=D(N{-}a)$. By \eqref{eq:telB},
\[
   2\,\mathcal S(N)=\sum_{a=1}^{N}\bigl[D(a)-D(a{-}1)\bigr]v(a)
   =D(N)v(N)-D(0)v(1)+\sum_{a=1}^{N-1}D(a)\bigl[v(a)-v(a{+}1)\bigr],
\]
the second equality being summation by parts (shift the index in the $D(a{-}1)$-sum). By
the values of $v$ above and $D(0)=1$, the last sum is $\mathcal D(N)$, so
\[
   \mathcal S(N)=\tfrac12\bigl[D(N)-N^{2H}+\mathcal D(N)\bigr].
\] 

\emph{Step 4 (collect).} Insert $T(N{+}1)=-\tfrac12D(N)$, $\Sigma(N)=\tfrac12[D(N)-1]$ and Step~3
into \eqref{eq:defect1}:
\[
   \sum_{h\le N}T(h)^2-\Phi(N)
   =\tfrac14D(N)\bigl[N^{2H}-(N{+}1)^{2H}\bigr]+\tfrac14\bigl[(N{+}1)^{2H}-N^{2H}\bigr]
     +\tfrac14D(N)+\tfrac14\mathcal D(N).
\]
Since $(N{+}1)^{2H}-N^{2H}=D(N)$, the right-hand side is
$-\tfrac14D(N)^2+\tfrac12D(N)+\tfrac14\mathcal D(N)$, which is the claim.
\end{proof}

The defect in Lemma~\ref{lem:trunc} contains the convolution $\mathcal D(N)$, the only term
whose size in $N$ is not explicit. The next lemma gives its leading asymptotics and the proof of
Proposition~\ref{prop:Lexists} uses it at $\varkappa=2H$.

\begin{lemma}[a convolution asymptotic]\label{lem:conv}
Let $\varkappa\in(0,1)$, $D(k)=(k+1)^{\varkappa}-k^{\varkappa}$ and
$\mathcal D(N)=\sum_{a=1}^{N-1}D(a)D(N-a)$. Then, with $\mathrm{B}$ the Beta function,
$\mathrm{B}(x,y)=\int_0^1t^{x-1}(1-t)^{y-1}\dd t=\Gamma(x)\Gamma(y)/\Gamma(x+y)$,
\[
   \mathcal D(N)=\varkappa^{2}\mathrm{B}(\varkappa,\varkappa)\,N^{2\varkappa-1}+O\!\left(N^{\varkappa-1}\right)
   \qquad(N\to\infty),
 \]
the constant in $O(\cdot)$ depending only on $\varkappa$.
  \end{lemma}

\begin{proof}
\emph{Idea.} Since $D(a)\approx\varkappa a^{\varkappa-1}$,
\[
   \mathcal D(N)\approx\varkappa^2\sum_{a}a^{\varkappa-1}(N-a)^{\varkappa-1}
   \approx\varkappa^2\int_0^Nx^{\varkappa-1}(N-x)^{\varkappa-1}\,\dd x
   =\varkappa^2\mathrm{B}(\varkappa,\varkappa)\,N^{2\varkappa-1}.
\]
To make this quantitative we split the sum at $N/2$, so that in each half one of the factors
$x^{\varkappa-1}$, $(N-x)^{\varkappa-1}$ stays away from its singularity.

\emph{Tools (facts about $D$).} Throughout, $C$ denotes a constant depending only on $\varkappa$. As
$x\mapsto x^\varkappa$ is concave with $(x^\varkappa)''=-\varkappa(1-\varkappa)x^{\varkappa-2}$, Taylor's
theorem gives, for $a\ge1$,
\[
   0<D(a)\le\varkappa a^{\varkappa-1},\qquad
   \eta(a):=D(a)-\varkappa a^{\varkappa-1},\quad |\eta(a)|\le\tfrac{\varkappa(1-\varkappa)}2\,a^{\varkappa-2},
\]
so $\sum_{a\ge1}|\eta(a)|<\infty$. Moreover $\sum_{a=1}^{Y}D(a)=(Y+1)^\varkappa-1\le CY^\varkappa$, by
telescoping.

\emph{Step 1 (split at $N/2$).} The summand of $\mathcal D(N)$ is symmetric under $a\mapsto N-a$, so,
with $A:=\lceil N/2\rceil-1$ the largest integer below $N/2$,
\[
   \mathcal D(N)=2\sum_{a=1}^{A}D(a)D(N-a)+\one_{N\ \rm even}D(N/2)^2,
\]
and the last term is at most $\varkappa^2(N/2)^{2\varkappa-2}=O(N^{2\varkappa-2})$.

\emph{Step 2 (replace the far factor).} For $a\le A$ we have $N-a>N/2$, so
$D(N-a)=\varkappa(N-a)^{\varkappa-1}+\eta(N-a)$ with $|\eta(N-a)|\le CN^{\varkappa-2}$. The error is
at most $2CN^{\varkappa-2}\sum_{a\le A}D(a)\le CN^{2\varkappa-2}$, so
\[
   2\sum_{a=1}^{A}D(a)D(N-a)=2\varkappa\sum_{a=1}^{A}D(a)(N-a)^{\varkappa-1}+O(N^{2\varkappa-2}).
\]

\emph{Step 3 (replace the near factor).} Write $D(a)=\varkappa a^{\varkappa-1}+\eta(a)$. As
$(N-a)^{\varkappa-1}\le(N/2)^{\varkappa-1}$ for $a\le A$, the error is at most
$2\varkappa(N/2)^{\varkappa-1}\sum_{a\ge1}|\eta(a)|=O(N^{\varkappa-1})$, so
\[
   2\varkappa\sum_{a=1}^{A}D(a)(N-a)^{\varkappa-1}
   =2\varkappa^{2}\sum_{a=1}^{A}a^{\varkappa-1}(N-a)^{\varkappa-1}+O(N^{\varkappa-1}).
\]

\emph{Step 4 (sum versus integral).} Put $g(x):=(N-x)^{\varkappa-1}$ and
$\varphi(x):=x^{\varkappa-1}g(x)$. On $[0,N/2]$, $g$ is increasing, $g\le(N/2)^{\varkappa-1}$ and
$|g'|\le(1-\varkappa)(N/2)^{\varkappa-2}$. Then
\[
   \sum_{a=1}^{A}\varphi(a)-\int_0^{A}\varphi=\sum_{a=1}^{A}\int_{a-1}^{a}\bigl[\varphi(a)-\varphi(x)\bigr]\dd x .
\]
The term $a=1$ contains the singularity at $0$. It is at most
$\varphi(1)+\int_0^1x^{\varkappa-1}g(x)\,\dd x\le(1+\tfrac1\varkappa)\,g(1)\le CN^{\varkappa-1}$. For
$a\ge2$ and $x\in[a-1,a]$, the mean value theorem applied to each factor gives
\[
   |\varphi(a)-\varphi(x)|\le|a^{\varkappa-1}-x^{\varkappa-1}|\,g(a)+x^{\varkappa-1}|g(a)-g(x)|
   \le(1-\varkappa)\Bigl[(a{-}1)^{\varkappa-2}(N/2)^{\varkappa-1}+(a{-}1)^{\varkappa-1}(N/2)^{\varkappa-2}\Bigr].
\]
Summing over $2\le a\le A$ with $\sum_{k\ge1}k^{\varkappa-2}=\zeta(2-\varkappa)<\infty$ and
$\sum_{k\le A}k^{\varkappa-1}\le CN^{\varkappa}$ gives $O(N^{\varkappa-1})+O(N^{2\varkappa-2})$,
which is $O(N^{\varkappa-1})$ as $\varkappa<1$.
Hence
\[
   \sum_{a=1}^{A}a^{\varkappa-1}(N-a)^{\varkappa-1}=\int_0^{A}x^{\varkappa-1}(N-x)^{\varkappa-1}\,\dd x
   +O(N^{\varkappa-1}).
\]

\emph{Step 5 (the integral).} Extending the integral from $A$ to $N/2$ costs $O(N^{2\varkappa-2})$,
since $N/2-A\le1$ and the integrand is $O(N^{2\varkappa-2})$ there. By the symmetry $x\mapsto N-x$
and the substitution $x=Ny$,
\[
   2\int_0^{N/2}x^{\varkappa-1}(N-x)^{\varkappa-1}\,\dd x=\int_0^{N}x^{\varkappa-1}(N-x)^{\varkappa-1}\,\dd x
   =\mathrm{B}(\varkappa,\varkappa)\,N^{2\varkappa-1}.
\]

\emph{Conclusion.} Combining Steps 1--5,
$\mathcal D(N)=\varkappa^2\mathrm{B}(\varkappa,\varkappa)N^{2\varkappa-1}+O(N^{\varkappa-1})+O(N^{2\varkappa-2})$,
and $O(N^{2\varkappa-2})\subset O(N^{\varkappa-1})$ because $\varkappa<1$.
\end{proof}

The next proposition combines the two lemmas into the expansion of $\Delta E_{\rm alt}$ and an
absolutely convergent series for $\kappa_{\rm alt}$.

\begin{proposition}[the rough-regime limit exists]\label{prop:Lexists}
Let $X$ have $d\ge2$ independent $\mathrm{fBm}(H)$ components, $H\in(\tfrac14,\tfrac12)$, and let $\delta=1/n$. Then
\[
   \Delta E_{\rm alt}\ =\ \kappa_{\rm alt}(H)\,\delta^{4H-1}+O\bigl(\delta^{2H}\bigr),
\]
so in particular $\kappa_{\rm alt}(H)=\lim_{\delta\to0}\delta^{-(4H-1)}\Delta E_{\rm alt}$ exists, and
\[
   \kappa_{\rm alt}(H)\ =\ H^2\,\zeta(2-4H)\ +\ \sum_{h\ge1}\Bigl[U(h)^2-H^2h^{4H-2}\Bigr],\qquad
   U(h)=-\tfrac14\bigl[(h+1)^{2H}-(h-1)^{2H}\bigr],
\]
the series converging absolutely (its terms are $O(h^{4H-4})$).
\end{proposition}
Since $\kappa_{\rm alt}(H)\ne0$ (Proposition~\ref{prop:kappasign}), the expansion has relative
error $O(\delta^{\varepsilon(H)})$ with $\varepsilon(H):=1-2H>0$. 

\begin{proof}[Proof of Proposition~\ref{prop:Lexists}]
\emph{Idea and plan.} $Q_{\rm alt}(n):=\delta^{-4H}E^n_{\rm alt}$ is a sum over ordered
quadruples of cells of products of two increment correlations $r(\cdot)$. Step~1 splits it
exactly into elementary sums (Step~2) and a sum $B_1(n)=\sum_{N<n}\Phi(N)$ (Step~3), whose
coefficient of $n$ is the constant term of $\Phi(N)$, found with Lemmas~\ref{lem:trunc}
and~\ref{lem:conv} (Step~4) and identified with the series of the statement (Step~5).

\emph{Notation and tools.} With $r$ and $T$ as above, $\rho_2=\sum_{h\ge1}r(h)^2$, $\Sigma(b)=\sum_{a=1}^{b}r(a)$,
$\mathfrak s=\sum_{b\ge1}r(b)\Sigma(b)$. For $H<\tfrac12$: $r(m)<0$ for $m\ge1$, $r(m)=c_Hm^{2H-2}+
O(m^{2H-4})$ with $c_H=H(2H-1)$, $\sum_m|r(m)|<\infty$, and $\sum_{m\in\ZZ}r(m)=0$, that is,
$\sigma:=\sum_{m\ge1}r(m)=-\tfrac12$, as $r(0)=1$. From $|r(m)|\le Cm^{2H-2}$ we get, for real $y\ge1$,
\begin{equation}\label{eq:rtailmom}
   \sum_{m\ge y}|r(m)|\le Cy^{2H-1},\qquad \sum_{m\le y}m\,|r(m)|\le Cy^{2H},
\end{equation}
with $C$ depending only on $H$. We use one elementary device repeatedly. If $|f(m)|\le Cm^{-\gamma}$
for $m\ge1$, with $\gamma>1$ and $\gamma\ne2$, then
\begin{equation}\label{eq:linsum}
   \sum_{m=1}^{n}(n-m)\,f(m)=n\sum_{m\ge1}f(m)+O\bigl(1+n^{2-\gamma}\bigr),
\end{equation}
because the left-hand side equals $n\sum_{m\ge1}f(m)-n\sum_{m>n}f(m)-\sum_{m\le n}mf(m)$, and the
last two sums are $O(n\cdot n^{1-\gamma})$ and $O(1+n^{2-\gamma})$.

\emph{Step 1 (an exact algebraic reduction).} By \eqref{eq:chenPL} for the word $1212$, and since
the two components are independent with $\E[\Delta^i_a\Delta^i_b]=\delta^{2H}r(b-a)$ for cells
$a,b\in\{0,\dots,n-1\}$, the only Wick pairing joins slots $1,3$ and slots $2,4$:
\[
   E^n_{\rm alt}=\sum_{c_1\le c_2\le c_3\le c_4}w(\mathbf c)\,
   \E[\Delta^1_{c_1}\Delta^1_{c_3}]\,\E[\Delta^2_{c_2}\Delta^2_{c_4}]
   =\delta^{4H}\sum_{c_1\le c_2\le c_3\le c_4}w(\mathbf c)\,r(c_3-c_1)\,r(c_4-c_2),
\]
where $w(\mathbf c)$ is the coefficient $c_{\mathbf a}$ of \eqref{eq:chenPL} at $\mathbf a=\mathbf c$:
$w(\mathbf c)=\prod_j1/m_j!$, the product running over the blocks of equal cells among
$c_1,\dots,c_4$ and $m_j$ being the size of the $j$-th block. Parametrise the ordered cell tuple
$c_1\le c_2\le c_3\le c_4$ by the gaps $h_1=c_2-c_1$, $h_2=c_3-c_2$, $h_3=c_4-c_3\ge0$; for fixed
gaps, $c_1$ takes the $(n-|h|)_+$ values $0,\dots,n-1-|h|$, where $|h|:=h_1+h_2+h_3$. Hence
\[
   Q_{\rm alt}(n)=\delta^{-4H}E^n_{\rm alt}=\sum_{h_1,h_2,h_3\ge0}(n-|h|)_+\,w(h)\,r(h_1{+}h_2)\,r(h_2{+}h_3),
\]
where $w(h):=w(\mathbf c)$ depends only on which gaps vanish. 
For $h_2=0$ the summand is $w_A(h_1,h_3)r(h_1)r(h_3)$ with $w_A=\tfrac12,\tfrac16,\tfrac1{24}$
according as none, one, or both of $h_1,h_3$ vanish. For $h_2\ge1$ the weight factorises,
$w=w_1(h_1)w_1(h_3)$, $w_1(0)=\tfrac12$, $w_1(\ge1)=1$. Passing to the lag variables $a=h_1+h_2=c_3-c_1$, $b=h_2+h_3=c_4-c_2$ of the two Wick pairs, which range
over $a,b\ge1$ with $1\le h_2\le a\wedge b$ and $|h|=a+b-h_2$, and sum over $h_2$ first. With
$w_1(a-h_2)w_1(b-h_2)=1-\tfrac12\one_{h_2=a}-\tfrac12\one_{h_2=b}+\tfrac14\one_{h_2=a=b}$ this gives
the exact identity
\begin{equation}\label{eq:Qaltsplit}\begin{gathered}
   Q_{\rm alt}(n)=A(n)+B_1(n)-B_2(n)+B_3(n),\\
   A(n)=\!\!\sum_{h_1,h_3\ge0}\!\!(n-h_1-h_3)_+\,w_A(h_1,h_3)\,r(h_1)r(h_3),\\
   B_1(n)=\!\!\sum_{a,b\ge1}\!\!r(a)r(b)\Lambda_n(a,b),\qquad
   \Lambda_n(a,b)=\!\!\sum_{h=1}^{a\wedge b}\!(n-a-b+h)_+,\\
   B_2(n)=\!\!\sum_{b\le n}\!(n-b)r(b)\Sigma(b),\qquad
   B_3(n)=\tfrac14\!\sum_{a\le n}\!(n-a)r(a)^2,
\end{gathered}\end{equation}
with $\Sigma(b)=\sum_{a=1}^br(a)$, and $A(n)$ the $h_2=0$ part.
The term $-B_2$ collects the two terms $-\tfrac12\one_{h_2=a}$ and $-\tfrac12\one_{h_2=b}$: the first
requires $a\le b$ (as $h_2\le a\wedge b$) and then $|h|=b$, the second requires $b\le a$ and then $|h|=a$, and by the symmetry $a\leftrightarrow b$
each gives $\tfrac12B_2(n)$. The term $\tfrac14\one_{h_2=a=b}$ gives $B_3(n)$, with $|h|=a$.

\emph{Step 2 (the elementary pieces).}\emph{The sum rule.} Counting which of $h_1,h_3$ vanish
and using $r(0)=1$,
\[
   \sum_{h_1,h_3\ge0}w_A(h_1,h_3)\,r(h_1)r(h_3)=\tfrac12\sigma^2+2\cdot\tfrac16\sigma+\tfrac1{24}
   =\tfrac18-\tfrac16+\tfrac1{24}=0 .
\]
\emph{The piece $A$.} Write $(n-h_1-h_3)_+=n-(h_1+h_3)$ on $h_1+h_3<n$. By the sum rule, $A(n)$ equals
$-n$ times the sum of $w_Ar(h_1)r(h_3)$ over $h_1+h_3\ge n$, minus the sum of
$(h_1+h_3)w_Ar(h_1)r(h_3)$ over $h_1+h_3<n$. On $h_1+h_3\ge n$ one of $h_1,h_3$ is at least $n/2$,
and by \eqref{eq:rtailmom} the first contribution is $O(n\cdot n^{2H-1}\sum_m|r(m)|)$ and the second is
$O(n^{2H}\sum_m|r(m)|)$. Hence
\[
   A(n)=O(n^{2H})=o(n).
\]
\emph{The pieces $B_2$, $B_3$.} Apply \eqref{eq:linsum}: to $f=r\Sigma$, with $|\Sigma(b)|\le\sum_m|r(m)|$
and $\gamma=2-2H\in(1,\tfrac32)$; and to $f=\tfrac14r^2$, with $\gamma=4-4H\in(2,3)$. This gives
\[
   B_2(n)=n\,\mathfrak s+O(n^{2H}),\qquad B_3(n)=\tfrac n4\rho_2+O(1).
\]

\emph{Step 3 (an exact one-dimensional reduction of $B_1$).}
$\Lambda_n(a,b)$ vanishes unless $a\vee b\le n-1$, so all sums here are finite. Exchanging the
$h$-sum with the $(a,b)$-sum and substituting $a=h+p$, $b=h+q$
(the condition $h\le a\wedge b$ becomes $p,q\ge0$, and $n-a-b+h=n-h-p-q$),
\[
   B_1(n)=\sum_{h\ge1}\ \sum_{p,q\ge0}r(h+p)\,r(h+q)\,(n-h-p-q)_+ .
\]
For an integer $k\ge1$ one has $(n-k)_+=\#\{N\in\NN:k\le N\le n-1\}$. Applying this with
$k=h+p+q$ and exchanging the sums,
\begin{equation}\label{eq:B1asPhi}
   B_1(n)=\sum_{N=1}^{n-1}\Phi(N),\qquad
   \Phi(N)=\!\!\sum_{\substack{h\ge1,\;p,q\ge0\\ h+p+q\le N}}\!\!r(h+p)\,r(h+q).
\end{equation}
This is the point of the reduction: summing over $N$ raises every power of $N$ by one, so the coefficient of $n^1$ in $B_1$ is the constant term of $\Phi$: if
$\Phi(N)=\sum_iA_iN^{\alpha_i}+c+O(N^{-\beta})$ with $0<\alpha_i<1$ and $\beta\in(0,1)$, then summing
$\sum_{N\le n-1}N^{\alpha}=\frac{n^{\alpha+1}}{\alpha+1}+O(n^{\alpha})$ gives
\begin{equation}\label{eq:B1exp}
   B_1(n)=\sum_i\frac{A_i}{\alpha_i+1}\,n^{\alpha_i+1}+c\,n+O\bigl(n^{1-\beta}+n^{\max_i\alpha_i}\bigr).
\end{equation}

\emph{Step 4 (the constant term of $\Phi$).} By Lemma~\ref{lem:trunc} and Lemma~\ref{lem:conv}
applied at $\varkappa=2H\in(\tfrac12,1)$, and since $D(N)=O(N^{2H-1})$ and $D(N)^2=O(N^{4H-2})$,
\begin{equation}\label{eq:Phidefect}
   \sum_{h\le N}T(h)^2-\Phi(N)=H^2\mathrm{B}(2H,2H)\,N^{4H-1}+O\bigl(N^{2H-1}\bigr).
\end{equation}
Here $\tfrac14\mathcal D(N)=\tfrac14(2H)^2\mathrm B(2H,2H)N^{4H-1}+O(N^{2H-1})$ by
Lemma~\ref{lem:conv}, and $\tfrac12D(N)-\tfrac14D(N)^2=O(N^{2H-1})$ because $4H-2<2H-1$.
 On the other hand $T(h)=-Hh^{2H-1}+O(h^{2H-2})$ gives $T(h)^2=H^2h^{4H-2}+O(h^{4H-3})$ with
$4H-3<-1$. Split
\[
   \sum_{h\le N}T(h)^2=H^2\sum_{h\le N}h^{4H-2}+\sum_{h\ge1}\bigl[T(h)^2-H^2h^{4H-2}\bigr]
   -\sum_{h>N}\bigl[T(h)^2-H^2h^{4H-2}\bigr].
\]
By (F3) of Lemma~\ref{lem:fp} in Appendix~\ref{app:fp}, applied with the exponent $4H-2\in(-1,0)$,
$\sum_{h\le N}h^{4H-2}=\frac{N^{4H-1}}{4H-1}+\zeta(2-4H)+O(N^{4H-2})$. The middle series converges
absolutely and the tail is $O(\sum_{h>N}h^{4H-3})=O(N^{4H-2})$. With (F1) of Lemma~\ref{lem:fp},
the constant is the finite part $\mathfrak f(T^2)$ of Definition~\ref{def:fp}:
 \begin{equation}\begin{gathered}\label{eq:T2partial}
   \sum_{h\le N}T(h)^2=\frac{H^2}{4H-1}\,N^{4H-1}+\mathfrak f(T^2)+O\bigl(N^{4H-2}\bigr),\\
   \mathfrak f(T^2)=H^2\zeta(2-4H)+\sum_{h\ge1}\bigl[T(h)^2-H^2h^{4H-2}\bigr].
\end{gathered}\end{equation}
Subtracting, and using $4H-2<2H-1$ for $H<\tfrac12$,
\begin{equation}\label{eq:Phifinal}
   \Phi(N)=\mathfrak b_H\,N^{4H-1}+\mathfrak f(T^2)+O\bigl(N^{2H-1}\bigr),\qquad
   \mathfrak b_H:=H^2\Bigl[\frac{1}{4H-1}-\mathrm{B}(2H,2H)\Bigr].
\end{equation}
This is an expansion of the form required in Step~3, with the single exponent
$\alpha=4H-1\in(0,1)$ and $\beta=1-2H$, so \eqref{eq:B1exp} gives
\begin{equation}\label{eq:B1final}
   B_1(n)=\frac{\mathfrak b_H}{4H}\,n^{4H}\;+\;\mathfrak f(T^2)\,n\;+\;O\bigl(n^{2H}\bigr),
\end{equation}
the error $O(n^{2H})$ dominating both $O(n^{4H-1})$ and $O(1)$ on $(\tfrac14,\tfrac12)$.
Finally $\mathrm B(2H,2H)=\Gamma(2H)^2/\Gamma(4H)$ and $\Gamma(4H)=(4H-1)\Gamma(4H-1)$ turn  the leading coefficient into
\[
   \frac{\mathfrak b_H}{4H}=\frac{H}{4}\Bigl[\frac{1}{4H-1}-\frac{\Gamma(2H)^2}{\Gamma(4H)}\Bigr]
   =\frac{H}{4(4H-1)}\Bigl(1-\frac{\Gamma(2H)^2}{\Gamma(4H-1)}\Bigr)=E^\infty_{\rm alt},
\]
by Proposition~\ref{prop:altcontinuum}. Inserting \eqref{eq:B1final} and Step~2 into
\eqref{eq:Qaltsplit},
\[
   Q_{\rm alt}(n)=E^\infty_{\rm alt}\,n^{4H}+\bigl[\mathfrak f(T^2)-\mathfrak s+\tfrac14\rho_2\bigr]n
   +O\bigl(n^{2H}\bigr),
\]
so that\, as $\Delta E_{\rm alt}=E^n_{\rm alt}-E_{\rm alt}=\delta^{4H}Q_{\rm alt}(n)-E^\infty_{\rm alt}$
(the continuum coordinate $E_{\rm alt}$ is $E^\infty_{\rm alt}$ of Proposition~\ref{prop:altcontinuum})
and $\delta=1/n$,
\begin{equation}\label{eq:Lexistsconst}
   \delta^{-(4H-1)}\Delta E_{\rm alt}=n^{-1}Q_{\rm alt}(n)-n^{4H-1}E^\infty_{\rm alt}
   =\mathfrak f(T^2)-\mathfrak s+\tfrac14\rho_2+O(n^{2H-1}).
\end{equation}

\emph{Step 5 (identification with the series of the statement).} Write $\mathfrak c:=\mathfrak f(T^2)$ for the constant Step~4
delivers. With $U=T-\tfrac12 r$ (this is the closed form of $U$ in the statement, since
$T(h)-\tfrac12r(h)=-\tfrac14[(h+1)^{2H}-(h-1)^{2H}]$) ,
$\sum U^2=\sum T^2-\sum rT+\tfrac14\rho_2$; and since $\Sigma(b)=\sigma-T(b+1)$, $T(b{+}1)=T(b)-r(b)$
and $\sigma=-\tfrac12$,
\[
   \mathfrak s=\sigma^2-\sum_b r(b)T(b)+\rho_2=\tfrac14-\sum_h r(h)T(h)+\rho_2 .
\]
To eliminate $\sum rT$ we telescope: $T(h)^2-T(h+1)^2=r(h)\bigl(2T(h)-r(h)\bigr)$, summed over $h\ge1$ with $T(h)\to0$, gives the exact identity
$2\sum_{h\ge1}r(h)T(h)=T(1)^2+\rho_2=\tfrac14+\rho_2$. Substituting,
\[
   \mathfrak c-\mathfrak s+\tfrac14\rho_2
   =\Bigl[\textstyle\sum U^2+\sum rT-\tfrac14\rho_2\Bigr]-\Bigl[\tfrac14-\sum rT+\rho_2\Bigr]+\tfrac14\rho_2
   =\textstyle\sum U^2+2\sum rT-\rho_2-\tfrac14=\sum U^2 ,
\]
all sums in the finite-part sense of Definition~\ref{def:fp} in Appendix~\ref{app:fp}. This is
legitimate because the constant of Step~4 is $\mathfrak f(T^2)$ by \eqref{eq:T2partial}, and
$\sum_hr(h)T(h)$ and $\rho_2$ converge absolutely, so (F1)--(F2) of Lemma~\ref{lem:fp} apply.
Finally $U(h)^2=H^2h^{4H-2}+O(h^{4H-4})$ (the $h^{2H-2}$ term of $U$ cancels, $U$ being a
central difference) and
$\sum_{h\le N}h^{4H-2}=\frac{N^{4H-1}}{4H-1}+\zeta(2-4H)+o(1)$ (valid since $2-4H\in(0,1)$), so the
finite part is the absolutely convergent expression of the statement:
\[
   \mathfrak f(T^2)-\mathfrak s+\tfrac14\rho_2=\mathfrak f(U^2)
   =H^2\zeta(2-4H)+\sum_{h\ge1}\bigl[U(h)^2-H^2h^{4H-2}\bigr]=\kappa_{\rm alt}(H).
\]

\emph{Conclusion.} By \eqref{eq:Lexistsconst} and Step~5,
$\delta^{-(4H-1)}\Delta E_{\rm alt}=\kappa_{\rm alt}+O(n^{2H-1})=\kappa_{\rm alt}+O(\delta^{1-2H})$,
and multiplying by $\delta^{4H-1}$ gives $\Delta E_{\rm alt}=\kappa_{\rm alt}\delta^{4H-1}+O(\delta^{2H})$.
\end{proof}

We also prove negativity directly from the representation of Proposition~\ref{prop:Lexists}, so that
Theorem~\ref{thm:roughrate} does not need the identification of Corollary~\ref{cor:regagree}.

\begin{proposition}[sign of $\kappa_{\rm alt}$, from the absolutely convergent representation]
\label{prop:kappasign}
For every $H\in(\tfrac14,\tfrac12)$ the series $S(H):=\sum_{h\ge1}[U(h)^2-H^2h^{4H-2}]$ has
nonnegative terms and
\[
   \kappa_{\rm alt}(H)=H^2\zeta(2-4H)+S(H)\ <\ 0 .
\]
\end{proposition}
\begin{proof}
\emph{Idea and plan.} Since $\zeta(2-4H)<0$ and $S(H)\ge0$, the claim is that the positive
series $S(H)$ is smaller than $H^2|\zeta(2-4H)|$. We bound $S(H)$ above by an explicit function
$S^+(H)$ (Steps~1--2), bound $|\zeta(2-4H)|$ below by Lemma~\ref{lem:zetabound} (Step~3), and
compare the two elementary bounds in the variable $u=4H-1$ (Step~4).

Write $\alpha:=2H-1\in(-\tfrac12,0)$, so that $U(h)=-\tfrac H2\int_{h-1}^{h+1}x^{\alpha}\dd x$
(the integral converges at $h=1$ because $\alpha>-1$).

\emph{Step 1 (termwise nonnegativity).} $x\mapsto x^{\alpha}$ is convex on $(0,\infty)$, so Jensen for the
uniform law on $[h-1,h+1]$ (whose mean is $h$) gives $\tfrac12\int_{h-1}^{h+1}x^{\alpha}\dd x\ge h^{\alpha}$, i.e.\
$|U(h)|\ge Hh^{\alpha}$ and $U(h)^2-H^2h^{2\alpha}\ge0$.

\emph{Step 2 (an explicit bound for $S$).} The $h=1$ term is $\tfrac{2^{4H}}{16}-H^2$. For $h\ge2$ put
$d_h:=|U(h)|-Hh^{\alpha}=\tfrac H2\int_{h-1}^{h+1}(x^{\alpha}-h^{\alpha})\dd x$. Taylor's
theorem with the Lagrange remainder gives, for $x\in[h-1,h+1]$ and some $\xi$ between $h$ and $x$,
\[
   x^{\alpha}-h^{\alpha}=\alpha h^{\alpha-1}(x-h)+\tfrac12\alpha(\alpha-1)\,\xi^{\alpha-2}(x-h)^2,
   \qquad 0<\alpha(\alpha-1)\,\xi^{\alpha-2}\le\alpha(\alpha-1)(h-1)^{\alpha-2},
\]
since $\xi\ge h-1\ge1$ and $\alpha-2<0$. Integrating, with  $\int_{h-1}^{h+1}(x-h)\dd x=0$ and $\int_{h-1}^{h+1}(x-h)^2\dd x=\tfrac23$,
\[
   d_h\le\tfrac{H\alpha(\alpha-1)}{6}(h-1)^{\alpha-2}\le\tfrac H8(h-1)^{\alpha-2},
\]
 using $\alpha(\alpha-1)\le\tfrac34$ for $\alpha\in(-\tfrac12,0)$. Also $|U(h)|\le H(h-1)^{\alpha}$
(the mean of the decreasing function $x^\alpha$ over $[h-1,h+1]$ is at most its value at $h-1$), so $|U(h)|+Hh^{\alpha}\le2H(h-1)^{\alpha}$, and as $d_h\ge0$ by Step~1,
\[
   U(h)^2-H^2h^{2\alpha}=d_h\bigl(|U(h)|+Hh^{\alpha}\bigr)\le\tfrac{H^2}{4}(h-1)^{2\alpha-2}
   =\tfrac{H^2}{4}(h-1)^{4H-4}\qquad(h\ge2).
\]
Summing over $h\ge2$, with $\sum_{h\ge2}(h-1)^{4H-4}=\zeta(4-4H)$, and adding the $h=1$ term,
\[
   S(H)\ \le\ S^+(H):=\frac{2^{4H}}{16}-H^2+\frac{H^2}{4}\,\zeta(4-4H).
\]

\emph{Step 3 (a lower bound for $|\zeta(2-4H)|$).} Lemma~\ref{lem:zetabound} at $s=2-4H\in(0,1)$,
where $1-s=4H-1$, gives $|\zeta(2-4H)|\ge Z(H):=\max\{\tfrac12,\ \tfrac1{4H-1}-\tfrac34\}$.

\emph{Step 4 (comparison).} We prove $S^+(H)<H^2Z(H)$ on $(\tfrac14,\tfrac12)$. Substitute
$u:=4H-1\in(0,1)$, so that $2^{4H}=2\cdot2^u$, $16H^2=(1+u)^2$, $4-4H=3-u$ and
$Z(H)=\max\{\tfrac12,\ \tfrac1u-\tfrac34\}$.

\emph{(a) Two convexity bounds.} For $s>1$ the map $x\mapsto x^{-s}$ is convex, so the midpoint
inequality $f(\tfrac{a+b}2)\le\frac{1}{b-a}\int_a^bf$ for convex $f$ (the left half of the
Hermite--Hadamard inequality) gives $k^{-s}\le\int_{k-1/2}^{k+1/2}x^{-s}\dd x$ for $k\ge2$. Summing over $k\ge2$,
\begin{equation}\label{eq:zetaub}
   \zeta(s)\ \le\ 1+\int_{3/2}^{\infty}x^{-s}\dd x\ =\ 1+\frac{(3/2)^{1-s}}{s-1},\qquad s>1 .
\end{equation}
Likewise $u\mapsto2^u$ and $u\mapsto(3/2)^u$ are convex, hence dominated on $[0,1]$ by their chords:
\begin{equation}\label{eq:chords}
   2^u\le1+u,\qquad (3/2)^u\le1+\tfrac u2\qquad(0\le u\le1).
\end{equation}

\emph{(b) An elementary upper bound for $S^+$.} By the first bound in \eqref{eq:chords},
$\frac{2^{4H}}{16H^2}=\frac{2\cdot2^u}{(1+u)^2}\le\frac{2}{1+u}$. By \eqref{eq:zetaub} with $s=3-u>1$
and $(3/2)^{u-2}=\tfrac49(3/2)^u\le\tfrac49\bigl(1+\tfrac u2\bigr)$ from \eqref{eq:chords},
$\tfrac14\zeta(4-4H)\le\tfrac14+\frac{2+u}{18(2-u)}$. Adding,
\begin{equation}\label{eq:Splusub}
   \frac{S^+(H)}{H^2}=\frac{2^{4H}}{16H^2}-1+\frac{\zeta(4-4H)}{4}\ \le\ \mathfrak T(u)-\tfrac34,
   \qquad \mathfrak T(u):=\frac{2}{1+u}+\frac{2+u}{18(2-u)} .
\end{equation}
It therefore suffices to prove $\mathfrak T(u)<Z+\tfrac34$ for $u\in(0,1)$, which we do in two cases whose
ranges meet at $u=\tfrac45$, i.e.\ at $H=0.45$.

\emph{(c) The case $u\ge\tfrac45$}, where we use only $Z\ge\tfrac12$ and must show
$\mathfrak T(u)<\tfrac54$. On $[0,1]$,
\[
   \mathfrak T'(u)=-\frac{2}{(1+u)^2}+\frac{2}{9(2-u)^2}\le-\frac{1}{2}+\frac{2}{9}=-\frac{5}{18}<0
\]
since $(1+u)^2\le4$ and $(2-u)^2\ge1$; hence $\mathfrak T$ decreases and
$\mathfrak T(u)\le\mathfrak T(\tfrac45)=\tfrac{10}9+\tfrac7{54}=\tfrac{67}{54}<\tfrac{5}4$, the last step being
$268<270$.

\emph{(d) The case $u\le\tfrac45$}, where we use $Z\ge\tfrac1u-\tfrac34$ and must show
$\mathfrak T(u)<\tfrac1u$. Since $\tfrac1u-\frac{2}{1+u}=\frac{1-u}{u(1+u)}$, the claim is that
\[
   R(u):=\tfrac1u-\mathfrak T(u)=\frac{1-u}{u(1+u)}-\frac{2+u}{18(2-u)}\ >\ 0,\qquad 0<u\le\tfrac45 .
\]
On $(0,1)$ the first summand is decreasing, being a positive decreasing numerator over a positive
increasing denominator, and $\frac{2+u}{18(2-u)}$ is increasing; so $R$ decreases, and
\[
   R(u)\ \ge\ R(\tfrac45)=\tfrac5{36}-\tfrac7{54}=\tfrac1{108}>0 .
\]

\emph{Conclusion.} In both cases $S(H)\le S^+(H)<H^2Z(H)\le H^2|\zeta(2-4H)|$ by Steps 2 and 3, and $\zeta(2-4H)<0$, so
$\kappa_{\rm alt}(H)=H^2\zeta(2-4H)+S(H)<0$. In case~(c), $Z-S^+/H^2\ge\tfrac12-(\tfrac{67}{54}-\tfrac34)
=\tfrac1{108}$, and in case~(d), $Z-S^+/H^2\ge R(u)\ge\tfrac1{108}$. This gives the explicit, if crude,
margin $\kappa_{\rm alt}(H)\le-\frac{H^2}{108}<-\frac{1}{1728}$.
\end{proof}

\begin{remark}[the Brownian value]\label{rem:kappahalf}
Setting $H=\tfrac12$ formally in Proposition~\ref{prop:Lexists} makes every term of $S$ vanish
($U(h)\equiv-\tfrac12$, $H^2h^{4H-2}\equiv\tfrac14$) and returns $\tfrac14\zeta(0)=-\tfrac18$. This is
the limit as $H\uparrow\tfrac12$, not the value at $H=\tfrac12$. The closed form of $T$ behind $U$ uses
$\sum_{m\ge1}r(m)=\lim_N\bigl(\tfrac12[(N+1)^{2H}-N^{2H}]-\tfrac12\bigr)=-\tfrac12$, valid for every
$H<\tfrac12$; at $H=\tfrac12$ one has $r(m)\equiv0$ for $m\ge1$, so this sum jumps to $0$ and the
representation does not apply. At $H=\tfrac12$ itself $\kappa_{\rm alt}=+\tfrac1{24}$ by
Theorem~\ref{thm:2H}(i), and $-\tfrac18$ enters Remark~\ref{rem:merger} as
$\lim_{H\uparrow1/2}\kappa_{\rm alt}$.
\end{remark}

\begin{remark}[numerical values of $\kappa_{\rm alt}$]\label{rem:kappanumbers}
The $\zeta$-series of Proposition~\ref{prop:kappaclosed} converges slowly, at rate $O(J^{-2H})$, since
$|\binom{2H}{j}|\asymp j^{-1-2H}$ while $\zeta(2j-4H)\to1$. Corollary~\ref{cor:geomkappa}
in Appendix~\ref{app:fp} rearranges it into a geometrically convergent form. The formula of
Proposition~\ref{prop:Lexists} gives
\[
\begin{array}{c|cccccc}
H & 0.30 & 0.35 & 0.40 & 0.45 & 0.48 & \to\tfrac12\\\hline
\kappa_{\rm alt}(H) & -0.341259 & -0.191518 & -0.146876 & -0.129644 & -0.125705 & -0.125000
\end{array}
\]
increasing to $-\tfrac18$ as $H\uparrow\tfrac12$ (a one-sided limit, see
Remark~\ref{rem:merger}).
\end{remark}

The next corollary combines Propositions~\ref{prop:Lexists}, \ref{prop:kappasign}
and~\ref{prop:adjexact} with the endpoint identity into the exact leading constant of
$\norm{\cW_4}$, the value $\cB(H)=2|\kappa_{\rm alt}(H)|$ of the Introduction.

\begin{corollary}[the exact rough-regime constant]\label{cor:kadj}
For $H\in(\tfrac14,\tfrac12)$, $M\ge4$ and $d$ independent fBm$(H)$ components, $\kappa_{\rm adj}=0$,
$\kappa_{\rm nest}=-\kappa_{\rm alt}$, and
\[
   \norm{\cW_4}=\sqrt{\tbinom d2}\;2\,|\kappa_{\rm alt}(H)|\,\delta^{4H-1}\bigl(1+O(\delta^{1-2H})\bigr).
\]
\end{corollary}

\begin{proof}
The three word-differences sum to zero and the adjacent one is of lower order, so
\eqref{eq:W4norm} is governed by $\Delta E_{\rm alt}$.
  By Proposition~\ref{prop:adjexact}, $\Delta E_{\rm adj}=O(\delta^{2H})=o(\delta^{4H-1})$, so
$\kappa_{\rm adj}=\lim_{\delta\to0}\delta^{-(4H-1)}\Delta E_{\rm adj}=0$. By
Lemma~\ref{lem:sumzero}, $\Delta E_{\rm nest}=-\Delta E_{\rm alt}-\Delta E_{\rm adj}$, so
$\kappa_{\rm nest}=-\kappa_{\rm alt}$, the limit $\kappa_{\rm alt}$ existing by
Proposition~\ref{prop:Lexists}. Expanding the square of $\Delta E_{\rm nest}$,
\[
   \sum_\bullet(\Delta E_\bullet)^2=(\Delta E_{\rm alt})^2+(\Delta E_{\rm adj})^2
   +(\Delta E_{\rm alt}+\Delta E_{\rm adj})^2=2(\Delta E_{\rm alt})^2(1+\chi+\chi^2),\qquad
   \chi:=\frac{\Delta E_{\rm adj}}{\Delta E_{\rm alt}} .
\]
Here $\chi=O(\delta^{2H}/\delta^{4H-1})=O(\delta^{1-2H})$, as
$\Delta E_{\rm alt}=\kappa_{\rm alt}\delta^{4H-1}(1+O(\delta^{1-2H}))$ with $\kappa_{\rm alt}<0$
(Propositions~\ref{prop:Lexists} and~\ref{prop:kappasign}). In particular $\Delta E_{\rm alt}\ne0$
for small $\delta$. Inserting this into \eqref{eq:W4norm},
\[
   \norm{\cW_4}=\sqrt{\tbinom d2}\;2\,|\Delta E_{\rm alt}|\sqrt{1+\chi+\chi^2}
   =\sqrt{\tbinom d2}\;2\,|\kappa_{\rm alt}|\,\delta^{4H-1}\bigl(1+O(\delta^{1-2H})\bigr).
\]
\end{proof}

The next theorem joins the upper bound of Proposition~\ref{prop:rough} to the lower bound of
Lemma~\ref{lem:altreduction}. It is the rough-regime half of the rate stated in the Introduction.

\begin{theorem}[the sharp rough-regime rate]\label{thm:roughrate}
Let $d\ge2$, let $X$ have independent $\mathrm{fBm}(H)$ components, and let $H\in(\tfrac14,\tfrac12)$
and $M\ge4$.
The level-$4$ weak error satisfies $\norm{\cW_4}\asymp\delta^{4H-1}$, so its squared norm is
$\asymp\delta^{2(4H-1)}$; precisely, at $d=2$, $\norm{\cW_4}=2|\kappa_{\rm alt}(H)|\,\delta^{4H-1}(1+o(1))$ with
$\kappa_{\rm alt}(H)<0$.
\end{theorem}

\begin{proof}[Proof of Theorem~\ref{thm:roughrate}]
The upper bound $\norm{\cW_4}=O(\delta^{4H-1})$ is Proposition~\ref{prop:rough}.
 For the lower bound,
$\Delta E_{\rm alt}=\kappa_{\rm alt}(H)\delta^{4H-1}(1+o(1))$ by Proposition~\ref{prop:Lexists}, with
$\kappa_{\rm alt}(H)<0$ by Proposition~\ref{prop:kappasign}. Both use the same absolutely
convergent expression for $\kappa_{\rm alt}$, so this step does not depend on the $\zeta$-series of
Proposition~\ref{prop:kappaclosed}. Lemma~\ref{lem:altreduction}, valid for every $d\ge2$, then gives
$\norm{\cW_4}\ge\sqrt2\,|\Delta E_{\rm alt}|=\sqrt2\,|\kappa_{\rm alt}(H)|\,\delta^{4H-1}(1+o(1))$. So $\norm{\cW_4}\asymp\delta^{4H-1}$. The exact asymptotic $\norm{\cW_4}=2|\kappa_{\rm alt}(H)|\delta^{4H-1}
(1+o(1))$ at $d=2$ (times $\sqrt{\binom d2}$ in general) is Corollary~\ref{cor:kadj}.
\end{proof}

The $o(1)$ holds at fixed $H$ as $\delta\to0$ and is not uniform: it degrades as $H\uparrow\tfrac12$,
as its Young-side counterpart does as $H\downarrow\tfrac12$ (Remark~\ref{rem:nonuniform}), at
$H=0.48$ the ratio $\delta^{-(4H-1)}\Delta E_{\rm alt}$ changes sign only at $n\approx1.2\cdot10^3$.

\begin{remark}[the one-dimensional augmented case]\label{rem:hp}
\citet[Lem.~A.1]{HagerPelizzari2026} bound the same discretization error for the time-augmented
one-dimensional fBm $Y=(t,X^H)$, $H\in(0,\tfrac12)$, by $C_{K,T,H}\,\Delta t^{2H}$ at every fixed
truncation level $K$. On $(\tfrac14,\tfrac12)$ this is smaller than the $\delta^{4H-1}$ of
Theorem~\ref{thm:roughrate}, and the difference comes from the alphabet. With one rough component
the signature coordinates are monomials in the increments, so a sup-norm bound on the interpolation
error of the covariance suffices. Here the $\delta^{4H-1}$ is carried by the alternating and nested
$(2,2)$-words (Remark~\ref{find:kappa}), that is, by the second moment of the L\'evy area between two
independent rough components, which has no analogue for $d=1$.
\end{remark}

\paragraph{Full truncation below the crossover.} The cell-by-cell bookkeeping of
Proposition~\ref{prop:rough} expands the macroscopic factor into cell increments and sums their
pairings in absolute value. This is harmless at level $4$, but from level $6$ on, macroscopic legs
that pair among themselves would be summed over a free pair of cells, at a cost $\delta^{2H-1}$
(Remark~\ref{rem:M6}). Instead we integrate by parts against the whole macroscopic factor and let
Chen's relation resum the undifferentiated legs into signatures of sub-intervals, whose moments are
bounded. Write $Z=(t,G)$ for the time-augmented Gaussian part, with time letter $0$ and
$e_0$ the corresponding basis vector,
and $Z^n=(t,G^n)$. We use three Gaussian facts, in the chaos frame set up before
Lemma~\ref{lem:cellproj}. (G1) If $F\in L^2(\cF_I)$ lies in the chaoses of order $\le q$, then
$F=\sum_{k\le q}I_k(f_k)$ with $f_k\in\frH(I)^{\odot k}$ (by \eqref{eq:cellchaos}), $D^jF\in L^2(\Omega;\frH(I)^{\otimes j})$, and
$\E\norm{D^jF}^2\le q^j\,\E F^2$, since $\E\norm{D^jI_k(f)}^2=\frac{k!}{(k-j)!}\,k!\norm f^2\le k^j\E I_k(f)^2$
and derivatives of different chaoses are orthogonal. (G2) On the chaoses of order $\le q$,
$\norm F_{L^s}\le(s-1)^{q/2}\norm F_{L^2}$ for $s\ge2$ \citep[Thm.~5.10]{Janson1997}. The same holds for
$\frH^{\otimes j}$-valued $F$ with scalar components in those chaoses, since for an orthonormal basis
$(e_\iota)$ Minkowski's inequality in $L^{s/2}$ gives
$\norm{\norm F^2}_{L^{s/2}}\le\sum_\iota\norm{\langle F,e_\iota\rangle}^2_{L^s}\le(s-1)^q\E\norm F^2$. (G3) For
$F=\sum_kI_k(f_k)$ and $\Phi$ in a finite sum of chaoses, Stroock's formula
$\Phi=\sum_kI_k(\E D^k\Phi/k!)$ \citep[see, e.g.,][]{Nualart2006} gives $\E[F\Phi]=\E F\,\E\Phi+\sum_{k\ge1}\langle f_k,\E D^k\Phi\rangle$, and
if $f_k\in\frH(I)^{\otimes k}$ the pairing only sees $P_I^{\otimes k}\E D^k\Phi$.

The next lemma bounds the moments of the signatures of $Z$ on sub-intervals, the objects into
which Chen's relation resums the macroscopic legs. The strong-error bounds of
\S\ref{subsec:strongbelow} also use it.

  \begin{lemma}[(C2) at all meshes controls all sub-intervals]\label{lem:subint}
Let $P\in\cP(H,L)$, $H\in(\tfrac14,1)$. For every word $w$ of length $m\le M$ in the letters of $Z$,
with $j$ time letters and $g=m-j$ letters of $G$, every $s<\infty$ and all $0\le r<r'\le1$,
\[
   \norm{\langle w,\Sig(Z)_{r,r'}\rangle}_{L^s}\ \le\ C(m,s,H,L)\,|r'-r|^{\,j+gH},
\]
and $\langle w,\Sig(Z)_{r,r'}\rangle$ lies in the Wiener chaoses of $G$ of order $\le g$ and parity $g$.
\end{lemma}
The proof uses only (C2) at $q=2$ on the dyadic meshes, and (C5) only for the existence of the
natural lift.

\begin{proof}
We prove the bound first on mesh cells, from (C2) and hypercontractivity (G2). Then on dyadic
intervals, which are unions of cells of different meshes, by Chen's relation. Then for general
endpoints by continuity and finally we add the time letters by integrating them out.

\emph{Cells.} The lift of the rescaled cell path is the $\delta^{-H}$-dilation of the restricted lift,
piecewise-linear interpolation commuting with affine time changes and dilations; so (C2) at $q=2$
gives $\norm{\langle v,\Sig(G)_{I}\rangle}_{L^2}\le L^{1/2}|I|^{|v|H}$ for every mesh cell $I$ and $|v|\le M$.
The coordinate is the limit in $L^2$ of homogeneous polynomials of degree $|v|$ in the increments of
$G$ on $I$: the lifts of the interpolants converge in probability \citep[Thm.~15.33(iv)]{FrizVictoir2010}, and
in every $L^q$ by uniform integrability, their moments being uniform in the mesh by
Lemma~\ref{lem:gaussmom} (with $K=L$ by (C5)); for $|v|\ge4$ add the continuity and the moment bound
of the Lyons extension \citep[Cor.~9.11(i), Thm.~9.5(ii)]{FrizVictoir2010}.
So it lies in the chaoses of order $\le|v|$ and parity $|v|$, and (G2) gives
$\norm{\langle v,\Sig(G)_I\rangle}_{L^s}\le C_{s,|v|}L^{1/2}|I|^{|v|H}$.

\emph{Dyadic sub-intervals, no time letters.} Let $r<r'$ be dyadic, $\ell=r'-r$. Then $[r,r']$ is the
union of consecutive dyadic intervals $J_1<\dots<J_N$, at most two of each length $2^{-i}\le\ell$, each
a cell of the mesh $2^{-i}$. By Chen,
$\langle w,\Sig(G)_{r,r'}\rangle=\sum_{p\le m}\sum_{w=v_1\cdots v_p}\sum_{i_1<\dots<i_p}\prod_{l=1}^p
\langle v_l,\Sig(G)_{J_{i_l}}\rangle$, the middle sum over deconcatenations into non-empty words. The
triangle inequality, H\"older with $p$ factors in $L^{ps}$, the cell bound, and dropping the order
constraint give
\[
   \norm{\langle w,\Sig(G)_{r,r'}\rangle}_{L^s}\ \le\ C\sum_{p\le m}\sum_{v_1\cdots v_p}\prod_{l=1}^p
   \Bigl(\sum_{i}|J_i|^{|v_l|H}\Bigr)\ \le\ C\,\ell^{mH},
\]
since $\sum_i|J_i|^{\alpha H}\le2\ell^{\alpha H}/(1-2^{-\alpha H})$ for $\alpha\ge1$. 

\emph{General endpoints.}The lift is a.s.\ continuous in $(r,r')$ at all levels (extension theorem), so for dyadic $r_k\to r$,
$r'_k\to r'$ Fatou's lemma transfers the bound. The convergence also holds in $L^2$, by uniform
integrability, so the chaos statement passes to the limit, the sum of the chaoses of order $\le g$
being closed in $L^2$.

\emph{Time letters.} Write $w=u_00u_10\cdots0u_j$ with $u_l$ words in the letters of $G$. For smooth
paths, integrating out the $j$ time slots gives
\[
   \langle w,\Sig(Z)_{r,r'}\rangle=\int_{r<\sigma_1<\dots<\sigma_j<r'}\ \prod_{l=0}^j\langle u_l,\Sig(G)_{\sigma_l,\sigma_{l+1}}\rangle
   \,\dd\sigma_1\cdots\dd\sigma_j,\qquad\sigma_0=r,\ \sigma_{j+1}=r',
\]
and both sides pass to the limit along piecewise-linear approximations pathwise: along a
subsequence the approximating lifts converge a.s.\ in the rough-path topology, so by continuity of the
extension and of the Young pairing with the bounded-variation path $t$
\citep[Cor.~9.11(i), Cor.~9.31]{FrizVictoir2010} all levels converge uniformly in $(\sigma,\sigma')$. Minkowski's
integral inequality and H\"older give $\norm{\cdot}_{L^s}\le C\int_{\rm simplex}\prod_l(\sigma_{l+1}-
\sigma_l)^{|u_l|H}\dd\sigma\le C\ell^{\,j+gH}$. The chaos statement passes through products and Bochner
integrals, $\sigma\mapsto$ integrand being continuous into $L^2$.
\end{proof}

\begin{remark}[two consequences for Definition~\ref{def:class}]\label{rem:C2simplify}
(i) The cell step shows that (C2) at $q=2$ implies (C2) at every $q$, with a constant depending on
$q$. (ii) Lemma~\ref{lem:subint} at levels $\le3$ gives $\norm{d(\mathbf G_s,\mathbf G_t)}_{L^q}\le
C_q|t-s|^H$ for the homogeneous distance and every $q$, so the Kolmogorov criterion
\citep[Thms.~A.10, A.12]{FrizVictoir2010}, used as in the proof of
Lemma~\ref{lem:gaussmom}, gives $\E\norm{\mathbf G}^r_{\pvar}\le C(r,p,H,L)$ for every $p\in(1/H,4)$,
and the drift $\theta te_1$ is absorbed by \citet[Thm.~9.33(i)]{FrizVictoir2010} with the ball--box comparison
\citep[Prop.~8.9, (8.6)]{FrizVictoir2010}. So the moment clause of
Definition~\ref{def:class} is implied by (C2), up to the value of the constant.
\end{remark}

The next lemma carries out the integration by parts against the macroscopic factor: it bounds
the Malliavin derivatives of that factor and turns them into a bound on its covariance with a variable
of one cell.

\begin{lemma}[integration by parts against the macroscopic factor]\label{lem:ibp}
Let $H\in(\tfrac14,\tfrac12)$ and $P\in\cP(H,L)$. Fix a cell $a$ and words $u,z$ in the letters of $Z$
with $|u|+|z|=r\le M$, and put
$\Phi=\langle u,\Sig(Z)_{0,t_a}\rangle\langle z,\Sig(Z^n)_{t_{a+1},1}\rangle$. For every $k\ge1$,
\[
   \E D^k\Phi=\sum_{\mathbf b}K_{\mathbf b},\qquad K_{\mathbf b}\in\textstyle\bigotimes_{l=1}^k\frH(I_{b_l}),
   \qquad \norm{K_{\mathbf b}}_{\frH^{\otimes k}}\le C\,\delta^{kH},
\]
the sum over $\mathbf b=(b_1,\dots,b_k)$ with all $b_l\ne a$, and $C=C(r,k,d,H,L)$ independent of
$a,\mathbf b,n$. Consequently $\norm{\Phi}_{L^s}\le C(s,r,d,H,L)$, and for
$F=\sum_{k\le q}I_k(f_k)\in L^2(\cF_{I_a})$,
\[
   \bigl|\E[F\Phi]-\E F\,\E\Phi\bigr|\ \le\ C\sum_{k=1}^q\norm{f_k}\,(L\,\delta^H)^k .
\]
\end{lemma}

\begin{proof}
We expand $\Phi$ by Chen's relation into factors living on single cells $b\ne a$ and
differentiate with the Leibniz rule. Each derivative falls on one cell and costs $\delta^H$,
and the undifferentiated factors resum into signatures of the gaps, with bounded moments. Finally
(G3) and the decorrelation condition of (C4) turn the derivative bound into the covariance bound.

\emph{Step 1 (Chen expansion and Leibniz rule).} Write $\Delta_b:=G_{t_{b+1}}-G_{t_b}$ for the
increment of $G$ over $I_b$, and $T^{(r)}$ for the tensor algebra over the letters of $Z$, truncated
at level $r$. By Chen, $\Sig(Z)_{0,t_a}=Y_0\otimes\dots\otimes Y_{a-1}$ with $Y_b=\Sig^{\le r}(Z)_{I_b}$, and
$\Sig(Z^n)_{t_{a+1},1}=Y_{a+1}\otimes\dots\otimes Y_{n-1}$ with $Y_b=\exp^{\le r}(\delta e_0+\Delta_b)$.
Each $Y_b$ is $\cF_{I_b}$-measurable with coordinates in the chaoses of order $\le r$
(Lemma~\ref{lem:subint}, for $b>a$ they are polynomials in $\Delta_b$), so by (G1) $D^jY_b$ takes values
in $\frH(I_b)^{\otimes j}\otimes T^{(r)}$. Thus $\langle u,\Sig(Z)_{0,t_a}\rangle=\sum_{u=u_0\cdots u_{a-1}}
\prod_c\langle u_c,Y_c\rangle$ (deconcatenation, empty words allowed) is a finite sum of finite products
of variables with moments of all orders, and the Leibniz rule applies. For fixed $\mathbf b$ let $B$
be the set of distinct cells in $\mathbf b$ and $j_c=\#\{l:b_l=c\}$. A term is non-zero only if the
words at the cells of $B$ are non-empty. For each choice of those words, the sum over the remaining deconcatenations factorises through the gaps between consecutive
cells of $B$, where Chen resums it into one coordinate of the signature of the gap: of $Z$ on the
prefix, of $Z^n$ on the suffix, each a node-to-node sub-interval. Hence $K_{\mathbf b}$ is a finite
sum, the number depending on $r,k$ only, of expectations of products of (i) at most $r+1$ gap
coordinates and (ii) the factors $D^{j_c}\langle u_c,Y_c\rangle$, $c\in B$.

\emph{Step 2 (moments of the two kinds of factors).}
(i) Prefix gap coordinates have $L^s$ norms $\le C$ by Lemma~\ref{lem:subint}. Suffix gap coordinates
are bounded pathwise by $C(1+\norm{\mathbf G^n}_{\pvar})^{r}$, $p\in(2\rho,4)$, by the extension theorem
and the translation estimate \citep[Thm.~9.33(i)]{FrizVictoir2010}, with moments of every order bounded
uniformly in $n$ by Lemma~\ref{lem:rhovarcontract} and (C5). (ii) A cell coordinate
$\langle u_c,Y_c\rangle$ with $j'$ time letters and $g'$ letters of $G$ lies in the chaoses of order
$\le g'$ and has $L^2$ norm $\le C\delta^{j'+g'H}$ (Lemma~\ref{lem:subint} for $c<a$, directly for $c>a$);
$D^{j_c}$ of it vanishes unless $j_c\le g'$, and by (G1)--(G2)
$\norm{\norm{D^{j_c}\langle u_c,Y_c\rangle}}_{L^s}\le C\delta^{g'H}\le C\delta^{j_cH}$. The norm of the
expectation of an elementary tensor is at most the expectation of the product of the norms, and
H\"older gives
\[
   \norm{K_{\mathbf b}}_{\frH^{\otimes k}}\ \le\ C\prod_{c\in B}\delta^{j_cH}\ =\ C\delta^{kH},
\]
as $\sum_{c\in B}j_c=k$. The case $k=0$ (no derivative, the same H\"older bound on the product of
the prefix and suffix coordinates) gives $\norm\Phi_{L^s}\le C$.

\emph{Step 3 (the covariance bound).}
$\Phi$ lies in a finite sum of chaoses, so (G3) applies: $\langle f_k,\E D^k\Phi\rangle=\sum_{\mathbf b}
\langle f_k,P^{\otimes k}_{I_a}K_{\mathbf b}\rangle$. On $\bigotimes_l\frH(I_{b_l})$ the operator
$P^{\otimes k}_{I_a}$ has norm at most $\prod_l\varrho_{\max}(I_a,I_{b_l})$ (Lemma~\ref{lem:cellproj},
applied to each tensor factor), so, summing over $\mathbf b$ with all $b_l\ne a$,
\[
   |\langle f_k,\E D^k\Phi\rangle|\ \le\ C\delta^{kH}\norm{f_k}\Bigl(\sum_{b\ne a}\varrho_{\max}(I_a,I_b)\Bigr)^k
   \ \le\ C\norm{f_k}(L\delta^H)^k
\]
by the decorrelation clause of (C4). Summing over $1\le k\le q$ and using
$\E[F\Phi]-\E F\,\E\Phi=\sum_{k\ge1}\langle f_k,\E D^k\Phi\rangle$ from (G3) gives the last claim.
\end{proof}

With Lemmas~\ref{lem:subint} and~\ref{lem:ibp} the upper bound of Proposition~\ref{prop:rough}
extends from level $4$ to every level. This is the rough-regime case of the rate stated in the
Introduction on the centred class, and Proposition~\ref{prop:drift5} repeats the argument with a drift.

  \begin{proposition}[full truncation in the rough regime]\label{prop:roughalllevels}
Let $H\in(\tfrac14,\tfrac12)$ and $P\in\cP_0(H,L)$. For every $m\le M$,
$\norm{\cW_m}\le C(M,d,H,L)\,\delta^{4H-1}$. Hence $\norm{\cW_{\le M}}_\hs=O(\delta^{4H-1})$ uniformly over
$\cP_0(H,L)$, for every fixed $M$.
\end{proposition}

\begin{proof}
As in Proposition~\ref{prop:rough}, we telescope cell by cell and show that each cell contributes
$O(\delta^{4H})$. What changes is the treatment of the macroscopic factor: it is never expanded into
cell increments, but enters through Lemma~\ref{lem:ibp}, so the number of its legs does not matter.

\emph{Step 1 (telescoping).} With $\Delta_a$ as in \eqref{eq:celltelescope}, write
$E_a:=\Sig(X)_{I_a}-\exp\Delta_a$ for the defect of
cell $a$. Then \eqref{eq:celltelescope} is
\[
   \Sig(X)-\Sig(X^n)=\sum_{a=0}^{n-1}\Sig(X)_{0,t_a}\otimes E_a\otimes\Sig(X^n)_{t_{a+1},1} .
\]
As $\cW_m=\E\pi_m\Sig(X^n)-\E\pi_m\Sig(X)$, the coordinate of $\cW_m$ at a word $w$ of length $m$ is
\[
   \langle w,\cW_m\rangle=-\sum_{a=0}^{n-1}\ \sum_{\substack{w=uvz\\ |v|=q\ge2}}\E[F\Phi],\qquad
   F:=\langle v,E_a\rangle,\quad \Phi:=\langle u,\Sig(X)_{0,t_a}\rangle\langle z,\Sig(X^n)_{t_{a+1},1}\rangle,
\]
since levels $0$ and $1$ of $E_a$ vanish. There are fewer than $m^2$ splittings $w=uvz$, and $\Phi$ is
as in Lemma~\ref{lem:ibp}, with no time letters.

\emph{Step 2 (the cell variable).} $F$ is
$\cF_{I_a}$-measurable, in the chaoses of order $\le q$ and parity $q$, with $\norm F_{L^2}\le C\delta^{qH}$
(Lemma~\ref{lem:subint}. The term $\langle v,\exp\Delta_a\rangle$ is a monomial of degree $q$ in
$\Delta_a$), so its chaos components $F=\sum_{k\le q}I_k(f_k)$ satisfy $\norm{f_k}\le C\delta^{qH}$.
\emph{Step 3 (one cell, by the length $q$ of $v$).} Lemma~\ref{lem:ibp} bounds
$|\E[F\Phi]-\E F\,\E\Phi|$ by $C\sum_{k\ge1}\norm{f_k}(L\delta^H)^k$, and $\norm\Phi_{L^2}\le C$.
\begin{itemize}[leftmargin=1.4em,topsep=2pt,itemsep=1pt]
\item $q=2$: $F$ is a cell area ($v=ii$ gives $F=0$), so $\E F=0$ and only $f_2$ survives. The term is
$\le C\delta^{2H}(L\delta^H)^2=CL^2\delta^{4H}$.
\item $q=3$: $\E F=0$ by parity, and only $f_1,f_3$ survive, the term is
$\le C\delta^{3H}(L\delta^H+L^3\delta^{3H})$.
\item $q=4$: $|\E F\,\E\Phi|\le C\delta^{4H}$, and $f_2,f_4$ contribute $\le C\delta^{4H}(L\delta^H)^2\le C\delta^{6H}$.
\item $q\ge5$: by Cauchy--Schwarz, $|\E[F\Phi]|\le\norm F_{L^2}\norm\Phi_{L^2}\le C\delta^{5H}$.
\end{itemize}

\emph{Conclusion.} The main term of each cell is of order $\delta^{4H}$, coming from $q=2$, $q=3$
($k=1$) and the mean term of $q=4$. All other terms are $O(\delta^{5H})\subset O(\delta^{4H})$. Summing
over the $n=\delta^{-1}$ cells and the fewer than $m^2$ splittings,
\[
   \norm{\cW_m}\ \le\ C(M,d,H,L)\,n\,\delta^{4H}\ =\ C(M,d,H,L)\,\delta^{4H-1}.
\]
\end{proof}

\begin{remark}[truncation level $M\ge6$: what the cell bookkeeping misses]\label{rem:M6}
Run group by group as in Proposition~\ref{prop:rough}, the level-$m$ groups $q=2$ and $q=3$ acquire,
for $m\ge6$, $m-4\ge2$ macroscopic legs beyond those pairing with the excess, and these pair among
themselves over a free pair of cells. Summed in absolute value, that pairing costs
\[
   \sum_{b_1,b_2}\bigl|\Cov(\Delta_{b_1},\Delta_{b_2})\bigr|
   =\delta^{2H}\sum_{b_1,b_2}|r(b_1-b_2)|\ \asymp\ \delta^{2H-1}\sum_{k\in\ZZ}|r(k)|,
\]
with $\sum_k|r(k)|<\infty$ for $H<\tfrac12$. This would only give $O(\delta^{6H-2})$, and
$6H-2<4H-1$ for every $H<\tfrac12$. The signed sum is
$\sum_{b_1,b_2}\Cov(\Delta_{b_1},\Delta_{b_2})=\Var(\sum_b\Delta_b)=O(1)$, and
Proposition~\ref{prop:roughalllevels} keeps it by never expanding the free legs: they stay inside
the expectation as sub-interval signatures, whose moments Lemma~\ref{lem:subint} bounds.

Numerically, though this is not proved, the order $\delta^{4H-1}$ is attained at level $6$, so
for $M\ge6$ the leading constant of $\norm{\cW_{\le M}}$ need not be that of level $4$: for $d=3$
and $H\in\{0.3,0.35,0.4\}$, exact evaluation (up to $n=512$) shows that the ratios of successive
differences of the words $121233$ and $122133$ to those of $1212$ settle at nonzero limits, and
that for $121211$ stays bounded away from zero with only a slight drift, whereas an order
$\delta^{6H-2}$ would make these ratios grow like $\delta^{2H-1}$.
\end{remark}

\paragraph{The drifted class below the crossover.} For $\theta\ne0$ the level-$3$ bias does not
vanish. Above the crossover it is settled by Proposition~\ref{prop:level3drift}, and all levels by
Proposition~\ref{prop:alllevelsdrift}, both through hypothesis (K2) of Lemma~\ref{lem:singleerror}, an $L^1$
bound on the kernel error $e$ on every cell product, the diagonal ones included. Below $\tfrac12$
this fails on the diagonal cells, where $\varrho$ is not locally integrable
(Remark~\ref{rem:C1scope}). We therefore avoid $\varrho$ altogether. At truncation $4$ every
drift contribution involves a single rough component, and then the covariance function controls the
bias, as in \citet[Lem.~A.1]{HagerPelizzari2026} (Remark~\ref{rem:hp}). This also yields the
constants. At every truncation level the rate follows by running the argument of
Proposition~\ref{prop:roughalllevels} on the time-augmented path.

The next lemma bounds the covariance of an increment of $G$ with a value of $G$ from the class
conditions alone. The proof of Proposition~\ref{prop:drift4} uses it to bound the interpolation error
$R_n-R$ of the covariance by $2L_V\delta^{2H}$.

  \begin{lemma}[a covariance bound on the class]\label{lem:Vbound}
Let $H\in(\tfrac14,\tfrac12)$ and $P\in\cP(H,L)$. Then for each component of the Gaussian part and all
$0\le s<t\le1$, $u\in[0,1]$, $|\Cov(G_t-G_s,G_u)|\le L_V|t-s|^{2H}$ with $L_V=C_HL$.
\end{lemma}

\begin{proof}
We bound the covariance of one mesh cell with a node by (C2) and (C4), and then pass to general
intervals through their dyadic decomposition. We write $G$ for one
component and $G_I:=G_t-G_s$ for $I=[s,t]$.

\emph{Step 1 (a cell against a node).} For a mesh $\delta=1/n$ and a cell $I_a$: (C2) at $q=2$ gives $\Var(G_{I_a})\le L\delta^{2H}$, and the
first condition of (C4) gives $\sum_{b\ne a}|\Cov(G_{I_a},G_{I_b})|\le\sum_{b\ne a}\int_{I_a\times I_b}|\varrho|
\le L\delta^{2H}$. So, as $G_0=0$ and $G_{t_b}=\sum_{c<b}G_{I_c}$,
\[
   |\Cov(G_{I_a},G_{t_b})|\ \le\ \Var(G_{I_a})+\sum_{c\ne a}|\Cov(G_{I_a},G_{I_c})|\ \le\ 2L\delta^{2H}
   \qquad\text{for every node }t_b .
\]

\emph{Step 2 (H\"older bound in $L^2$).} An interval $[u,v]$  with dyadic endpoints is a
union of dyadic cells, at most two of each length $2^{-k}\le|v-u|$, each a mesh cell for $n=2^k$. By
Minkowski and $\Var(G_{I})\le L|I|^{2H}$ on each such cell,
\[
   \norm{G_v-G_u}_{L^2}\ \le\ 2\sqrt L\sum_{k:\,2^{-k}\le|v-u|}2^{-kH}\ \le\ K_H\sqrt L\,|v-u|^H,\qquad
   K_H:=\frac{2}{1-2^{-H}},
\]
and by $L^2$-continuity this holds for all $u,v$.

\emph{Step 3 (a dyadic cell against any value).} For a dyadic cell $J$ and any $u$, let $t_b$
be the node below $u$ of the mesh of width $|J|$ and split $G_u=G_{t_b}+(G_u-G_{t_b})$. Step~1 bounds
the first covariance by $2L|J|^{2H}$, and Cauchy--Schwarz with Step~2 bounds the second by
$\sqrt L|J|^H\cdot K_H\sqrt L|J|^H$. So
\[
   |\Cov(G_J,G_u)|\ \le\ (2+K_H)\,L\,|J|^{2H}.
\]

\emph{Step 4 (summing over the dyadic decomposition).} Summing Step~3 over the dyadic
decomposition of $[s,t]$, at most two cells of each length $2^{-k}\le|t-s|$, gives the claim with
$C_H=2(2+K_H)/(1-2^{-2H})$, first for dyadic $s,t$ and then for all by $L^2$-continuity.
\end{proof}

At truncation $4$ every drift contribution is an iterated integral $\int\!\!\int_{s<t}w\,\dd G^k_s\dd G^k_t$
of one Gaussian component against a smooth weight $w$. The next lemma writes its expectation as a
linear functional of the covariance, bounded by the sup norm. Propositions~\ref{prop:drift4}
and~\ref{prop:driftconst} use it.

  \begin{lemma}[quadratic functionals of one rough component]\label{lem:quadfun}
Let $w\in C^2(T)$, $T=\{0\le s\le t\le1\}$, $\omega(t)=w(t,t)$, and let $Z$ be a continuous real path
with $Z_0=0$ and finite $p$-variation. Then
\begin{multline*}
   \int\!\!\int_{s<t}w\,\dd Z_s\dd Z_t=\tfrac12\omega(1)Z_1^2-\tfrac12\!\int_0^1\!\omega'Z^2
   -Z_1\!\int_0^1\!\partial_1w(s,1)Z_s\dd s\\
   +\!\int_0^1\!\partial_1w(t,t)Z_t^2\dd t+\!\int\!\!\int_{s<t}\!\partial_1\partial_2w\,Z_sZ_t\,\dd s\,\dd t,
\end{multline*}
the left side being, for polynomial $w$, a combination of coordinates of the canonical lift of
$(t,Z_t)$ ($s^a(u-s)^b(1-u)^c\leftrightarrow a!\,b!\,c!\,\langle0^az0^bz0^c\rangle$, $0$ the time
letter), and for general $w$ defined by the right side. For a centred Gaussian $Z$ with covariance $Q$ the expectation is a linear functional
$\mathcal L_w(Q)$ with $|\mathcal L_w(Q)|\le4\norm w_{C^2(T)}\sup|Q|$, where $\norm w_{C^2(T)}$ is the
largest sup norm on $T$ of $w$ and of its partial derivatives of order $\le2$.
\end{lemma}

\begin{proof}
We prove the identity for smooth $Z$ by two integrations by parts, extend it by continuity, and
take expectations term by term.

\emph{Smooth $Z$.} For smooth $Z$ integrate by parts first in $s$: $\int_0^tw(s,t)\dot Z_s\dd s=\omega(t)Z_t-Y_t$,
with $Y_t:=\int_0^t\partial_1w(s,t)Z_s\dd s$, as $Z_0=0$. Then integrate in $t$:
$\int_0^1\omega Z\dot Z\,\dd t=\tfrac12\omega(1)Z_1^2-\tfrac12\int_0^1\omega'Z^2$, and, $Y$ being of
class $C^1$ with $Y_0=0$ and $\dot Y_t=\partial_1w(t,t)Z_t+\int_0^t\partial_1\partial_2w(s,t)Z_s\dd s$,
$\int_0^1Y\dot Z\,\dd t=Y_1Z_1-\int_0^1\dot YZ\,\dd t$. Subtracting the second from the first gives the
identity.

\emph{Extension.} Both sides are continuous under uniform convergence with bounded $p$-variation (the
left by continuity of the Young pairing, \citet[Ch.~9]{FrizVictoir2010}), so the identity extends by
density.
\emph{Expectation.} Taking expectations term by term, with $Q(s,t)=\E[Z_sZ_t]$,
\[
\begin{aligned}
   \mathcal L_w(Q)={}&\tfrac12\omega(1)Q(1,1)-\tfrac12\int_0^1\omega'(t)Q(t,t)\,\dd t
   -\int_0^1\partial_1w(s,1)Q(s,1)\,\dd s\\
   &+\int_0^1\partial_1w(t,t)Q(t,t)\,\dd t+\int\!\!\int_{s<t}\partial_1\partial_2w(s,t)\,Q(s,t)\,\dd s\,\dd t .
\end{aligned}
\]
As $|\omega|\le\norm w_{C^2}$, $|\omega'|=|(\partial_1w+\partial_2w)(t,t)|\le2\norm w_{C^2}$ and the
triangle $\{s<t\}$ has area $\frac{1}{2}$, the five terms are bounded by $\frac{1}{2},1,1,1,\frac{1}{2}$ times $\norm w_{C^2}\sup|Q|$,
which sum to $4$.
\end{proof}

The next proposition combines Lemmas~\ref{lem:Vbound} and~\ref{lem:quadfun}; Main result~2 of
the Introduction and the discussion of centring after Proposition~\ref{prop:level4} use it.

\begin{proposition}[the drifted rough regime at truncation $4$]\label{prop:drift4}
Let $H\in(\tfrac14,\tfrac12)$, $M\ge4$ and $P\in\cP(H,L)$ with drift $\theta$. Then $\cW_1=\cW_2=0$,
$\norm{\cW_3}\le C|\theta|L_V\delta^{2H}$ and $\norm{\cW_4(P)-\cW_4(P^{(0)})}\le C\theta^2L_V\delta^{2H}$,
with $C=C(d)$ and $L_V$ as in Lemma~\ref{lem:Vbound},
where $P^{(0)}\in\cP_0(H,L')$, $L'=L'(L,M,H,d)$, is the law of the Gaussian
part $G$. Hence $\norm{\cW_{\le4}(P)}=\norm{\cW_4(P^{(0)})}+O(\delta^{2H})=O(\delta^{4H-1})$
uniformly over $\cP(H,L)$, and for the drifted independent-fBm law the level-$\le4$ constant of
Corollary~\ref{cor:kadj} is unchanged. For that law
$\cW_3=\theta\delta^{2H}(-\frac{1}{12},\frac{1}{6},-\frac{1}{12})+O(\delta^{2H+1}+\delta^2)$ on the words
$1kk,k1k,kk1$ ($k\ge2$), all other level-$3$ words having zero bias.
\end{proposition}
The conditions (C1)--(C5) concern $G$ only, so $P^{(0)}$ satisfies them with the same $L$.
The moment condition of Definition~\ref{def:class} for $G$ holds by
Remark~\ref{rem:C2simplify}(ii) with a constant depending on $L,M,H,d$, and $L'$ is the larger of
$L$ and this constant. With Proposition~\ref{prop:level3drift}, the
level-$3$ formula holds on all of $(\tfrac14,1)$. At level $5$ the drift terms are again
$O(\delta^{4H-1})$, and for drifted independent fBm of exactly that order (proved for $d\ge3$;
Propositions~\ref{prop:drift5} and~\ref{prop:attain5}).

\begin{proof}
We expand every coordinate in powers of $\theta$. A parity argument removes all terms with an
odd number of slots of some component $G^i$, and the surviving drift terms are quadratic functionals
of one component, which Lemmas~\ref{lem:quadfun} and~\ref{lem:Vbound} bound by $O(\delta^{2H})$.

\emph{Step 1 (expansion in $\theta$ and parity).}
$\Sig(\theta te_1+G)$ is the image of the signature of $(t,G)$ under the linear map sending the time
letter to $\theta e_1$, and likewise for the interpolant, the drift $\theta te_1$ being linear,
hence interpolated exactly. So, the
$\theta^j$ coefficient of each coordinate is a sum of coordinates of $\Sig(t,G)$ with $j$ time letters.
Such a coordinate with an odd number of $G^i$-slots for some $i$ has mean zero for $X$ and $X^n$: the
law of $G$ is invariant under $G^i\mapsto-G^i$ and the lift, a limit of interpolations, is
equivariant.

\emph{Step 2 (levels $1$ to $3$).} Level $1$ is node-exact. At level $2$, $\langle ii\rangle$ is endpoint-determined and
$\langle ik\rangle$ has an odd $G$-count in every $\theta$-power. At level $3$ the survivors are one
time letter and a $G^k$-pair, $k\ge2$ (the word $111$ is endpoint-determined): integrating the time
slot out gives $\int\!\!\int_{s<u}w\,\dd G^k\dd G^k$ with $w\in\{s,u-s,1-u\}$, whose bias is
$\mathcal L_w(R_n-R)$ by Lemma~\ref{lem:quadfun}, $R$ and $R_n$ being the covariances of $G^k$ and
$G^{n,k}$. With $\lambda=(s-t_a)/\delta$ for $s\in I_a$ and $\varepsilon_s=G^n_s-G_s=\lambda(G_{t_{a+1}}-G_s)
-(1-\lambda)(G_s-G_{t_a})$, $R_n-R=\E[\varepsilon_sG^n_u]+\E[G_s\varepsilon_u]$, and Lemma~\ref{lem:Vbound}
bounds each term by $L_V\delta^{2H}$ ($G^n_u$ being a convex combination of two values of $G$), so
\[
   \sup|R_n-R|\ \le\ 2L_V\delta^{2H},\qquad |\mathcal L_w(R_n-R)|\ \le\ 8\norm w_{C^2(T)}L_V\delta^{2H}.
\]

\emph{Step 3 (level $4$).} At level $4$ the $\theta^0$
part is $\cW_4(P^{(0)})$, the $\theta^1,\theta^3$ parts have odd $G$-counts, the $\theta^4$ part is
exact, and the $\theta^2$ part has two time letters and one $G^k$-pair, $k\ge2$, i.e.\ weights
$w=s^{j_1}(u-s)^{j_2}(1-u)^{j_3}/(j_1!j_2!j_3!)$, $j_1+j_2+j_3=2$, again $O(L_V\delta^{2H})$ by
Step~2. As
$2H>4H-1$, the drift contributions are lower order. With $\theta^2\le L^2$ and
Proposition~\ref{prop:rough} applied to $P^{(0)}$, this gives the bound on $\norm{\cW_{\le4}(P)}$.

\emph{Step 4 (the level-$3$ constants for fBm).} By Step~2 the biases at $1kk,k1k,kk1$ are
$\theta$ times the biases of $\int\!\!\int_{s<u}w\,\dd G^k\dd G^k$ with $w=s,\,u-s,\,1-u$, for which
$(\partial_1w-\partial_2w)(t,t)=1,-2,1$. Proposition~\ref{prop:driftconst}, whose proof does not use
the present one, gives $\theta\delta^{2H}(-\frac1{12},\frac16,-\frac1{12})+O(\delta^{2H+1}+\delta^2)$.
\end{proof}

From level $5$ on we run the argument of Proposition~\ref{prop:roughalllevels} on the time-augmented
path $Z=(t,G)$; no hypothesis beyond the class is needed.

\begin{proposition}[the drifted rough regime at every truncation level]\label{prop:drift5}
Let $H\in(\tfrac14,\tfrac12)$ and $P\in\cP(H,L)$ with drift $\theta$. Then
$\norm{\cW_m(P)}\le C(M,d,H,L)\,\delta^{4H-1}$ for every $m\le M$, uniformly over $\cP(H,L)$.
\end{proposition}

\begin{proof}
As in Proposition~\ref{prop:roughalllevels}, now for the time-augmented path $Z=(t,G)$. What
changes is that the cell excess may contain time letters, and each of them gains a factor $\delta$.
As in Step~1 of the proof of Proposition~\ref{prop:drift4}, the coefficient of $\theta^j$ in $\cW_m(P)$ is a sum of
at most $2^m$ $Z$-biases at words with $j$ time letters, and $|\theta|^j\le(1\vee L)^m$. Run the telescoping
\eqref{eq:celltelescope} for $Z$, with $E_a=\Sig(Z)_{I_a}-\exp(\delta e_0+\Delta_a)$
and $\Delta_a=G_{t_{a+1}}-G_{t_a}$, and split each
term by the number $j_E$ of time letters in the excess word $v$ ($|v|=q$, $g=q-j_E$ letters of $G$).

\emph{No time letter in the excess, $j_E=0$.}
The proof of Proposition~\ref{prop:roughalllevels} applies word for word, Lemma~\ref{lem:ibp}
allowing time letters in the prefix and suffix words. Each term is $O(\delta^{4H})$.

\emph{At least one time letter in the excess, $j_E\ge1$.} $g=0$ gives $F=0$, time being interpolated exactly. Otherwise $F$ lies in the chaoses of
order $\le g$ and parity $g$ with $\norm F_{L^2}\le C\delta^{j_E+gH}$ (Lemma~\ref{lem:subint}), and by
Lemma~\ref{lem:ibp}: for $g=1$, $\E F=0$ and the $k=1$ term is $\le C\delta^{j_E+H}L\delta^H$; for $g=2$,
$|\E F\,\E\Phi|+\norm{f_2}L^2\delta^{2H}\le C\delta^{j_E+2H}$; for $g=3$, $\E F=0$ and $f_1,f_3$ give
$\le C\delta^{j_E+3H}(L\delta^H+L^3\delta^{3H})$; for $g\ge4$, Cauchy--Schwarz gives $\le C\delta^{j_E+4H}$.
Each is $\le C\delta^{1+2H}$ per cell, $O(\delta^{2H})=o(\delta^{4H-1})$ after summation.

\emph{Conclusion.} Summing over the $n$ cells, the finitely many words and the
$\theta$-powers gives $\norm{\cW_m(P)}\le C(M,d,H,L)\,\delta^{4H-1}$.
\end{proof}

\noindent
So the terms with a time letter in the cell excess are of lower order, while those with time letters
only in the macroscopic factor are not: they can attain $\delta^{4H-1}$ (Proposition~\ref{prop:attain5}).

\begin{proposition}[the level-$5$ drift bias is attained]\label{prop:attain5}
For the drifted independent-fBm law with $d\ge3$, drift $\theta e_1$ and $H\in(\tfrac14,\tfrac12)$, the
coordinate of $\cW_5$ at the word $12323$ equals
$\theta\,\frac{\kappa_{\rm alt}(H)}{2}\,\delta^{4H-1}\bigl(1+O(\delta^{1-2H})\bigr)$, with
$\kappa_{\rm alt}(H)<0$ (Proposition~\ref{prop:kappasign}). In particular
$\norm{\cW_5}\ge c|\theta|\delta^{4H-1}$.
\end{proposition}

\begin{proof}
The drift enters once, as a time letter in front of the word $2323$. Splitting the remaining
interval at the next node, the bias becomes a weighted sum over cells of the level-$4$ differences
$\Delta E_{\rm alt}$ at the coarser meshes $1/\ell$, $\ell<n$, each evaluated by
Proposition~\ref{prop:Lexists}.

\emph{Step 1 (reduction to one time letter).}
Only the $\theta^1$ term survives, the letter $1$ occurring once and the $\theta^0$ term having $G^1$
once. It is $\theta$ times the bias of $\langle0\,2323\rangle=\int_0^1\langle2323,\Sig(G)_{t,1}\rangle\dd t$,
likewise for the interpolant.

\emph{Step 2 (split at the next node).} For $t\in I_c$ put $u=t_{c+1}$ and $\lambda=(t-t_c)/\delta$. Then
$\Sig(G^n)_{t,1}=\exp((1-\lambda)\Delta_c)\otimes\Sig(G^n)_{u,1}$ exactly, $\Delta_c=G_{t_{c+1}}-G_{t_c}$.
By Chen's relation,
\[
   \langle2323,\Sig(G^n)_{t,1}\rangle=\sum_{2323=w'w''}\langle w',\exp((1-\lambda)\Delta_c)\rangle\,
   \langle w'',\Sig(G^n)_{u,1}\rangle ,
\]
the sum over the five deconcatenations, $w'$ being the part on the short piece $[t,u]$. As $u$ is a node, stationarity of
increments and self-similarity give
\[
   \E\langle2323,\Sig(G^n)_{u,1}\rangle=(1-u)^{4H}\bigl(E_{\rm alt}+\Delta E_{\rm alt}^{[\ell]}\bigr),\qquad
   \E\langle2323,\Sig(G)_{t,1}\rangle=(1-t)^{4H}E_{\rm alt},
\]
with
$\Delta E_{\rm alt}^{[\ell]}$ the value of $\Delta E_{\rm alt}$ at mesh $1/\ell$, $\ell=n-c-1$
($\Delta E^{[0]}_{\rm alt}:=0$), and $E_{\rm alt}$ the continuum value of the word $1212$, equal to that
of $2323$.

\emph{Step 3 (the short piece is of lower order).} The terms with
$|w'|\ge1$ letters on the short piece are $O(\delta^{2H})$ uniformly in $t$. For $|w'|=1$, $w'=2$, and
Gaussian integration by parts in the increments (the case $k=1$ of (G3)) gives
\[
   \E[\Delta^2_c\Psi]=\sum_{b>c}\Cov(\Delta_c,\Delta_b)\,\E\partial_{\Delta^2_b}\Psi,\qquad
   \sum_{b}|\Cov(\Delta_c,\Delta_b)|\le C\delta^{2H}.
\]
Here $\Psi$ is a coordinate of $\Sig(G^n)_{u,1}$, and
by Chen $\partial_{\Delta^2_b}\Psi$ is a finite sum of products of a coordinate of $\Sig(G^n)_{u,t_b}$,
a polynomial of degree at most $2$ in $\Delta_b$ with bounded coefficients, and a coordinate of
$\Sig(G^n)_{t_{b+1},1}$. By Cauchy--Schwarz, H\"older and the uniform moments of
Lemmas~\ref{lem:rhovarcontract} and~\ref{lem:gaussmom},
$|\E\partial_{\Delta^2_b}\Psi|\le\norm{\partial_{\Delta^2_b}\Psi}_{L^2}\le C$ uniformly in $b,c,n$. So the
term $|w'|=1$ is $O(\delta^{2H})$. For the others, for $|w'|=2$ the bound is
$(\sum_b|\Cov(\Delta_c,\Delta_b)|)^2$, for $|w'|=3$ it is $\Var(\Delta_c)|\Cov(\Delta_c,G_1-G_u)|$, for
$|w'|=4$ it is $\E[(\Delta^2_c)^2(\Delta^3_c)^2]$, all $O(\delta^{4H})$.

\emph{Step 4 (sum over cells).} Integrating over $t$ and subtracting the continuum value,
the bias of $\langle0\,2323\rangle$ is
\[
   \sum_{c=0}^{n-1}\delta(1-t_{c+1})^{4H}\Delta E^{[n-c-1]}_{\rm alt}+O(\delta)+O(\delta^{2H}),
\]
the $O(\delta)$ being
the Riemann error of $\int(1-t)^{4H}E_{\rm alt}\dd t$ and the $O(\delta^{2H})$ coming from Step~3. By Proposition~\ref{prop:Lexists} and $\kappa_{\rm alt}\ne0$ (Proposition~\ref{prop:kappasign}),
$\Delta E^{[\ell]}_{\rm alt}=\kappa_{\rm alt}\ell^{1-4H}(1+O(\ell^{2H-1}))$ for $\ell\ge\ell_0$, and each of the finitely many
$\ell<\ell_0$ has a finite value, so the bound holds for all $\ell\ge1$ with a larger constant. With
$1-t_{c+1}=\ell\delta$, so that $\delta(\ell\delta)^{4H}\ell^{1-4H}=\delta^{4H+1}\ell$, and with
$\sum_{\ell<n}\ell=\frac{n(n-1)}{2}$, $\sum_{\ell<n}\ell^{2H}=O(n^{2H+1})$, the sum is
\[
   \kappa_{\rm alt}\delta^{4H+1}\sum_{\ell=1}^{n-1}\ell\,\bigl(1+O(\ell^{2H-1})\bigr)
   =\frac{\kappa_{\rm alt}}{2}\delta^{4H-1}\bigl(1+O(\delta^{1-2H})\bigr).
\]

\emph{Conclusion.} The main term is $\theta\frac{\kappa_{\rm alt}}{2}\delta^{4H-1}$. The remainders
$O(\delta)$ and $O(\delta^{2H})$ of Step~4 are $O(\delta^{2-4H})$ and
$O(\delta^{1-2H})$ relative to $\delta^{4H-1}$, and $1-2H<2-4H$, so both are absorbed in the
relative error $O(\delta^{1-2H})$.
\end{proof}

The next proposition evaluates the leading constant of the drift biases of
Proposition~\ref{prop:drift4} for fBm, on both sides of the crossover. It gives the level-$3$ constants
stated there and the level-$4$ $\theta^2$ constants of Remark~\ref{rem:dscope}.

  \begin{proposition}[the drift constant, for fBm]\label{prop:driftconst}
Let $G$ be fBm$(H)$, $H\in(\tfrac14,1)$, and $w\in C^2(T)$. Then
\[
   \E\!\int\!\!\int_{s<t}w\,\dd G^n_s\dd G^n_t-\E\!\int\!\!\int_{s<t}w\,\dd G_s\dd G_t
   =-\frac{\delta^{2H}}{12}\int_0^1(\partial_1w-\partial_2w)(t,t)\,\dd t+O(\delta^{2H+1}+\delta^2).
\]
\end{proposition}
This is the $\Phi^{\rm s}=-\Phi^{\rm t}=-\frac{1}{12}$ structure of Lemma~\ref{lem:singleerror}, now on
both sides of the crossover.

\begin{proof}
We need the covariance error $Q:=R_n-R$ on the diagonal, on the line $u=1$, and in the bulk. Only the
diagonal gives a term of order $\delta^{2H}$: there the interpolant loses the variance
$\lambda(1-\lambda)\delta^{2H}$, whose cell average is $\frac{1}{6}\delta^{2H}$. Everything else is
controlled by the linear interpolation errors of $t^{2H}$ and $(1-t)^{2H}$ and by the bilinear
interpolation error of $|u-s|^{2H}$ off the diagonal band.

\emph{Step 1 (the covariance error).}
By Lemma~\ref{lem:quadfun} the left side is $\mathcal L_w(Q)$ with $Q=R_n-R$, a sum of five terms, and
$R_n=(I_n\otimes I_n)R$, $I_n$ the piecewise-linear interpolation at the nodes as in
Proposition~\ref{prop:adjexact}, and $\lambda=\lambda(t)=(t-t_a)/\delta$ for $t\in I_a$. Write $R(s,u)=\frac{1}{2}(f(s)+f(u)-g(u-s))$ with $f(x)=x^{2H}$, $g(x)=|x|^{2H}$.
Then
\[
   Q(s,u)=\tfrac12\bigl(\epsilon_f(s)+\epsilon_f(u)\bigr)-\tfrac12\epsilon^{\rm B}_g(s,u),
\]
with $\epsilon_f=I_nf-f$ and
$\epsilon^{\rm B}_g$ the bilinear interpolation error of $g(u-s)$.

\emph{Step 2 (the linear interpolation errors).}
 $\int|\epsilon_f|=O(\delta^{2H+1}+\delta^2)$: cell $0$ contributes
$O(\delta^{2H+1})$ and the others at most $\delta^3\sum_{a\ge1}|f''(a\delta)|$, which is $O(\delta^{2H+1})$ for
$H<\tfrac12$ and $O(\delta^2)$ for $H>\tfrac12$ ($\epsilon_f=0$ at $H=\tfrac12$). Likewise for $\tilde f(s)=(1-s)^{2H}$.

\emph{Step 3 (the line $u=1$).} $Q(1,1)=0$. At the node $u=1$ bilinear
interpolation reduces to linear interpolation, so $Q(s,1)=\frac{1}{2}(\epsilon_f(s)-\epsilon_{\tilde f}(s))$,
and the term of $\mathcal L_w$ on this line is $O(\int|\epsilon_f|+\int|\epsilon_{\tilde f}|)$ by Step~2.

\emph{Step 4 (the diagonal: the main term).} On
the diagonal the bilinear interpolant of $g$ equals $2\lambda(1-\lambda)\delta^{2H}$ (the two
off-diagonal corners of $I_a\times I_a$ carry $g=\delta^{2H}$, each with weight $\lambda(1-\lambda)$), so
$Q(t,t)=\epsilon_f(t)-\lambda(1-\lambda)\delta^{2H}$, and the two diagonal terms of $\mathcal L_w$ combine to
\[
\begin{gathered}
   \int_0^1\varphi(t)\,Q(t,t)\,\dd t=-\delta^{2H}\int_0^1\varphi\,\lambda(1-\lambda)\,\dd t+O\Bigl(\int|\epsilon_f|\Bigr),\\
   \varphi(t):=\partial_1w(t,t)-\tfrac12\omega'(t)=\tfrac12(\partial_1w-\partial_2w)(t,t),
\end{gathered}
\]
$\varphi$ being of class $C^1$. On
each cell $\int_{I_a}\varphi\lambda(1-\lambda)=\frac{1}{6}\int_{I_a}\varphi+O(\delta^2\norm{\varphi'}_\infty)$, since
$\int_{I_a}\lambda(1-\lambda)=\frac{\delta}{6}$ and $\varphi$ varies by at most $\delta\norm{\varphi'}_\infty$ on
$I_a$. Summing over the $n$ cells, the diagonal terms equal
$-\frac{\delta^{2H}}{12}\int_0^1(\partial_1w-\partial_2w)(t,t)\,\dd t+O(\delta^{2H+1}+\int|\epsilon_f|)$, the main term of the display.

\emph{Step 5 (the double integral).} In the double integral $\epsilon_f$ contributes $O(\int|\epsilon_f|)$; on
cell products at lag $\ge2$, $|\epsilon^{\rm B}_g|\le C\delta^2|u-s|^{2H-2}$, integrating to $O(\delta^{2H+1})$
for $H<\tfrac12$ and $O(\delta^2)$ for $H>\tfrac12$ (zero at $H=\tfrac12$). On the $O(n)$ cell products at
lag $\le1$, $|\epsilon^{\rm B}_g|\le C\delta^{2H}$ on area $\delta^2$ each, $O(\delta^{2H+1})$.

\emph{Conclusion.} The main term comes from Step~4; Steps~2, 3 and~5 contribute
$O(\delta^{2H+1}+\delta^2)$, which is $o(\delta^{2H})$. The constants
depend on $\norm w_{C^2}$ and $H$. They stay bounded as $H\to\tfrac12$: the second derivatives
$f''$ and $g''$ carry the factor $2H(2H-1)$, which compensates the growth of the lag sums of
Steps~2 and~5 near the crossover: $|2H-1|\sum_{a=1}^na^{2H-2}\le2$ for $H<\tfrac12$ and
$\le n^{2H-1}$ for $H>\tfrac12$, so the remainder is at most $C(\delta^{2H+1}+\delta^2)$ with $C$
depending on $\norm w_{C^2}$ and bounded for $H$ in compact subsets of $(\tfrac14,1)$.
\end{proof}

\noindent
This gives the level-$3$ constants of Proposition~\ref{prop:drift4} and the level-$4$ $\theta^2$
constants $(-\frac{1}{24},\frac{1}{12},-\frac{1}{12},0,\frac{1}{12},-\frac{1}{24})\theta^2\delta^{2H}$ on
$11kk,1k1k,1kk1,k11k,k1k1,kk11$, on all of $(\tfrac14,1)$.

\section{The strong error}\label{sec:strong}

We now compare the weak error with the strong error
$\eta_\delta=\bigl(\E\norm{\Sig^{\le M}(X^n)-\Sig^{\le M}(X)}_\hs^2\bigr)^{1/2}$ of \S\ref{sec:setup}. On the
Young range we determine its order for fBm drivers; in mean square it lies a full power of $\delta$
above the squared weak error. Below the crossover we bound it uniformly over the class, with an
exponent that cannot be improved.

\subsection{The Young range}\label{subsec:strongyoung}

Throughout this subsection $\tilde A_a$ denotes the L\'evy area over $[a,a+1]$ of a pair
$(B^1,B^2)$ of independent $\mathrm{fBm}(H)$, and $\rho_A(m):=\Cov(\tilde A_0,\tilde A_m)$, $m\in\ZZ$. By
stationarity of the increments and self-similarity, the L\'evy areas $A_a$ of $(X^1,X^2)$ over the
cells $I_a$ of an independent-fBm driver satisfy $\Cov(A_a,A_b)=\delta^{4H}\rho_A(|a-b|)$.
The proof of Proposition~\ref{prop:strong} sums these covariances over all pairs of cells. The
next lemma gives their decay in the lag $m$, so that the lag series $\sum_m\rho_A(m)$ and
$\sum_m|m|\rho_A(m)$ that control it converge absolutely.

\begin{lemma}[decay of the cross-cell area correlations]\label{lem:rhoAasym}
For $H\in(\tfrac14,1)$ and $c_H=H(2H-1)$,
\[
   \rho_A(m)\ =\ \frac{c_H^2\,(2H-2)}{72}\,m^{4H-6}\,\bigl(1+O(m^{-1})\bigr)
   \qquad (m\to\infty).
\]
In particular $\rho_A(m)<0$ for all large $m$ when $H\ne\tfrac12$, and $\sum_m|\rho_A(m)|$ and
$\sum_m|m|\,|\rho_A(m)|$ converge.
\end{lemma}

\begin{proof}
\emph{Idea.} We write $\rho_A(m)$ as a four-fold integral of two increment kernels against two
sign kernels, and Taylor-expand the increment kernels in the displacements inside the cells. The sign
kernels are antisymmetric, so they annihilate the terms of order $0$ and $1$. The leading term is of
order $2$.

\emph{Step 1 (the covariance as a four-fold integral).} For cells at lag $|a-b|\ge1$ the kernel $\varrho(s,s')\varrho(t,t')$, $\varrho(u,v)=c_H|u-v|^{2H-2}$, is
integrable on $(I_a\times I_b)^2$, its only singularities being corners. For the interpolants on a
finer mesh the cell areas are double sums of products of increments of $X^i$ and $X^j$, so Isserlis'
theorem gives their covariance exactly. Letting the mesh go to zero, by dominated convergence and the
$L^2$ convergence of the areas (Lemma~\ref{lem:L1conv}, after rescaling time by $1/n$),
\begin{gather*}
   \Cov(A_a,A_b)=\tfrac14\!\iint_{I_a}\!\iint_{I_b}\!\mathrm{sgn}(t-s)\,\mathrm{sgn}(t'-s')\,
   \varrho(s,s')\,\varrho(t,t')\,\dd s\,\dd t\,\dd s'\,\dd t' .
\end{gather*}
 \emph{Step 2 (Taylor expansion at lag $m\ge2$).} Take the unit cells $[0,1]$ and $[m,m+1]$, $m\ge2$,
write $s=u_1$, $t=u_2$, $s'=m+v_1$, $t'=m+v_2$ with $u=(u_1,u_2)$, $v=(v_1,v_2)\in[0,1]^2$, and
$\varrho(u,v)=g(|u-v|)$ with $g(r)=c_Hr^{2H-2}$. Then $\varrho(s,s')=g(m+x)$ and $\varrho(t,t')=g(m+y)$
with the displacements $x:=v_1-u_1$, $y:=v_2-u_2\in[-1,1]$. The function $g$ is smooth on
$[m-1,m+1]\subset[1,\infty)$, with $|g^{(j)}|\le C_jm^{2H-2-j}$ there, so Taylor's formula to second
order in each factor gives
\[
   g(m+x)\,g(m+y)\ =\ g^2+g g'\,(x+y)+g'^2\,xy+\tfrac12\,g g''\,(x^2+y^2)+R_3,
   \qquad |R_3|\le C\,m^{4H-7},
\]
with $g,g',g''$ evaluated at $m$. Every term of $R_3$ carries at least three derivatives in total.

\emph{Step 3 (the sign kernels kill orders $0$ and $1$).} Write $\sigma(u):=\mathrm{sgn}(u_2-u_1)$. The
tools are the three moments
\[
   \int_{[0,1]^2}\sigma(u)\,\dd u=0,\qquad\int_{[0,1]^2}\sigma(u)\,u_1\,\dd u=-\tfrac16,\qquad
   \int_{[0,1]^2}\sigma(u)\,u_2\,\dd u=\tfrac16 ;
\]
the second is $\int_0^1u_1[(1-u_1)-u_1]\,\dd u_1=\tfrac12-\tfrac23$, and the third follows by exchanging
$u_1$ and $u_2$, which flips the sign of $\sigma$. A monomial in $(u,v)$ integrates against
$\sigma(u)\sigma(v)$ to the product of its $u$-integral and its $v$-integral, so by the first moment it
contributes only if it contains a $u$-variable and a $v$-variable. This removes $g^2$ and $gg'(x+y)$
entirely. In $xy=v_1v_2-v_1u_2-u_1v_2+u_1u_2$ only $-v_1u_2-u_1v_2$ survives, and in
$x^2=v_1^2-2u_1v_1+u_1^2$, $y^2=v_2^2-2u_2v_2+u_2^2$ only the mixed terms survive:
\[
   \iint\sigma(u)\sigma(v)\,xy\ =\ 2\cdot\tfrac16\cdot\tfrac16=\tfrac2{36},\qquad
   \iint\sigma(u)\sigma(v)\,(x^2+y^2)\ =\ -2\bigl(\tfrac16\bigr)^2-2\bigl(\tfrac16\bigr)^2=-\tfrac4{36}.
\]
With the prefactor $\tfrac14$ of the four-fold integral,
\[
   \rho_A(m)\ =\ \tfrac14\Bigl(\tfrac2{36}\,g'^2-\tfrac12\cdot\tfrac4{36}\,g g''\Bigr)+O(m^{4H-7})
   \ =\ \tfrac1{72}\bigl(g'(m)^2-g(m)g''(m)\bigr)+O(m^{4H-7}).
\]

\emph{Conclusion.} From $g'(r)=c_H(2H-2)r^{2H-3}$ and $g''(r)=c_H(2H-2)(2H-3)r^{2H-4}$,
\[
   g'^2-g g''\ =\ c_H^2(2H-2)\,m^{4H-6}\bigl[(2H-2)-(2H-3)\bigr]\ =\ c_H^2(2H-2)\,m^{4H-6},
\]
which is the main term of the statement. The remainder $O(m^{4H-7})$ is the factor $1+O(m^{-1})$ when
$c_H\ne0$, and at $H=\tfrac12$ both vanish ($g\equiv0$, independent increments).
Only smoothness of $\varrho$ off the diagonal is used, so the asymptotic holds on all
of $(\tfrac14,1)$. Finally $4H-6<-1$ and $4H-5<-1$ for $H<1$, so both series converge.
\end{proof}

So $|\rho_A(m)|\asymp m^{4H-6}$, two powers faster than the naive $m^{4H-4}$, which ignores the
antisymmetry of the area and would not be summable for $H\ge\tfrac34$.

The lower bound in Proposition~\ref{prop:strong} needs the long-run variance of the areas to be
strictly positive.

\begin{lemma}[the area long-run variance is strictly positive]\label{lem:lrv}
Let $H\in(\tfrac12,1)$. The long-run variance $\sigma^2_\infty:=\sum_{m\in\ZZ}\rho_A(m)$ is finite and
strictly positive. Explicitly,
\[
   \sigma^2_\infty\ \ge\ 2\pi\,C_H^2\!\int_\RR\frac{(\lambda-\sin\lambda)^2}{\lambda^4}\,|\lambda|^{2-4H}\,\dd\lambda\ >\ 0,
\]
with $C_H:=\Gamma(2H+1)\sin(\pi H)/(2\pi)>0$, the constant for which
$c_H|u|^{2H-2}=C_H\int_\RR e^{\mathrm iu\lambda}|\lambda|^{1-2H}\,\dd\lambda$ for $u\ne0$ (an improper
integral).
\end{lemma}
\begin{proof}
\emph{Idea.} The sequence $(\rho_A(m))_{m\in\ZZ}$ is the sequence of Fourier coefficients of a
nonnegative spectral density $\tilde f$ on $[-\pi,\pi]$, and the Fej\'er means of $\sum_m\rho_A(m)$ show that
$\sigma^2_\infty$ is at least $2\pi$ times the value of $\tilde f$ near frequency $0$. We compute $\tilde f$ from
the harmonizable representation of fBm and evaluate its lower bound at $0$ in closed form.

Finiteness is the absolute summability of Lemma~\ref{lem:rhoAasym}.
\emph{Step 1 (the areas as double Wiener integrals).} For
positivity use the harmonizable representation
$\dd B^k_s=\int_\RR e^{\mathrm is\lambda}\sqrt{\hat\varrho(\lambda)}\,\dd W^k(\lambda)\,\dd s$, $k=1,2$,
with independent
complex Gaussian noises $\E|\dd W^k(\lambda)|^2=\dd\lambda$ and spectral density
$\hat\varrho(\lambda)=C_H
|\lambda|^{1-2H}$, so that $\varrho(u,v)=\int_\RR e^{\mathrm i(u-v)\lambda}\hat\varrho(\lambda)\,\dd\lambda$. Precisely, $\int_I\int_J\varrho=\int_\RR\widehat{\one_I}\,\overline{\widehat{\one_J}}\,\hat\varrho\,\dd\lambda$
for bounded intervals $I,J$, with $\widehat{\one_I}(\lambda)=\int_Ie^{\mathrm is\lambda}\dd s$. Since the
integral defining $C_H$ converges only as an improper integral, one inserts $e^{-\varepsilon|\lambda|}$:
then Fubini applies, the damped kernel $2C_H\Gamma(2-2H)\,\mathrm{Re}[(\varepsilon-\mathrm i(u-v))^{2H-2}]$
is bounded by $2C_H\Gamma(2-2H)|u-v|^{2H-2}$, and $|\widehat{\one_I}(\lambda)|\le\min(|I|,2/|\lambda|)$, so
dominated convergence on both sides gives the identity as $\varepsilon\downarrow0$. The covariance
determines the law of the increments, so this realisation may be used.
The value of $C_H$ in the statement, the spectral constant of fractional Gaussian noise,
follows from the classical formula $\int_\RR e^{\mathrm iu\lambda}|\lambda|^{-\alpha}\dd\lambda=2\Gamma(1-\alpha)
\sin(\tfrac{\pi\alpha}2)|u|^{\alpha-1}$, $\alpha=2H-1\in(0,1)$, and the reflection formula for
$\Gamma$. Then
\begin{gather*}
   \tilde A_a=\iint_{\RR^2}e^{\mathrm ia(\lambda+\xi)}\,\Upsilon(\lambda,\xi)\sqrt{\hat\varrho(\lambda)\hat\varrho(\xi)}\,
   \dd W^1(\lambda)\,\dd W^2(\xi),\\
   \Upsilon(\lambda,\xi)=\tfrac12\!\iint_{[0,1]^2}\!\mathrm{sgn}(t-s)\,e^{\mathrm i(s\lambda+t\xi)}\,\dd s\,\dd t,
\end{gather*}
the noises being Hermitian, $\overline{\dd W^k(\lambda)}=\dd W^k(-\lambda)$, so $\tilde A_a$ is
real.
\emph{Step 2 (the spectral density of $\rho_A$).} The It\^o isometry for independent
noises gives

\[
  \rho_A(m)=\iint_{\RR^2}e^{-\mathrm im(\lambda+\xi)}
  |\Upsilon(\lambda,\xi)|^2\hat\varrho(\lambda)\hat\varrho(\xi)\,\dd\lambda\,\dd\xi .
\]
Pushing forward under $\omega=\lambda+\xi$ and periodizing (valid since $m\in\ZZ$, so that
$e^{-\mathrm im\omega}$ is $2\pi$-periodic in $\omega$),
$\rho_A(m)=\int_{-\pi}^{\pi}e^{-\mathrm im\theta}\tilde f(\theta)\,\dd\theta$ with nonnegative spectral
density $\tilde f(\theta)=\sum_{k\in\ZZ}q(\theta+2\pi k)\ge q(\theta)$, where
$q(\omega)=\int_\RR|\Upsilon(\lambda,\omega-\lambda)|^2\hat\varrho(\lambda)\hat\varrho(\omega-\lambda)
\,\dd\lambda\ge0$.

\emph{Step 3 (Fej\'er means: $\sigma^2_\infty\ge2\pi q(0)$).} For every $\lambda\ne0$ the integrand of $q(\omega)$ is nonnegative and continuous in
$\omega$ at $\omega=0$, so Fatou's lemma gives $\liminf_{\omega\to0}q(\omega)\ge q(0)$. With the Fej\'er kernel
$K_N(\theta)=\sum_{|m|<N}(1-\tfrac{|m|}N)e^{-\mathrm im\theta}=\tfrac1N\bigl(\sin(N\theta/2)/\sin(\theta/2)\bigr)^2\ge0$,
$\int_{-\pi}^\pi K_N=2\pi$, which
concentrates at $0$, absolute summability gives
\[
   \sigma^2_\infty=\lim_{N\to\infty}\sum_{|m|<N}\bigl(1-\tfrac{|m|}N\bigr)\rho_A(m)
   =\lim_{N\to\infty}\int_{-\pi}^\pi K_N\,\tilde f\ \ge\ \liminf_{N\to\infty}\int_{|\theta|<\eta}K_N\,q
   \ \ge\ 2\pi\inf_{|\theta|<\eta}q(\theta)
\]
for every $\eta>0$, and letting $\eta\downarrow0$, $\sigma^2_\infty\ge2\pi q(0)$, because
$\inf_{|\theta|<\eta}q(\theta)\to\min\bigl(q(0),\liminf_{\theta\to0}q(\theta)\bigr)=q(0)$.

\emph{Step 4 (evaluation of $q(0)$).} With $u=t-s$, whose density on $\{s<t\}$ is $1-u$,
\[
   \Upsilon(\lambda,-\lambda)=\tfrac12\int_0^1(1-u)\bigl(e^{-\mathrm iu\lambda}-e^{\mathrm iu\lambda}\bigr)\dd u
   =-\mathrm i\int_0^1(1-u)\sin(u\lambda)\,\dd u=\mathrm i\,\frac{\sin\lambda-\lambda}{\lambda^2},
\]
hence
$|\Upsilon(\lambda,-\lambda)|^2=(\lambda-\sin\lambda)^2/\lambda^4>0$ for $\lambda\ne0$ and
\[
   q(0)=\int_\RR|\Upsilon(\lambda,-\lambda)|^2\hat\varrho(\lambda)^2\,\dd\lambda
   =C_H^2\!\int_\RR\frac{(\lambda-\sin\lambda)^2}{\lambda^4}\,|\lambda|^{2-4H}\,\dd\lambda .
\]
The integrand is nonnegative and integrable, $\asymp|\lambda|^{4-4H}$ near $0$ ($4-4H>0$) and
$\asymp|\lambda|^{-4H}$ near $\infty$ ($4H>1$), so $q(0)\in(0,\infty)$ and
$\sigma^2_\infty\ge2\pi q(0)>0$.
\end{proof}

\begin{remark}[the size of $\sigma^2_\infty$]\label{rem:lrvsize}
The correlations $\rho_A(m)$, $m\ge1$, are negative for large $m$ (Lemma~\ref{lem:rhoAasym}) and,
numerically, for $1\le m\le60$ at $H=0.51,\dots,0.99$ (by direct evaluation of the four-fold integral;
the code is in the companion repository). So the cross-cell terms reduce
$\sigma^2_\infty$ below $\rho_A(0)=\Var(\tilde A_0)$, and the question is how large the reduction is.
Summing the exactly evaluated $\rho_A(m)$ gives
\[
\begin{array}{c|ccc}
H & 0.6 & 0.7 & 0.8\\\hline
\rho_A(0)=2\bigl(E^\infty_{\rm adj}-E^\infty_{\rm nest}\bigr) & 0.1678893 & 0.1087955 & 0.0635642\\
2\sum_{m\ge1}\rho_A(m) & -0.0013655 & -0.0042364 & -0.0067924\\
\sigma^2_\infty & \mathbf{0.166524} & \mathbf{0.104559} & \mathbf{0.056772}
\end{array}
\]
The single-cell term is in closed form: expanding
$\Var(\tilde A_0)=\tfrac14\E[(\langle12,\Sblk\rangle-\langle21,\Sblk\rangle)^2]$ by the shuffles
$(12)\shuffle(12)=2\cdot1212+4\cdot1122$ and $(12)\shuffle(21)=2\cdot1221+1212+2121+2\cdot2112$ gives
$\rho_A(0)=2(E^\infty_{\rm adj}-E^\infty_{\rm nest})$ exactly, with the two values from
Proposition~\ref{prop:altcontinuum}. The cross-cell sum is small and dominated by $m=1$, so the
cancellation is mild and the leading constant of $\Var(\sum_aA_a)$  is of the order of the
single-cell area variance. The explicit lower bound of Lemma~\ref{lem:lrv} is $0.61$, $0.56$
and $0.49$ times $\sigma^2_\infty$ at $H=0.6,0.7,0.8$.
\end{remark}

The next proposition gives the exact order of the strong error for fBm drivers on the Young
range.

\begin{proposition}[the strong exponent]\label{prop:strong}
Let $H\in(\tfrac12,1)$, $M\ge2$ and $d\ge2$. Let $X$ be the independent-component fBm$(H)$ driver and $X^n$ its piecewise-linear interpolant on the
uniform mesh $\delta=1/n$. The strong error satisfies
\[
  \eta_\delta=\bigl(\E\norm{\Sig^{\le M}(X^n)-\Sig^{\le M}(X)}_\hs^2\bigr)^{1/2}
  \ \asymp\ \delta^{\,2H-\frac{1}{2}}\qquad(\delta\to0),
 \]
with constants depending only on $H$, $M$ and $d$.
\end{proposition}

For $M=1$ or
$d=1$ the defect vanishes identically, and at $H=\tfrac12$ the rate is the classical $\delta^{1/2}$.

\smallskip\noindent\textbf{What the exponent yields.} For $M\ge4$ the weak error satisfies
$\norm{\cW_{\le M}}\asymp\delta^{2H}$ (Theorem~\ref{thm:2H} and Proposition~\ref{prop:alllevels}), so
\[
  \frac{\norm{\cW_{\le M}}_\hs^2}{\eta_\delta^2}\ \asymp\ \frac{\delta^{4H}}{\delta^{4H-1}}\ =\ \delta
  \qquad\text{at every }H\in(\tfrac12,1):
\]
the squared weak error lies one full mesh power below the mean-square strong error. The exponent is stated on $(\tfrac12,1)$ only, since
Lemma~\ref{lem:lrv}, on which the lower bound rests, uses $\hat\varrho(\lambda)=C_H|\lambda|^{1-2H}$
with $2-2H\in(0,1)$. Below the crossover see \S\ref{subsec:strongbelow}.
\begin{proof}[Proof of Proposition~\ref{prop:strong}]
We prove the lower bound. The upper bound,
\begin{equation}\label{eq:strongupper}
  \E\norm{\Sig^{\le M}(X)-\Sig^{\le M}(X^n)}_\hs^2\ \le\ C(H,M,d)\,\delta^{4H-1},
\end{equation}
is \citet[Thm.~1.1(i)]{HuberStrong}. The coordinate $12$ of the defect at level~$2$ is the sum of the $n$
cell areas. Its variance is a double sum of the covariances $\delta^{4H}\rho_A(|a-b|)$, which we
evaluate with Lemma~\ref{lem:rhoAasym} and bound below with Lemma~\ref{lem:lrv}.

  \emph{Reduction to the level-$2$ area.} By the shuffle relation
$\langle ij,\Sblk\rangle+\langle ji,\Sblk\rangle=X^iX^j$ the symmetric part of level $2$ is
endpoint-determined, hence pathwise identical for $X$ and $X^n$ (equal nodal endpoints). Level $1$ is
exact at the nodes. As the Hilbert--Schmidt norm is the $\ell^2$ sum over levels,
\[
   \eta_\delta^2\ \ge\ \Bigl\lVert\sum_{a=0}^{n-1} D_a\Bigr\rVert_{L^2}^2=\Var\Bigl(\sum_{a=0}^{n-1} D_a\Bigr),
\]
where, by \eqref{eq:areaChen} in \S\ref{subsec:strongbelow} with $a=0$, $b=n$,
$\sum_a D_a$ is the total level-$2$ area defect in the coordinate $12$ and $D_a:=A_a(X)-A_a(X^n)=A_a(X)$ the within-cell
L\'evy area ($A_a(X^n)=0$ on straight-line cells). We have $\E D_a=0$, since the law of $X$ is
invariant under $X^1\mapsto-X^1$ (the components are centred and independent) and this map changes the
sign of $D_a$, and $\Cov(D_a,D_b)=\delta^{4H}\rho_A(|a-b|)$ by the scaling relation stated before
Lemma~\ref{lem:rhoAasym}.

\emph{Summation.} By Lemma~\ref{lem:rhoAasym} both $\sum_m|\rho_A(m)|$ and $\sum_m|m|\,|\rho_A(m)|$
converge. Exactly $n-|m|$ ordered pairs $(a,b)\in\{0,\dots,n-1\}^2$ have $b-a=m$, so
\[
   \Var\Bigl(\sum_{a=0}^{n-1} D_a\Bigr)=\delta^{4H}\sum_{a,b=0}^{n-1}\rho_A(|a-b|)
   =\delta^{4H}\!\!\sum_{|m|<n}\!(n-|m|)\rho_A(m).
\]
With $\sigma^2_\infty=\sum_{m\in\ZZ}\rho_A(m)$,
\[
   \sum_{|m|<n}(n-|m|)\rho_A(m)\ =\ n\,\sigma^2_\infty\ -\ n\!\sum_{|m|\ge n}\!\rho_A(m)\ -\ \sum_{|m|<n}|m|\,\rho_A(m)
   \ =\ n\,\sigma^2_\infty+O(1):
\]
the tail term is at most $n\sum_{|m|\ge n}|\rho_A(m)|\lesssim n\cdot n^{4H-5}\to0$, and the last term is
at most $\sum_{m}|m|\,|\rho_A(m)|<\infty$. With $n=1/\delta$,
$\Var(\sum_{a<n}D_a)=\delta^{4H}(n\sigma^2_\infty+O(1))=\sigma^2_\infty\,\delta^{4H-1}(1+o(1))$.
As $\sigma^2_\infty>0$ by Lemma~\ref{lem:lrv}, this is two-sided, with no threshold at
$H=\tfrac34$ and no logarithmic factor. Hence
$\eta_\delta\gtrsim\delta^{2H-1/2}$ at every $H\in(\tfrac12,1)$.

The $o(1)$ is not uniform in $H$: it is controlled by $\sum_{|m|\ge n}|\rho_A(m)|\lesssim
n^{4H-5}$ and $\sum_{|m|<n}|m||\rho_A(m)|$, both of which degrade as $H\uparrow1$. So the two-sided
statement holds at fixed $H\in(\tfrac12,1)$, while the exponent $2H-\tfrac12$ is uniform there.
\end{proof}

Level~$2$ supplies the lower bound, but it need not carry the leading term of $\eta_\delta^2$. The
level-$3$ defect contains $\sum_a X_{0,t_a}\otimes D_a$, which is built from the same cell areas. So
$\sigma^2_\infty$ is the constant of the lower bound. We do not claim that it is the constant of
$\eta_\delta^2$.

\begin{remark}[the upper bound: proof strategy]\label{rem:strongupper}
The natural route to $\eta_\delta=O(\delta^{2H-1/2})$ is the Chen telescoping of
Proposition~\ref{prop:rough}, $\Sig(X)-\Sig(X^n)=\sum_aG^{<a}\otimes E_a\otimes\exp\Delta^{>a}$,
whose level-$2$ excess is the mean-zero area $A_a$. Its macroscopic factors $G^{<a}$,
$\exp\Delta^{>a}$ are random and correlated with the areas, and Cauchy--Schwarz pair by pair gives
only the mean-square bound $\delta^{4H-2}$. The Young continuity estimate of the signature map also
gives at best the exponent $2H-1$, i.e.\ the mean-square bound $\delta^{4H-2}$: the extra
$\delta^{1/2}$ comes from a central-limit cancellation in $\sum_aA_a$, which continuity does not see.
\\citet{HuberStrong} instead writes the defect exactly as one multiple integral of $X$
against the deterministic kernel $\one_{\Delta_M}-P_n^{\otimes M}\one_{\Delta_M}$, $P_n$ being cell
averaging, so the macroscopic factors become bounded deterministic cofactors, and it sums over cells
before bounding. The fully factorised Wick pattern is controlled by $\rho_A$ through
Lemma~\ref{lem:rhoAasym}, and the other patterns are of smaller order.
\end{remark}

\smallskip\noindent\textbf{Comparison with the bound of \citet{Schenck2026}.} \cite[Proposition~2.7]{Schenck2026}
charges the discretization to the strong $L^2$ error, bounded by $\delta^{2H}$ in mean square on
$(\tfrac12,1)$ for his independent-coordinate fractional Ornstein--Uhlenbeck drivers. For fBm drivers
this is loose in two ways. (a) The mean-square strong error is of exact order $\delta^{4H-1}$
(Proposition~\ref{prop:strong}), and $4H-1>2H$ there. For his fOU class we do not prove the matching
rate. (b) The weak error is a mean, and its squared norm is of order $\delta^{4H}$, because for
centred laws with independent components the expected L\'evy area vanishes. Point (b) alone gives
an improvement by the factor $\delta^{2H}$, without going through the strong error. On
$(\tfrac14,\tfrac12)$ the exponent $4H-1$ is, up to an arbitrary $\varepsilon>0$, his mean-square
upper bound for fOU drivers (his Proposition~2.7(ii)), and at level~$2$ for fBm the upper bound of
\citet[Thm.~2]{NeuenkirchTindelUnterberger2010}.

\subsection{Below the crossover}\label{subsec:strongbelow}

For $H\le\tfrac12$ we do not determine the strong error exactly, but we bound it uniformly over the
class, up to an arbitrarily small loss in the exponent. The bound is a direct consequence of the
convergence-rate theorem of \citet{FrizRiedel2014}. We include the short derivation because it is
uniform over the class and handles the drift. Lemma~\ref{lem:stronglower} below shows that its
exponent cannot be improved.

\begin{proposition}[a class-uniform strong rate]\label{prop:strongunif}
Let $H\in(\tfrac14,\tfrac12]$, $\rho=1/(2H)$, $M\ge1$ and $0<\vartheta<2H-\tfrac12$. There is $C=C(H,M,L,d,\vartheta)$
such that $\eta_\delta(P)\le C\delta^{\vartheta}$ for every $P\in\cP(H,L)$ and every mesh $\delta$.
\end{proposition}
\begin{proof}
  \emph{Idea.} The theorem of \citet{FrizRiedel2014} bounds the rough-path distance between the lifts
of $G^n$ and $G$ by a power of $\sup_s\norm{G_s-G^n_s}_{L^2}$, and this supremum is of order $\delta^H$.
The drift is added by a translation, and Lemma~\ref{lem:liprough} passes to the signature. Throughout, $S_N(\mathbf y)$ is the Lyons lift of a rough path
$\mathbf y$ to the step-$N$ group, and $T_h\mathbf y$ is the translate of $\mathbf y$ by a
bounded-variation path $h$ \citep[\S9.4.6]{FrizVictoir2010}.

\emph{Step 1 (the $L^2$ distance of the interpolant).} By Lemma~\ref{lem:subint} (a single letter of $G$),
$\norm{G_{r,r'}}_{L^2}\le C|r'-r|^H$ for $0\le r<r'\le1$, uniformly over the class. For $s$ in a cell $[t_a,t_{a+1}]$,
$G_s-G^n_s=(1-\lambda)G_{t_a,s}-\lambda G_{s,t_{a+1}}$ with $\lambda=(s-t_a)/\delta\in[0,1]$, so
\begin{equation}\label{eq:supmod}
   \sup_{0\le s\le1}\norm{G_s-G^n_s}_{L^2}\ \le\ C\delta^H .
\end{equation}

\emph{Step 2 (the rough-path distance).}
 The process $G$ is centred, continuous and Gaussian with independent components, and $V_\rho(R_G)\le L$ by (C5),
with $\rho\in[1,2)$; $G^n$ is its piecewise-linear interpolation on the uniform mesh. Choose $\nu>\rho$ with
$\tfrac1\nu+\tfrac1\rho>1$ and $H(1-\rho/\nu)\ge\vartheta$. This is possible because $H(1-\rho/\nu)$ increases to
$H(2-\rho)=2H-\tfrac12$ as $\nu$ increases to $\rho/(\rho-1)$ (to $\tfrac12$ as $\nu\to\infty$ if $\rho=1$), and
$\rho/(\rho-1)>\rho$ because $\rho<2$. Fix $q>2\nu$ and
put $N:=\lfloor q\rfloor$. The hypothesis of \citet[Thm.~1]{FrizRiedel2014} is a bound
$K\ge V_\rho(R_G;[0,1]^2)$ on the grid $\rho$-variation of their Definition~1, which (C5) gives with
$K=L$. That theorem yields, for every $r\ge1$,
\[
   \bignormL{\rho_{q\text{-var}}\bigl(S_N(G^n),S_N(\mathbf G)\bigr)}{L^r}\ \le\ C\,r^{N/2}
   \sup_{0\le s\le1}\norm{G^n_s-G_s}_{L^2}^{1-\rho/\nu}\ \le\ C_r\,\delta^{\vartheta},
\]
(their $\gamma$ is our $\nu$) with $C=C(q,\rho,\nu,L,N,d)$ and $\mathbf G$ the natural lift, the limit of the lifts of the piecewise-linear
approximations \citep[Thm.~15.33(iv)]{FrizVictoir2010}. The last step is \eqref{eq:supmod}, with
$\delta^{H(1-\rho/\nu)}\le\delta^\vartheta$ as $\delta\le1$.

\emph{Step 3 (the drift).} Write $\mathbf x^1:=S_N(G^n)$ and $\mathbf x^2:=S_N(\mathbf G)$, geometric $q$-rough paths.
The lifts $S_N(X^n)$ and $S_N(\mathbf X)$ are their translates by the same path $h:=\theta te_1$
\citep[Rem.~9.34, Cor.~9.35]{FrizVictoir2010}. Let $\varepsilon:=\rho_{q\text{-var}}(\mathbf x^1,\mathbf x^2)$ and take the control
$\omega:=\sum_{j=1,2}\hnorm{\mathbf x^j}^q_{q\text{-var};[\cdot,\cdot]}+|h|_{\onevar;[\cdot,\cdot]}+\sum_{i\le N}
(\rho^{(i)}_{q\text{-var};[\cdot,\cdot]}(\mathbf x^1,\mathbf x^2)/\varepsilon)^{q/i}$, where
$\rho^{(i)}_{q\text{-var};[s,t]}(\mathbf x^1,\mathbf x^2):=\norm{\pi_i(\mathbf x^1)-\pi_i(\mathbf x^2)}_{(q/i)\text{-var};[s,t]}$ is
the level-$i$ term of \eqref{eq:rhopvar}. Each summand is superadditive, so $\omega$ is a control. This is a variant (normalised by $\varepsilon$) of the control in the proof of \citet[Prop.~8.7]{FrizVictoir2010}.
With the distance controlled by $\omega$,
$\rho_{q,\omega}(\mathbf x,\mathbf y):=\max_{i\le N}\sup_{s<t}|\pi_i(\mathbf x_{s,t})-\pi_i(\mathbf y_{s,t})|/\omega(s,t)^{i/q}$,
we have $\rho_{q,\omega}(\mathbf x^1,\mathbf x^2)\le\varepsilon$ (on $\{\varepsilon=0\}$ there is nothing to prove), so
\citet[Thm.~9.33(ii)]{FrizVictoir2010} gives $\rho_{q,\omega}(T_h\mathbf x^1,T_h\mathbf x^2)\le C\varepsilon$, and by superadditivity
\begin{equation*}
   \rho_{q\text{-var}}(T_h\mathbf x^1,T_h\mathbf x^2)\ \le\ C\varepsilon\,(1+\omega(0,1))^{N/q},\qquad
   \omega(0,1)\ \le\ \hnorm{\mathbf x^1}^q_{q\text{-var}}+\hnorm{\mathbf x^2}^q_{q\text{-var}}+L+N ,
\end{equation*}
using $|h|_{\onevar}=|\theta|\le L$ and $\rho^{(i)}_{q\text{-var};[0,1]}(\mathbf x^1,\mathbf x^2)\le\varepsilon$.

\emph{Step 4 (the signature).} Lemma~\ref{lem:liprough}, applied
with $q$ in place of $p$ and $\bar N$ in place of $N$, gives
\[
   \norm{\hSblk_\delta-\Sblk}_\hs\ \le\ C\,\bar N^{M-1}\rho_{q\text{-var}}(T_h\mathbf x^1,T_h\mathbf x^2)
   \ \le\ C\,\bar N^{M-1}(1+\omega(0,1))^{N/q}\,\varepsilon,
\]
$\bar N$ being $1$ plus the larger
homogeneous norm of the two translates.

\emph{Step 5 (moments and conclusion).} These norms have moments of every order, bounded over the class and
uniformly in the mesh. Indeed, fix $p\in(2\rho,4)$ with $p\le q$. The homogeneous $q$-variation of $S_N(\mathbf y)$ is at
most its homogeneous $p$-variation, which is at most $C\hnorm{\mathbf y}_{\pvar}$ \citep[Thm.~9.5]{FrizVictoir2010}.
Lemma~\ref{lem:gaussmom}, with $K=L$ by (C5), bounds the moments of $\hnorm{S_3(G^n)}_{\pvar}$ and of
$\hnorm{\mathbf G}_{\pvar}$ uniformly in the mesh and over the class, and the translation adds at most
$C(1+L)$ \citep[Thm.~9.33(i)]{FrizVictoir2010}. So $Z:=\bar N^{M-1}(1+\omega(0,1))^{N/q}$ lies in
$L^4$, with a norm bounded over the class and uniformly in the mesh, and H\"older's inequality gives
\[
   \eta_\delta(P)^2\ \le\ C\,\E\bigl[Z^2\varepsilon^2\bigr]\ \le\ C\,\norm{Z}_{L^4}^2\,\norm{\varepsilon}_{L^4}^2,
   \qquad\text{so}\qquad
   \eta_\delta(P)\ \le\ C\norm{\varepsilon}_{L^4}\ \le\ C\delta^{\vartheta}
\]
by Step~2 with $r=4$. The constants depend only on $H$, $M$, $L$, $d$ and $\vartheta$, through the
choices of $\nu$, $q$, $N$ and $p$.
\end{proof}

The lower bound of Lemma~\ref{lem:stronglower} below rests on an exact identity for
the level-$2$ interpolation error.

\paragraph{The level-$2$ interpolation error.} Let $\mathbf x$ be a geometric rough path with level-$1$
increments $x_{s,t}$ and L\'evy area $A_{s,t}$ (the antisymmetric part of its level-$2$ component), write
$A_c:=A_{t_c,t_{c+1}}$ for its cell areas, and $\mathbf x^n$ for the canonical lift of the interpolant.
Level-$2$ Chen is algebraic, so no regularity enters. For $s<u<t$ the antisymmetric part of
Chen's relation at level~$2$ reads $A_{s,t}=A_{s,u}+A_{u,t}+\tfrac12\bigl(x_{s,u}\otimes x_{u,t}-x_{u,t}\otimes x_{s,u}\bigr)$.
Iterated over the cells of $[t_a,t_b]$, with $x_c:=x_{t_c,t_{c+1}}$, it gives
\[
   A_{t_a,t_b}\ =\ \sum_{c=a}^{b-1}A_c\ +\ \tfrac12\!\sum_{a\le c<c'<b}\bigl(x_c\otimes x_{c'}-x_{c'}\otimes x_c\bigr).
\]
The cross terms depend only on the cell increments, which $\mathbf x$ and $\mathbf x^n$ share.
  The interpolant is straight inside each cell, so its area there
vanishes and subtracting the two instances of this identity gives, for node-aligned $t_a<t_b$,
\begin{equation}\label{eq:areaChen}
   A^n_{t_a,t_b}\ =\ A_{t_a,t_b}\ -\ \sum_{c=a}^{b-1}A_c .
\end{equation}
So the level-$2$ interpolation error over $[0,1]$ is exactly $\sum_cA_c$: every mesh effect at
level~$2$ passes through this sum.
The rate $2H-\tfrac12$ is sharp, as \citet{FrizRiedel2014} note through the L\'evy area of fBm. The
following lemma is a short self-contained proof of that sharpness for the interpolant.

Both the next lemma and Lemma~\ref{lem:opsigns} use the spectral form of the fBm covariance: with
$a_H>0$ such that $|x|^{2H}=a_H\int_\RR(1-\cos x\xi)|\xi|^{-1-2H}\dd\xi$, the covariance
$C(s,t):=\E[B_sB_t]=\tfrac12(|s|^{2H}+|t|^{2H}-|t-s|^{2H})$ of an fBm$(H)$ $B$ is
\begin{equation}\label{eq:fbmspectral}
   C(s,t)=\tfrac{a_H}2\int_\RR\Re\bigl[(e^{\mathrm is\xi}-1)(e^{-\mathrm it\xi}-1)\bigr]|\xi|^{-1-2H}\dd\xi ,
\end{equation}
because $\Re\bigl[(e^{\mathrm is\xi}-1)(e^{-\mathrm it\xi}-1)\bigr]=(1-\cos s\xi)+(1-\cos t\xi)-(1-\cos(t-s)\xi)$.

  \begin{lemma}[the strong error of fBm below the crossover]\label{lem:stronglower}
Let $H\in(\tfrac14,\tfrac12]$, $M,d\ge2$, and let $X$ be fBm$(H)$ with independent components. The
L\'evy areas $A_a$ of $(X^1,X^2)$ over the cells $I_a$ satisfy, for every mesh $\delta=1/n$,
\[
   \Bigl\lVert\sum_{a<n}A_a\Bigr\rVert_{L^2}\ \ge\ c\,\delta^{2H-1/2},\qquad c=c(H)>0 ,
\]
and $\sum_{a<n}A_a$ is the coordinate $12$ of $\Sig(X)-\Sig(X^n)$; hence $\eta_\delta\ge c\,\delta^{2H-1/2}$.
\end{lemma}

\begin{proof}
\emph{Idea.} After rescaling, it suffices to show that the sum $Y$ of the $n$ unit-cell areas has
variance at least of order $n$. Conditionally on $B^2$, the (approximated) sum is a linear functional
of $B^1$, whose variance has the spectral form \eqref{eq:specstep} with weight $|\xi|^{1-2H}$.
The total spectral mass is of order $n$, and the spectral density is at most $C_1n$, so most of
the mass sits where the weight is bounded below.

\emph{Step 1 (reduction to unit cells).} By \eqref{eq:areaChen}, which uses only Chen's relation at level~$2$ and so holds for every $H>\tfrac14$, the coordinate $12$ of $\Sig(X)-\Sig(X^n)$ is $\sum_{a<n}A_a$, with $A_a$ the L\'evy area of
$(X^1,X^2)$ over the $a$-th cell; so $\eta_\delta\ge\norm{\sum_aA_a}_{L^2}$. By self-similarity, $(A_a)_{a<n}$ has the
law of $(\delta^{2H}\tilde A_a)_{a<n}$, where $\tilde A_a$ is the L\'evy area over $[a,a+1]$ of a pair $(B^1,B^2)$ of
independent fBm$(H)$. The $\tilde A_a$ are centred (the law is invariant under $B^1\mapsto-B^1$, which changes their sign), so $\norm{\sum_aA_a}_{L^2}^2=\delta^{4H}\Var(Y)$ for
$Y:=\sum_{a<n}\tilde A_a$, and it suffices to show $\Var(Y)\ge c_0\,n$ with $c_0=c_0(H)>0$. Then
$c=\sqrt{c_0}$, since $\delta^{4H}n=\delta^{4H-1}$. At $H=\tfrac12$ the $\tilde A_a$
are i.i.d.\ with variance $\tfrac14$, so let $H<\tfrac12$.

\emph{Step 2 (a spectral formula) .} By \eqref{eq:fbmspectral}, applied to $B=B^1$, for a step function $f=\sum_jf_j\one_{[u_j,u_{j+1})}$ and $\hat f(\xi):=\int f(t)e^{\mathrm it\xi}\dd t$,
using $\sum_jf_j\bigl(e^{\mathrm iu_{j+1}\xi}-e^{\mathrm iu_j\xi}\bigr)=\mathrm i\xi\hat f(\xi)$ and bilinearity,
\begin{equation}\label{eq:specstep}
   \E\Bigl(\sum_jf_jB_{u_j,u_{j+1}}\Bigr)^2=\frac{a_H}{2}\int_\RR|\hat f(\xi)|^2|\xi|^{1-2H}\dd\xi .
\end{equation}

\emph{Step 3 (approximation and conditioning).} Let $Y_k$ be the same sum of areas for the piecewise-linear interpolants of
$(B^1,B^2)$ on the mesh $1/k$, $k\in\NN$, which refines the integer mesh. The interpolant on the
integer mesh of this fine interpolant is that of $(B^1,B^2)$, so by \eqref{eq:areaChen}, applied to both
paths on $[0,n]$, $Y-Y_k$ is the coordinate $12$ of the difference of the signatures over $[0,n]$ of
$(B^1,B^2)$ and of its fine interpolant. This tends to $0$ in $L^2$ as $k\to\infty$ (Lemma~\ref{lem:L1conv}, after rescaling $[0,n]$ to $[0,1]$), so $\Var(Y)=\lim_k\Var(Y_k)$. For a smooth path on a cell
$I_a=[a,a+1]$, integration by parts gives $A=\tfrac12x^1_{I_a}x^2_{I_a}-\int_{I_a}x^2_{a,s}\,\dd x^1_s$. For the
interpolants the integral equals $\sum_jB^1_{J_j}\tfrac12(B^2_{a,u_j}+B^2_{a,u_{j+1}})$, the sum over the fine intervals
$J_j=[u_j,u_{j+1}]\subset I_a$ (on $J_j$ both components are linear, so the integrand is the average
of its endpoint values). Hence $Y_k=\sum_jB^1_{J_j}\,\psi_k|_{J_j}$ for the step function
\[
   \psi_k:=\sum_{a<n}\ \sum_{J_j\subset I_a}\one_{J_j}\Bigl(\tfrac12B^2_{I_a}-\tfrac12\bigl(B^2_{a,u_j}+B^2_{a,u_{j+1}}\bigr)\Bigr),
\]
which is a function of $B^2$. As $B^1$ is independent of $B^2$, the law of total variance and
\eqref{eq:specstep} applied conditionally on $B^2$ give
\[
   \Var(Y_k)\ \ge\ \E[\Var(Y_k\mid B^2)]\ =\ \frac{a_H}{2}\int_\RR\E|\hat\psi_k(\xi)|^2|\xi|^{1-2H}\dd\xi .
\]

\emph{Step 4 (two bounds on $\E|\hat\psi_k|^2$).} (i) \emph{Total mass.} $\E\norm{\psi_k}^2_{L^2}=n\,v_k$ with $v_k:=\E\int_0^1(\tfrac12B^2_{0,1}-m_k(s))^2\dd s$, where $m_k(s)$
is the average of $B^2_{0,\cdot}$ at the two fine nodes around $s$ (the $n$ cells contribute equally, by
stationarity of the increments). By continuity of $B^2$, $v_k\to
v:=\E\int_0^1(\tfrac12B^2_{0,1}-B^2_{0,s})^2\dd s=\tfrac1{2H+1}-\tfrac14>0$. The integrand is
$-\tfrac14+\tfrac12s^{2H}+\tfrac12(1-s)^{2H}$. So $v_k\ge v/2$ for $k\ge k_0$, and by Plancherel
 \[
   \int_\RR\E|\hat\psi_k|^2\,\dd\xi=2\pi\,\E\norm{\psi_k}^2_{L^2}\ \ge\ \pi v\,n\qquad(k\ge k_0).
\]
(ii) \emph{A uniform bound.} $\hat\psi_k(\xi)=\sum_ae^{\mathrm ia\xi}\hat g_a(\xi)$ with $g_a:=\psi_k(a+\cdot)\one_{[0,1)}$. Each value $g_a(t)$ is a
sum of three increments of $B^2$ over subintervals of $I_a$ with coefficients $\pm\tfrac12$, so $\norm{g_a(t)}_{L^2}\le\tfrac32$,
and $g_a(t)=\int\varphi_{a,t}\,\dd B^2$ for a step function $\varphi_{a,t}$ supported in $I_a$ with $|\varphi_{a,t}|\le\tfrac12$. For $|a-b|=m\ge2$
the supports are at distance $\ge m-1\ge1$, where the increment covariance of $B^2$ has the density
$c_H|x-y|^{2H-2}$, $|c_H|=H(1-2H)\le1$, so
\[
   |\E[g_a(t)g_b(s)]|=H(1-2H)\Bigl|\iint\varphi_{a,t}(x)\varphi_{b,s}(y)|x-y|^{2H-2}\dd x\dd y\Bigr|\le(m-1)^{2H-2}.
\]
For $|a-b|\le1$ we use $|\E[g_a(t)g_b(s)]|\le(\tfrac32)^2$. Fewer than $3n$ ordered pairs $(a,b)$ have
$|a-b|\le1$, and at most $2n$ have $|a-b|=m$ for each $m\ge2$. Hence, for every $\xi$,
\[
   \E|\hat\psi_k(\xi)|^2\ \le\ \sum_{a,b}\int_0^1\!\!\int_0^1|\E[g_a(t)g_b(s)]|\,\dd s\,\dd t
   \ \le\ \Bigl(3\cdot\tfrac94+2\sum_{m\ge2}(m-1)^{2H-2}\Bigr)n\ \le\ C_1n,
\]
with $C_1:=7+2\sum_{m\ge2}(m-1)^{2H-2}<\infty$, finite because $2-2H>1$.

\emph{Step 5 (cutting off low frequencies).}
For $\lambda>0$, since $1-2H>0$, the weight is at least $\lambda^{1-2H}$ on $\{|\xi|\ge\lambda\}$, and by
(ii) the mass of $\E|\hat\psi_k|^2$ on $\{|\xi|<\lambda\}$ is at most $2\lambda C_1n$. With (i),
\[
   \int_\RR\E|\hat\psi_k|^2|\xi|^{1-2H}\dd\xi\ \ge\ \lambda^{1-2H}\int_{|\xi|\ge\lambda}\E|\hat\psi_k|^2\dd\xi
   \ \ge\ \lambda^{1-2H}(\pi v-2\lambda C_1)\,n .
\]
Taking $\lambda=\pi v/(4C_1)$ and using Step~3 gives, for $k\ge k_0$,
\[
   \Var(Y_k)\ \ge\ c_0\,n,\qquad c_0:=\frac{a_H}{2}\cdot\frac{\pi v}{2}\Bigl(\frac{\pi v}{4C_1}\Bigr)^{1-2H}>0 ,
\]
with $c_0$ independent of $n$ and $k$, and letting $k\to\infty$ gives $\Var(Y)\ge c_0n$, which proves the claim.
  \end{proof}

So the exponent in Proposition~\ref{prop:strongunif} is optimal up to an arbitrarily
small loss: for fBm below the crossover, $c\,\delta^{4H-1}\le\eta_\delta^2\le C_\vartheta\delta^{2\vartheta}$
for every $\vartheta<2H-\tfrac12$. Up to this loss, the squared weak error, of order
$\delta^{2(4H-1)}$ (Theorem~\ref{thm:roughrate}), is the square of $\eta_\delta^2$.
\section{The weak error cannot be removed using the mesh data}\label{sec:mesh}

No procedure based on the mesh data can reduce the weak error below the order $\delta^{\gamma(H)}$
uniformly over the class. For each mesh we exhibit two laws in the class whose mesh data have the same
law but whose expected signatures differ by order $\delta^{\gamma(H)}$. The second law adds to
fBm a small Gaussian ``ripple'' that vanishes at every mesh node (Construction~\ref{con:bridge}).
Proposition~\ref{prop:adm} shows that both laws lie in the class and have the same mesh data, and
Propositions~\ref{prop:sep}, \ref{prop:sepbrownian} and~\ref{prop:seprough} bound their separation
on $(\tfrac12,1)$, at $H=\tfrac12$ and on $(\tfrac14,\tfrac12)$.

\begin{construction}[the tied-down ripple]\label{con:bridge}
Fix $H\in(\tfrac14,1)$, $\epsilon>0$ and a block width $\delta_0=1/n_0$. Let
$C(u,v)=\tfrac12(u^{2H}+v^{2H}-|u-v|^{2H})$ be the covariance of $\mathrm{fBm}(H)$, let $Z$ be an
$\mathrm{fBm}(H)$, and call $\varpi(u):=Z_u-uZ_1$, $u\in[0,1]$, the tied-down fBm. Its covariance
is
\[
   B_H(u,v)=C(u,v)-u\,C(v,1)-v\,C(u,1)+uv ,
\]
and $\varpi(0)=\varpi(1)=0$. (It is not the fBm bridge, i.e.\ $Z$ conditioned on $Z_1=0$, whose
covariance is $C(u,v)-C(u,1)C(v,1)$.) Let $G$ be an $\RR^d$-valued process with independent
$\mathrm{fBm}(H)$ components, and let $(\varpi^i_a)_{1\le i\le d,\ 0\le a<n_0}$ be i.i.d.\ copies of
$\varpi$, independent of $G$. With the blocks $J_a:=[a\delta_0,(a+1)\delta_0]$ define the ripple
$R=(R^1,\dots,R^d)$ by $R^i_t:=\delta_0^H\varpi^i_a((t-a\delta_0)/\delta_0)$ for $t\in J_a$. It is
continuous because $\varpi$ vanishes at both ends. Let $P_0$ be the law of $G$ and call the law
$P_1^{(\delta_0)}$ of $G+\epsilon R$ the tied-down perturbation of $P_0$. At mesh $\delta$ we use the pair $P_0$, $P_1:=P_1^{(\delta)}$, whose blocks
are the mesh cells.
\end{construction}

\begin{proposition}[admissible, with the same mesh data]\label{prop:adm}
There is $L_0=L_0(H,\epsilon)$ (depending also on the fixed $M$ and $d$) such that for every $L\ge L_0$:
\emph{(a)} $P_0\in\cP_0(H,L)$ and $P_1^{(\delta_0)}\in\cP_0(H,L)$ for every $\delta_0$.
\emph{(b)} At mesh $\delta$, the mesh data $\cD=(X_{t_j})_{j\le n}$ have the same law under $P_0$ and
under $P_1=P_1^{(\delta)}$.
\emph{(c)} $\mu(P_0)$ and $\mu(P_1)$ agree through level $3$.
\end{proposition}
\begin{proof}
\emph{(a)} \emph{Centred:} $P_1^{(\delta_0)}$ is the law of a sum of independent centred Gaussian
processes, so it is centred Gaussian, with $\theta=0$. \emph{Independent components:} each component
$i$ carries its own family $(\varpi^i_a)_a$, independent of the other components and of $G$, so the
components of $G+\epsilon R$ are independent, and they share one increment covariance because the
$\varpi^i_a$ are identically distributed. \emph{Normalised:} $t=1$ is a block node, so
$R_1=0$ and $\Var(G^i_1+\epsilon R^i_1)=\Var(G^i_1)=1$. Hypotheses (C1)--(C5) and the moment condition of Definition~\ref{def:class} are verified for $P_0$
and $P_1^{(\delta_0)}$, at every mesh and with one $L_0(H,\epsilon)$ for all $\delta_0$, in
Remark~\ref{rem:C1C2}.

\emph{(b)} Since $\varpi_a(0)=\varpi_a(1)=0$, the ripple vanishes at every mesh node. So the nodal
vectors coincide pathwise under the coupling, and the laws of the mesh data agree.

\emph{(c)} By (a) both laws lie in $\cP_0(H,L)$, so Proposition~\ref{prop:level4} gives
$\cW_m(P_j)=0$ for $m\le3$ and $j=0,1$, that is,
$\mu(P_j)_m=\E_{P_j}\bigl[\Sig^{\le M}(X^n)\bigr]_m$ at mesh $\delta$. The interpolant $X^n$ is a
function of the mesh data, whose law is the same under $P_0$ and $P_1$ by (b), so the two right-hand
sides agree. Hence $\mu(P_1)-\mu(P_0)$ vanishes through level $3$.
\end{proof}

On the Young range the separation is computed to leading order with the single-error asymptotic
of \S\ref{sec:bias}.

\begin{proposition}[separation: exact leading constant for $H>\tfrac12$]\label{prop:sep}
Fix $H\in(\tfrac12,1)$, $\epsilon>0$ and $d\ge2$, and let
\[
   \kappa_{\rm sep}(H)=\frac{1}{6}-\frac{2}{(2H+1)(2H+2)}\ <\ 0\qquad\text{for all }H\in(\tfrac12,1),
\]
which vanishes only at $H=1$. For every pair of letters $i\ne j$, the differences
$\Delta\mu_\bullet:=\mu(P_1)_\bullet-\mu(P_0)_\bullet$ at the adjacent, alternating and nested words
($iijj$, $ijij$, $ijji$, and likewise with $i,j$ exchanged) are
\begin{gather*}
   \Delta\mu_{\rm adj}=-\tfrac12\epsilon^2\kappa_{\rm sep}(H)\,\delta^{2H}\bigl(1+o(1)\bigr),\qquad
   \Delta\mu_{\rm alt}=\epsilon^2\kappa_{\rm sep}(H)\,\delta^{2H}\bigl(1+o(1)\bigr),\\
   \Delta\mu_{\rm nest}=-\tfrac12\epsilon^2\kappa_{\rm sep}(H)\,\delta^{2H}\bigl(1+o(1)\bigr).
\end{gather*}
Consequently, for $M\ge4$,
\[
   \norm{\mu(P_0)-\mu(P_1)}_\hs\ \ge\ \sqrt{\tbinom d2}\,\sqrt3\,|\kappa_{\rm sep}(H)|\,\epsilon^2\,\delta^{2H}
   +o(\delta^{2H}),
\]
the right-hand side being the Hilbert--Schmidt norm of the level-$4$ block of $\mu(P_1)-\mu(P_0)$,
with equality for $4\le M\le5$. In particular $\norm{\mu(P_0)-\mu(P_1)}_\hs\asymp\delta^{2H}$.
\end{proposition}
\begin{proof}
For the leading order we expand the level-$4$ coordinates in $\epsilon^2$. The terms linear in the ripple
density have the single-error form of Lemma~\ref{lem:singleerror}, with the ripple density in place
of the interpolation error, and the term quadratic in it is of lower order.

\emph{Step 1 (upper bound).} By Proposition~\ref{prop:adm}(b) the two laws have the same
mesh-data law, hence a common mean $\bar\mu:=\E_{P_0}[\Sig^{\le M}(X^n)]=\E_{P_1}[\Sig^{\le M}(X^n)]$
of the interpolant signature. Write $\mu_j:=\mu(P_j)$. Both laws lie in $\cP_0(H,L)$ at every mesh
(Proposition~\ref{prop:adm}(a)), and the bound of Proposition~\ref{prop:alllevels} is uniform over
$\cP_0(H,L)$, so $\norm{\mu_0-\bar\mu}=O(\delta^{2H})$ and $\norm{\mu_1-\bar\mu}=O(\delta^{2H})$.
Hence $\norm{\mu_1-\mu_0}=O(\delta^{2H})$.

\emph{Step 2 (expansion in $\epsilon^2$).} Since $R$ is independent of $G$, the
increment-covariance densities add, $\varrho_1=\varrho_0+\epsilon^2\varrho_R$, with
$\varrho_0=\varrho$ the fBm density and $\varrho_R$ the ripple density. After the Wick pairing,
each level-$4$ $(2,2)$-coordinate $E_\bullet$ is the integral over the ordered simplex
$\{t_1<t_2<t_3<t_4\}$ of a product $\varrho(t_A)\varrho(t_B)$, one covariance factor for each Wick
pair, $t_A$ and $t_B$ being the time pairs of the two Wick pairs (the representation of the
$(2,2)$-coordinates stated before \eqref{eq:coinc}). The expansion
\[
   \varrho_1(t_A)\varrho_1(t_B)-\varrho(t_A)\varrho(t_B)
   =\epsilon^2\bigl[\varrho_R(t_A)\varrho(t_B)+\varrho(t_A)\varrho_R(t_B)\bigr]
   +\epsilon^4\varrho_R(t_A)\varrho_R(t_B)
\]
gives
\[
   \Delta\mu_\bullet:=\mu(P_1)_\bullet-\mu(P_0)_\bullet
   =\epsilon^2\!\int_{\rm simplex}\!\bigl[\varrho_R(t_A)\varrho(t_B)+\varrho(t_A)\varrho_R(t_B)\bigr]
   +\epsilon^4\!\int_{\rm simplex}\!\varrho_R(t_A)\varrho_R(t_B) ,
\]
the single-error structure of Lemma~\ref{lem:singleerror} with the interpolation error $e$ replaced
by $\varrho_R$, plus a pure-ripple term.
\emph{Step 3 (the ripple density has zero cell mean).}  On cell $I_a$ the ripple density is
$\varrho_R(s,t)=\delta^{2H-2}b_H(\tfrac{s-t_a}\delta,\tfrac{t-t_a}\delta)$ with
\[
   b_H(u,v)=\partial_u\partial_v B_H(u,v)
   =c_H|u-v|^{2H-2}-f(u)-f(v)+1,\qquad f(u):=H\bigl(u^{2H-1}+(1-u)^{2H-1}\bigr),
\]
and $\varrho_R=0$ across cells (the ripples of different cells are independent).
$b_H\in L^1([0,1]^2)$ as $2H-2>-1$. The covariance $B_H$ vanishes on
$\partial[0,1]^2$ ($\varpi(0)=\varpi(1)=0$), so $\int_{[0,1]^2}b_H=B_H\big|_{\rm corners}=0$, hence
(the two halves of the square carry equal integrals, $b_H$ being symmetric)
\[
   \int_{0<u<v<1}b_H\ =\ \tfrac12\int_{[0,1]^2}b_H\ =\ 0 :
\]
the ripple has zero cell mean, which is (K1) of Lemma~\ref{lem:singleerror}.

\emph{Step 4 (the pure-ripple term is $o(\delta^{2H})$).} The zero cell mean also
controls the pure-ripple term. Each pair of the pairing carries a block-diagonal factor, so both of
its variables lie in one cell. For the alternating and nested pairings the variables of one pair
enclose a variable of the other pair, so, cells being intervals, the ordering $t_1<t_2<t_3<t_4$ puts
all four variables into one cell. For the adjacent pairing the two pairs may lie in different cells
$a<b$, and then the term factorises into
$\bigl(\int_{(I_a\times I_a)\cap\{s<t\}}\varrho_R\bigr)\bigl(\int_{(I_b\times I_b)\cap\{s<t\}}\varrho_R\bigr)=0$.
What remains is a sum of $n$ single-cell contributions. Rescaling the four variables of cell
$I_a$ to the unit cell, each is $\delta^4(\delta^{2H-2})^2=\delta^{4H}$ times a pairing integral of
$b_H\otimes b_H$ over the unit simplex, which is at most $\norm{b_H}_{L^1}^2$. With $n=\delta^{-1}$
cells,
\[
   \epsilon^4\Bigl|\int_{\rm simplex}\varrho_R(t_A)\varrho_R(t_B)\Bigr|
   \ \le\ \epsilon^4\norm{b_H}_{L^1}^2\,\delta^{4H-1}\ =\ o(\delta^{2H}),
\]
because $4H-1>2H$ for $H>\tfrac12$. This is the reason $H=\tfrac12$ is excluded here.

\emph{Step 5 (the single-ripple terms).} As $f(u)=f(1-u)$, $b_H$ is invariant under the cell
reflection $(u,v)\mapsto(1-v,1-u)$, so Lemma~\ref{lem:reflect} applies through its case~\emph{(ii)},
and (K2) holds as in Remark~\ref{rem:singleerrorscope}(ii): block-diagonality leaves only the $m=0$
cell, and the self-similarity of the ripple makes the rescaling exact. The weights $w$ of
Lemma~\ref{lem:weights} satisfy Hypothesis~(W), so Lemma~\ref{lem:singleerror} applies to each
single-ripple term, with remainder $O(\delta^{4H-1}+\delta^2\log\tfrac1\delta)=o(\delta^{2H})$. Its
constants are $\Phi^{\rm s}_R:=\phi^{\rm s}(0)$ and $\Phi^{\rm t}_R:=\phi^{\rm t}(0)$, and
Lemma~\ref{lem:reflect} gives $\Phi^{\rm s}_R=-\Phi^{\rm t}_R$. Together with
\[
   \Phi^{\rm t}_R-\Phi^{\rm s}_R=\int_{0<u<v<1}b_H(u,v)\,(v-u)\,\dd u\,\dd v=:\kappa_{\rm sep}(H)
\]
this gives $\Phi^{\rm s}_R=-\Phi^{\rm t}_R=-\tfrac12\kappa_{\rm sep}(H)$.

\emph{Step 6 (assembly over the words, and the norm).} By Lemma~\ref{lem:singleerror} with the
weights of Lemma~\ref{lem:weights}, whose diagonal integrals $\int_0^1\partial_\bullet w(x,x)\,\dd x$
are $0$ or $\pm\tfrac12$, each single-ripple term equals
$\delta^{2H}\bigl[\Phi^{\rm s}_R\int_0^1\partial_sw(x,x)\,\dd x+\Phi^{\rm t}_R\int_0^1\partial_tw(x,x)\,\dd x\bigr]$
to leading order. The values of these integrals for the six weights are: for the adjacent word,
$(0,-\tfrac12)$ (error on $\{1,2\}$) and $(\tfrac12,0)$ (error on $\{3,4\}$). For the alternating
word, $(-\tfrac12,\tfrac12)$ for both error pairs. For the nested word, $(0,0)$ (error on $\{1,4\}$)
and $(\tfrac12,-\tfrac12)$ (error on $\{2,3\}$). Summing over the two error pairs of each word, the
adjacent and nested words give $\tfrac12(\Phi^{\rm s}_R-\Phi^{\rm t}_R)=-\tfrac12\kappa_{\rm sep}$
and the alternating word gives $\Phi^{\rm t}_R-\Phi^{\rm s}_R=\kappa_{\rm sep}$, so
 $\Delta\mu=\epsilon^2 \kappa_{\rm sep}(H)\,(-\tfrac12,1,-\tfrac12)\,\delta^{2H}+o(\delta^{2H})$, the three
displayed coordinates, and, at $d=2$,
$\norm{\Delta\mu}=\epsilon^2|\kappa_{\rm sep}|\sqrt{2(\tfrac14+1+\tfrac14)}\,\delta^{2H}
=\sqrt3\,\epsilon^2|\kappa_{\rm sep}|\,\delta^{2H}$. For general $d$ each of the $\binom d2$ pairs of
letters contributes the same, as in the count \eqref{eq:W4norm}. The other level-$4$ words
agree: $\langle iiii,\Sig\rangle=(X^i_1)^4/24$ depends only on the endpoint, whose law is the same,
and a word in which some letter occurs an odd number of times has mean zero under both centred laws
with independent components. Levels $1$--$3$ agree by Proposition~\ref{prop:adm}(c), and level~$5$
vanishes for both centred laws, whence the equality for $M\le5$.

\emph{Step 7 (the constant).} With $c_H=H(2H-1)$ and the elementary beta integrals
$\int_{u<v}(v-u)^{2H-1}=\tfrac1{2H(2H+1)}$,
$\int_{u<v}(v-u)u^{2H-1}=\int_{u<v}(v-u)(1-v)^{2H-1}=\tfrac1{2H(2H+1)(2H+2)}$,
$\int_{u<v}(v-u)v^{2H-1}=\int_{u<v}(v-u)(1-u)^{2H-1}=\tfrac1{2(2H+2)}$ and $\int_{u<v}(v-u)=\tfrac16$,
\[
   \kappa_{\rm sep}(H)=\frac{2H-1}{2(2H+1)}-\frac{2H^2+H+1}{(2H+1)(2H+2)}+\frac{1}{6}=\frac{1}{6}-\frac{2}{(2H+1)(2H+2)} .
\]
Here $\kappa_{\rm sep}(H)=0\iff(2H+1)(2H+2)=12\iff(2H+5)(H-1)=0\iff H=1$. Since
$\kappa_{\rm sep}(\tfrac12)=-\tfrac16<0$ and $\kappa_{\rm sep}$ is continuous, $\kappa_{\rm sep}<0$
on $[\tfrac12,1)$.
\end{proof}

At $H=\tfrac12$ the cells decouple exactly, and the separation can be computed in closed form for
every mesh.

\begin{proposition}[the Brownian separation]\label{prop:sepbrownian}
At $H=\tfrac12$ the pair $(P_0,P_1)$ of Construction~\ref{con:bridge} satisfies, for every $n$ and
every $\epsilon>0$,
\[
   \mu(P_1)_\bullet-\mu(P_0)_\bullet\ =\ \delta\,\Dsep_\bullet(\epsilon),\qquad
   \Dsep(\epsilon)=\Bigl(\tfrac{x}{12},\ \tfrac{x^{2}}{24},\ -\tfrac{x}{12}-\tfrac{x^{2}}{24}\Bigr),
   \qquad x:=\epsilon^2 ,
\]
for $(\mathrm{adj},\mathrm{alt},\mathrm{nest})$ and every pair of letters; hence, at $d=2$ and for
$M\ge4$,
$\norm{\mu(P_0)-\mu(P_1)}_\hs\ge\delta\bigl(2\sum_\bullet \Dsep_\bullet^2\bigr)^{1/2}
=\tfrac{\epsilon^2}{6}\,\delta\,(1+O(\epsilon^2))\ne0$ for every $\epsilon>0$, with equality for
$4\le M\le5$. For general $d\ge2$ the right-hand side is multiplied by $\sqrt{\binom d2}$. In particular
$\norm{\mu(P_0)-\mu(P_1)}_\hs\asymp\delta$ at $H=\tfrac12$.
\end{proposition}
\begin{proof}
We reduce to a single cell with Chen's identity and the independence of the cells, and then
compute the unit-cell coefficients from the explicit covariance of the perturbed process.

\emph{Step 1 (exact cell decoupling).} At $H=\tfrac12$ the driver has independent increments and the ripple is
independent across cells, so the cell paths are independent under both $P_0$ and $P_1$, and Chen's
identity gives $\mu(P_j)=\bigotimes_a\E[G^{(j)}_a]$, where $G^{(j)}_a$ is the signature, truncated
at level $M$, of the path on the cell $I_a$ under $P_j$. The two laws have the same nodal laws
(Proposition~\ref{prop:adm}(b)) and, being centred Gaussian with independent components, the same
level-$1,2,3$ cell means. So every cross-cell product in the level-$4$ coefficient agrees, and the
level-$4$ difference is the sum of the $n$ single-cell differences. By Brownian scaling each is
$\delta^{2}$ times the unit-cell one, which gives exactly $n\delta^{2}\Dsep_\bullet=\delta \Dsep_\bullet$.

\emph{Step 2 (the unit cell).} On $[0,1]$ the perturbed process is $Y_t=\sqrt{1+x}\,\beta_t+t\,Z$ with $\beta$
a standard Brownian bridge and $Z\sim\cN(0,1)$ independent, i.e.\ covariance
$C(s,t)=(1+x)\min(s,t)-x\,st$ and increment-covariance measure
$\varrho_1=(1+x)\,\delta_{\{s=t\}}-x$. (The Brownian path of the cell is $\beta+tZ$, and the
ripple adds an independent bridge $\epsilon\beta'$; $\beta+\epsilon\beta'$ has the law of
$\sqrt{1+x}\,\beta$.) Let $E_\bullet(x)$ be the unit-cell level-$4$ coordinates of $Y$. After the
Wick pairing, $E_\bullet(x)$ is the integral over the ordered simplex of $\varrho_1(t_A)\varrho_1(t_B)$,
which we evaluate with the diagonal convention \eqref{eq:coinc}: each coincident adjacent pair
carries the factor $\tfrac12$, consistent with the Brownian value $E_{\rm adj}(0)=\tfrac18$. For the
adjacent word, with Wick pairs $\{1,2\}$ and $\{3,4\}$, the product splits into a term with two
diagonal factors, two terms with one diagonal and one constant factor, and a term with two constant
factors:
\[
   E_{\rm adj}(x)=(1+x)^2\cdot\tfrac12\cdot\tfrac12\cdot\tfrac12
   \ -\ 2\cdot(1+x)\tfrac12\cdot x\cdot\tfrac16\ +\ x^2\cdot\tfrac1{24},
\]
the factors $\tfrac12$ coming from the coincident pairs and from the area of $\{s<t\}$, $\tfrac16$
from the volume of $\{s<t_3<t_4\}$ (or $\{t_1<t_2<t\}$), and $\tfrac1{24}$ from the volume of the
ordered simplex. The alternating and nested words are expanded in the same way. Subtracting the
values at $x=0$ gives
  \[
   \Dsep_{\rm adj}=\tfrac{(1+x)^2-1}{8}-\tfrac{x(1+x)}{6}+\tfrac{x^2}{24}=\tfrac{x}{12},\qquad
   \Dsep_{\rm alt}=\tfrac{x^2}{24},\qquad
   \Dsep_{\rm nest}=-\tfrac{x(1+x)}{12}+\tfrac{x^2}{24}=-\tfrac{x}{12}-\tfrac{x^2}{24},
\]
using in $\Dsep_{\rm alt}$ and $\Dsep_{\rm nest}$ that the terms containing $\delta_{\{t_1=t_3\}}$ or
$\delta_{\{t_2=t_4\}}$ are degenerate on the ordered simplex. The three sum to zero
(Lemma~\ref{lem:sumzero}).
\end{proof}

\begin{remark}[the Brownian separation is not the limit from above]\label{rem:sepmerger}
Proposition~\ref{prop:sepbrownian} gives the $\epsilon^2$-coefficient
$(\tfrac1{12},0,-\tfrac1{12})$, whereas $\lim_{H\downarrow1/2}\kappa_{\rm sep}(H)(-\tfrac12,1,-\tfrac12)
=(\tfrac1{12},-\tfrac16,\tfrac1{12})$. The difference $(0,\tfrac16,-\tfrac16)$ is the phenomenon of
Remark~\ref{rem:merger} for the interpolation bias: at $H=\tfrac12$ the Young contribution of
Lemma~\ref{lem:singleerror} and the rough single-cell contribution have the same order $\delta$ and
add. So the Brownian separation constant is $\epsilon^2/6$ to leading order in $\epsilon$, not the
$\epsilon^2/(2\sqrt3)$ suggested by $\kappa_{\rm sep}(\tfrac12)=-\tfrac16$. The mechanism: for the alternating word the weight
$w=\tfrac12[(1-s)^{2H}-(1-t)^{2H}-(t-s)^{2H}]$ tends to $0$ pointwise as $H\downarrow\tfrac12$
while its diagonal derivative tends to $\pm\tfrac12$, a boundary layer that the first-order
expansion of Lemma~\ref{lem:singleerror} does not see uniformly.
\end{remark}

\subsection*{The level-4 coordinates from the covariance alone, and the rough-regime separation}

To extend the lower bound to the rough regime we use one observation: for a centred Gaussian driver
with i.i.d.\ components, the three level-$4$ $(2,2)$-coordinates are explicit covariance
functionals, given by absolutely convergent integrals for every $H\in(\tfrac14,1)$. No
renormalisation and no analytic continuation in $H$ is needed, because the shuffle identity
collapses the inner double integral before the singular kernel is paired with itself.
Lemma~\ref{lem:cov4} is used in Propositions~\ref{prop:seprough} and~\ref{prop:adjexact} and in
Remark~\ref{rem:noAC}.

\begin{lemma}[the level-4 coordinates as covariance functionals]\label{lem:cov4}
Let $X=(X^1,X^2)$ have i.i.d.\ centred Gaussian components with covariance $C(s,t)$, variance
$v(t)=C(t,t)$ and increment variance $D(s,t)=\E[(X_t-X_s)^2]=v(s)+v(t)-2C(s,t)$. Write $\partial_1D$ and
$\partial_2D$ for the partial derivatives of $D$ in its first and second argument. Assume
\begin{itemize}[leftmargin=2.1em,topsep=2pt,itemsep=1pt]
\item[\emph{(D1)}] $D$ is symmetric and continuous on $[0,1]^2$ with $D(x,x)=0$. $v$ and $C(1,\cdot)$
are absolutely continuous on $[0,1]$. Off the diagonal $\partial_1D$ and $\partial_2D$ exist and are
locally integrable, $\partial_1\partial_2D=\partial_2\partial_1D$ a.e.\ there, and for each fixed $t$
both $s\mapsto D(s,t)$ and $s\mapsto\partial_2D(s,t)$ are absolutely continuous on compact
subintervals of each component of $\{s\ne t\}$.
\item[\emph{(D2)}] $\partial_1D\,\partial_2D\in L^1([0,1]^2)$ and $D\,\partial_1\partial_2D\in L^1([0,1]^2)$.
\item[\emph{(D3)}] \emph{(diagonal vanishing)} $\Xi(\tau)\to0$ as $\tau\downarrow0$, where
\begin{equation}\label{eq:diagbdry}
   \Xi(\tau)\ :=\ \int_0^{1-\tau} D(s,s+\tau)\,\bigl[\partial_2D-\partial_1D\bigr](s,s+\tau)\,\dd s .
\end{equation}
\item[\emph{(D4)}] \emph{(mesh hypothesis)} $X$ is a centred
member of $\cP_0(H,L)$ for some $H\in(\tfrac14,1)$, $\sup_n\norm{D_n}_\infty<\infty$, and the
increment covariance splits as $\varrho=\varrho_*+g$ with
$\sup_n|D_n\,\E[\varrho_*\,|\,\cF_n]|\le K|t-s|^{4H-2}$  off the diagonal band $\{|t-s|\le\delta\}$
for a constant $K$, and $g\in L^1([0,1]^2)$, where $D_n$ is the increment variance of the
piecewise-linear interpolant $X^n$ and $\cF_n$ is the $\sigma$-field on $[0,1]^2$ generated by the
cell products $I_a\times I_b$, under Lebesgue measure. Off the diagonal cells,
$\varrho^n=\E[\varrho\,|\,\cF_n]$ is the cell average of $\varrho$.
\end{itemize}
Then
\[
   E_{\rm adj}=\tfrac12\int_0^1 v(t)\Bigl[\tfrac{\dd}{\dd t}C(1,t)-\tfrac12 v'(t)\Bigr]\dd t,
   \qquad
   E_{\rm nest}=\tfrac18 D(0,1)^2+\tfrac18\int_{[0,1]^2}\partial_1D\,\partial_2D\,\dd s\,\dd t,
\]
and $E_{\rm alt}=\tfrac18v(1)^2-E_{\rm adj}-E_{\rm nest}$ by Lemma~\ref{lem:sumzero}.
\end{lemma}

The membership condition of (D4) is used for both coordinates, through Lemma~\ref{lem:L1conv}. Its other
conditions are used only for the nested coordinate. The split $\varrho=\varrho_*+g$ cannot be dropped. For the tied-down perturbation with
$H<\tfrac12$, the ripple kernel \eqref{eq:bHsplit} contains the term
$\epsilon^2\delta_0^{2H-2}\bigl(1-f(u)-f(v)\bigr)$, which blows up like $u^{2H-1}$ along the
lines through the block nodes, also far from the diagonal. So $\varrho$ itself does not satisfy
the domination in (D4). This term is the integrable part $g$ of the split
(Remark~\ref{rem:cov4hyp}, the paragraph on (D4)).

\begin{remark}[the diagonal hypothesis and the threshold $H>\tfrac14$]\label{rem:cov4rel}
For a driver with stationary increments, $D(s,t)=\phi(|t-s|)$ and
$\Xi(\tau)=(1-\tau)\bigl(\phi^2\bigr)'(\tau)$, so (D3) is equivalent to $\phi(\tau)\phi'(\tau)\to0$.
For $\phi(\tau)=\tau^{2H}$ we get $\Xi(\tau)=4H(1-\tau)\tau^{4H-1}$, which tends to $0$ if and only if
$H>\tfrac14$. This is where the threshold $H>\tfrac14$ comes from.
\end{remark}
\begin{proof}[Proof of Lemma~\ref{lem:cov4}]
Both coordinates are first computed for the interpolant $X^n$, whose kernels are bounded, and
then passed to the limit with Lemma~\ref{lem:L1conv}. For the nested word the limit is an integral of
$D\,\partial_1\partial_2D$, which an integration by parts off a diagonal strip of width $\tau$ turns
into the stated form. (D3) removes the boundary term on the strip as $\tau\downarrow0$.

\emph{Step 1 (the adjacent coordinate).} We compute through the interpolants. The
variance $\sigma_a^2:=D(t_a,t_{a+1})$ of the increment over $I_a$ plays the role of the fBm cell
variance $\delta^{2H}$. For a
continuous piecewise-linear $Y$ with $Y_0=0$ the shuffle gives
$\langle1122,\Sig(Y)\rangle=\tfrac12\int_0^1(Y^1_t)^2(Y^2_1-Y^2_t)\,\dd Y^2_t$ (Riemann--Stieltjes), and
$\E[(Y_1-Y_t)\dot Y_t]=\frac{\dd}{\dd t}(\E[Y_1Y_t]-\tfrac12\E Y_t^2)$ inside each cell. So, by
$X^1\perp X^2$, $E_{\rm adj}(X^n)=\tfrac12\int_0^1v_n\,\psi_n'$ with $v_n(t)=\E(X^n_t)^2$ and
$\psi_n(t)=\E[X^n_1X^n_t]-\tfrac12v_n(t)$. With $\lambda=(t-t_a)/\delta$ on $I_a$,
$X^n_t=(1-\lambda)X_{t_a}+\lambda X_{t_{a+1}}$ gives
$v_n=I_nv-\lambda(1-\lambda)\sigma_a^2$ and $\psi_n=I_n\psi+\tfrac12\lambda(1-\lambda)\sigma_a^2$, where
$\psi:=C(1,\cdot)-\tfrac12v$ and $I_n$ is piecewise-linear interpolation at the nodes. Hence
$\psi_n'=(I_n\psi)'+\tfrac12(1-2\lambda)\sigma_a^2/\delta$ on $I_a$. We multiply out
$\tfrac12\int_0^1v_n\psi_n'$ and use, on each cell,
$\int_{I_a}(1-2\lambda)\,I_nv\,\dd t=-\tfrac\delta6\bigl(v(t_{a+1})-v(t_a)\bigr)$,
$\int_{I_a}\lambda(1-\lambda)\,\dd t=\tfrac\delta6$ and
$\int_{I_a}\lambda(1-\lambda)(1-2\lambda)\,\dd t=0$. This gives
\[
   E_{\rm adj}(X^n)=\tfrac12\int_0^1I_nv\,(I_n\psi)'
   +\tfrac1{24}\sum_a\sigma_a^2\bigl(v(t_a)-v(t_{a+1})\bigr)
   -\tfrac1{12}\sum_a\sigma_a^2\bigl(\psi(t_{a+1})-\psi(t_a)\bigr).
\]
The first term equals $\tfrac12\sum_a(\psi(t_{a+1})-\psi(t_a))\tfrac12(v(t_a)+v(t_{a+1}))$, a
Riemann--Stieltjes sum, and tends to $\tfrac12\int_0^1v\,\dd\psi$, $v$ being continuous and $\psi$
absolutely continuous by (D1). The two other terms are bounded by
$\max_a\sigma_a^2\,(\mathrm{TV}(v)+\mathrm{TV}(\psi))\to0$ because $D$ is continuous with $D(x,x)=0$.
Hence $E_{\rm adj}(X^n)\to\tfrac12\int_0^1v(t)\,\dd\psi(t)$.
On the other hand $E_{\rm adj}(X^n)\to E_{\rm adj}$ by Lemma~\ref{lem:L1conv}, $X$ being a member of the
class by (D4). This is the stated formula, $\psi'(t)=\frac{\dd}{\dd t}C(1,t)-\tfrac12v'(t)$.

\emph{Step 2 (the nested coordinate of the interpolants, and its limit). } By the shuffle,
$\int_{t_1<t_2<t_3<t_4}\dd X^2_{t_2}\dd X^2_{t_3}=\tfrac12(X^2_{t_4}-X^2_{t_1})^2$, whence
$E_{\rm nest}(X^n)=\tfrac12\iint_{s<t}D_n\,\varrho^n$ for the piecewise-linear interpolant $X^n$,
with $D_n$ its increment variance and $\varrho^n$ its increment covariance. $\varrho^n$ is a bounded density,
so the interchange is Fubini's. Split at the diagonal band $\{|t-s|\le\delta\}$. There, by (C2)
(part of (D4)), $|D_n|\lesssim\delta^{2H}$ and $|\varrho^n|\lesssim\delta^{2H-2}$ on a set of measure
$\asymp\delta$, so the band contributes $O(\delta^{4H-1})$, which vanishes since $H>\tfrac14$. Off
the band $\varrho^n\to\varrho$ and $D_n\to D$ pointwise, and the family $(D_n\varrho^n)_n$ is
uniformly integrable: by (D4) its $\varrho_*$ part is dominated by the integrable $|t-s|^{4H-2}$,
and its $g$ part is $D_n\E[g|\cF_n]$ with $\sup_n\norm{D_n}_\infty<\infty$. The family
$\{\E[g\,|\,\mathcal G]\}$ is uniformly integrable over any family of sub-$\sigma$-fields when
$g\in L^1$ (the cell $\sigma$-fields are not nested in $n$), and $\E[g|\cF_n]\to g$ a.e.\ by
Lebesgue differentiation. Vitali's theorem gives the passage to the limit.

\emph{Step 3 (the limit as an integral of $D\,\partial_1\partial_2D$).} The convergence
$E_{\rm nest}(X^n)\to E_{\rm nest}$ is Lemma~\ref{lem:L1conv}, valid for every member of the class.
The limit is therefore $E_{\rm nest}=\tfrac12\int_{\{s<t\}}D\,\varrho$, where
$\varrho=\partial_1\partial_2C$ off the diagonal. Off the diagonal
$\partial_1\partial_2D=-2\partial_1\partial_2C$, since the terms $v(s)$ and $v(t)$ of $D$ have no
mixed derivative. Since both $D$ and $\partial_1\partial_2D$ are symmetric,
$E_{\rm nest}=\tfrac14\int_{[0,1]^2}D\,\varrho=-\tfrac18\int_{[0,1]^2}D\,\partial_1\partial_2D$.
For $H\le\tfrac12$ the off-diagonal density of $\partial_1\partial_2D$ is not integrable across the
diagonal, and the distribution is its finite part. But $D$ vanishes on the diagonal, and the second
clause of (D2) says that the product is in $L^1$, so the finite part is the ordinary integral and
\begin{equation}\label{eq:nestlimit}
   \int_{[0,1]^2}D\,\partial_1\partial_2D\ =\ \lim_{\tau\downarrow0}\ \int_{\{|t-s|>\tau\}}
   D\,\partial_1\partial_2D
\end{equation}
by dominated convergence, both sides being absolutely convergent.

\emph{Step 4 (integration by parts off the diagonal).} Fix $\tau\in(0,1)$ and integrate by parts in the first variable on each of the two components of
$\{|t-s|>\tau\}$, which (D1) allows. On the upper piece, for fixed $t\in(\tau,1)$,
\[
  \int_0^{t-\tau}\!\! D\,\partial_1\partial_2D\,\dd s
  =\bigl[D\,\partial_2D\bigr]_{s=0}^{s=t-\tau}-\int_0^{t-\tau}\!\!\partial_1D\,\partial_2D\,\dd s ,
\]
and $\int_\tau^1 D(0,t)\,\partial_2D(0,t)\,\dd t=\tfrac12\bigl[D(0,1)^2-D(0,\tau)^2\bigr]$, while
the inner boundary contributes $\int_0^{1-\tau}D(s,s+\tau)\,\partial_2D(s,s+\tau)\,\dd s$. The lower
piece is its mirror image: the outer boundary $s=1$ gives
$\tfrac12\bigl[D(1-\tau,1)^2-D(0,1)^2\bigr]$ and the inner boundary
$-\int_0^{1-\tau}D(s,s+\tau)\,\partial_1D(s,s+\tau)\,\dd s$, the symmetry of $D$ turning
$\partial_2D(t+\tau,t)$ into $\partial_1D(t,t+\tau)$. Collecting, and with $\Xi$ as in
\eqref{eq:diagbdry},
\[
  \int_{\{|t-s|>\tau\}}\!\! D\,\partial_1\partial_2D
  \ =\ \Xi(\tau)-D(0,1)^2+\tfrac12D(0,\tau)^2+\tfrac12D(1-\tau,1)^2
  -\int_{\{|t-s|>\tau\}}\!\!\partial_1D\,\partial_2D .
\]
Let $\tau\downarrow0$: the two square terms vanish by \emph{(D1)}, the last integral converges to
$\int_{[0,1]^2}\partial_1D\,\partial_2D$ by dominated convergence under \emph{(D2)}, and the
diagonal term vanishes by \emph{(D3)}. Hence \eqref{eq:nestlimit} equals
$-D(0,1)^2-\int_{[0,1]^2}\partial_1D\,\partial_2D$, and multiplying by $-\tfrac18$ (Step~3) gives the stated $E_{\rm nest}$.
\end{proof}

Remark~\ref{rem:cov4hyp} in Appendix~\ref{app:verif} checks the hypotheses for the two laws to
which Lemma~\ref{lem:cov4} is applied.

\begin{remark}[the continuum coordinates without analytic continuation]\label{rem:noAC}
Lemma~\ref{lem:cov4} reproduces Proposition~\ref{prop:altcontinuum} on all of $(\tfrac14,1)$ by direct
evaluation. For fBm, $v(t)=t^{2H}$, $\tfrac{\dd}{\dd t}C(1,t)=H[t^{2H-1}+(1-t)^{2H-1}]$ and
$\int_{[0,1]^2}\partial_1D\,\partial_2D=-4H^2\int_{[0,1]^2}|t-s|^{4H-2}=-\tfrac{2H}{4H-1}$, giving
$E^\infty_{\rm nest}=\tfrac18-\tfrac{2H}{8(4H-1)}=\tfrac{2H-1}{8(4H-1)}$ and, after one Beta integral,
$E^\infty_{\rm adj}=\tfrac H2\int_0^1t^{2H}(1-t)^{2H-1}\dd t=\tfrac H2 \mathrm{B}(2H{+}1,2H)=\tfrac{H}{4(4H-1)}\Gamma(2H)^2/\Gamma(4H-1)$. The
alternating value follows from Lemma~\ref{lem:sumzero}. So for every $H\in(\tfrac14,1)$ these
formulas are values of absolutely convergent covariance integrals, and the proof of
Proposition~\ref{prop:altcontinuum} needs no analytic continuation in $H$.
\end{remark}

The proof of Proposition~\ref{prop:seprough} writes the leading coefficient of the nested
separation through two functionals $W$ and $\Theta$ of tied-down fBm. The next lemma gives their signs,
which make that coefficient strictly negative.

\begin{lemma}[two spectral identities, and the signs of $W$ and $\Theta$]\label{lem:opsigns}
Let $\Bop$, $\Cop$ be the integral operators on $L^2[0,1]$ with kernels $B_H$ (the unit
covariance of tied-down fBm) and $C$ (the fBm$(H)$ covariance), $H\in(\tfrac14,1)$, and set
\begin{gather*}
   \Theta(H):=\!\int_{[0,1]^2}\!\partial_uB_H(u,v)\,\partial_vB_H(u,v)\,\dd u\dd v,\\
   W(H):=\!\int_{[0,1]^2}\!\partial_uB_H(u,v)\,\mathrm{sgn}(v-u)|v-u|^{2H-1}\dd u\dd v .
\end{gather*}
Let $\Bop e_k=\lambda_ke_k$ and $\Cop f_m=\mu_mf_m$ be the (nonnegative) eigensystems,
and put $A_{kl}:=\int_0^1e_k'e_l$, $P_{km}:=\int_0^1e_k'f_m$. Then all $e_k$ lie in $H^1[0,1]$ and
vanish at $0$ and at $1$, the two double series below converge, and
\[
   \boxed{\ \Theta=-\sum_{k,l}\lambda_k\lambda_l\,A_{kl}^{2}\ <\ 0,
   \qquad
   W=\frac{1}{H}\sum_{k,m}\lambda_k\mu_m\,P_{km}^{2}\ >\ 0 .\ }
\]
\end{lemma}
\begin{proof}
Both functionals are integrals of products of first derivatives of covariance kernels. We expand
each kernel in the eigenbasis of its operator. The eigenfunctions of $\Bop$ vanish at both endpoints,
so an integration by parts makes the matrix $A$ antisymmetric and turns $Q$ into $-P$. Each double
series then becomes a sum of squares with a definite sign.

\emph{Step 0 (regularity of the eigenfunctions of $\Bop$).} Apart from
$H\,\mathrm{sgn}(v-u)|v-u|^{2H-1}$, the function $\partial_uB_H$ equals
$H(1-v)u^{2H-1}-vH(1-u)^{2H-1}-C(v,1)+v$. This is unbounded for $H<\tfrac12$, but its endpoint
singularities are square-integrable under the same condition as the diagonal one. So
$\partial_uB_H\in L^2([0,1]^2)$ exactly when $4H-2>-1$, i.e.\ $H>\tfrac14$, the threshold at
which the L\'evy area exists. Then $g\mapsto\partial(\Bop g)$ is Hilbert--Schmidt on $L^2$ and
$\Bop g\in H^1[0,1]$ for $g\in L^2$. Hence $e_k=\lambda_k^{-1}\Bop e_k$ lies in $H^1$, and
$e_k(0)=\lambda_k^{-1}\int_0^1B_H(0,v)e_k(v)\dd v=0$, likewise $e_k(1)=0$, because the tied-down
fBm covariance $B_H$ vanishes on the boundary. (Only $\mathrm{Ran}(\Bop)$ is used. The range of
$\Bop^{1/2}$ is not contained in $H^1$ for $H<\tfrac12$.)

\emph{Step 1 (the operators $\Bop$ and $\Cop$ are injective).} By \eqref{eq:fbmspectral}, for a finite signed
measure $\nu$ on $[0,1]$ and $Z$ an fBm$(H)$,
\[
   \Var\int_0^1Z\,\dd\nu=\tfrac{a_H}2\int_\RR\Bigl|\int(e^{\mathrm is\xi}-1)\,\nu(\dd s)\Bigr|^2|\xi|^{-1-2H}\dd\xi .
\]
If this vanishes, the Fourier transform of $\nu$, which is continuous, equals its value at $0$
everywhere, so $\nu$ is a multiple of $\delta_0$. For $\Bop$: $B_H$ is the covariance of the
tied-down fBm $\varpi(u)=Z_u-uZ_1$, so $\Bop g=0$ gives
$0=\langle g,\Bop g\rangle=\Var\int_0^1Z\,\dd\nu$ for the signed measure
$\nu(\dd s)=g(s)\dd s-(\int_0^1ug)\,\delta_1(\dd s)$. As $\nu(\{0\})=0$, $\nu=0$, so $g=0$ almost
everywhere. For $\Cop$ the same argument applies with $\nu(\dd s)=g(s)\dd s$, which has no atom.
Hence all $\lambda_k$ and $\mu_m$ are positive, and $(e_k)$ and $(f_m)$ are complete in $L^2[0,1]$.

\emph{Step 2 ($\Theta$).} By Mercer, $B_H(u,v)=\sum_k\lambda_ke_k(u)e_k(v)$ uniformly.
$\partial\Bop$ being Hilbert--Schmidt, the differentiated series
$\partial_uB_H(u,v)=\sum_k\lambda_ke_k'(u)e_k(v)$ converges in $L^2([0,1]^2)$. By symmetry of $B_H$,
$\partial_vB_H(u,v)=(\partial_uB_H)(v,u)=\sum_l\lambda_le_l(u)e_l'(v)$, also in $L^2([0,1]^2)$, so the
product integrates term by term:
\[
   \Theta=\sum_{k,l}\lambda_k\lambda_l\Bigl(\int_0^1e_k'e_l\Bigr)\Bigl(\int_0^1e_ke_l'\Bigr)
   =\sum_{k,l}\lambda_k\lambda_l\,A_{kl}A_{lk},
\]
absolutely convergent by Cauchy--Schwarz and Bessel's inequality
$\sum_lA_{kl}^2=\sum_l\langle e_k',e_l\rangle^2\le\norm{e_k'}^2$:
\[
   \sum_{k,l}|\lambda_kA_{kl}|\,|\lambda_lA_{lk}|
   \ \le\ \sum_k\lambda_k^2\norm{e_k'}^2\ =\ \sum_k\norm{\partial(\Bop e_k)}^2\ <\ \infty ,
\]
the last sum being the squared Hilbert--Schmidt norm of $g\mapsto\partial(\Bop g)$ (Step~0), as
$(e_k)$ is complete (Step~1). As
$e_k,e_l\in H^1$ vanish at both endpoints, $A_{kl}+A_{lk}=[e_ke_l]_0^1=0$, i.e.\ $A_{lk}=-A_{kl}$, and
$\Theta=-\sum_{k,l}\lambda_k\lambda_lA_{kl}^2\le0$. Equality would force $A_{kl}=0$ whenever
$\lambda_k\lambda_l>0$. By Step~1 all $\lambda_l>0$ and $(e_l)$ is complete in
$L^2[0,1]$, so $e_k'\perp e_l$ for all $l$ gives $e_k'=0$. So $e_k$ is constant, hence $e_k\equiv0$ by the
boundary values, contradicting $\norm{e_k}=1$. Thus $\Theta<0$.

\emph{Step 3 ($W$).} Write $\mathrm{sgn}(v-u)|v-u|^{2H-1}=\tfrac1{2H}\partial_vD_0(u,v)$ with
$D_0(u,v)=|v-u|^{2H}=v_0(u)+v_0(v)-2C(u,v)$, $v_0(t)=t^{2H}$. The $v_0'(v)$ part contributes
$\tfrac1{2H}\int_0^1v_0'(v)\bigl[B_H(1,v)-B_H(0,v)\bigr]\dd v=0$, the covariance $B_H$ vanishing on
the boundary. Hence
\[
   W=-\frac{1}{H}\int_{[0,1]^2}\partial_uB_H(u,v)\,\partial_vC(u,v)\,\dd u\,\dd v .
\]
Expanding both kernels as in Step~2, $\partial_vC(u,v)=\sum_m\mu_mf_m(u)f_m'(v)$ in $L^2([0,1]^2)$
(the same computation gives $\partial_uC\in L^2([0,1]^2)$ for $H>\tfrac14$), so
\[
   \int_{[0,1]^2}\partial_uB_H\,\partial_vC
   =\sum_{k,m}\lambda_k\mu_m\Bigl(\int_0^1e_k'f_m\Bigr)\Bigl(\int_0^1e_kf_m'\Bigr)
   =\sum_{k,m}\lambda_k\mu_m\,P_{km}Q_{km},
\]
$Q_{km}:=\int_0^1e_kf_m'$. Integration by parts gives $P_{km}+Q_{km}=[e_kf_m]_0^1=0$, since $e_k$
vanishes at the endpoints (no boundary behaviour of $f_m$ is needed). Thus $Q_{km}=-P_{km}$ and
$W=\tfrac1H\sum_{k,m}\lambda_k\mu_mP_{km}^2\ge0$. Equality would force
$P_{km}=0$ for every $m$, so $e_k'\perp\overline{\mathrm{Ran}\,\Cop}=L^2$ ($\Cop$ is injective
by Step~1), whence $e_k'=0$ and $e_k\equiv0$ as in Step~2. So $W>0$.

\emph{Step 4 (convergence of the series for $W$).} It converges absolutely by the same Cauchy--Schwarz and Bessel bound as in Step~2, now
with $\sum_mP_{km}^2\le\norm{e_k'}^2$ and $\sum_kQ_{km}^2\le\norm{f_m'}^2$, which bound it by the
product of the Hilbert--Schmidt norms of $g\mapsto\partial(\Bop g)$ and $g\mapsto\partial(\Cop g)$.
\end{proof}

\begin{remark}[values of $W$ and $\Theta$]\label{rem:WTheta}
The $\mp H\,\mathrm{sgn}(u-v)|u-v|^{2H-1}$ term of $\partial_uB_H$ gives both constants explicit
singular parts, $W=\tfrac1{2(4H-1)}+W_{\rm reg}$ and
$\Theta=-\tfrac{H}{2(4H-1)}+\Theta_{\rm reg}$, with absolutely convergent remainders. Numerically,
\[
\begin{array}{c|cccccc}
H & 0.28 & 0.30 & 0.35 & 0.40 & 0.45 & \tfrac12\\\hline
W(H) & 4.2459 & 2.5315 & 1.1947 & 0.7231 & 0.4800 & \tfrac13\\
\Theta(H) & -1.0866 & -0.6566 & -0.3165 & -0.1918 & -0.1249 & -\tfrac1{12}
\end{array}
\]
and both diverge at $H=\tfrac14$. The Brownian values $W(\tfrac12)=\tfrac13$,
$\Theta(\tfrac12)=-\tfrac1{12}$ are those required for consistency with Proposition~\ref{prop:sepbrownian}
(Remark~\ref{rem:roughsepbrownian}), and Lemma~\ref{lem:opsigns} reproduces
$\Theta(\tfrac12)=-\tfrac{16}{\pi^4}\sum_{k+l\ \mathrm{odd}}(l^2-k^2)^{-2}=-\tfrac1{12}$.
\end{remark}

The next proposition gives the lower bound of Corollary~\ref{cor:undetectable} on the rough range
$(\tfrac14,\tfrac12)$. It uses the nested word, whose separation is of the larger order
$\delta^{4H-1}$ there (Remark~\ref{rem:roughsepbrownian}).

\begin{proposition}[separation in the rough regime]\label{prop:seprough}
Let $H\in(\tfrac14,\tfrac12)$ and let $(P_0,P_1)$ be the pair of Construction~\ref{con:bridge}
with $x:=\epsilon^2>0$ fixed. Then
\[
   \mu(P_1)_{\rm nest}-\mu(P_0)_{\rm nest}
   =\Bigl[-x\,H\,W(H)+\tfrac{x^2}{2}\,\Theta(H)\Bigr]\,\delta^{4H-1}\;+\;O(\delta^{2H}),
\]
and by Lemma~\ref{lem:opsigns} the bracket is strictly negative for every $x>0$. Consequently
$\norm{\mu(P_0)-\mu(P_1)}_\hs\gtrsim\delta^{4H-1}$ for $M\ge4$.
\end{proposition}
\begin{proof}
By Lemma~\ref{lem:cov4} the nested coordinate is a quadratic functional of the increment variance, and the
ripple adds $xD_R$ to the increment variance. We expand in $x$ and evaluate the two terms by rescaling each cell
exactly to the unit cell, where they become $W$ and $\Theta$.

\emph{Step 1 (reduction to the increment variances).} Both laws are centred Gaussian with i.i.d.\ components and satisfy (D1)--(D4)
(Remark~\ref{rem:cov4hyp}(a),(b)), so Lemma~\ref{lem:cov4} applies to each. Their
increment variances satisfy $D_1=D_0+xD_R$, $D_R(s,t)=v_R(s)+v_R(t)-2\Gamma_\delta(s,t)$,
$v_R(t)=\Var(R_t)$, with $\Gamma_\delta=\Cov(R_s,R_t)$ block diagonal:
$\Gamma_\delta(s,t)=\delta^{2H}B_H(\tfrac{s-t_a}\delta,\tfrac{t-t_a}\delta)$ on $I_a$, zero
across cells. As $R$ vanishes at every node, $D_1(0,1)=D_0(0,1)$, so the first term
of Lemma~\ref{lem:cov4} cancels. The cross terms are equal by symmetry of $D_0,D_R$, so, with
$\Delta_{\rm nest}:=\mu(P_1)_{\rm nest}-\mu(P_0)_{\rm nest}$,
\[
   \Delta_{\rm nest}=\frac{x}{4}\int_{[0,1]^2}\partial_1D_R\,\partial_2D_0
   +\frac{x^2}{8}\int_{[0,1]^2}\partial_1D_R\,\partial_2D_R .
\]

\emph{Step 2 (the $x^2$ term).} The $x^2$ term equals $\tfrac{x^2}{2}\Theta\,\delta^{4H-1}$ exactly. Expand
$\partial_1D_R=v_R'(s)-2\partial_s\Gamma_\delta$ and $\partial_2D_R=v_R'(t)-2\partial_t\Gamma_\delta$:
\begin{align*}
   \int_{[0,1]^2}\partial_1D_R\,\partial_2D_R
   &=\Bigl(\int_0^1v_R'\Bigr)^2
   -2\int_0^1v_R'(s)\bigl[\Gamma_\delta(s,1)-\Gamma_\delta(s,0)\bigr]\dd s\\
   &\quad-2\int_0^1v_R'(t)\bigl[\Gamma_\delta(1,t)-\Gamma_\delta(0,t)\bigr]\dd t
   +4\int_{[0,1]^2}\partial_s\Gamma_\delta\,\partial_t\Gamma_\delta .
\end{align*}
The first three terms vanish, $v_R$ and $\Gamma_\delta(s,\cdot)$ vanishing at every node. In the
last, $\Gamma_\delta$ is block diagonal, and on $I_a\times I_a$ the substitution
$s=t_a+\delta u$, $t=t_a+\delta v$ gives $\partial_s\Gamma_\delta=\delta^{2H-1}\partial_uB_H(u,v)$,
$\partial_t\Gamma_\delta=\delta^{2H-1}\partial_vB_H(u,v)$ and $\dd s\,\dd t=\delta^2\dd u\,\dd v$, so
each of the $n=\delta^{-1}$ cells contributes exactly $\delta^{4H}\Theta$ and
$4\int\int\partial_s\Gamma_\delta\,\partial_t\Gamma_\delta=4\sum_a\delta^{4H}\Theta=4\delta^{4H-1}
\Theta$.

\emph{Step 3 (the $x$ term).} Split $\partial_1D_R$ again. Integrating the $v_R'$ part in
$t$, and then by parts in $s$ (with $v_R(0)=v_R(1)=0$ and $D_0(s,1)-D_0(s,0)=(1-s)^{2H}-s^{2H}$),
turns it into $\int_0^1v_R'(s)[D_0(s,1)-D_0(s,0)]\dd s=2H\int_0^1v_R(s)[s^{2H-1}+(1-s)^{2H-1}]\dd s
=O(\delta^{2H})$, since $v_R=O(\delta^{2H})$. The $-2\partial_s\Gamma_\delta$ part is block diagonal
and rescales exactly, with $\partial_tD_0(s,t)=2H\,\mathrm{sgn}(t-s)|t-s|^{2H-1}
=2H\delta^{2H-1}\mathrm{sgn}(v-u)|v-u|^{2H-1}$ on $I_a\times I_a$, to
$-2\sum_a2H\delta^{4H}W=-4H\,W\,\delta^{4H-1}$.

\emph{Conclusion.} Multiplying Step~2 by $\tfrac{x^2}8$ and Step~3 by $\tfrac x4$, the
main terms are $\tfrac{x^2}2\Theta\,\delta^{4H-1}$ and $-xHW\delta^{4H-1}$, and the only error term is
the $v_R'$ part of the $x$ term, $O(\delta^{2H})$, which is $o(\delta^{4H-1})$ because $2H>4H-1$ for
$H<\tfrac12$. Collecting gives the result. For the last assertion, the level-$4$ block of
$\mu(P_0)-\mu(P_1)$ contains the two nested words $1221$ and $2112$, so
$\norm{\mu(P_0)-\mu(P_1)}_\hs\ge\sqrt2\,|\Delta_{\rm nest}|\gtrsim\delta^{4H-1}$.
\end{proof}

\begin{remark}[why the nested word, and consistency at $H=\tfrac12$]\label{rem:roughsepbrownian}
Lemma~\ref{lem:cov4} also gives the adjacent coordinate in closed form. Since
$\Gamma_\delta(1,t)\equiv0$ ($t=1$ is a node),
\[
   \mu(P_1)_{\rm adj}-\mu(P_0)_{\rm adj}=\frac{\epsilon^2}{2}\int_0^1v_R(t)\,\partial_tC_0(1,t)\,\dd t
   \ >\ 0 ,
\]
exactly, with both factors positive. This is $\asymp\delta^{2H}$ for every $H$, with leading constant
$\tfrac{\epsilon^2}{2}\int_0^1B_H(u,u)\dd u=-\tfrac12\kappa_{\rm sep}(H)\epsilon^2$. It recovers
Proposition~\ref{prop:sep} in the Young regime and $\Dsep_{\rm adj}=\epsilon^2/12$ at $H=\tfrac12$.
Below the crossover $\delta^{2H}=o(\delta^{4H-1})$, so the adjacent word cannot carry the rough-regime
lower bound, while the nested word does. At $H=\tfrac12$ the exponents merge and
Proposition~\ref{prop:seprough} becomes exact: $W(\tfrac12)=\tfrac13$,
$\Theta(\tfrac12)=-\tfrac1{12}$ and the $\delta^{2H}$ term is
$\tfrac x4\cdot2\int_0^1v_R=\tfrac{x\delta}{12}$, so
$\Delta_{\rm nest}=\delta(\tfrac x{12}-\tfrac x6-\tfrac{x^2}{24})=-\delta(\tfrac x{12}+\tfrac{x^2}{24})
=\delta \Dsep_{\rm nest}(\epsilon)$, exactly Proposition~\ref{prop:sepbrownian}.
\end{remark}

\subsection*{Consequence: no mesh-data procedure removes the weak error uniformly}

Proposition~\ref{prop:adm} and the separation bounds give the following statement directly.

\begin{corollary}[a lower bound for mesh-data procedures]\label{cor:undetectable}
Fix $M\ge4$, $d\ge2$, $H\in(\tfrac14,1)$ and $\epsilon>0$, let $L\ge L_0(H,\epsilon)$ as in
Proposition~\ref{prop:adm}, and let $P_0,P_1$ be the pair of Construction~\ref{con:bridge} at mesh
$\delta$. Let $T$ be a measurable function of the mesh data $\cD$ with values in $\TM$ and
$\E_{P_0}\norm T_\hs<\infty$. Then
\[
   \sup_{P\in\cP_0(H,L)}\bigl\lVert\E_{P}[T]-\mu(P)\bigr\rVert_\hs\ \ge\
   \max_{j\in\{0,1\}}\ \bigl\lVert\E_{P_j}[T]-\mu(P_j)\bigr\rVert_\hs\ \ge\
   \tfrac12\norm{\mu(P_1)-\mu(P_0)}_\hs\ \ge\ c\,\delta^{\gamma(H)}
\]
for all sufficiently fine meshes $\delta$ (for every mesh at $H=\tfrac12$), with $c=c(H,\epsilon,d)>0$
independent of $\delta$ and of $T$.
\end{corollary}

\begin{proof}
The first inequality holds because $P_0,P_1\in\cP_0(H,L)$ (Proposition~\ref{prop:adm}(a)). By
Proposition~\ref{prop:adm}(b), $\cD$ has the same law under $P_0$ and $P_1$, so
$\E_{P_0}[T]=\E_{P_1}[T]=:m$. By the triangle inequality,
$\norm{\mu(P_1)-\mu(P_0)}_\hs\le\norm{m-\mu(P_0)}_\hs+\norm{m-\mu(P_1)}_\hs$, which is at most twice the
maximum. This is the second inequality. For the third: for $H\in(\tfrac12,1)$,
Proposition~\ref{prop:sep} gives $\norm{\mu(P_0)-\mu(P_1)}_\hs\ge\sqrt3|\kappa_{\rm sep}(H)|\epsilon^2\delta^{2H}+o(\delta^{2H})$
with $\kappa_{\rm sep}(H)\ne0$; at $H=\tfrac12$, Proposition~\ref{prop:sepbrownian} gives
$\norm{\mu(P_0)-\mu(P_1)}_\hs\ge\delta\bigl(2\sum_\bullet\Dsep_\bullet(\epsilon)^2\bigr)^{1/2}$ with a nonzero
constant, for every mesh, and for $H\in(\tfrac14,\tfrac12)$, Proposition~\ref{prop:seprough} gives
$\norm{\mu(P_0)-\mu(P_1)}_\hs\gtrsim\delta^{4H-1}$. In each case the exponent is $\gamma(H)$.
\end{proof}

Taking for $T$ the signature $\Sig^{\le M}(X^n)$ of the interpolant, whose mean is the same under both
laws, shows that $\max_j\norm{\cW_{\le M}(P_j)}_\hs\ge c\,\delta^{\gamma(H)}$, so the upper bounds of
\S\S\ref{sec:bias}--\ref{sec:rough} are attained within $\cP_0(H,L)$. More generally, no correction of
the weak error computed from the mesh data can reduce it below the order $\delta^{\gamma(H)}$
uniformly over the class: any such corrected quantity is again a function $T$ of $\cD$. The statement
is about the class, since $P_1$ is built for the mesh at hand (Remark~\ref{rem:meshtied}).

\section{Discussion and open questions}\label{sec:discussion}

\paragraph{Summary.} For Gaussian drivers with independent components the weak error of the
piecewise-linear interpolant is of order $\delta^{\gamma(H)}$, $\gamma(H)=\min(2H,4H-1)$, uniformly over
the class and at every truncation level, and for $M\ge4$ this order is attained. On the centred class
the weak error first appears at level~$4$, where we compute its leading constant for independent
fBm. For the even levels $6,\dots,M$ we prove only the upper bound $O(\delta^{\gamma(H)})$. The strong error is larger: for fBm, in mean square, a full power of
$\delta$ larger on $(\tfrac12,1)$, and below the crossover, up to an arbitrarily small loss in the
exponent, the square root of the weak error. Finally, no procedure based on the mesh data removes the
weak error uniformly over the class. In the estimation
problem of \citet{HuberII} the weak error is the bias of the block average, and these results
determine the discretization part of the minimax rate there.

\paragraph{The constants are fixed-$H$ limits.} The leading constants hold as $\delta\to0$ at fixed
$H$, not uniformly near the crossover, where the Young and rough contributions add
(Remarks~\ref{rem:nonuniform}, \ref{rem:merger} and~\ref{rem:sepmerger}, and the comment after
Theorem~\ref{thm:roughrate}).

\paragraph{Open questions.}

\emph{(i) Non-Gaussian drivers.} Centring and independent components already give $\cW_3=0$ and make
$\cW_4$ a functional of the component covariances (Proposition~\ref{prop:nongauss}). Without
independence, vanishing mixed third moment functions suffice for $\cW_3=0$. Non-Gaussianity first
enters the bias at level~$5$, through third-order moment functions of single components, and not at
all at odd levels for symmetric components. Open is the size of $\cW_5$ for asymmetric laws such as
independent Rosenblatt processes.

\emph{(ii) Constants beyond level~$4$ for $H<\tfrac12$.} The rate $\norm{\cW_{\le M}}=O(\delta^{4H-1})$ holds
at every $M$ (Propositions~\ref{prop:roughalllevels}, \ref{prop:drift5}), but $\cB(H)$ is only the
level-$4$ constant. The full-truncation constant also collects, on $\cP_0$, the even levels $6,\dots,M$ (of the
same order at level~$6$ numerically; Remark~\ref{rem:M6}) and, over $\cP$, the drift terms at levels
$\ge5$, which attain $\delta^{4H-1}$ for drifted independent fBm (proved for $d\ge3$;
Proposition~\ref{prop:attain5}). It is not evaluated
here.

\emph{(iii) The strong error beyond fBm.} The exact order $\delta^{2H-1/2}$ of
Proposition~\ref{prop:strong} is proved for fBm drivers on $(\tfrac12,1)$. We do not prove the
matching rate for other drivers, such as fractional Ornstein--Uhlenbeck processes, which are natural
candidates for the class but which we have not verified to belong to it. Below the
crossover the class-uniform bound of Proposition~\ref{prop:strongunif} loses an arbitrarily small
power of $\delta$ against the lower bound of Lemma~\ref{lem:stronglower}.

\emph{(iv) Scope.} $M$ is fixed, and we do not treat the asymptotics as $M\to\infty$. The rates hold
for every $d\ge2$ (and trivially at $d=1$), while the leading constants of \S\ref{sec:bias} are computed
at $d=2$ and scale as in
Remark~\ref{rem:dscope}. The lower limit $H=\tfrac14$ comes from the construction used throughout:
\citet{CoutinQian2002} obtain the geometric lift as an $L^2$ limit of piecewise-linear approximations
on that range. Open are correlated fractional drivers, where Lemmas~\ref{lem:weights}
and \ref{lem:momentid} fail and extra Wick pairings survive. For branched (non-geometric) signatures
\citep{Gubinelli2010,HairerKelly2015} we expect that the level-$4$ improvement does not carry over,
since the symmetric level-$2$ bracket is then no longer determined by the endpoint, so that a bias may
reappear at level~$2$. We have not proved this.

\clearpage
\appendix

\section{Verification of the class conditions}\label{app:verif}

This appendix checks two things: that the class of Definition~\ref{def:class} contains, at every mesh,
the three laws the paper uses; and that the hypotheses of Lemma~\ref{lem:cov4} hold at the two laws
to which that lemma is applied.

\begin{remark}[verification of (C1)--(C5) and of the moment condition. This proves Proposition~\ref{prop:adm}(a)] \label{rem:C1C2}
We check (C1)--(C5), at every mesh
$\delta=1/n$, for the three laws used in this paper: the independent-fBm law $P_0$, its tied-down
perturbation $P_1^{(\delta_0)}$ at every block width $\delta_0=1/n_0$ (Construction~\ref{con:bridge})
and the drifted law $P_\theta$, the law of $\theta te_1$ plus a path of $P_0$. All constants depend on
$H$, $\epsilon$, $d$ and the fixed truncation level $M$ only, and in the drifted case also on $\theta$. $P_\theta$ has the Gaussian part of $P_0$, so only its moment condition needs
a word. We write $h:=\delta/\delta_0$ and $J$ for the blocks, and say that a cell straddles a
block node if the node lies in its interior.

\emph{The ripple kernel.} By pinning and the independence of the blocks, the increment covariance of
the ripple $R$ vanishes across blocks. On a block $J$, in block coordinates $u,v\in[0,1]$, it has for
$H\ne\tfrac12$ the density
\begin{equation}\label{eq:bHsplit}
   \varrho_R\ =\ \underbrace{c_H|s-t|^{2H-2}}_{\text{the fBm kernel, restricted}}
   \ +\ \delta_0^{2H-2}g(u,v),\qquad g:=1-f(u)-f(v),
\end{equation}
with $f(u)=H\bigl(u^{2H-1}+(1-u)^{2H-1}\bigr)$. This is $\delta_0^{2H-2}b_H(u,v)$
(Proposition~\ref{prop:sep}). At $H=\tfrac12$ the first term is the
diagonal measure $\dd s\,\delta_{\{s=t\}}$ and $g\equiv-1$ (the Brownian bridge). The law
$P_1^{(\delta_0)}$ has the kernel $\varrho_1=\varrho_0+\epsilon^2\varrho_R$, $\varrho_0$ that of fBm. We use
three facts about $f$: $\int_0^1f=1$. For an interval $I\subseteq[0,1]$ of length $\ell$,
$\int_If\le\ell^{2H}$ if $H<\tfrac12$ (as $u^{2H-1}$ is decreasing) and $\int_If\le2H\ell$ if $H>\tfrac12$,
and $\int_0^1|f-P_hf|\le Ch$ for $H>\tfrac12$ ($f$ has integrable derivative) and $\le Ch^{2H}$ for
$H<\tfrac12$ (the two end cells contribute at most $4h^{2H}$, since $\int_I|f-P_hf|\le2\int_If\le2h^{2H}$ on each, and the cell at distance $kh$ from the nearer end
at most $h^2\sup|f'|\le Ch^{2H}k^{2H-2}$), $P_h$ being averaging over the cells of width $h$. For
$H\ge\tfrac12$ the sections of the ripple kernel are bounded:
$\int_J|\varrho_R|(s,\dd t)\le C\delta_0^{2H-1}$, since the restricted fBm kernel integrates to at most
$2H\delta_0^{2H-1}$ and $|g|\le1+4H$ (total variation $2$ at $H=\tfrac12$).

\emph{(C1).} For $P_0$ and $H\ne\tfrac12$, the off-diagonal cell pairs ($a\ne b$) contribute
$O(\delta^{2H-1})$ by Remark~\ref{rem:C1scope}, which treats the fBm kernel on all of $(\tfrac14,1)$;
this is (C1) for $H<\tfrac12$. For $H>\tfrac12$ the $n$ diagonal cells add at most
$\sum_a2\int_{I_a^2}\varrho=2n\delta^{2H}=2\delta^{2H-1}$, because cell averaging at most doubles the
$L^1$ norm on a cell and $\int_{I_a^2}\varrho=\delta^{2H}$. At $H=\tfrac12$ the total variation of $e$
is at most $2$. For $P_1^{(\delta_0)}$, cell averaging is linear, so it suffices to bound the cell-averaging error
$e_R$ of $\varrho_R$ by $C\delta^{2H-1}$ in $L^1$ (in total variation at $H=\tfrac12$; over
off-diagonal cells for $H<\tfrac12$). 

\emph{The case $\delta\ge\delta_0$.} Then
$\norm{e_R}\le2\norm{\varrho_R}$ over the relevant cells. For $H>\tfrac12$,
$\norm{\varrho_R}_{L^1}=n_0\delta_0^{2H}\norm{b_H}_{L^1}=C\delta_0^{2H-1}\le C\delta^{2H-1}$, at
$H=\tfrac12$ the total variation is $2$. For $H<\tfrac12$ a block meets two different cells only if it
straddles a mesh node, so at most $n$ blocks contribute, each at most $C\delta_0^{2H}\le C\delta^{2H}$ on
off-diagonal cell pairs, in total $C\delta^{2H-1}$.

\emph{The case $\delta<\delta_0$.} Split the cell pairs
into three kinds. (i) Both cells inside one block: there $e_R$ is $\delta_0^{2H}$ times the error of
averaging $b_H$ over the (possibly shifted) lattice of width $h$ in block coordinates. Its fBm part
is $O(h^{2H-1})$ by the bound for $P_0$ just proved, applied in block units at mesh width $h$
(over off-diagonal cell pairs if $H<\tfrac12$), the shift of the lattice does not matter, as the fBm
kernel is translation invariant, and its $g$ part is at most $2\int_0^1|f-P_hf|=O(h^{2H-1})$, as $h<1$. So each block contributes
$C\delta_0^{2H}h^{2H-1}=C\delta_0\delta^{2H-1}$, and the $n_0$ blocks $C\delta^{2H-1}$. (ii) Cells in
different blocks, neither straddling: $\varrho_R=0$ there. (iii) Pairs with a straddling cell $I_a$, of
which there are at most $n_0$. For $H\ge\tfrac12$ the row and column of $I_a$ carry
$\norm{e_R}\le4\delta\sup_s\int|\varrho_R|(s,\dd t)\le C\delta\delta_0^{2H-1}$, in total
$Cn_0\delta\delta_0^{2H-1}=C\delta\delta_0^{2H-2}\le C\delta^{2H-1}$. For $H<\tfrac12$ the off-diagonal
pairs of the row and the column of $I_a$ carry at most $4\sum_{b\ne a}\int_{I_a\times I_b}|\varrho_R|\le C\delta^{2H}$ (the
kernel is symmetric) by
the kernel clause of (C4), verified below, in total $Cn_0\delta^{2H}\le C\delta^{2H-1}$.

\emph{(C2).} For $P_0$ the rescaled cell path $\delta^{-H}(G_{t_a+\delta\cdot}-G_{t_a})$ has the law of
fBm on $[0,1]$, by stationarity of the increments and self-similarity. For $H>\tfrac12$ its
signature is the Young one, bounded by $c_{M,p}(1\vee\norm{\cdot}_{\pvar})^M$, and the path
$p$-variation has moments of every order by the Gaussian Kolmogorov criterion
\citep[Thm.~A.10]{FrizVictoir2010}: for $\alpha\in(1/p,H)$, $\norm z_{\pvar}\le\norm z_{\alpha\text{-H\"ol}}$
on $[0,1]$. For $H\le\tfrac12$ Lemma~\ref{lem:gaussmom} applies, as $V_\rho(R_{\rm fBm})<\infty$
\citep[Prop.~15.5]{FrizVictoir2010}, and the extension theorem \citep[Thm.~9.5]{FrizVictoir2010} bounds
every level. For $P_1^{(\delta_0)}$ and $H>\tfrac12$ we use the increment variance of the ripple,
\begin{equation}\label{eq:Rvar}
   \Var(R^i_t-R^i_s)\ \le\ 2|t-s|^{2H}\qquad(0\le s<t\le1).
\end{equation}
Within a block this is $\delta_0^{2H}D_B\le2|t-s|^{2H}$ with $D_B$ of Remark~\ref{rem:cov4hyp}(b).
If a block node lies between $s$ and $t$, then $R_s$ and $R_t$ are independent and
$\Var(R_t-R_s)=v_R(s)+v_R(t)$, with $v_R(t)\le{\rm dist}(t,\text{nodes})^{2H}\le|t-s|^{2H}$ (same
remark), and likewise for $s$. So the rescaled cell path of $R$ has increment variance at most
$2|u-v|^{2H}$, uniformly in $a$, $\delta$ and $\delta_0$, and the Gaussian Kolmogorov criterion bounds
every moment of its $p$-variation. With the $P_0$ part and Minkowski's inequality, this gives (C2).
For $H\le\tfrac12$ we argue through the covariance. Each component's covariance has grid
$\rho$-variation, $\rho=\frac{1}{2H}$, with $V_\rho(R;[s,t]^2)^\rho\le C|t-s|$. For fBm this is
\citet[Prop.~15.5]{FrizVictoir2010}, as $2H\rho=1$. The ripple is block-independent and vanishes at
the block nodes, so a grid rectangle $[t_i,t_{i+1}]\times[t_j,t_{j+1}]$ receives a non-zero
contribution only from the at most two blocks containing an endpoint of $[t_i,t_{i+1}]$ in their
interior; hence $V_\rho(R_R;[s,t]^2)^\rho\le2^{\rho-1}\sum_JV_\rho(R_J;(J\cap[s,t])^2)^\rho$. Each block
covariance is $\delta_0^{2H}$ times the fBm covariance plus rank-one terms with factors $u$ and $C(u,1)$
in block coordinates, so $V_\rho(R_J;I^2)^\rho\le C|I|+C\delta_0(|I|/\delta_0)^{(1+2H)\rho}\le C|I|$ for
$I\subset J$, and $\sum_JC|J\cap[s,t]|=C|t-s|$. Grid $\rho$-variation is a seminorm, so
$V_\rho(R_{P_1};[s,t]^2)\le C|t-s|^{2H}$, and the covariance of the rescaled cell process has
$V_\rho\le C$ uniformly in $a$, $\delta$ and $\delta_0$. Lemma~\ref{lem:gaussmom} and the extension
theorem give (C2).

\emph{(C3), $H\ge\tfrac12$.} For $P_0$,
$\int_0^1\varrho(s,v)\dd v=H(s^{2H-1}+(1-s)^{2H-1})\le2H$ for $H>\tfrac12$, and $1$ for Brownian motion. For
$P_1^{(\delta_0)}$ the ripple adds $\epsilon^2\!\int|\varrho_R|(s,\dd v)\le C\epsilon^2\delta_0^{2H-1}\le
C\epsilon^2$. (C3) is stated in total variation because the ripple density changes sign.

\emph{(C4), $H<\tfrac12$.} \emph{The kernel condition.} For $P_0$,
$\sum_{b\ne a}\int_{I_a\times I_b}|\varrho_0|\le2\delta^{2H}\sum_{m\ge1}\bar r(m)$, with
$\bar r(m)=|c_H|\int_{[0,1]\times[m,m+1]}|u-v|^{2H-2}\asymp m^{2H-2}$ summable. For $P_1^{(\delta_0)}$,
$|\varrho_1|\le|\varrho_0|+\epsilon^2|\varrho_R|$. The first term of \eqref{eq:bHsplit} is dominated by
$|\varrho_0|$, and for a cell $I_a$ and a block $J$ meeting it, with $\tilde I$ the block coordinates of
$I_a\cap J$ and $\ell=|\tilde I|$,
$\sum_{b\ne a}\int_{(I_a\cap J)\times(I_b\cap J)}\delta_0^{2H-2}|g|\le
\delta_0^{2H}\int_{\tilde I\times[0,1]}|g|\le\delta_0^{2H}(2\ell+\textstyle\int_{\tilde I}f)\le3\delta_0^{2H}\ell^{2H}$.
If $\delta<\delta_0$, then $\ell\le h$ and a cell meets at most two blocks. If $\delta\ge\delta_0$, only the
at most two blocks straddling an end of $I_a$ contribute, with $\ell\le1$. Either way the sum is at most
$6\delta^{2H}$, and $\sum_{b\ne a}\int_{I_a\times I_b}|\varrho_1|\le\bigl(2(1+\epsilon^2)\sum_m\bar r(m)+6\epsilon^2\bigr)
\delta^{2H}$.

\emph{The decorrelation condition.} For $P_0$ it is Proposition~\ref{prop:cancorr}:
$\varrho_{\max}(I_a,I_b)\le C(H)\bar r(|a-b|)$. For $P_1^{(\delta_0)}$ write $G_1=G_0+\epsilon R$, $R$
independent of $G_0$. An element of the first chaos of $G_1$ on $I_a$ is $\xi=\xi_0+\epsilon\xi_R$ with
$\xi_0$, $\xi_R$ in the first chaoses of $G_0$, $R$ on $I_a$ and $\norm\xi^2=\norm{\xi_0}^2+\epsilon^2\norm{\xi_R}^2$
(the map $\sum_kc_kG_{1,J_k}\mapsto(\sum_kc_kG_{0,J_k},\sum_kc_kR_{J_k})$ preserves this identity and
extends to the closures). Cauchy--Schwarz in $\RR^2$ gives
$\varrho^{G_1}_{\max}\le\varrho^{G_0}_{\max}+\varrho^{R}_{\max}$, and the same argument over the independent
blocks gives $\varrho^R_{\max}(I_a,I_b)\le\max_J\varrho^{R^J}_{\max}(I_a\cap J,I_b\cap J)$, over the
blocks $J$ meeting both cells, $R^J$ the ripple of block $J$. If $\delta\ge\delta_0$, a block meets two
different cells only if it straddles a mesh node, so for each $a$ at most two $b\ne a$ contribute, each
at most $1$. Let $\delta<\delta_0$. A cell meets at most two blocks, and as $\varrho_{\max}$ is invariant
under scaling we work in block coordinates, where $R^J$ is the tied-down fBm $\varpi=Z-uZ_1$ and the
cells become pieces $I\subseteq[0,1]$ of length at most $h<1$. For a step function $\varphi$ on $I$,
$\varpi(\varphi):=\int\varphi\,\dd\varpi=Z(\varphi)-(\int\varphi)Z_1$ with $Z(\varphi):=\int\varphi\,\dd Z$.
Three facts hold. (a) By \eqref{eq:cancorrspec} and the Poincar\'e inequality of
Remark~\ref{rem:cancorrsharp}, rescaled to length $|I|$ (as $Z(\varphi)\overset d=|I|^HZ(\tilde\varphi)$ for
the rescaled $\tilde\varphi$), $\norm\varphi_{L^2(I)}\le C(H)|I|^{1/2-H}\norm{Z(\varphi)}_{L^2}$, hence
$|\int\varphi|\le C(H)h^{1-H}\norm{Z(\varphi)}_{L^2}$. (b) $\Cov(Z(\varphi),Z_1)=\int\varphi f$ with $f$ as
in \eqref{eq:bHsplit}, since $\partial_uC(u,1)=f(u)$. As $4H-2>-1$,
$\norm f_{L^2(I)}\le C(H)h^{2H-1/2}\kappa_I$ with $\kappa_I:=(1+x_I/h)^{2H-1}+(1+(1-y_I)/h)^{2H-1}\le2$ for
$I=[x_I,y_I]$, so $|\Cov(Z(\varphi),Z_1)|\le C(H)h^H\kappa_I\norm{Z(\varphi)}_{L^2}$. (c) By
Proposition~\ref{prop:cancorr}, $\varrho^Z_{\max}(I,I')\le C(H)\bar r(m)$ for pieces at lag $m\ge2$ of the
(shifted) mesh lattice, pieces being subsets of lattice cells and $\varrho_{\max}$ monotone under
inclusion. By (a) and (b), $\norm{\varpi(\varphi)}^2\ge\norm{Z(\varphi)}^2-2|\int\varphi|\,|\Cov(Z(\varphi),Z_1)|
\ge(1-C'h)\norm{Z(\varphi)}^2\ge\tfrac12\norm{Z(\varphi)}^2$ for $h\le h_0(H)$, and expanding
$\Cov(\varpi(\varphi),\varpi(\psi))$ in the same way,
\[
   \varrho^{\varpi}_{\max}(I,I')\ \le\ 2\bigl[\varrho^Z_{\max}(I,I')+Ch(\kappa_I+\kappa_{I'})+Ch^{2-2H}\bigr].
\]
Summing over the at most $1/h+2$ pieces $I'\ne I$ of a block: $\sum\varrho^Z_{\max}\le C(1+\sum_m\bar r(m))$,
$\sum h\kappa_I\le6$, $\sum_{I'}h\kappa_{I'}\le4h\sum_{k\le1/h+1}(1+k)^{2H-1}\le Ch^{1-2H}\le C$, and
$\sum h^{2-2H}\le3h^{1-2H}\le3$. For $h_0<h<1$ a block holds at most $1/h_0+2$ pieces, each with
$\varrho_{\max}\le1$. Hence $\sup_a\sum_{b\ne a}\varrho^{G_1}_{\max}(I_a,I_b)\le C(H)$, uniformly in $\delta$
and $\delta_0$.

\emph{(C5).} Throughout $\rho=1\vee1/(2H)$ as in \eqref{eq:rhoH}, and the object is the covariance
$R_G$, not the density. Write $R_{\rm fBm}$ for the fBm covariance. For $P_0$, $V_\rho(R_{\rm fBm})<\infty$ is classical and underlies
the lift \citep{CoutinQian2002,FrizVictoir2010}. $P_\theta$ has the same kernel, the drift $\theta te_1$ being
deterministic.

For $P_1^{(\delta_0)}$ the key is pinning, $\varpi_a(0)=\varpi_a(1)=0$ (block-diagonality alone
would not do: for independent restarted fBm on each block the grid $\{0,1\}$ already gives
$V_\rho^\rho=\delta_0^{(2H-1)\rho}\to\infty$ for $H<\tfrac12$). By pinning, every rectangular
increment of $R$ meets at most two blocks. Write $J_a$ for the $a$-th block, of width $\delta_0$ (the cells of
Construction~\ref{con:bridge} at width $\delta_0$. The cells $I_a$ of Definition~\ref{def:class} have
width $\delta$). Let $R_R=\sum_aR^{(a)}$, with $R^{(a)}$ the block-$a$ ripple covariance extended by
zero, and let $\Delta^a$ be its increment over the rectangle at hand. At most two $\Delta^a$ are
nonzero, so $|\Delta|^\rho\le2^{\rho-1}\sum_a|\Delta^a|^\rho$ and
\[
   V_\rho^\rho\bigl(R_R;[0,1]^2\bigr)\ \le\ 2^{\rho-1}\sum_aV_\rho^\rho\bigl(R^{(a)}\bigr)
   \ =\ 2^{\rho-1}\,\delta_0^{\,2H\rho-1}\,V_\rho^\rho(R_\varpi) .
\]
The inequality holds because one grid serves every block at once,
$\sum_{\rm cells}\sum_a|\Delta^a|^\rho=\sum_a\sum_{\rm cells}|\Delta^a|^\rho
\le\sum_aV_\rho^\rho(R^{(a)})$. $V_\rho^\rho$ itself is not additive. The equality holds because,
again by pinning, $V_\rho^\rho(R^{(a)};[0,1]^2)=V_\rho^\rho(R^{(a)};J_a^2)$, and each block carries a
copy of the covariance $R_\varpi=B_H$ of tied-down fBm scaled by $\delta_0^{2H}$, and there are
$n_0=1/\delta_0$ blocks, so $\sum_a\delta_0^{2H\rho}=\delta_0^{2H\rho-1}$. Here
$V_\rho(R_\varpi)<\infty$ follows from
$R_\varpi(u,v)=R_{\rm fBm}(u,v)-v\,R_{\rm fBm}(u,1)-u\,R_{\rm fBm}(v,1)+uv$
(Construction~\ref{con:bridge}; here $R_{\rm fBm}=C$): the three corrections are rank one with unit total variation in
each factor, so $V_\rho(R_\varpi)\le V_\rho(R_{\rm fBm})+3$. The exponent $2H\rho-1$ equals
$2H-1\ge0$ for $H\ge\tfrac12$, where $\rho=1$, and vanishes for $H<\tfrac12$, where
$\rho=\tfrac1{2H}$. So $\delta_0^{2H\rho-1}\le1$, and $2^{\rho-1}V_\rho^\rho(R_\varpi)\le2^{\rho-1}\bigl(V_\rho(R_{\rm fBm})+3\bigr)^\rho$ is uniform in
$\delta_0$, as Remark~\ref{rem:meshtied} requires.
The covariance of $P_1^{(\delta_0)}$ is $R_{\rm fBm}+\epsilon^2R_R$, and grid $\rho$-variation is a
seminorm, so $V_\rho(R_{P_1})\le V_\rho(R_{\rm fBm})+\epsilon^2V_\rho(R_R)$: (C5) holds for
$P_1^{(\delta_0)}$ at every $\delta_0$, with one constant depending on $H$ and $\epsilon$ only.

\emph{The moment condition.} It is imposed at $p=p_\delta(H)$ in the inhomogeneous norm, and for
$P_1^{(\delta_0)}$ we need one constant for every $\delta_0$ (Remark~\ref{rem:meshtied}). Above the
crossover, $p\in(1/H,2)$, the lift is canonical and $\norm{\mathbf X}_{\pvar}$ is the path $p$-variation.
For $P_0$, $\E|G_t-G_s|^2=|t-s|^{2H}$, and for $P_1^{(\delta_0)}$ the same holds with the constant $2$ for
the ripple by \eqref{eq:Rvar}. The Gaussian Kolmogorov criterion, with $\norm z_{\pvar}\le
\norm z_{\alpha\text{-H\"ol}}$ for $\alpha\in(1/p,H)$, bounds every moment uniformly in $\delta_0$, and
Minkowski's inequality combines the two parts. At or below the crossover, $p_\delta(H)=\tfrac12(1/H+4)\in(2\rho,4)$.
Lemma~\ref{lem:gaussmom}, with the bound $V_\rho(R_{P_1})\le V_\rho(R_{\rm fBm})+\epsilon^2V_\rho(R_R)$ of
(C5), uniform in $\delta_0$, bounds every moment of the homogeneous norm of the natural lift, and the
ball--box estimate converts it to the inhomogeneous norm. Young cross-integration of $G$ against $R$ is
not needed. 

\emph{The drift.} For $P_\theta$ write $X=\theta te_1+G$ with $G\sim P_0$. Above the
crossover $\lfloor p\rfloor=1$, and the drift adds its $p$-variation $|\theta|$, so
$\norm{\mathbf X}_{\pvar}\le|\theta|+\norm G_{\pvar}\le(1+|\theta|)(1+\norm G_{\pvar})$. At or below
the crossover $\mathbf X$ is the translate of $\mathbf G$ by the drift (Definition~\ref{def:class}).
The drift has finite $1$-variation $|\theta|$, so the translation estimate
\citep[Thm.~9.33(i)]{FrizVictoir2010} gives $\hnorm{\mathbf X}_{\pvar}\le C(p)\bigl(\hnorm{\mathbf
G}_{\pvar}+|\theta|\bigr)$, and the ball--box comparison \citep[Prop.~8.9]{FrizVictoir2010} bounds the
inhomogeneous norm by $C(p,d)\bigl(1+\hnorm{\mathbf X}_{\pvar}\bigr)^{\lfloor p\rfloor}$. In both cases
\[
   \norm{\mathbf X}_{\pvar}\ \le\ C(p,d)\,(1+|\theta|)^{\lfloor p\rfloor}\,
   \bigl(1+\hnorm{\mathbf G}_{\pvar}\bigr)^{\lfloor p\rfloor},
\]
with $\hnorm{\mathbf G}_{\pvar}=\norm G_{\pvar}$ above the crossover. All moments of
$\hnorm{\mathbf G}_{\pvar}$ are finite, by Lemma~\ref{lem:gaussmom} at or below the crossover and by
the Gaussian Kolmogorov criterion above it. Hence
 $\E\norm{\mathbf X}^{2M}_{\pvar}\le C_{\rm dr}(1+|\theta|)^{2M\lfloor p\rfloor}$ with $C_{\rm dr}=C_{\rm dr}(H,M,d)\ge1$, so
$P_\theta\in\cP(H,L)$ whenever $L\ge\max\{L_{P_0},\,C_{\rm dr}(1+|\theta|)^{2M\lfloor p\rfloor}\}$, where
$L_{P_0}$ is the constant of (C1)--(C5) for $P_0$, which $P_\theta$ shares (this also gives
$|\theta|\le L$); in particular for $|\theta|\le1$ and $L\ge\max\{L_{P_0},2^{2M\lfloor p\rfloor}C_{\rm dr}\}$.
Taking $L_0(H,\epsilon)$ larger than all constants above for $P_0$ and $P_1^{(\delta_0)}$ proves Proposition~\ref{prop:adm}(a).
\end{remark}
\begin{remark}[the diagonal hypotheses at the two laws the lemma is applied to]\label{rem:cov4hyp}
We check \emph{(D1)}--\emph{(D4)} at the two laws to which Lemma~\ref{lem:cov4} is applied.
Throughout, $\varsigma(t)$ denotes the distance from $t$ to the nearest node of the mesh or block
grid in play, and $\varsigma_\mp(t)$ the distances to the two flanking nodes (not to be confused with
the $\rho$ of \eqref{eq:rhoH}).

\emph{(D1) and (D2).} These are immediate, because every increment variance below is $|t-s|^{2H}$ plus a term
whose gradient is bounded by an integrable power. In case~\emph{(b)}, $\partial_1D_1$ blows up like
$\varsigma^{2H-1}$ along the block-node lines, at any distance from the diagonal. This singularity is
integrable, $2H-1>-1$, and $\partial_2D_R(\cdot,t)$ is absolutely continuous across each node, with
mixed derivative $O(\varsigma^{2H-1})$ there. This gives the second condition of \emph{(D2)}, which for
\emph{(a)} reduces to $|t-s|^{4H-2}\in L^1$, i.e.\ $H>\tfrac14$.

\emph{(D4).} Besides membership in the class (Remark~\ref{rem:C1C2}), (D4) has three further conditions. The first, $\sup_n\norm{D_n}_\infty<\infty$, holds at both
laws because $X^n_t$ is a convex combination of two nodal values, so
$\norm{X^n_t}_{L^2}\le\sup_t\norm{X_t}_{L^2}$ and $D_n\le4\sup_tv(t)$. The second, the domination of
$\varrho_*$, holds because off the band $D_n(s,t)\lesssim|t-s|^{2H}$ (triangle inequality with
$\norm{X^n_t-X_t}_{L^2}\lesssim\delta^H\le|t-s|^H$), while the cell average of $\varrho_*$ is
$\lesssim|t-s|^{2H-2}$, the corner singularity of adjacent cells being integrable since $2H-2>-2$.
In case~\emph{(b)} this needs $D_R\lesssim|s-t|^{2H}$, which follows from $\Lambda\ge0$:
$D_B=(v-u)^{2H}-(v-u)\Lambda+(v-u)^2\le2(v-u)^{2H}$ within a block, and across blocks
$D_R=v_R(s)+v_R(t)\le2|t-s|^{2H}$ by the bound on $v_R$ below. The third condition is a split
$\varrho=\varrho_*+\tilde g$ with $\tilde g\in L^1$; we write $\tilde g$ for the integrable part
called $g$ in (D4), to keep it apart from the $g$ of \eqref{eq:bHsplit}. For \emph{(a)} take
$\varrho_*=\varrho$, $\tilde g=0$.
For \emph{(b)} take $\varrho_*=\varrho_{\rm fBm}+\epsilon^2c_H|s-t|^{2H-2}\one_{\text{same block}}$
and $\tilde g:=\epsilon^2\delta_0^{2H-2}g(u,v)$ on each block square $J\times J$, in block
coordinates $u,v$, and $\tilde g:=0$ across blocks, with $g=1-f(u)-f(v)$ and $f$ as in
\eqref{eq:bHsplit}. This is the split \eqref{eq:bHsplit}. The global fBm
kernel must be in $\varrho_*$ because $\varrho$ in \emph{(D4)} is the whole increment covariance of
the law and $c_H|s-t|^{2H-2}\notin L^1$ for $H<\tfrac12$. The function $\tilde g$ is integrable because
$2H-1>-1$, with a $\delta_0$-dependent constant, which suffices since $\delta_0$ is fixed here.

\emph{(D3).} This is the only condition which needs work. $\Xi$ is a quadratic form in $D$.

\emph{(a) $P_0$, the independent-\textnormal{fBm} law.} $D_0(s,t)=|t-s|^{2H}$, so
$\Xi(\tau)=4H(1-\tau)\tau^{4H-1}\to0$ exactly on $H>\tfrac14$.

\emph{(b) $P_1=P_0+\epsilon R$, the pair of Construction~\ref{con:bridge}}, at a fixed block
width $\delta_0=1/n_0$, which is the grid the nodes below refer to, and $x=\epsilon^2$, so that
$D_1=D_0+xD_R$. Two facts about the ripple suffice.
First, $\varpi_a$ is the tied-down fBm $Z_u-uZ_1$ of Construction~\ref{con:bridge}, for
which $B_H(u,u)-u^{2H}=u\bigl[u-u^{2H}-1+(1-u)^{2H}\bigr]\le0$  (for $H\le\tfrac12$ because
$u^{2H}\ge u$ and $(1-u)^{2H}\le1$, for $H\ge\tfrac12$ because $(1-u)^{2H}\le1-u\le u^{2H}+1-u$),
while $B_H(u,u)=B_H(1-u,1-u)$ by the reversal $Z_1-Z_{1-\cdot}\overset{d}=Z$, so
$v_R(t)=\Var(R_t)\le\min\{\varsigma_-(t),\varsigma_+(t)\}^{2H}\le\varsigma(t)^{2H}$, and
$|v_R'(t)|\le C(H)\bigl[\varsigma_-(t)^{2H-1}+\varsigma_+(t)^{2H-1}\bigr]$ from the explicit
$v_R(t)=\delta_0^{2H}B_H(u,u)$, $u:=(t-u_a)/\delta_0$ the position within block $J_a=[u_a,u_a+\delta_0]$.
Second, within one block,
$D_R(s,t)=\delta_0^{2H}D_B\bigl(\tfrac{s-u_a}{\delta_0},\tfrac{t-u_a}{\delta_0}\bigr)$ with
\begin{gather*}
   D_B(u,v)=(v-u)^{2H}-(v-u)\Lambda(u,v)+(v-u)^2,\\
   \Lambda(u,v)=v^{2H}-u^{2H}+(1-u)^{2H}-(1-v)^{2H}\in[0,2] ,
\end{gather*}
so that $D_R(s,s+\tau)=\tau^{2H}+O(\delta_0^{2H-1}\tau)$ and
$[\partial_2-\partial_1]D_R(s,s+\tau)=4H\tau^{2H-1}+O(\delta_0^{2H-1})+O(\delta_0^{-1}\tau\,\Pi)$ with
$\Pi:=\varsigma_-(s)^{2H-1}+\varsigma_+(s)^{2H-1}+\varsigma_-(s+\tau)^{2H-1}+\varsigma_+(s+\tau)^{2H-1}$, whose
$s$-integral over $[0,1]$ is $O(\delta_0^{2H-1})$. Where $s$ and $s+\tau$ straddle a node the ripple
covariance vanishes and $D_R(s,s+\tau)=v_R(s)+v_R(s+\tau)\le2\tau^{2H}$, both points lying
within $\tau$ of that node, while $[\partial_2-\partial_1]D_R=v_R'(s+\tau)-v_R'(s)$ has
$\int|v_R'|\le C(H)\tau^{2H}/H$ over the straddling interval. Hence the $n_0$ straddling
intervals contribute $O(n_0\tau^{4H})$ and the block interiors $4H(1+x)^2\tau^{4H-1}$, so
\[
   \Xi[D_1](\tau)\ =\ 4H(1+\epsilon^2)^2\,\tau^{4H-1}\ +\ O_{\delta_0}\bigl(\tau^{2H}\bigr)
   \ \longrightarrow\ 0\qquad(H>\tfrac14).
\]
Here $n_0$ is fixed while $\tau\downarrow0$.
\end{remark}

\section{The finite-part functional and the two regularisations}\label{app:fp}

Proposition~\ref{prop:Lexists} evaluates $\kappa_{\rm alt}$ as an absolutely convergent series with
a single $\zeta(2-4H)$ term, and Proposition~\ref{prop:kappaclosed} as a $\zeta$-series obtained by
continuing three individually divergent power sums. We prove that the two agree. Both are values of
one linear functional, the constant term in the large-$N$ expansion of a partial sum. Once this
functional is set up, linearity identifies them, and the $\zeta$-values come from Euler--Maclaurin.
Throughout $H\in(\tfrac14,\tfrac12)$, so that
\begin{equation}\label{eq:expwindow}
   4H+1\in(2,3),\quad 4H\in(1,2),\quad 4H-1\in(0,1),\quad 4H-2\in(-1,0),\quad 4H-3\in(-2,-1);
\end{equation}
none of these exponents is an integer, so no logarithmic term arises below. Both endpoints are
excluded: $4H-2\downarrow-1$ as $H\downarrow\tfrac14$ (the pole of $\zeta(2-4H)$, matching
$\kappa_{\rm alt}\to-\infty$), and $4H-2\uparrow0$ as $H\uparrow\tfrac12$.

\begin{definition}[finite part]\label{def:fp}
Let $\mathfrak F$ be the set of functions $\varphi:\NN\to\RR$ for which there are an integer
$p\ge0$, exponents $\alpha_1>\cdots>\alpha_p>0$ and reals $A_1,\dots,A_p,c$ with
\[
   \varphi(N)=\sum_{i=1}^{p}A_iN^{\alpha_i}+c+o(1)\qquad(N\to\infty,\ N\in\NN).
\]
Write $\mathrm{fp}(\varphi):=c$. For $f:\NN\to\RR$ put $S_f(N):=\sum_{h\le N}f(h)$ and, whenever
$S_f\in\mathfrak F$, set $\mathfrak f(f):=\mathrm{fp}(S_f)$.
\end{definition}

The constant $c$ is unique, because a finite sum of distinct positive powers of $N$ plus a
constant is $o(1)$ only if all its coefficients vanish. A linear combination of two representations
is again one, after merging equal exponents, so $\mathfrak F$ is a vector space, and by uniqueness
$\mathrm{fp}$ is linear on it. The restriction to integer $N$ matters: for $f(h)\asymp h^\sigma$,
$\sigma>0$, the step function $S_f$ has jumps that grow without bound, so no such expansion holds
for real $N$.

The next lemma collects the four rules by which we compute $\mathfrak f$. The proof of
Proposition~\ref{prop:Lexists} uses (F1)--(F3), and Proposition~\ref{prop:regagree} below uses all four.

\begin{lemma}[the finite-part toolkit]\label{lem:fp}
\emph{(F1)} If $\sum_h|f(h)|<\infty$ then $S_f\in\mathfrak F$ and $\mathfrak f(f)=\sum_{h\ge1}f(h)$.

\noindent\emph{(F2)} $\mathfrak f$ is linear on its domain.

\noindent\emph{(F3)} For $-1<\sigma<3$, $\sigma\notin\ZZ_{\ge0}$,
\[
   \sum_{h\le N}h^{\sigma}=\frac{N^{\sigma+1}}{\sigma+1}+\frac{N^{\sigma}}{2}+\zeta(-\sigma)
   +\frac{\sigma}{12}\,N^{\sigma-1}+O(N^{\sigma-3}),
\]
so $S_{h\mapsto h^\sigma}\in\mathfrak F$ and $\mathfrak f(h\mapsto h^{\sigma})=\zeta(-\sigma)$.

\noindent\emph{(F4)} Let $\beta\in\RR\setminus\ZZ$, let $(c_i)_{i\ge0}$ satisfy
$C:=\sum_{i\ge0}|c_i|<\infty$, and put $i_0:=\max(0,\lfloor\beta\rfloor+1)$. Then for integers
$N\ge2$ the series $\varphi(N):=\sum_{i\ge0}c_iN^{\beta-i}$ converges absolutely and
\[
   \varphi(N)=\sum_{0\le i<i_0}c_iN^{\beta-i}+o(1),
\]
a finite sum of strictly positive powers. In particular $\varphi\in\mathfrak F$ and
$\mathrm{fp}(\varphi)=0$.
\end{lemma}
\begin{proof}
(F1) and (F2) are immediate from Definition~\ref{def:fp}. For (F3), Euler--Maclaurin for the
Riemann zeta function gives, for every $s\ne1$,
$\sum_{h\le N}h^{-s}=\zeta(s)+\frac{N^{1-s}}{1-s}+\frac{N^{-s}}{2}-\frac{s}{12}N^{-s-1}
+O(N^{-s-3})$; substituting $s=-\sigma$ yields the display. The surviving exponents
$\sigma+1,\sigma,\sigma-1$ are all nonzero because $\sigma\notin\ZZ_{\ge0}$, and the
remaining terms are $o(1)$ because $\sigma<3$. For (F4), $\beta-i_0<0$ by construction, and $\beta-i>0$ strictly for $i<i_0$ because
$\beta\notin\ZZ$. The remaining terms satisfy
$\bigl|\sum_{i\ge i_0}c_iN^{\beta-i}\bigr|\le N^{\beta-i_0}\sum_{i\ge i_0}|c_i|\le CN^{\beta-i_0}=o(1)$,
since $N^{\beta-i}\le N^{\beta-i_0}$ for $i\ge i_0$ and $N\ge1$.
\end{proof}

The hypothesis $\sigma\notin\ZZ_{\ge0}$ in (F3) is needed: at $\sigma=0$ the finite part of
$\sum_{h\le N}1=N$ is $0$, not $\zeta(0)=-\tfrac12$. Only $\sigma=4H\in(1,2)$ and
$\sigma=4H-2\in(-1,0)$ occur below.

Let $r(m)=\tfrac12\bigl[|m+1|^{2H}-2|m|^{2H}+|m-1|^{2H}\bigr]$ be
the fGn autocovariance, $T(h):=\sum_{m\ge h}r(m)$ and $U(h):=T(h)-\tfrac12r(h)$ for $h\ge1$. For
$H<\tfrac12$ the sum defining $T$ telescopes, since $(m+1)^{2H}-m^{2H}\to0$, and $U$ is the function
of Proposition~\ref{prop:Lexists}:
\[
   T(h)=\sum_{m\ge h}r(m)=-\tfrac12\bigl[h^{2H}-(h-1)^{2H}\bigr],\qquad
   U(h)=T(h)-\tfrac12r(h)=-\tfrac14\bigl[(h+1)^{2H}-(h-1)^{2H}\bigr].
\]
The next proposition computes the finite part of $\sum_hU(h)^2$ as a $\zeta$-series. It is the
identity behind Proposition~\ref{prop:kappaclosed}, through Corollary~\ref{cor:regagree}.

\begin{proposition}[the $\zeta$-series is the finite part]\label{prop:regagree}
For every $H\in(\tfrac14,\tfrac12)$ the partial sums of $U^2$ lie in $\mathfrak F$ and
\[
   \mathfrak f\bigl(h\mapsto U(h)^2\bigr)
   \;=\;-\frac{1}{16}+\frac{1}{8}\sum_{j\ge1}\Bigl|\tbinom{2H}{j}\Bigr|\,\zeta(2j-4H).
\]
\end{proposition}
\begin{proof}
We split $U^2$ into three power-type terms, take the finite part of each with Lemma~\ref{lem:fp},
and add. The two pure powers give $\zeta(-4H)$ by (F3). In the mixed term the binomial series turns
$(h^2-1)^{2H}$ into a series of powers, each of which (F3) evaluates.

\emph{Step 1 (a three-term splitting).} Squaring the closed form for $U$,
\begin{equation}\label{eq:U2split}
   16\,U(h)^2=(h+1)^{4H}+(h-1)^{4H}-2\,(h^2-1)^{2H}\qquad(h\ge1),
\end{equation}
the last two terms vanishing at $h=1$. Take the finite part of each summand and add by (F2).

\emph{Step 2 (the two pure powers).} By (F3) with $\sigma=4H\in(1,2)$,
\begin{align*}
    \sum_{h\le N}(h+1)^{4H}&=\sum_{k\le N+1}k^{4H}-1
   =\frac{(N{+}1)^{4H+1}}{4H+1}+\frac{(N{+}1)^{4H}}{2}+\bigl(\zeta(-4H)-1\bigr)\\
   &+\ \frac{4H}{12}(N{+}1)^{4H-1}+O(N^{4H-3}).
\end{align*}
By (F4) and $\beta\notin\ZZ$ (see \eqref{eq:expwindow}), each
$(N+1)^{\beta}=\sum_{i\ge0}\binom\beta iN^{\beta-i}$ with $\beta\in\{4H+1,4H,4H-1\}$ lies in
$\mathfrak F$ with finite part $0$ (the coefficients are summable, $|\binom\beta i|\asymp i^{-1-\beta}$
with $\beta>0$),
$O(N^{4H-3})=o(1)$. Hence $\mathfrak f\bigl(h\mapsto(h+1)^{4H}\bigr)=\zeta(-4H)-1$. Identically,
$\sum_{h\le N}(h-1)^{4H}=\sum_{k\le N-1}k^{4H}$, with $(N-1)^\beta=\sum_{i\ge0}\binom\beta i(-1)^iN^{\beta-i}$ in (F4), gives
$\mathfrak f\bigl(h\mapsto(h-1)^{4H}\bigr)=\zeta(-4H)$.

\emph{Step 3 (the mixed term).} For $k\ge2$ the binomial series
$(k^2-1)^{2H}=k^{4H}(1-k^{-2})^{2H}=\sum_{j\ge0}(-1)^j\binom{2H}{j}k^{4H-2j}$ converges absolutely.
Separating $j=0,1$,
\[
   (k^2-1)^{2H}=k^{4H}-2H\,k^{4H-2}+G(k),\qquad
   G(k):=\sum_{j\ge2}(-1)^j\tbinom{2H}{j}k^{4H-2j}=O(k^{4H-4}),
\]
and $4H-4<-1$, so $\sum_{k\ge2}|G(k)|<\infty$. Since $4H-2\in(-1,0)$, (F3) applies to both power
sums, and with (F1)--(F2),
\[
   F:=\mathfrak f\bigl(h\mapsto(h^2-1)^{2H}\bigr)
   =\bigl(\zeta(-4H)-1\bigr)-2H\bigl(\zeta(2-4H)-1\bigr)+\sum_{k\ge2}G(k),
\]
the two subtracted $1$'s being the omitted $k=1$ terms $1^{4H}=1^{4H-2}=1$. Fubini applies because
$\sum_{j\ge2}\bigl|\binom{2H}{j}\bigr|\bigl(\zeta(2j-4H)-1\bigr)<\infty$ (the second factor
is $O(2^{-2j})$), and gives
$\sum_{k\ge2}G(k)=\sum_{j\ge2}(-1)^j\binom{2H}{j}\bigl(\zeta(2j-4H)-1\bigr)$. Since
$\sum_{j\ge0}(-1)^j\binom{2H}{j}=(1-1)^{2H}=0$, the constants $-1$ telescope away across all
$j\ge0$ and
\begin{equation}\label{eq:Fseries}
   F=\sum_{j\ge0}(-1)^j\tbinom{2H}{j}\,\zeta(2j-4H),
\end{equation}
absolutely convergent since $\bigl|\binom{2H}{j}\bigr|\asymp j^{-1-2H}$ and $\zeta(2j-4H)\to1$.

\emph{Step 4 (assembly).} Applying $\mathfrak f$ to \eqref{eq:U2split} and using (F2),
$16\,\mathfrak f(U^2)=\bigl(\zeta(-4H)-1\bigr)+\zeta(-4H)-2F$, and the $j=0$ term of
\eqref{eq:Fseries} cancels both copies of $\zeta(-4H)$:
\[
   16\,\mathfrak f(U^2)=-1-2\sum_{j\ge1}(-1)^j\tbinom{2H}{j}\zeta(2j-4H).
\]
For $2H\in(0,1)$ and $j\ge1$ the product $2H(2H-1)\cdots(2H-j+1)$ has exactly $j-1$ negative
factors, so $\mathrm{sign}\binom{2H}{j}=(-1)^{j-1}$ and $(-1)^j\binom{2H}{j}=-\bigl|\binom{2H}{j}
\bigr|$. The display follows.
\end{proof}

The next corollary shows that the values of $\kappa_{\rm alt}$ given by
Propositions~\ref{prop:Lexists} and~\ref{prop:kappaclosed} are both $\mathfrak f(U^2)$.

\begin{corollary}[$\kappa_{\rm alt}$ has one value]\label{cor:regagree}
For every $H\in(\tfrac14,\tfrac12)$,
\[
   \underbrace{H^2\zeta(2-4H)+\sum_{h\ge1}\Bigl[U(h)^2-H^2h^{4H-2}\Bigr]}_{\text{Proposition~\ref{prop:Lexists}}}
   \;=\;\underbrace{-\frac{1}{16}+\frac{1}{8}\sum_{j\ge1}\Bigl|\tbinom{2H}{j}\Bigr|\zeta(2j-4H)}
   _{\text{Proposition~\ref{prop:kappaclosed}}} .
\]
\end{corollary}
\begin{proof}
Both sides are $\mathfrak f(U^2)$. On the left write
$U(h)^2=H^2h^{4H-2}+\bigl[U(h)^2-H^2h^{4H-2}\bigr]$. The bracket is $O(h^{4H-4})$, hence absolutely
summable, so (F1)--(F3) give
$\mathfrak f(U^2)=H^2\zeta(2-4H)+\sum_{h\ge1}[U(h)^2-H^2h^{4H-2}]$, which is the left-hand side. The right-hand side is Proposition~\ref{prop:regagree}.
\end{proof}

The $\zeta$-series converges slowly (Remark~\ref{rem:kappanumbers}). The next corollary rewrites it
so that it converges geometrically. Proposition~\ref{prop:roughwall} uses this form.

\begin{corollary}[a geometrically convergent form]\label{cor:geomkappa}
For every $H\in(\tfrac14,\tfrac12)$,
\[
   \kappa_{\rm alt}(H)=\frac{1}{16}+\frac{1}{8}\sum_{j\ge1}\Bigl|\tbinom{2H}{j}\Bigr|
   \bigl(\zeta(2j-4H)-1\bigr).
\]
\end{corollary}

\begin{proof}
By Proposition~\ref{prop:Lexists} and Corollary~\ref{cor:regagree},
$\kappa_{\rm alt}=-\tfrac1{16}+\tfrac18\sum_{j\ge1}b_j\zeta(2j-4H)$ with $b_j:=|\binom{2H}{j}|$. For
$j\ge1$ the product $2H(2H-1)\cdots(2H-j+1)$ has exactly $j-1$ negative factors, so
$(-1)^j\binom{2H}{j}=-b_j$, and $0=(1-1)^{2H}=\sum_{j\ge0}(-1)^j\binom{2H}{j}=1-\sum_{j\ge1}b_j$. Hence
$\sum_{j\ge1}b_j=1$, and, both series converging absolutely,
\[
   \kappa_{\rm alt}=-\tfrac1{16}+\tfrac18\sum_{j\ge1}b_j+\tfrac18\sum_{j\ge1}b_j\bigl(\zeta(2j-4H)-1\bigr)
   =\tfrac1{16}+\tfrac18\sum_{j\ge1}b_j\bigl(\zeta(2j-4H)-1\bigr).
\]
\end{proof}

The
$j=1$ term is negative, as $\zeta(2-4H)<-\tfrac12$ (Lemma~\ref{lem:zetabound}). Every term with $j\ge2$ has $2j-4H>2$ and obeys
$0<\zeta(s)-1\le2^{-s}\bigl(1+\tfrac2{s-1}\bigr)$, so the tail decays geometrically: the error after
$J$ terms is $O(4^{-J})$, against the $O(J^{-2H})$ of the raw series (Remark~\ref{rem:kappanumbers}).

\clearpage
\bibliographystyle{plainnat}
\bibliography{references}

@article{Schenck2026,
  author  = {Schenck, Bryson},
  title   = {Finite-Sample Bounds for Expected Signature Estimation under Weak Dependence},
  journal = {arXiv preprint arXiv:2605.20541}, year = {2026}}

@article{Lyons1998,
  author = {Lyons, Terry J.}, title = {Differential equations driven by rough signals},
  journal = {Revista Matem\'atica Iberoamericana}, volume = {14}, number = {2},
  pages = {215--310}, year = {1998}}

@book{LyonsCaruanaLevy2007,
  author = {Lyons, Terry J. and Caruana, Michael and L\'evy, Thierry},
  title = {Differential Equations Driven by Rough Paths},
  series = {Lecture Notes in Mathematics}, volume = {1908},
  publisher = {Springer}, address = {Berlin}, year = {2007}}

@book{FrizVictoir2010,
  author = {Friz, Peter K. and Victoir, Nicolas B.},
  title = {Multidimensional Stochastic Processes as Rough Paths: Theory and Applications},
  series = {Cambridge Studies in Advanced Mathematics}, volume = {120},
  publisher = {Cambridge University Press}, year = {2010}}

@article{CassLittererLyons2013,
  author = {Cass, Thomas and Litterer, Christian and Lyons, Terry},
  title = {Integrability and tail estimates for {G}aussian rough differential equations},
  journal = {The Annals of Probability}, volume = {41}, number = {4},
  pages = {3026--3050}, year = {2013}}

@article{CoutinQian2002,
  author = {Coutin, Laure and Qian, Zhongmin},
  title = {Stochastic analysis, rough path analysis and fractional {B}rownian motions},
  journal = {Probability Theory and Related Fields}, volume = {122}, number = {1},
  pages = {108--140}, year = {2002}}

@article{FrizRiedel2014,
  author = {Friz, Peter K. and Riedel, Sebastian},
  title = {Convergence rates for the full {G}aussian rough paths},
  journal = {Annales de l'Institut Henri Poincar\'e, Probabilit\'es et Statistiques},
  volume = {50}, number = {1}, pages = {154--194}, year = {2014}}

@article{NeuenkirchTindelUnterberger2010,
  author = {Neuenkirch, Andreas and Tindel, Samy and Unterberger, J\'er\'emie},
  title = {Discretizing the fractional {L}\'evy area},
  journal = {Stochastic Processes and their Applications}, volume = {120}, number = {2},
  pages = {223--254}, year = {2010}}

@article{Gubinelli2010,
  author = {Gubinelli, Massimiliano}, title = {Ramification of rough paths},
  journal = {Journal of Differential Equations}, volume = {248}, number = {4},
  pages = {693--721}, year = {2010}}

@article{HairerKelly2015,
  author = {Hairer, Martin and Kelly, David},
  title = {Geometric versus non-geometric rough paths},
  journal = {Annales de l'Institut Henri Poincar\'e, Probabilit\'es et Statistiques},
  volume = {51}, number = {1}, pages = {207--251}, year = {2015}}

@article{CassFerrucci2024,
  author = {Cass, Thomas and Ferrucci, Emilio},
  title = {On the {W}iener chaos expansion of the signature of a {G}aussian process},
  journal = {Probability Theory and Related Fields}, volume = {189}, number = {3},
  pages = {909--947}, year = {2024}}

@article{BaudoinCoutin2007,
  author = {Baudoin, Fabrice and Coutin, Laure},
  title = {Operators associated with a stochastic differential equation driven by fractional {B}rownian motions},
  journal = {Stochastic Processes and their Applications}, volume = {117}, number = {5},
  pages = {550--574}, year = {2007}}

@book{Janson1997,
  author = {Janson, Svante}, title = {Gaussian Hilbert Spaces},
  series = {Cambridge Tracts in Mathematics}, volume = {129},
  publisher = {Cambridge University Press}, year = {1997}}

@article{KolmogorovRozanov1960,
  author = {Kolmogorov, Andrey N. and Rozanov, Yuri A.},
  title = {On strong mixing conditions for stationary {G}aussian processes},
  journal = {Theory of Probability and its Applications}, volume = {5}, number = {2},
  pages = {204--208}, year = {1960}}

@book{Reutenauer1993,
  author = {Reutenauer, Christophe}, title = {Free {L}ie Algebras},
  series = {London Mathematical Society Monographs, New Series}, volume = {7},
  publisher = {Oxford University Press}, year = {1993}}

@article{ChevyrevOberhauser2022,
  author = {Chevyrev, Ilya and Oberhauser, Harald},
  title = {Signature Moments to Characterize Laws of Stochastic Processes},
  journal = {Journal of Machine Learning Research}, volume = {23}, number = {176},
  pages = {1--42}, year = {2022}}

@article{Levin2013,
  author = {Levin, Daniel and Lyons, Terry and Ni, Hao},
  title = {Learning from the Past, Predicting the Statistics for the Future, Learning an Evolving System},
  journal = {arXiv preprint arXiv:1309.0260}, year = {2013}}

@incollection{Lyons2014,
  author = {Lyons, Terry},
  title = {Rough Paths, Signatures and the Modelling of Functions on Streams},
  booktitle = {Proceedings of the International Congress of Mathematicians (ICM), Seoul},
  publisher = {Kyung Moon SA}, pages = {163--184}, year = {2014}}

@article{HagerPelizzari2026,
  author  = {Hager, Paul P. and Pelizzari, Luca},
  title   = {Expected signatures via partial integration, coordinate change and symmetrization},
  journal = {arXiv preprint arXiv:2607.29534}, year = {2026}}

@article{FrizVictoir2011,
  author  = {Friz, Peter K. and Victoir, Nicolas},
  title   = {A note on higher dimensional $p$-variation},
  journal = {Electronic Journal of Probability},
  volume  = {16},
  pages   = {1880--1899},
  year    = {2011},
  doi     = {10.1214/EJP.v16-951}
}

@book{AdlerTaylor2007,
  author    = {Adler, Robert J. and Taylor, Jonathan E.},
  title     = {Random Fields and Geometry},
  series    = {Springer Monographs in Mathematics},
  publisher = {Springer},
  address   = {New York},
  year      = {2007}
}

@unpublished{HuberII, author={Huber, Florian}, title={Optimal estimation of expected signatures under weak dependence: minimax rates and an ergodic central limit theorem}, note={in preparation}, year={2026}}

@article{DeyaNeuenkirchTindel2012,
  author  = {Deya, Aur{\'e}lien and Neuenkirch, Andreas and Tindel, Samy},
  title   = {A {M}ilstein-type scheme without {L}{\'e}vy area terms for {SDE}s driven by fractional {B}rownian motion},
  journal = {Ann. Inst. Henri Poincar{\'e} Probab. Stat.},
  volume  = {48},
  number  = {2},
  pages   = {518--550},
  year    = {2012},
  doi     = {10.1214/10-AIHP392},
  eprint  = {1001.3344},
  archivePrefix = {arXiv}
}

@article{ChevyrevLyons2016,
  author  = {Chevyrev, Ilya and Lyons, Terry},
  title   = {Characteristic functions of measures on geometric rough paths},
  journal = {Ann. Probab.},
  volume  = {44},
  number  = {6},
  pages   = {4049--4082},
  year    = {2016},
  doi     = {10.1214/15-AOP1068},
  eprint  = {1307.3580},
  archivePrefix = {arXiv}
}

@article{Passeggeri2020,
  author  = {Passeggeri, Riccardo},
  title   = {On the signature and cubature of the fractional {B}rownian motion for {$H>\frac12$}},
  journal = {Stochastic Process. Appl.},
  volume  = {130},
  number  = {3},
  pages   = {1226--1257},
  year    = {2020},
  doi     = {10.1016/j.spa.2019.04.013},
  eprint  = {1609.07352},
  archivePrefix = {arXiv}
}

@article{BayerFrizRiedelSchoenmakers2016,
  author  = {Bayer, Christian and Friz, Peter and Riedel, Sebastian and Schoenmakers, John},
  title   = {From rough path estimates to multilevel {M}onte {C}arlo},
  journal = {SIAM J. Numer. Anal.},
  volume  = {54},
  number  = {3},
  pages   = {1449--1483},
  year    = {2016},
  doi     = {10.1137/140995209},
  eprint  = {1305.5779},
  archivePrefix = {arXiv}
}

@article{HuLiuNualart2016,
  author  = {Hu, Yaozhong and Liu, Yanghui and Nualart, David},
  title   = {Rate of convergence and asymptotic error distribution of {E}uler approximation schemes for fractional diffusions},
  journal = {Ann. Appl. Probab.},
  volume  = {26},
  number  = {2},
  pages   = {1147--1207},
  year    = {2016},
  doi     = {10.1214/15-AAP1114},
  eprint  = {1408.6471},
  archivePrefix = {arXiv}
}

@article{LiuTindel2019,
  author  = {Liu, Yanghui and Tindel, Samy},
  title   = {First-order {E}uler scheme for {SDE}s driven by fractional {B}rownian motions: the rough case},
  journal = {Ann. Appl. Probab.},
  volume  = {29},
  number  = {2},
  pages   = {758--826},
  year    = {2019},
  doi     = {10.1214/17-AAP1374},
  eprint  = {1703.03625},
  archivePrefix = {arXiv}
}

@book{Nualart2006,
  author = {Nualart, David},
  title = {The {M}alliavin Calculus and Related Topics},
  edition = {2nd}, series = {Probability and its Applications},
  publisher = {Springer}, address = {Berlin}, year = {2006}}

@unpublished{HuberStrong,
  author = {Huber, Florian},
  title  = {The mean-square rate of piecewise-linear signature approximation for fractional {B}rownian motion in the {Y}oung regime},
  note   = {Preprint},
  year   = {2026}}

\end{document}